# Higher order influence functions and minimax estimation of nonlinear functionals

**James Robins**[1], **Lingling Li**[1], **Eric Tchetgen**[1] **and Aad van der Vaart**[2]

*Harvard School of Public Health and Vrije Universiteit Amsterdam*

**Abstract:** We present a theory of point and interval estimation for nonlinear functionals in parametric, semi-, and non-parametric models based on higher order influence functions (Robins (2004), Section 9; Li et al. (2004), Tchetgen et al. (2006), Robins et al. (2007)). Higher order influence functions are higher order U-statistics. Our theory extends the first order semiparametric theory of Bickel et al. (1993) and van der Vaart (1991) by incorporating the theory of higher order scores considered by Pfanzagl (1990), Small and McLeish (1994) and Lindsay and Waterman (1996). The theory reproduces many previous results, produces new non-$\sqrt{n}$ results, and opens up the ability to perform optimal non-$\sqrt{n}$ inference in complex high dimensional models. We present novel rate-optimal point and interval estimators for various functionals of central importance to biostatistics in settings in which estimation at the expected $\sqrt{n}$ rate is not possible, owing to the curse of dimensionality. We also show that our higher order influence functions have a multi-robustness property that extends the double robustness property of first order influence functions described by Robins and Rotnitzky (2001) and van der Laan and Robins (2003).

## 1. Introduction

Over the past 3 years, we have developed a theory of point and interval estimation for nonlinear functionals $\psi(F)$ in parametric, semi-, and non-parametric models based on higher order likelihood scores and influence functions that applies equally to both $\sqrt{n}$ and non-$\sqrt{n}$ problems (Robins [16], Section 9, Li et al. [9], Tchetgen et al. [21], Robins et al. [18]). The theory reproduces results previously obtained by the modern theory of non-parametric inference, produces many new non-$\sqrt{n}$ results, and most importantly opens up the ability to perform non-$\sqrt{n}$ inference in complex high dimensional models, such as models for the estimation of the causal effect of time varying treatments in the presence of time varying confounding and informative censoring. See Tchetgen et al. [22] for examples of the latter.

Higher order influence functions are higher order U-statistics. Our theory extends the first order semiparametric theory of Bickel et al. [3] and van der Vaart [25] by incorporating the theory of higher order scores and Bhattacharrya bases considered by Pfanzagl [11], Small and McLeish [20] and Lindsay and Waterman [8].

[1]Harvard School of Public Health, Department of Biostatistics and Epidemiology, 677 Huntington Ave, Boston, MA 02115, USA, e-mail: robins@hsph.harvard.edu; lingling_li@hphc.org; etchetge@hsph.harvard.edu

[2]Vrije Universiteit Amsterdam, Department of Mathematics, e-mail: aad@cs.vu.nl

*AMS 2000 subject classifications:* 62G05, 62C20.

*Keywords and phrases:* influence functions, nonparametric, minimax, robust inference, U-statistics.

The purpose of this paper is to demonstrate the scope and flexibility of our methodology by deriving rate-optimal point and interval estimators for various functionals that are of central importance to biostatistics. We now describe some of these functionals. We suppose we observe i.i.d copies of a random vector $O = (Y, A, X)$ with unknown distribution $F$ on each of $n$ study subjects. In this paper, we largely study non-parametric models that place no restrictions on $F$, other than bounds on both the $L_p$ norms and on the smoothness of certain density and conditional expectation functions. The variable $X$ represents a random vector of baseline covariates such as age, height, weight, hematocrit, and laboratory measures of lung, renal, liver, brain, and heart function. $X$ is assumed to have compact support and a density $f_X(x)$ with respect to the Lebesgue measure in $R^d$, where, in typical applications, $d$ is in the range 5 to 100. $A$ is a binary treatment and $Y$ is a response, higher values of which are desirable. Then, in the absence of confounding by additional unmeasured factors, the functional $\psi(F) = E\{E[Y|A=1,X]\} - E\{E[Y|A=0,X]\}$ is the mean effect of treatment in the total study population. Our results for $E\{E[Y|A=1,X]\} - E\{E[Y|A=0,X]\}$ follow from results for the functional $\psi(F) = E\{E[Y|A=1,X]\}$ based on data $(AY, A, X)$ rather than $(Y, A, X)$. If $Y$ is missing for some study subjects, and $A$ is now the indicator that takes the value 1 when $Y$ is observed and zero otherwise, then the functional $E\{E[Y|A=1,X]\}$ is the marginal mean of $Y$ under the missing at random assumption that the probability $P[A=0|X,Y] = P[A=0|X]$ that $Y$ is missing does not depend on the unobserved $Y$.

Returning to data $O = (Y, A, X)$, the functional

$$\psi(F) = E\{\text{cov}(Y, A|X)\}/E[\text{var}\{A|X\}]$$
$$= E[w(X)\{E[Y|A=1,X] - E[Y|A=0,X]\}],$$

with $w(X) = \text{var}\{A|X\}/E[\text{var}\{A|X\}]$ is the variance weighted average treatment effect. Our results for $E\{\text{cov}(Y,A|X)\}/E[\text{var}\{A|X\}]$ are derived from results for the functionals $\psi(F) = E\{\text{cov}(Y,A|X)\}$ and $\psi(F) = E\left[\{E(Y|X)\}^2\right]$.

We note that Robins and van der Vaart's [19] construction of an adaptive confidence set for a regression function $E(Y|X=x)$ depended on being able to construct a confidence interval for $\psi(F) = E\left[\{E(Y|X)\}^2\right]$. They constructed an interval for $E\left[\{E(Y|X)\}^2\right]$ when the marginal distribution of $X$ was known. In this paper, we construct a confidence interval for $E\left[\{E(Y|X)\}^2\right]$ when the marginal of $X$ is unknown and, in Section 5, use it to obtain an adaptive confidence set for $E(Y|X=x)$.

The functional $E\{\text{cov}(Y,A|X)\}$ is the functional $E\{\text{var}(Y|X)\}$ in the special case in which $Y = A$ w.p.1. Minimax estimation of $\text{var}(Y|X)$ has recently been discussed by Wang et al. [27] and Cai et al. [6] in the setting of non-random $X$.

The function $\gamma(x) = E[Y|A=1, X=x] - E[Y|A=0, X=x]$ is the effect of treatment on the subgroup with $X = x$. It is important to estimate the function $\gamma(x)$, in addition to the average treatment effect in the total population, because treatment should be given, since beneficial, to those subjects with $\gamma(x) > 0$ but withheld, since harmful, from subjects with $\gamma(x) < 0$. We show that one can obtain adaptive confidence sets for $\gamma(x)$ if one can set confidence intervals for the functional $\psi(F) = E\left[\gamma(X)^2\right]$. We construct intervals for $E\left[\gamma(X)^2\right]$ under the additional assumption that the data $O = (Y, A, X)$ came from a randomized trial. In a randomized trial, in contrast to an observational study, the randomization probabilities, $P(A = 1|X) = E(A|X)$ are known by design. We plan to report confidence intervals for $E\left[\gamma(X)^2\right]$ with $E(A|X)$ unknown elsewhere.

All of the above functionals $\psi(F)$ have a positive semiparametric information bound (SIB) and thus a (first order) efficient influence function with a finite variance. In fact all the functionals $\psi(F)$ have efficient influence function

(1.1) $\qquad IF(b(F), p(F), \psi(F)) \equiv if(O, b(X, F), p(X, F), \psi(F)),$

where $b(x, F), p(x, F)$ are functions of certain conditional expectations, and, for any $b^*(x), p^*(x),$

$$E_F[IF(b^*, p^*, \psi(F))] = E_F[h_1(O)\{b^*(X) - b(X; F)\}\{p^*(X) - p(X; F)\}]$$

where $h_1(O)$ is a known function. We refer to functionals in our class as doubly-robust to indicate that $IF(b(F), p(F), \psi(F))$ continues to have mean zero when either (but not both) $p(F)$ is misspecified as $p^*$ or $b(F)$ is misspecified as $b^*$. The functions $b(x, F), p(x, F), IF(O, b(X, F), p(X, F), \psi(F))$, and $h_1(O)$ differ depending on the functional $\psi(F)$ of interest.

As the functionals $\psi(F)$ are all closely related, we shall use $E\{\text{cov}(Y, A|X)\}$ as a prototype in this introduction. For $\psi(F) \equiv E\{\text{cov}(Y, A|X)\}$, $b(X; F) = E_F(Y|X)$, $p(X; F) = E_F(A|X)$,

$$IF(b(F), p(F), \psi(F)) = \{Y - b(X; F)\}\{A - p(X; F)\} - \psi(F),$$

and $h_1(O) \equiv 1$.

Whenever a functional $\psi(F)$ has a non-zero SIB, given sufficiently stringent bounds on $L_p$ norms and on smoothness, it is possible to use the estimated first order influence function to construct regular estimators and honest asymptotic confidence intervals whose width shrinks at the usual parametric rate of $n^{-1/2}$. We recall that, by definition, regular estimators are $n^{1/2}$-consistent. When $X$ is high dimensional, the a priori smoothness restrictions on $p(X; F)$ and $b(X; F)$ necessary for point or interval estimators of $E\{\text{cov}(Y, A|X)\}$ to achieve the parametric rate of $n^{-1/2}$ are so severe as to be substantively implausible. As a consequence, we replace the usual approach based on first order influence functions by one based on higher order influence functions.

To provide quantitative results, we require a measure of the maximal possible complexity (e.g. smoothness) of $p(\cdot; F)$ and $b(\cdot; F)$ believed substantively plausible. We use Hölder balls for concreteness, although our methods extend to other measures of complexity. A function $h(\cdot)$ lies in the Hölder ball $H(\beta, C)$, with Hölder exponent $\beta > 0$ and radius $C > 0$, if and only if $h(\cdot)$ is bounded in supremum norm by $C$ and all partial derivatives of $h(x)$ up to order $\lfloor\beta\rfloor$ exist, and all partial derivatives of order $\lfloor\beta\rfloor$ are Lipschitz with exponent $(\beta - \lfloor\beta\rfloor)$ and constant $C$. We make the assumption that $b(\cdot, F), p(\cdot, F)$ lie in given Hölder balls $H(\beta_b, C_b), H(\beta_p, C_p)$. Furthermore, it turns out we must also make assumptions about the complexity of the function $g(X; F) \equiv E_F[h_1(O)|X] f_X(X)$, which we take to lie in a given $H(\beta_g, C_g)$. For $\psi(F) = E\{\text{cov}(Y, A|X)\}, g(X; F) = f_X(X)$.

Using higher order influence functions, we construct regular estimators and honest (i.e uniform over our model) asymptotic confidence intervals for functionals $\psi(F)$ in our class whose width shrinks at the usual parametric rate of $n^{-1/2}$ whenever $\beta/d \equiv \frac{\beta_b + \beta_p}{2}/d > 1/4$ and $\beta_g > 0$. This result cannot be improved on, since even when $g(x)$ is known a priori, $\beta/d > 1/4$ is necessary for a regular estimator to exist.

When $\beta/d \leq 1/4$ and $g(x)$ is known a priori, we have shown using arguments similar to those of Birge and Massart [5] that the minimax rate of convergence for an estimator and minimax rate of shrinkage of a confidence interval is $n^{-\frac{4\beta/d}{4\beta/d+1}} \geq n^{-\frac{1}{2}}$. When $g(x)$ is unknown, we construct point and interval estimators with this same rate of $n^{-\frac{4\beta/d}{4\beta/d+1}}$ whenever

$$(1.2) \qquad \beta_g/d > \beta/d \frac{2(\Delta+1)(1-4\beta/d)}{(\Delta+2)(1+4\beta/d) - 4(\beta/d)(1-4\beta/d)(\Delta+1)},$$

where $\Delta = \left|\frac{\beta_p}{\beta_b} - 1\right|$. For example, if $\Delta = 0$, $\beta/d = 1/8$, we require $\beta_g/d$ exceed $1/22$ to achieve the rate $n^{-\frac{4\beta/d}{4\beta/d+1}}$. When the previous inequality does not hold and $\Delta = 0$, we have constructed, in a yet unpublished paper, estimators that converge at rate

$$(1.3) \qquad \log(n) \, n^{-\frac{1}{2} + \frac{\beta_g/d}{1+2\beta_g/d} \frac{(m^*+1)^2}{2\beta/d}}, \text{ with}$$

$$m^* \equiv \left\lceil \left( \left[\frac{\beta}{d}\left(4\frac{\beta}{d} + \left(1 - 4\frac{\beta}{d}\right)\frac{1+2\beta_g/d}{\beta_g/d}\right)\right]^{1/2} - (1+2\beta/d) \right) \right\rceil.$$

We conjecture that this rate is minimax, up to log factors. In this paper, however, the estimators we construct are inefficient when the previous inequality fails to hold, converging at rates less than the conjectured minimax rate of Equation (1.3).

Let us return to the case where $Y = A$ w.p.1. Then $\psi(F) = E\{\text{var}(Y|X)\}$ and $p(\cdot) = b(\cdot)$ so $\Delta = 0$. Now, for fixed $\beta$, Equation (1.3) converges to $\log(n) \, n^{-2\beta/d}$ as $\beta_g \to 0$, which agrees (up to a log factor) with the minimax rate of $n^{-2\beta/d}$ given by Wang et al. [27] and Cai et al. [6] under the semiparametric homoscedastic model $\text{var}(Y|X) = \sigma^2$ with equal-spaced non-random $X$. This result might suggest that $X$ being random rather than equal-spaced can result in faster rates of convergence only when the density of $X$ has some smoothness, as quantified here by $\beta_g > 0$. But this suggestion is not correct. Recall that we obtained the rate $\log(n) \, n^{-2\beta/d}$ for $\psi(F) = E\{\text{var}(Y|X)\}$ as $\beta_g \to 0$ under a non-parametric model. In Section 4, we construct a simple estimator of $\sigma^2$ under the homoscedastic model with $X$ random with unknown density that, for $\beta/d < 1/4, \beta < 1$, and without smoothness restrictions on $f_X(x)$, converges at the rate $n^{-\frac{4\beta/d}{4\beta/d+1}}$, which is faster than the equal-spaced non-random minimax rate of $n^{-2\beta/d}$.

The paper is organized as follows. In Section 2, we define the higher order (estimation) influence functions of a functional $\psi(F)$ for $F$ contained in a model $\mathcal{M}$ and prove two fundamental theorems – the extended information equality theorem and the efficient estimation influence function theorem. Further, in the context of a parametric model whose dimension increases with sample size, we outline why estimators based on higher order influence can outperform those based on first order influence functions in high-dimensional models. In Section 3, we introduce the class of functionals we study in the remainder of the paper and describe their importance in biostatistics. The theory of Section 2, however, is not directly applicable to these functionals because they have first order but not higher order influence functions. We show that higher order influence functions fail to exist precisely because the Dirac delta function is not an element of the Hilbert space $L_2$ of square integrable functions. We describe two approaches to overcoming this difficulty. The first approach is based on approximating the Dirac delta function by a projection operator onto a subspace of $L_2$ of dimension $k(n)$, where $k(n)$ can be as large as the square of the sample size $n$. The second approach is based on approximating the functional $\psi(F)$ by a truncated functional $\widetilde{\psi}_{k(n)}(F)$. The truncated functional has influence functions of all orders, is equal to $\psi(F)$ if either a $k(n)$ dimensional working parametric model (with $k(n) < n^2$) for the function $b(\cdot)$ or the function $p(\cdot)$ in Equation (1.1) is correct, and remains close to $\psi(F)$ even if both working models are misspecified. We then use higher order influence function based estimators of $\widetilde{\psi}_{k(n)}(F)$ as estimators of $\psi(F)$. These estimators $\widehat{\psi}_{m,k(n)}$ are asymptotically normal with variance and bias for $\psi(F)$ depending both on the choice of the dimension $k(n)$ of the working models and on the order $m$ of the influence function of $\widetilde{\psi}_{k(n)}(F)$. We show that these same estimators $\widehat{\psi}_{m,k(n)}$ can also be obtained under the approximate Dirac delta function approach. We derive the optimal estimator $\widehat{\psi}_{m_{opt},k_{opt}(n)}(\beta_b, \beta_p, \beta_g)$ in the class as a function of the Hölder balls in which the functions $b$, $p$, and $g$ are assumed to lie. Finally we conclude Section 3 by showing that the estimators $\widehat{\psi}_{m,k(n)}$ have a multi-robustness property that extends the double-robustness property of the first order influence function estimator $\widehat{\psi}_1$.

In Section 4, we consider whether the estimators $\widehat{\psi}_{m_{opt},k_{opt}(n)}(\beta_b, \beta_p, \beta_g)$ are rate-minimax. We show that whenever $\beta/d \equiv \frac{\beta_b + \beta_p}{2}/d > 1/4$ and $\beta_g > 0$, $\widehat{\psi}_{m_{opt},k_{opt}(n)}(\beta_b, \beta_p, \beta_g)$ is not only rate minimax but is semiparametric efficient. Further, by letting the order $m = m(n)$ of the U-statistic depend on sample size, we construct a single estimator $\widehat{\psi}_{m(n),k(n)}$ that is semiparametric efficient for all $\beta/d > 1/4$ even when $g(\cdot)$ cannot be estimated at an algebraic rate. We show, however, that when $\beta/d < 1/4$, $\widehat{\psi}_{m_{opt},k_{opt}(n)}(\beta_b, \beta_p, \beta_g)$ does not in general converge at the minimax rate. In Section 4.1, however, we construct a new estimator $\widehat{\psi}^{eff}_{\mathcal{K},J}(\beta_g, \beta_b, \beta_p)$ that converges at the minimax rate of $n^{-\frac{4\beta/d}{4\beta/d+1}}$ whenever Eq. (1.2) holds. In Section 5, we use the results obtained in earlier sections to construct adaptive confidence intervals for a regression function $E[Y|X=x]$ when the marginal of $X$ is unknown and for the treatment effect function and optimal treatment regime in a randomized clinical trial. In Section 6.1, we discuss how to obtain higher order U-statistic point estimators and confidence intervals for functionals $\tau(F)$ that are implicitly defined as the solution to an equation $\psi(\tau, F) = 0$. In Section 6.2, we define higher order testing influence functions and efficient scores and describe their relationship to the higher order estimation influence functions and efficient influence functions of Section 2. Finally, in Section 6.3, we discuss the relationship between the higher order U-statistic point estimators of an implicitly defined functional $\tau(F)$ and higher order testing influence functions.

Before proceeding, several additional comments are in order. In this paper, we investigate the asymptotic properties of our higher order U-statistic point and interval estimators. The reader is referred to Li et al. [9] for an investigation of the finite sample properties of our procedures through simulation. Furthermore due to space limitations we only provide proofs for selected theorems. Proofs of the remaining theorems can be found in an accompanying technical report. In addition, precise regularity conditions are sometimes omitted from both the statements and the proofs of various theorems. This reflects the fact that the goal of this paper is to provide a broad overview of our theory as it currently stands.

Different subject matter experts will clearly disagree as to the maximum possible complexity of $p(x; F)$, $b(x; F)$ and $g(x; F)$. Thus it is important to have methods that adapt to the actual smoothness of these functions. Elsewhere, we plan to provide point estimators that optimally adapt to unknown smoothness. In contrast to point estimators, however, for honest confidence intervals, the degree of possible adaption to unknown smoothness is small. Therefore we propose that an analyst should report a mapping from a priori smoothness assumptions encoded in the exponents and radii of Hölder balls (or in other measures of complexity) to the associated optimal $1 - \alpha$ honest confidence intervals proposed in this paper. Such a mapping is finally only useful if substantive experts can approximately quantify their informal opinions concerning the shape and wiggliness of $p$, $b$, and $g$ using the measure of complexity on offer by the analyst. It is an open question which, if any, complexity measure is suitable for this purpose.

Finally, most of our mathematical results concern rates of convergence. We offer only a few results on the constants in front of those rates. This is not because the constant is less important than the rate in predicting how a proposed procedure will perform in the moderate sized samples occurring in practice. Rather, at present, we do not possess the mathematical tools necessary to obtain useful results concerning constants. A more extended discussion of the issue is found in Section 3 of Li et al. [9].

In the following, we use $X_n \asymp Y_n$ to mean $X_n = O_p(Y_n)$ and $Y_n = O_p(X_n)$; $X_n \sim Y_n$ to mean $\frac{X_n}{Y_n} \xrightarrow{P} 1$; and $X_n \gg Y_n$ ($X_n \ll Y_n$) to respectively mean $\frac{Y_n}{X_n} \xrightarrow{P} 0$ $\left(\frac{X_n}{Y_n} \xrightarrow{P} 0\right)$ as $n \to \infty$.

## 2. Theory of higher order influence functions

Given $n$ i.i.d observations $\mathbf{O} \equiv \mathbf{O}_n \equiv \{O_i, i = 1, \ldots, n\}$ from a model

$$\mathcal{M}(\Theta) = \{F(\cdot; \theta), \theta \in \Theta\},$$

we consider inference on a nonlinear functional $\psi(\theta)$. In general, $\psi(\theta)$ can be infinite dimensional but for now we only consider the one dimensional case. In the following all quantities can depend on the sample size $n$, including the support of $O$, the parameter space $\Theta$, and the functional $\psi(\theta)$. We generally suppress the dependence on $n$ in the notation. We will be particularly interested in models in which the parameter $\theta$ is infinite dimensional and $\theta$, $\Theta$, and $\psi(\cdot)$ do not depend on $n$. We also briefly discuss models in which subvectors of $\theta$ are finite-dimensional parameters whose dimension $k(n) = n^{1+\rho}$ increases as power $1 + \rho$ (often $\rho > 0$) of $n$ and thus $\theta_n, \Theta_n$, and $\psi_n(\cdot)$ depend on $n$.

Our first task is to define higher order influence functions. Before proceeding we recall some facts about $U$−statistics. Consider a function $b_m(o_1, o_2, \ldots, o_m) \equiv b(o_1, o_2, \ldots, o_m)$ where we often suppress $b$'s subscript $m$. For integers $i_1, i_2, \ldots, i_m$ lying in $\{1, \ldots, n\}$, we define

$$B_{m,i_1,\ldots,i_m} \equiv b_m(O_{i_1}, O_{i_2}, \ldots, O_{i_m}) \equiv b(O_{i_1}, O_{i_2}, \ldots, O_{i_m})$$

and

$$\mathbb{V}_n[b_m] \equiv \frac{(n-m)!}{n!} \sum_{i_1 \neq i_2 \ldots \neq i_m} B_{m,i_1,\ldots,i_m}.$$

In an abuse of notation, we will consider the following expressions to be equivalent

$$\mathbb{V}_n[B_m] \equiv \mathbb{V}_n[B_{m,,i_1,\ldots,i_m}] \equiv \mathbb{V}_n[b_m].$$

Thus $\mathbb{V}_n[b_m]$ is an $m^{\text{th}}$ order U-statistic with kernel $b_m(o_1, o_2, \ldots, o_m)$. We do not assume that $b_m(o_1, o_2, \ldots, o_m)$ is symmetric. We will write $\mathbb{V}_n[B_m]$ as $\mathbb{B}_{n,m}$. So, suppressing the dependence on $n$, $\mathbb{B}_m \equiv \mathbb{V}[B_m]$.

Any $\mathbb{B}_m$ has a unique (up to permutation) decomposition $\mathbb{B}_m = \sum_{s=1}^{m} \mathbb{D}_s^{(b)}(\theta)$ under any $F(.;\theta)$ as a sum of degenerate U-statistics $\mathbb{D}_s^{(b)}(\theta)$, where degeneracy of $\mathbb{D}_s^{(b)}(\theta)$ means that $D_s^{(b)}(\theta) = d_s^{(b)}(O_{i_1}, O_{i_2}, \ldots, O_{i_s}; \theta)$ satisfies

$$E_\theta\left[d_s^{(b)}\left(o_{i_1}, \ldots, o_{i_{l-1}}, O_{i_l}, o_{i_{l+1}}, \ldots, o_{i_s}; \theta\right)\right] = 0, l = 1, \ldots, s,$$

where upper and lower case letters, respectively, denote random variables and their possible realizations.

Let $\mathcal{U}_m(\theta)$ be the Hilbert space of all $U$-statistics of order $m$ with mean zero and finite variance with inner product defined by covariances with respect to the $n$-fold product measure $F^n(\cdot;\theta)$. Note that any $U$-statistic $\mathbb{B}_s$ of order $s$, $s < m$, is also an $m^{\text{th}}$ order $U$- statistic with $\mathbb{D}_l^{(b)}(\theta)$ identically zero for $m \geq l > s$.

Since any two degenerate $U$- statistics of different orders are uncorrelated, the $\mathcal{U}_m(\theta)$-Hilbert space projection of $\mathbb{B}_m$ on $\mathcal{U}_l(\theta)$ is $\sum_{s=1}^{l} \mathbb{D}_s^{(b)}(\theta)$ for $l < m$. Thus a $U$-statistic $\mathbb{B}_m$ is degenerate $\Leftrightarrow \mathbb{B}_m = \mathbb{D}_m^{(b)}(\theta) \Leftrightarrow \Pi_\theta[\mathbb{B}_m|\mathcal{U}_{m-1}(\theta)] = 0 \Leftrightarrow \mathbb{B}_m \in \mathcal{U}_{m-1}(\theta)^{\perp_{m,\theta}}$, where $\Pi_\theta[\cdot|\cdot] \equiv \Pi_{\theta,m}[\cdot|\cdot]$ is the projection operator of the Hilbert space $\mathcal{U}_m(\theta)$ (with the dependence on $m$ suppressed when no ambiguity can arise) and, for any linear subspace $\mathcal{R}$ of $\mathcal{U}_m(\theta)$, $\mathcal{R}^{\perp_{m,\theta}}$ is its orthocomplement in the Hilbert space $\mathcal{U}_m(\theta)$. Given any $\mathbb{B}_m = \mathbb{V}[B_m]$, $\mathbb{D}_m^{(b)}(\theta)$ is explicitly given by $\mathbb{V}[d_{m,\theta}\{B_m\}]$ where $d_{m,\theta}$ maps $B_m \equiv b(O_{i_1}, O_{i_2}, \ldots, O_{i_m})$ to

$$(2.1) \qquad d_{m,\theta}\{B_m\} = b(O_{i_1}, O_{i_2}, \ldots, O_{i_m})$$
$$+ \sum_{t=0}^{m-1} (-1)^{m-t} \sum_{i_{r_1} \neq i_{r_2} \cdots \neq i_{r_t}} E_\theta\left(b(O_{i_1}, O_{i_2}, \ldots, O_{i_m}) | O_{i_{r_1}}, O_{i_{r_2}}, \ldots, O_{i_{r_t}}\right).$$

Given a function $g(\zeta)$, $\zeta \equiv \{\zeta_1, \ldots, \zeta_r\}^T$, define for $m = 0, 1, 2, \ldots$,

$$g_{\backslash \bar{l}_m}(\zeta) \equiv g_{\backslash l_1, \ldots, l_m}(\zeta) \equiv \frac{\partial^m g(\zeta)}{\partial \zeta_{l_1} \ldots \partial \zeta_{l_m}}$$

with $l_s \in \{1, \ldots, r\}$, where the $\backslash$ symbol denotes differentiation by the variables occurring to its right and the overbar $\bar{l}_m$ denotes the vector $(l_1, \ldots, l_m)$. Given a sufficiently smooth $r$-dimensional parametric submodel $\widetilde{\theta}(\zeta)$ mapping $\zeta \in R^r$ injectively into $\Theta$, define for $\theta$ in the range of $\widetilde{\theta}(\cdot)$, $\psi_{\backslash \bar{l}_m}(\theta) \equiv \left(\psi \circ \widetilde{\theta}\right)_{\backslash l_1, \ldots, l_m}(\zeta)|_{\zeta = \widetilde{\theta}^{-1}\{\theta\}}$ and $f_{\backslash \bar{l}_m}(\mathbf{O}_n; \theta) \equiv \left(f \circ \widetilde{\theta}\right)_{\backslash l_1, \ldots, l_m}(\zeta)|_{\zeta = \widetilde{\theta}^{-1}\{\theta\}}$, where $f(\mathbf{O}_n; \theta) \equiv \prod_i f(O_i; \theta)$ is the density of $\mathbf{O}_n$ with respect to a dominating measure. That is $\psi_{\backslash \bar{l}_m}(\theta)$ and $f_{\backslash \bar{l}_m}(\mathbf{O}_n, \theta)$ are higher order derivatives of $\psi(\cdot)$ and $f(\mathbf{O}_n; \cdot)$ under a parametric submodel $\widetilde{\theta}(\zeta)$, where the model $\widetilde{\theta}$ has been suppressed in the notation. An $s^{\text{th}}$ order score associated with the submodel $\widetilde{\theta}(\zeta)$ is defined to be

$$\widetilde{\mathbb{S}}_{s,\bar{l}_s}(\theta) \equiv f_{\backslash \bar{l}_s}(\mathbf{O}_n; \theta) / f(\mathbf{O}_n; \theta),$$

where $\widetilde{\mathbb{S}}_{s,\bar{l}_s}(\theta)$ is a U-statistic of order $s$. To understand why $\widetilde{\mathbb{S}}_{s,\bar{l}_s}(\theta)$ is a $U$-statistic we provide formulae for an arbitrary score $\widetilde{\mathbb{S}}_{s,\bar{l}_s}(\theta)$ in terms of the subject specific scores

$$S_{l_1 \ldots l_m, j}(\theta) \equiv f_{/l_1 \ldots l_m, j}(O_j; \theta) / f_j(O_j; \theta)$$

$j = 1, \ldots, n$ for $s = 1, 2, 3$. Suppressing the $\theta$-dependence, results in Waterman and Lindsay [8] imply

$$\widetilde{\mathbb{S}}_{1,\overline{l}_1} = \sum_j S_{l_1,j},$$

$$\widetilde{\mathbb{S}}_{2,\overline{l}_2} = \sum_j S_{l_1 l_2,j} + \sum_{s \neq j} S_{l_1,j} S_{l_2,s},$$

$$\widetilde{\mathbb{S}}_{3,\overline{l}_3} = \sum_j S_{l_1 l_2 l_3, j} + \sum_{s \neq j} S_{l_1 l_2, j} S_{l_3, s} + S_{l_3 l_2, j} S_{l_1, s} + S_{l_1 l_3, j} S_{l_2, s}$$
$$+ \sum_{s \neq j \neq t} S_{l_1,j} S_{l_2,s} S_{l_3,t}.$$

Note these equations express each $\widetilde{\mathbb{S}}_{m,\overline{l}_m}$ as a sum of degenerate U-statistics. We now define a $m^{\text{th}}$ order estimation influence function $\mathbb{IF}_{m,\psi(\cdot)}(\theta) \equiv \mathbb{IF}_{m,\psi}(\theta) \equiv \mathbb{IF}_m(\theta)$ for $\psi(\theta)$ where we suppress the dependence on $\psi$ when no ambiguity will arise.

**Definition 2.1.** A U-statistic $\mathbb{IF}_m(\theta)$ of order $m$ and finite variance is said to be an $m^{\text{th}}$ order estimation influence function for $\psi(\theta)$ if (i) $E_\theta[\mathbb{IF}_m(\theta)] = 0$, $\theta \in \Theta$ and (ii) for $s = 1, 2, \ldots, m$ and every suitably smooth and regular (see Appendix) $r$ dimensional parametric submodel $\widetilde{\theta}(\zeta)$, $r = 1, 2, \ldots, m$,

$$\psi_{\backslash \overline{l}_s}(\theta) = E_\theta\left[\mathbb{IF}_m(\theta) \widetilde{\mathbb{S}}_{s,\overline{l}_s}(\theta)\right].$$

Estimation influence functions need not always exist, but when they do they are useful for deriving point estimators of $\psi$ with small bias and for deriving confidence interval estimators centered on an estimate of $\psi$. We will generally refer to estimation influence functions simply as influence functions. We remark that $\mathbb{IF}_m(\theta)$ is an influence function under the above definition if and only if it is one under the modified version in which the dimension of the parametric submodel $\widetilde{\theta}(\zeta)$ is unrestricted. A key result is the following theorem which is related to results of Small and McLeish [20].

**Theorem 2.2** (Extended information equality theorem). *Given a $m^{\text{th}}$ order influence function $\mathbb{IF}_m(\theta)$, for any smooth, regular submodel $\widetilde{\theta}(\zeta)$ and $s \leq m$,*

$$\partial^s E_\theta\left[\mathbb{IF}_m\left(\widetilde{\theta}(\zeta)\right)\right] / \partial \zeta_{l_1} \cdots \partial \zeta_{l_s}|_{\zeta = \widetilde{\theta}^{-1}\{\theta\}} = -\psi_{\backslash \overline{l}_s}(\theta)$$

*Thus, if the functionals $E_\theta[\mathbb{IF}_m(\theta^*)]$ and $-[\psi(\theta^*) - \psi(\theta)]$ have bounded Fréchet derivatives with respect to $\theta^*$ to order $m + 1$ for a norm $||\cdot||$,*

$$E_\theta[\mathbb{IF}_m(\theta + \delta\theta)] = -[\psi(\theta + \delta\theta) - \psi(\theta)] + O\left(||\delta\theta||^{m+1}\right)$$

*since the functions $E_\theta[\mathbb{IF}_m(\theta^*)]$ and $-[\psi(\theta^*) - \psi(\theta)]$ of $\theta^*$ have the same Taylor expansion around $\theta$ up to order $m$.*

The proof is in the Appendix. Define the $m^{\text{th}}$ order tangent space $\Gamma_m(\theta)$ at $\theta$ for the model $\mathcal{M}(\Theta)$ to be the subspace of $\mathcal{U}_m(\theta)$ formed by taking the closed linear span of all scores of order $m$ or less as we vary over all regular parametric submodels $\widetilde{\theta}(\varsigma)$ (whose range includes $\theta$) of our model $\mathcal{M}(\Theta)$. We say a model is (locally) nonparametric for $m^{\text{th}}$ order inference if $\Gamma_m(\theta) = \mathcal{U}_m(\theta)$.

Given any $m^{\text{th}}$ order estimation influence function $\mathbb{IF}_m(\theta)$, define the $m^{\text{th}}$ order efficient estimation influence function to be

$$\mathbb{IF}_m^{eff}(\theta) = \Pi_\theta \left[ \mathbb{IF}_m(\theta) | \Gamma_m(\theta) \right],$$

where $\Pi_\theta [\cdot|\cdot] \equiv \Pi_{\theta,m}[\cdot|\cdot]$ is the $\mathcal{U}_m(\theta)$-projection operator. In the appendix, we prove the following:

**Theorem 2.3** (Efficient estimation influence function theorem)**.**

1. $\mathbb{IF}_m^{eff}(\theta)$ *is unique in the sense that for any two* $m^{th}$ *order influence functions*

$$\Pi_\theta \left[ \mathbb{IF}_m^{(1)}(\theta) | \Gamma_m(\theta) \right] = \Pi_\theta \left[ \mathbb{IF}_m^{(2)}(\theta) | \Gamma_m(\theta) \right] \quad a.s.$$

2. $\mathbb{IF}_m^{eff}(\theta)$ *is a* $m^{th}$ *order estimation influence function and has variance less than or equal to any other* $m^{th}$ *order estimation influence function.*
3. $\mathbb{IF}_m(\theta)$ *is a* $m^{th}$ *order estimation influence function if and only if*

$$\mathbb{IF}_m(\theta) \in \left\{ \mathbb{IF}_m^{eff}(\theta) + \mathbb{U}_m(\theta) ; \mathbb{U}_m(\theta) \in \Gamma_m^{\perp_{m,\theta}}(\theta) \right\}$$

*where* $\Gamma_m^{\perp_{m,\theta}}(\theta)$ *is the ortho-complement of* $\Gamma_m(\theta)$ *in* $\mathcal{U}_m(\theta)$.
4. *If* $\mathbb{IF}_m(\theta)$ *exists then* $\mathbb{IF}_s^{eff}(\theta)$ *exists for* $s < m$ *and* $\Pi_\theta \left[ \mathbb{IF}_m(\theta) | \Gamma_s(\theta) \right] = \mathbb{IF}_s^{eff}(\theta)$.
5. *If the model* $\mathcal{M}(\Theta)$ *is (locally) nonparametric, then*

   (a) *there is at most one* $m^{th}$ *order estimation influence function* $\mathbb{IF}_m(\theta)$ *for* $\psi(\theta)$,
   
   (b)
   $$\mathbb{IF}_m(\theta) = \mathbb{IF}_{m-1}(\theta) + \mathbb{IF}_{mm}(\theta)$$
   
   *where*
   $$\mathbb{IF}_{m-1}(\theta) = \Pi_{m,\theta} \left[ \mathbb{IF}_m(\theta) | \mathcal{U}_{m-1}(\theta) \right]$$
   
   *and* $\mathbb{IF}_{mm}(\theta)$ *is a degenerate* $m^{th}$ *order U-statistic and thus*
   
   $$E_\theta \left[ \mathbb{IF}_{m-1}(\theta) \mathbb{IF}_{mm}(\theta) \right] = 0.$$
   
   (c) (i) *Suppose, for a given* $m \geq 2$, $\mathbb{IF}_{m-1}(\theta)$ *exists and a kernel* $if_{m-1,m-1}(o_{i_1}, \ldots, o_{i_{m-1}}; \theta)$ *of* $\mathbb{IF}_{m-1,m-1}(\theta)$ *has a first order influence function with kernel* $if_{1,if_{m-1,m-1}(o_{i_1},\ldots,o_{i_{m-1}};\cdot)}(O_{i_m}; \theta)$ *for all* $o_{i_1}, \ldots, o_{i_{m-1}}$ *in a set* $\mathcal{O}_{m-1}$ *which has probability 1 under* $F^{(m-1)}(\cdot, \theta)$. *Then* $\mathbb{IF}_m(\theta)$ *exists and*

(2.2) $$m\mathbb{IF}_{m,m}(\theta) = \mathbb{V}\left( d_{m,\theta} \left[ if_{1,if_{m-1,m-1}(O_{i_1},\ldots,O_{i_{m-1}};\cdot)}(O_{i_m}; \theta) \right] \right)$$

*where the operator* $d_{m,\theta}$ *is given in Equation (2.1).*

(ii) *Conversely, if* $\mathbb{IF}_m$ *exists then the symmetric kernel* $if_{m-1,m-1}^{sym}(o_{i_1}, \ldots, o_{i_{m-1}}; \theta)$ *of* $\mathbb{IF}_{m-1,m-1}(\theta)$ *has a first order influence function for all* $o_{i_1}, \ldots, o_{i_{m-1}}$ *in a set* $\mathcal{O}_{m-1}$ *which has probability 1 under* $F^{(m-1)}(\cdot, \theta)$. *Further*

$$m^{-1} d_{m,\theta} \left[ if_{1,if_{m-1,m-1}^{sym}(O_{i_1},\ldots,O_{i_{m-1}};\cdot)}(O_{i_m}; \theta) \right] = if_{m,m}^{sym}(O_{i_1}, \ldots, O_{i_m}; \theta).$$

**Remark 2.4.** Pfanzagl [11] previously proved part 5.c(i) for $m = 2$. Our theorem offers a generalization of his result. Note, in part (i) of 5(c), we can always take the kernel to be the symmetric kernel.

**Remark 2.5.** Provided one knows how to calculate first order influence functions, one can obtain $\mathbb{IF}_2(\theta), \ldots, \mathbb{IF}_m(\theta)$ recursively using part (5.c). An example of such a calculation is given in Section 3.2.2 below. Thus part (5.c) has the interesting implication that even though higher order influence functions are defined in terms of their inner products with higher order scores $\widetilde{\mathbb{S}}_{m,\bar{l}_m}$, nevertheless, in (locally) nonparametric models, one can derive all the higher order influence functions of a functional $\psi(\theta)$ without even knowing how to compute the scores $\widetilde{\mathbb{S}}_{m,\bar{l}_m}$ for any $m > 1$. In fact, one need not even be aware of the structure of the scores $\widetilde{\mathbb{S}}_{m,\bar{l}_m}$ in terms of the subject-specific higher order scores $S_{l_1\ldots l_s,j}(\theta)$. In contrast, in parametric or semiparametric models whose tangent space $\Gamma_m(\theta)$ does not equal the set $\mathcal{U}_m(\theta)$ of all $m^{\text{th}}$ order $U-statistics$, one can often (but not always) still obtain an inefficient influence $\mathbb{IF}_m(\theta)$ by applying part (5.c) of the Theorem. However, calculation of the efficient influence function $\mathbb{IF}_m^{eff}(\theta) = \Pi_\theta [\mathbb{IF}_m(\theta) | \Gamma_m(\theta)]$ by projection generally requires explicit knowledge of the scores $\widetilde{\mathbb{S}}_{m,\bar{l}_m}$ to derive $\Gamma_m(\theta)$. For this reason, it can be considerably more difficult to analyze certain parametric models (with dimension increasing with sample size) than to analyze (locally) nonparametric models. We will consider derivation of and projections onto $\Gamma_m(\theta)$ in a forthcoming paper. In the current paper, however, we do calculate $IF_2^{eff}(\theta)$ in one model that is not (locally) nonparametric so as to provide some sense of the issues that arise. Specifically in Section 4, we calculate $IF_2^{eff}(\theta)$ for a truncated version of the functional $E\left[\{E[Y|X]\}^2\right]$ in a model that assumes the marginal distribution of $X$ is known.

**Remark 2.6** (Implications of Theorem 2.3 for the variance of unbiased estimators). Suppose we have $n$ iid draws $\mathbf{O} = (O_1, \ldots, O_n)$ from $F(o;\theta), \theta \in \Theta$, and a U-statistic $\widehat{\psi}_m$ of order $m \leq n$ with $\text{var}_\theta\left[\widehat{\psi}_m\right] < \infty$ for $\theta \in \Theta$ satisfying $E_\theta\left[\widehat{\psi}_m\right] = \psi(\theta)$ for all $\theta \in \Theta$. That is, $\widehat{\psi}_m$ is unbiased for $\psi(\theta)$. We will use Theorem (2.3) to generalize a number of well-known results on minimum variance unbiased estimation to arbitrary models.

By $E_\theta\left[\widehat{\psi}_m\right] = \psi(\theta)$, we immediately conclude that, viewing $\widehat{\psi}_m$ as a $k^{\text{th}}$ order U-statistic, $\widehat{\psi}_m - \psi(\theta)$ is a $k^{th}$ order estimation influence function for $\psi(\theta)$ for $n \geq k \geq m$. By Theorem 2.3, $\text{var}_\theta\left[\widehat{\psi}_m\right] \geq \text{var}_\theta\left[\mathbb{IF}_m^{eff}(\theta)\right]$. We refer to $\text{var}_\theta\left[\mathbb{IF}_m^{eff}(\theta)\right]$ as the $m^{\text{th}}$ order Bhattacharyya variance bound at $\theta$ for the parameter $\psi(\theta)$ in model $\mathcal{M}(\Theta)$, as this definition, in a precise analogy to Bickel et al. [3]'s generalization of the Cramer-Rao variance bound, generalizes Bhattacharyya's [2] variance bound to arbitrary semi- and non- parametric models. Indeed our first order Bhattacharyya bound is precisely Bickel et al.'s [3] generalization of the Cramer-Rao variance bound.

We shall refer to an $m^{\text{th}}$ order U-statistic estimator $\widehat{\psi}_m$ as $m^{\text{th}}$ order "unbiased locally efficient" at $\theta^*$ for $\psi(\theta)$ in model $\mathcal{M}(\Theta)$ if it is unbiased for $\psi(\theta)$ under the model with variance at $\theta^*$ equal to the $m^{\text{th}}$ order Bhattacharyya bound at $\theta^*$. If $\widehat{\psi}_m$ is "unbiased locally efficient" at $\theta^*$ for all $\theta^* \in \Theta$, we say it is 'unbiased globally efficient'. By Theorem 2.3, $\text{var}_\theta\left[\mathbb{IF}_k^{eff}(\theta)\right] \geq \text{var}_\theta\left[\mathbb{IF}_m^{eff}(\theta)\right]$ for $n \geq k > m$. As a consequence if an $m^{\text{th}}$ order 'unbiased locally efficient' estimator $\widehat{\psi}_{m,eff}$ exists at $\theta^*$ then, for $n \geq k \geq m$, $\mathbb{IF}_k^{eff}(\theta^*) = \mathbb{IF}_m^{eff}(\theta^*)$ so the $m^{\text{th}}$ and $k^{\text{th}}$ order Bhattacharyya bounds are equal at $\theta^*$ and $\widehat{\psi}_{m,eff}$ is also $k^{\text{th}}$ order 'unbiased locally efficient' at $\theta^*$.

From the fact that for an unbiased estimator $\widehat{\psi}_m$, $\widehat{\psi}_m - \psi(\theta)$ is an $m^{\text{th}}$ order influence function, we conclude that the variance of $\widehat{\psi}_m$ attains the bound $\text{var}_{\theta^*}\left[\mathbb{IF}_m^{eff}(\theta^*)\right]$ at $\theta^*$ if and only if $\widehat{\psi}_m - \psi(\theta^*) = \mathbb{IF}_m^{eff}(\theta^*)$. It follows that $\widehat{\psi}_m$ is 'unbiased globally efficient' if and only if $\widehat{\psi}_m - \psi(\theta) = \mathbb{IF}_m^{eff}(\theta)$ for all $\theta \in \Theta$. We thus have proved the following theorem in the $\Rightarrow$ direction. The $\Leftarrow$ direction is immediate.

**Theorem 2.7.** *In a model $\mathcal{M}(\Theta)$, there exists an $m^{th}$ order unbiased globally efficient U-statistic estimator of $\psi(\theta)$, if and only if, for all $\theta \in \Theta$, $\mathbb{IF}_m^{eff}(\theta) + \psi(\theta)$ is a function $\widehat{\psi}_{m,eff}$ of the data $\mathbf{O}$, not depending on $\theta$. In that case, $\widehat{\psi}_{m,eff}$ is the unique unbiased globally efficient estimator.*

In a locally nonparametric model all unbiased $m^{\text{th}}$ order estimators are unbiased globally efficient, as there is a unique $m^{\text{th}}$ order influence function. For example, the usual unbiased estimator $\widehat{\sigma}^2 = \sum_{i=1}^n \left\{X_i - \sum_{j=1}^n X_j/n\right\}^2 / (n-1)$ of the variance of a random variable $X$ is a second order U-statistic and thus is a $k^{\text{th}}$ order unbiased globally efficient U-statistic for $k \geq 2$ in the locally nonparametric model consisting of all distributions under which $\widehat{\sigma}^2$ has a finite variance.

In Section 4 we use the results from this remark to compare the relative efficiencies of competing rate-optimal unbiased estimators in a model which is not locally nonparametric.

We now describe the main heuristic idea behind using higher order influence functions. Technical details are suppressed. Consider the estimator

$$(2.3) \qquad \widehat{\psi}_m = \psi\left(\widehat{\theta}\right) + \mathbb{IF}_m^{eff}\left(\widehat{\theta}\right)$$

based on a sample size $n$, where $\widehat{\theta}$ is an initial rate optimal estimator of $\theta$ from a separate independent training sample. That is we assume that our actual sample size is $N$ and we randomly split the $N$ observations into two samples: an analysis sample of size $n$ and a training sample of size $N - n$ where $(N - n)/N = c^*$, $1 > c^* > 0$. We obtain our initial estimate $\widehat{\theta}$ from the training sample data. Sample splitting has no effect on optimal rates of convergence, although in the form described here does affect 'constants'. Throughout the paper, we derive the properties of our estimators conditional on the data in the training sample. In a later section, we describe how one can sometimes obtain an optimal constant by choosing $(N - n)/N = N^{-\epsilon}, \epsilon > 0$ rather than $c^*$.

**Remark 2.8.** Note that sample splitting is avoided in most statistical applications by using modern "empirical process theory" to prove that 'plug-in' estimators such as $\widehat{\psi}_m = \left\{\psi(\theta) + \mathbb{IF}_m^{eff}(\theta)\right\}_{\theta=\widehat{\theta}}$ that estimate $\theta$ from the same sample used to calculate $\mathbb{IF}_m^{eff}(\cdot)$ have nice statistical properties. However empirical process theory is not applicable in our setting because we are interested in function classes whose size (entropy) is so large that they fail to be Donsker. For this reason we initially believed that explicit sample splitting would be difficult to avoid in our methodology. However, in Robins et al. [18], we describe a new method that effectively allows one to use all the data for estimator construction.

Expanding and evaluating conditionally on the training sample (or equivalently on $\widehat{\theta}$), we find by Theorem 2.2 that the conditional bias

$$E_\theta\left[\widehat{\psi}_m - \psi(\theta)|\widehat{\theta}\right] = \psi\left(\widehat{\theta}\right) - \psi(\theta) + E_\theta\left[\mathbb{IF}_m^{eff}\left(\widehat{\theta}\right)|\widehat{\theta}\right]$$

is $O_p\left(||\widehat{\theta} - \theta||^{m+1}\right)$ which decreases with $m$ provided $||\widehat{\theta} - \theta|| < 1$.

In Theorem 3.22 below, we show that if

$$sup_{o \in \mathcal{O}}\left|f\left(o;\widehat{\theta}\right) - f(o;\theta)\right| \to 0$$

as $||\widehat{\theta} - \theta|| \to 0$, where $f(o;\theta)$ is the density of $O$ under $\theta$ and $\mathcal{O}$ has probability one under all $\theta \in \Theta$, then

$$\text{var}_\theta\left[\widehat{\psi}_m|\widehat{\theta}\right] \equiv \text{var}_\theta\left[\mathbb{IF}_m^{eff}\left(\widehat{\theta}\right)|\widehat{\theta}\right] = \text{var}_{\widehat{\theta}}\left[\mathbb{IF}_m^{eff}\left(\widehat{\theta}\right)\right]\left(1 + O_p(||\widehat{\theta} - \theta||)\right)$$

Now, by Theorem 2.3, $\text{var}_{\widehat{\theta}}\left[\mathbb{IF}_m^{eff}\left(\widehat{\theta}\right)\right]$ increases with $m$. Further, $\text{var}_{\widehat{\theta}}\left[\mathbb{IF}_1^{eff}\left(\widehat{\theta}\right)\right] \asymp 1/n$, since, conditional on $\widehat{\theta}$, $\mathbb{IF}_1^{eff}\left(\widehat{\theta}\right)$ is the sample average of *iid* random variables.

To proceed further we shall need to be more explicit about the model $\mathcal{M}(\Theta)$. For now, we consider finite-dimensional parametric models whose dimension $k(n)$ increases with sample size. That is $\theta \equiv \theta_n$ depends on $n$ and the dimension of $\Theta \equiv \Theta_n$ is $k(n)$. Suppose $k(n) \asymp n^\gamma, \gamma \geq 0$. Let $\widehat{\theta}_n$ be the maximum likelihood estimator of $\theta$. If $k(n)$ increases slower than the sample size (i.e., $\gamma < 1$), then, a) under regularity conditions, $||\widehat{\theta}_n - \theta_n|| = O_p\left(\{k(n)/n\}^{1/2}\right) = O_p\left(n^{-\frac{1}{2}(1-\gamma)}\right)$ with $||\cdot||$ the usual Euclidean norm in $R^{k(n)}$; and b) $\text{var}_{\widehat{\theta}}\left[\mathbb{IF}_m^{eff}\left(\widehat{\theta}\right)\right]$, although increasing with $m$, remains order $1/n$; as a consequence, if $m$ is chosen greater than the solution $m^*$ to $n^{-\frac{m^*+1}{2}(1-\gamma)} = n^{-1/2}$, the bias of $\widehat{\psi}_m$ will be $o_p\left(n^{-1/2}\right)$, the rate of convergence will be the usual parametric rate of $n^{-1/2}$, and thus, for $n$ sufficiently large, the squared bias of $\widehat{\psi}_m$ will be less than the variance. As a consequence, as discussed in Section 3.2.5, we can construct honest (i.e uniform over $\theta_n \in \Theta_n$) asymptotic confidence intervals centered at $\widehat{\psi}_{m^*}$ with width of order $n^{-1/2}$. Here is a concrete example.

**Example.** Suppose $O = (Y, X)$ with $Y$ Bernoulli and with $X$ having a density with respect to the uniform measure $\mu(\cdot)$ on the unit cube $[0,1]^d$ in $R^d$. Suppose $\psi = E\left[(E[Y|X])^2\right]$. Let $\{z_l(\cdot)\} \equiv \{z_l(x); 1, 2, \ldots\}$ be a countable, linearly independent, sequence of either spline, polynomial, or compact wavelet basis functions dense in $L_2(\mu)$. Set $\overline{z}_k(x) = (z_1(x), \ldots, z_k(x))^T$. We assume

$$E(Y|X=x) \in \left\{\begin{array}{c} b\left(x; \overline{\eta}_{k^*(n)}\right) \equiv \left[1 + \exp\left(-\overline{\eta}_{k^*(n)}^T \overline{z}_{k^*(n)}(x)\right)\right]^{-1}; \\ \overline{\eta}_{k^*(n)} \in \mathcal{N}_{k^*(n)} \end{array}\right\},$$

$$f_X(x) \in \left\{\begin{array}{c} f\left(x; \overline{\omega}_{k^{**}(n)}\right) \equiv c\left(\overline{\omega}_{k^{**}(n)}\right) \exp\left[\overline{\omega}_{k^{**}(n)}^T \overline{z}_{k^{**}(n)}(x)\right]; \\ \overline{\omega}_{k^{**}(n)} \in \mathcal{W}_{k^{**}(n)} \end{array}\right\},$$

where $c\left(\overline{\omega}_{k^{**}(n)}\right)$ is a normalizing constant and $\mathcal{N}_{k^*(n)}$ and $\mathcal{W}_{k^{**}(n)}$ are open bounded subsets of $R^{k^*(n)}$ and $R^{k^{**}(n)}$. Hence, $\Theta_n = \mathcal{N}_{k(n)} \times \mathcal{W}_{k(n)}$ has dimension $k(n) = k^*(n) + k^{**}(n)$ and $\psi(\theta) = \psi_n(\theta_n) = \int b^2\left(x; \overline{\eta}_{k^*(n)}\right) f\left(x; \overline{\omega}_{k^{**}(n)}\right) d\mu(x)$.

He [7] and Portnoy [12] prove that, under regularity conditions, $||\widehat{\theta}_n - \theta_n|| = O_p\left(\{k(n)/n\}^{1/2}\right)$ when $k(n) = n^\gamma \ll n$. Below we shall see that $\text{var}_{\widehat{\theta}}\left[\mathbb{IF}_m^{eff}\left(\widehat{\theta}\right)|\widehat{\theta}\right] \asymp 1/n$ for $n^\gamma \ll n$.

Consider next models whose dimension $k(n) \asymp n^\gamma$ increases faster than $n$ (i.e., $\gamma > 1$). In such models, the MLE $\widehat{\theta}_n$ is generally inconsistent and indeed there may exist no consistent estimator of $\theta_n$. In that case, $||\widehat{\theta}_n - \theta_n||$ fails to be $o_p(1)$ and the conditional bias $E_\theta\left[\widehat{\psi}_m - \psi(\theta)|\widehat{\theta}\right]$ may not decrease with $m$. In order to guarantee consistent estimators of $\theta_n$ exist, it is necessary to place further a priori restrictions on the complexity of $\Theta_n$. Typical examples of complexity-reducing assumptions would be an $\epsilon$-sparseness assumption that only $k(n)^\epsilon, 0 < \epsilon < 1$, of the $k(n)$ parameters are non-zero or a smoothness assumption that specifies that the rate of decrease of the $j^{th}$ component of $\theta_n$ is equal to $1/j$ raised to a given (positive) power. Even after imposing such complexity-reducing assumptions, $\psi(\theta) \equiv \psi_n(\theta_n)$ may not be estimable at rate $n^{-1/2}$.

For instance consider the previous example but now with $\gamma^*$ and $\gamma^{**}$ exceeding 1, so $k^{**}(n) = n^{\gamma^{**}} \gg n$, $k^*(n) = n^{\gamma^*} \gg n$ and $k(n) = k^{**}(n) + k^*(n) \asymp n^\gamma \gg n$ with $\gamma = \max(\gamma^{**}, \gamma^{**})$. Consider the norms

$$\left\|\overline{\eta}_{k^*(n)}\right\| = \left\{\int b^2\left(x; \overline{\eta}_{k^*(n)}\right) d\mu(x)\right\}^{1/2},$$

$$\left\|\overline{\omega}_{k^{**}(n)}\right\|_p = \left\{\int f\left(x; \overline{\omega}_{k^{**}(n)}\right)^p d\mu(x)\right\}^{1/p} \text{ and}$$

$$\|\theta\|_p = \left\|\overline{\eta}_{k^*(n)}\right\| + \left\|\overline{\omega}_{k^{**}(n)}\right\|_p.$$

Suppose, under a particular smoothness assumption, optimal rate estimators $\widehat{\overline{\eta}}_{k^*(n)}$ and $\widehat{\overline{\omega}}_{k^{**}(n)}$ of $\overline{\eta}_{k^*(n)}$ and $\overline{\omega}_{k^{**}(n)}$ satisfy $\left\|\widehat{\overline{\eta}}_{k^*(n)} - \overline{\eta}_{k^*(n)}\right\| = O_p(n^{-\gamma_\eta})$ and $\left\|\widehat{\overline{\omega}}_{k^{**}(n)} - \overline{\omega}_{k^{**}(n)}\right\|_p = O_p(n^{-\gamma_\omega})$ for some $\gamma_\eta > 0, \gamma_\omega > 0$ and all $p \geq 2$. Hence, $||\widehat{\theta} - \theta||_p = O_p(\max\{n^{-\gamma_\eta}, n^{-\gamma_\omega}\})$. For $\gamma > 1$, based on arguments given later, we expect that $\text{var}_{\widehat{\theta}}\left[\widehat{\psi}_m - \psi(\theta)|\widehat{\theta}\right] \asymp \frac{n^{(\gamma-1)(m-1)}}{n}$ and

$$E_\theta\left[\widehat{\psi}_m - \psi(\theta)|\widehat{\theta}\right] = O_p\left(\left\|\widehat{\overline{\eta}}_{k^*(n)} - \overline{\eta}_{k^*(n)}\right\|^2 \left\|\widehat{\overline{\omega}}_{k^{**}(n)} - \overline{\omega}_{k^{**}(n)}\right\|_{m-1}^{m-1}\right)$$
$$= O_p\left(n^{-2\gamma_\eta - (m-1)\gamma_\omega}\right)$$
$$= O_p\left(||\widehat{\theta} - \theta||_{m-1}^{m+1}\right).$$

To find the estimator $\widehat{\psi}_{m_{best}}$ in the class $\widehat{\psi}_m$ with optimal rate of convergence, let $m^* = 1 + \frac{1 - 4\gamma_\eta}{(\gamma-1) + 2\gamma_\omega}$ be the value of $m$ that equates the order $n^{-4\gamma_\eta - 2(m-1)\gamma_\omega}$ of the squared bias and the order $\frac{n^{(\gamma-1)(m-1)}}{n}$ of the variance. Then $m_{best} = \lfloor m^* \rfloor$ if the order $n^{-4\gamma_\eta - 2(m-1)\gamma_\omega} + n^{(\gamma-1)(m-1)-1}$ of the mean squared error at $\lfloor m^* \rfloor$ is less than or equal to that at $\lceil m^* \rceil$. Otherwise, $m_{best} = \lceil m^* \rceil$. The rate of convergence of $\widehat{\psi}_{m_{best}}$ will often be slower than $n^{-1/2}$. Note $m_{best} = 1$ whenever $\gamma > 2$, regardless of $\gamma_\eta$ and $\gamma_\omega$.

By using the estimator $\widehat{\psi}_{\lceil m^* \rceil}$ rather than $\widehat{\psi}_{m_{best}}$, we can guarantee that the variance asymptotically dominates bias and construct honest (i.e uniform over $\theta_n \in \Theta_n$) asymptotic confidence intervals centered at $\widehat{\psi}_{\lceil m^* \rceil}$. Of course, the sample size $n$ at which, for all $\theta_n \in \Theta_n$, the finite sample coverage of the intervals discussed above is close to the asymptotic (i.e. nominal) coverage is generally unknown and could be very large. For this reason, a better, but unfortunately as yet technically out of reach, approach to confidence interval construction is discussed in Section 3.2.5.

In contrast to the case of parametric models of increasing dimension, in the infinite dimensional models which we consider in the following section, the functionals $\psi(\theta)$ of interest have first order influence functions $\mathbb{IF}_1(\theta)$ but do not have higher order influence functions. As a consequence, an initial 'truncation' step is needed before we can apply the approach outlined in the preceding paragraph.

Finally, even in the case of parametric models with $k(n) \gg n$ and complexity reducing assumptions imposed, , when the minimax rate for estimation of $\psi(\theta)$ is slower than $n^{-1/2}$, the optimal estimator $\widehat{\psi}_{m_{best}}$ in the class $\widehat{\psi}_m$ will generally not be rate minimax. See Section 3.2.6 and Sections 4.1.1 for additional discussion.

## 3. Inference for a class of doubly robust functionals

### 3.1. The class of functionals

In this Section we consider models in which the parameter $\theta$ is infinite dimensional and $\theta$, $\Theta$, and $\psi(\cdot)$ do not depend on $n$. We make the following three assumptions (Ai)–(Aiii):

(Ai) The data $O$ includes a vector $X$, where, for all $\theta \in \Theta$, the distribution of $X$ is supported on the unit cube $[0,1]^d$ ( or more generally a compact set) in $R^d$ and has a density $f(x)$ with respect to the Lebesgue measure. Further $\Theta = \Theta_1 \times \Theta_2$ where $\theta_1 \in \Theta_1$ governs the marginal law of $X$ and $\theta_2 \in \Theta_2$ governs the conditional distribution of $O|X$.

(Aii) The parameter $\theta_2$ contains components $b = b(\cdot)$ and $p = p(\cdot)$, $b : [0,1]^d \to \mathcal{R}$ and $p : [0,1]^d \to \mathcal{R}$, such that the functional $\psi(\theta)$ of interest has a first order influence function $\mathbb{IF}_{1,\psi}(\theta) = \mathbb{V}[IF_{1,\psi}(\theta)]$, where

$$IF_{1,\psi}(\theta) = H(b,p) - \psi(\theta), \tag{3.1}$$

with $H(b,p) \equiv h(O, b(X), p(X))$

$$\equiv b(X) p(X) h_1(O) + b(X) h_2(O) + p(X) h_3(O) + h_4(O) \tag{3.2}$$

$$\equiv BPH_1 + BH_2 + PH_3 + H_4,$$

and the known functions $h_1(\cdot), h_2(\cdot), h_3(\cdot), h_4(\cdot)$ do not depend on $\theta$.

(Aiii)  (a) $\Theta_{2b} \times \Theta_{2p} \subseteq \Theta_2$ where $\Theta_{2b}$ and $\Theta_{2p}$ are the parameter spaces for the functions $b$ and $p$. Furthermore the sets $\Theta_{2b}$ and $\Theta_{2p}$ are dense in $L_2(F_X(x))$ at each $\theta_1^* \in \Theta_1$.

or

(b) $b^*(\cdot) = p^*(\cdot)$, $h_3(O) = h_2(O)$ w.p.1, and $\Theta_{2b} \subseteq \Theta_2$ is dense in $L_2(F_X(x))$ at each $\theta_1^* \in \Theta_1$.

**Remark 3.1.** (Aiiib) can be viewed as a special case of (Aiiia) as discussed in Example 1a below, so we need only prove results under assumption (Aiiia).

Assumptions (Ai)–(Aiii) have a number of important implications that we summarize in a Theorem and two Lemmas.

**Theorem 3.2** (Double-robustness). *Assume (Ai)–(Aiii) hold, and recall p and b are elements of $\theta$. Then*

$$E_\theta\left[H\left(b\ ,p^*\right)\right] = E_\theta\left[H\left(b^*,p\right)\right] = E_\theta\left[H\left(b,p\right)\right] = \psi\left(\theta\right)$$

*for all $(p^*, b^*) \in \Theta_{2p} \times \Theta_{2b}, \theta \in \Theta$.*

*Proof.* $E_\theta\left[H\left(b^*,p\right)\right] - E_\theta\left[H\left(b,p\right)\right] = E_\theta\left[\{H_1 p\left(X\right) + H_2\}\{b\left(X\right) - b^*\left(X\right)\}\right]$ and $E_\theta\left[H\left(b\ ,p^*\right)\right] - E_\theta\left[H\left(b,p\right)\right] = E_\theta\left[\{H_1 b\left(X\right) + H_3\}\{p\left(X\right) - p^*\left(X\right)\}\right]$. The theorem then follows from part 1) of the following lemma. □

Theorem 3.2 states that $H\left(\cdot,\cdot\right)$ has mean $\psi\left(\theta\right)$ under $F\left(\cdot;\theta\right)$ even when $p$ is misspecified as $p^*$ or $b$ is misspecified as $b^*$. We refer to the functional $\psi\left(\theta\right)$ as doubly robust because of this property. The next lemma shows that $H\left(b^*, p^*\right)$ is not unbiased if both $b$ and $p$ are simultaneously misspecified. That is, $E_\theta\left[H\left(b^*, p^*\right)\right] \neq \psi\left(\theta\right)$.

**Lemma 3.3.** *Assume (Ai)–(Aiii) hold. Then*

1. $E_\theta\left[\{H_1 B + H_3\}|X\right] = E_\theta\left[\{H_1 P + H_2\}|X\right] = 0$
2. $E_\theta\left[H\left(b^*, p^*\right)\right] - E_\theta\left[H\left(b,p\right)\right] = E_\theta\left[\left(B - B^*\right)\left(P - P^*\right) H_1\right]$
   and $\psi\left(\theta\right) \equiv E_\theta\left[H\left(b,p\right)\right] = E_\theta\left[-BPH_1 + H_4\right]$

*Proof.* Part (1): By assumptions (Ai) and (Aiiia) we have paths $\widetilde{\theta}_l\left(t\right), l = 1, 2, \ldots,$ in our model with $\widetilde{\theta}_l\left(0\right) = \theta$ and $p_l\left(t\right) = p_l\left(x;t\right) = p\left(x\right) + tc_l\left(x\right), b_l\left(x;t\right) = b\left(x\right), F_l\left(x;t\right) = F\left(x\right)$ for $l = 1, 2, \ldots,$ where the sequence $c_l\left(\cdot\right)$ is dense in $L_2\left[F_X\left(x\right)\right]$. Let $S_l\left(\theta\right)$ be the score for path $\widetilde{\theta}_l\left(t\right)$ at $t = 0$. Then by $\psi\left(\widetilde{\theta}_l\left(t\right)\right) = E_{\widetilde{\theta}_l(t)}\left[H\left(b, p_l\left(t\right)\right)\right]$

$$d\psi\left(\widetilde{\theta}_l\left(t\right)\right)/dt_{|t=0} = E_\theta\left[\{H_1 B + H_3\} c_l\left(X\right)\right]$$
$$+ E_\theta\left[H\left(b,p\right) S_l\left(\theta\right)\right].$$

By $\mathbb{IF}_{1,\psi}\left(\theta\right) = H\left(b,p\right) - \psi\left(\theta\right),$

$$d\psi\left(\widetilde{\theta}_l\left(t\right)\right)/dt_{|t=0} = E_\theta\left[H\left(b,p\right) S_l\right].$$

Thus $E\left[\{H_1 B + H_3\} c_l\left(X\right)\right] = 0$. But $\{c_l\left(\cdot\right)\}$ is dense in $L_2\left[F_0\left(X\right)\right]$ so

$$E\left[H_1 B + H_3 | X\right] = 0.$$

An analogous argument proves $E_\theta\left[\{H_1 P + H_2\}|X\right] = 0$. Part (2): $E_\theta\left[H\left(b^*, p^*\right)\right] - E_\theta\left[H\left(b,p\right)\right] =$

$$E_\theta\left[\left(B^* P^* - BP\right) H_1 + \left(B^* - B\right) H_2 + \left(P^* - P\right) H_3\right]$$
$$= E_\theta\left[\left(B^* P^* - BP\right) H_1 - \left(B^* - B\right) PH_1 - \left(P^* - P\right) BH_1\right]$$
$$= E_\theta\left[\left(B - B^*\right)\left(P - P^*\right) H_1\right],$$

where the second equality is by part 1). Choosing $P^* = B^* = 0$ w.p.1 completes the proof of the theorem since then $E_\theta\left[H\left(b^*, p^*\right)\right] = E_\theta\left[H_4\right]$. □

Below we will need the following partial converse to Lemma 3.3.

**Lemma 3.4.** *Let $\Theta_{2b}, \Theta_{2p}, \Theta_1$ and $\Theta$ and $H(b,p)$ be as defined in (Ai)–(Aiiia). Suppose that*

$$E_\theta\left[\{H_1 B + H_3\}|X\right] = E_\theta\left[\{H_1 P + H_2\}|X\right] = 0$$

*and $\psi(\theta) = E_\theta[H(b,p)]$. Then $\mathbb{V}[H(b,p) - \psi(\theta)]$ is the first order influence function of $\psi(\theta)$.*

*Proof.* The influence function of the functional $E_\theta[H(b^*, p^*)]$ for known functions $b^*, p^*$ is $\mathbb{V}[H(b^*, p^*) - E_\theta[H(b^*, p^*)]]$. Thus by the linearity of first order influence functions, the Lemma is true if and only if for each $\theta_0 \in \Theta$, the functional $\tau(b,p) = E_{\theta_0}[H(b,p)]$ with $\theta_0$ fixed has influence function equal to 0 w.p.1 at $(b,p) = (b_0, p_0) \subset \theta_0$. That the influence function is equal to 0 follows from the fact that, under the assumptions of the Lemma, for sets $\{c_l(\cdot)\}$ and $\{d_l(\cdot)\}$ dense in $L_2[F_0(X)]$,

$$dE_{\theta_0}[H(b_0(X) + tc_l(X), p_0(X) + td_l(X))]/dt_{|t=0}$$
$$= E_\theta[\{H_1 b_0(X) + H_3\} d_l(X)] + E_\theta[\{H_1 p_0(X) + H_2\} c_l(X)] = 0. \quad \square$$

Results of Ritov and Bickel [14] and Robins and Ritov [15] imply it is not possible to construct honest asymptotic confidence intervals for $\psi(\theta)$ whose width shrinks to 0 as $n \to \infty$ if $b(\cdot)$ and $p(\cdot)$ are too rough. Therefore we also place a priori bounds on their roughness. Our bounds will be based on the following definition.

**Definition 3.5.** A function $h(\cdot)$ with domain $[0,1]^d$ is said to belong to a Hölder ball $H(\beta, C)$, with Hölder exponent $\beta > 0$ and radius $C > 0$, if and only if $h(\cdot)$ is uniformly bounded by $C$, all partial derivatives of $h(\cdot)$ up to order $\lfloor \beta \rfloor$ exist and are bounded, and all partial derivatives $\nabla^{\lfloor \beta \rfloor}$ of order $\lfloor \beta \rfloor$ satisfy

$$\sup_{x, x+\delta x \in [0,1]^d} \left|\nabla^{\lfloor \beta \rfloor} h(x+\delta x) - \nabla^{\lfloor \beta \rfloor} h(x)\right| \leq C||\delta x||^{\beta - \lfloor \beta \rfloor}.$$

We note that the $L_p, 2 < p < \infty$ and $L_\infty$ rates of convergence for estimation of a marginal density or conditional expectation $h(\cdot) \in H(\beta, C)$ are $O\left(n^{-\frac{\beta}{2\beta+d}}\right)$ and $O\left(\left(\frac{n}{\log n}\right)^{-\frac{\beta}{2\beta+d}}\right)$ respectively. We refer to an estimator attaining these rates as rate optimal.

We make the following fourth assumption:

(Aiv) We assume $b(\cdot), p(\cdot),$ and $g(\cdot)$ lie in given Hölder balls $H(\beta_b, C_b), H(\beta_p, C_p), H(\beta_g, C_g)$ where

(3.3) $$g(x) \equiv E\{H_1|X=x\} f(x).$$

Furthermore we assume $g(X) > \sigma_g > 0$ w.p.1. Finally we assume, as can always be arranged by a suitable choice of estimator, that the initial training sample estimators $\widehat{b}(.), \widehat{p}(.),$ and $\widehat{g}(\cdot)$ are rate optimal, have more than $\max\{\beta_b, \beta_g, \beta_p\}$ derivatives, and have $L_\infty$ norm bounded by a constant $c_\infty$. Further $\inf_{x \in [0,1]^d} \widehat{g}(x) > \sigma_g$. The reason for the restrictions on $g(\cdot)$ will become clear below.

The restrictions (Ai)–(Aiv) are the only restrictions common to all functionals and models in the class. Additional model and/or functional specific restrictions will be given below.

To motivate our interest in such a class of functionals and models we provide a number of examples. In each case, one can use Lemma 3.4 to verify that the influence function of $\psi(\theta)$ is as given. All but Examples 3 and 4 are examples of (locally) nonparametric models.

**Example 1.** Suppose $O=(A,Y,X)$ with $A$ and $Y$ univariate random variables.

**Example 1a** (Expected product of conditional expectations). Let $\psi(\theta) = E_\theta[p(X)b(X)]$ where $b(X) = E_\theta[Y|X], p(X) = E_\theta[A|X]$. In this model

$$IF_{1,\psi}(\theta) = p(X)b(X) - \psi(\theta)$$
$$+ p(X)\{Y - b(X)\} + b(X)\{A - p(X)\}$$

so $H_1 = -1, H_2 = A, H_3 = Y, H_4 = 0$.

We also consider the special case of this model where $A = Y$ with probability one (w.p.1). Then, as in assumption (Aiiib), $b(X) = p(X)$ w.p.1, $H_2 = H_3$ w.p.1. Then $\psi(\theta) = E_\theta[b^2(X)]$. In Section 5, we show how our confidence interval for $E_\theta[b^2(X)]$ can be used to obtain an adaptive confidence interval for the regression function $b(\cdot)$.

**Example 1b** (Expected conditional covariance).

$$\psi(\theta) = E_\theta[AY] - E_\theta[p(X)b(X)] = E_\theta[\text{cov}_\theta\{Y,A|X\}]$$

has influence function

$$AY - \{p(X)b(X) + p(X)\{Y - b(X)\} + b(X)\{A - p(X)\}\} - \psi(\theta),$$

so $H_1 = 1, H_2 = -A, H_3 = -Y, H_4 = AY$.

Example 1c below shows that a confidence interval and point estimators for $E_\theta[\text{cov}_\theta\{Y,A|X\}]$ can be used to obtain confidence intervals and point estimator for the variance weighted average treatment effect in an observational study.

**Example 1c** (Variance-weighted average treatment effect). Suppose, in an observational study, $O = \{Y^*, A, X\}$, $A$ is a binary treatment taking values in $\{0,1\}$, $Y^*$ is a univariate response and $X$ is a vector of pretreatment covariates. Consider the parameter $\tau(\theta)$ given by:

$$(3.4) \qquad \tau(\theta) = \frac{E_\theta[\text{cov}_\theta(Y^*, A|X)]}{E_\theta[\text{var}_\theta(A|X)]} = \frac{E_\theta[\text{cov}_\theta(Y^*, A|X)]}{E_\theta[\pi(X)\{1 - \pi(X)\}]},$$

where $\pi(X) = pr(A = 1|X)$ is often referred to as the propensity score. We are interested in $\tau(\theta)$ for several reasons. First, in the absence of confounding by unmeasured factors, $\tau(\theta)$ is the variance-weighted average treatment effect since $\tau(\theta)$ can be rewritten as $E_\theta[w_\theta(X)\gamma(X;\theta)]$ where $w_\theta(X) = \frac{\text{var}_\theta(A|X)}{E_\theta[\text{var}_\theta(A|X)]}$ and

$$\gamma(x;\theta) = E_\theta(Y^*|A=1, X=x) - E_\theta(Y^*|A=0, X=x)$$

is the average conditional treatment effect at level $x$ of the covariates. Second, under the semiparametric model

$$(3.5) \qquad \gamma(X;\theta) = \upsilon(\theta) \ w.p.1$$

that assumes the treatment effect does not depend on $X$, $\tau(\theta) = \upsilon(\theta)$. In Remark 4.2, we briefly consider inference on $\tau(\theta)$ under model (3.5). However since the model (3.5) may not hold and therefore the parameter $\upsilon(\theta)$ may be undefined, our main goal is to make inference on $\tau(\theta)$ without imposing (3.5).

Now if for any $\tau \in R$, we define $\psi(\tau, \theta)$ to be

$$\psi(\tau, \theta) = E_\theta \left[ \{Y^*(\tau) - E_\theta(Y^*(\tau)|X)\} \{A - E_\theta(A|X)\} \right],$$

with $Y^*(\tau) = Y^* - \tau A$, it is easy to verify that $\tau(\theta)$ may also be characterized as the solution $\tau = \tau(\theta)$ to the equation $\psi(\tau, \theta) = 0$. Thus inference on $\tau(\theta)$ is easily obtained from inference on $\psi(\tau, \theta)$. In particular a $(1 - \alpha)$ confidence set for $\tau(\theta)$ is the set of $\tau$ such that a $(1 - \alpha)$ CI interval for $\psi(\tau, \theta)$ contains 0. Therefore, with no loss of generality, we consider the construction of a $(1-\alpha)$ CI for $\psi(\widetilde{\tau}, \theta)$ for a fixed value $\tau = \widetilde{\tau}$, and write $Y = Y^*(\widetilde{\tau})$ and $\psi(\theta) = \psi(\widetilde{\tau}, \theta)$. Thus $\psi(\theta) = E_\theta[\text{cov}_\theta \{Y, A|X\}]$ and we are in the setting of Example 1b.

In Section 6, we show the rates at which the width of the confidence sets for $\psi(\widetilde{\tau}, \theta)$ and for $\tau(\theta)$ shrink with $n$ are equal.

**Example 2a** (Missing at random). Suppose $O = (AY, A, X)$ where $Y$ is an outcome that is not always observed, $A$ is the binary missingness indicator, $X$ is a $d$-dimensional vector of always observed continuous covariates, and let $b(X) = E(Y|A=1, X)$, $\pi(X) = P(A=1|X)$ be the propensity score, and $p(X) = 1/\pi(X)$. We suppose $\pi(X) > \sigma > 0$ and define

(3.6) $$\psi(\theta) = E_\theta \left[ \frac{AY}{\pi(X)} \right] = E_\theta [b(X)].$$

Interest in $\psi(\theta)$ lies in the fact that $\psi(\theta)$ is the marginal mean of $Y$ under the missing (equivalently, coarsening) at random (MAR) assumption that $P(A=1|X,Y) = \pi(X)$. In this model $IF_{1,\psi}(\theta) = Ap(X)(Y - b(X)) + b(X) - \psi(\theta)$ so $H_1 = -A, H_2 = 1, H_3 = AY, H_4 = 0$.

Note that if one has assumed a priori that $f_X(\cdot)$ and $p(X)$ lay in Hölder balls with respective exponents $\beta_{f_X}$ and $\beta_p$, then $\beta_g$ would be $\min(\beta_{f_X}, \beta_p)$, since $g(X) = -f_X(X)/p(X)$.

**Example 2b** (Missing not-at random). Consider again the setting of Example 2a but we no longer assume MAR. Rather we assume

$$P(A = 1|X, Y) = \{1 + \exp\{-[\gamma(X) + \alpha Y]\}\}^{-1}$$

may depend on $Y$, where now $\gamma(X)$ is an unknown function and $\alpha$ is a known constant (to be later varied in a sensitivity analysis). In this case the marginal mean of $Y$ is given by $\psi(\theta) = E_\theta(AY[1 + \exp\{-[\gamma(X) + \alpha Y]\}])$. Robins and Rotnitzky [17] proved this model places no restrictions on $F(o)$ and derived

$$IF_{1,\psi}(\theta) = A\{1 + \exp\{-\alpha Y\} p(X)\} \{Y - b(X)\} + b(X) - \psi(\theta)$$

where, now,

$$b(X) = E[Y \exp\{-\alpha Y\}|A = 1, X]/E[\exp\{-\alpha Y\}|A = 1, X],$$

and $p(X) = \exp\{-\gamma(X)\}$. Thus

$$H_1 = -\exp\{-\alpha Y\} A, \ H_2 = \{1 - A\}, H_3 = AY \exp\{-\alpha Y\},$$

and $H_4 = AY$. When $\alpha = 0$ this provides an alternate parametrization of Example 2a.

**Example 3** (Marginal structural models and the average treatment effect). Consider the set-up of Example 1c including the non-identifiable assumption of no unmeasured confounders, except now $A$ is discrete with possibly many levels and $f(a|X) > \delta > 0$ w.p.1. A marginal structural model assumes $E_{f_X}\{E_\theta(Y^*| A = a, X)\} = d(a, \upsilon(\theta))$, where $d(a, \upsilon)$ is a known function and $\upsilon(\theta)$ is an unknown vector parameter of dimension $d^*$. When $A$ is dichotomous with $a \in \{0, 1\}$ and $d(a, \upsilon) = \upsilon_1 + \upsilon_2 a$, then $\upsilon_2(\theta)$ is the average treatment effect parameter. Let $f^*(a)$ be any density with the same support as $A$ and let $s^*(a)$ be a $d^*$-vector function, both chosen by the analyst. Then $\upsilon(\theta)$ is identified as the (assumed) unique value of $\upsilon$ satisfying

$$\psi_\upsilon(\theta) \equiv E_\theta\left[s(O, A, \upsilon)\frac{f^*(A)}{f(A|X)}\right] = 0,$$

where $s(O, a, \upsilon) = \{Y^* - d(a, \upsilon)\}s^*(a)$. Thus a $(1-\alpha)$ confidence set for $\upsilon(\theta)$ is the set of vectors $\upsilon$ such that a $(1-\alpha)$ CI for $\psi_\upsilon(\theta)$ contains 0. Therefore, with no loss of generality, we consider the construction of a $(1-\alpha)$ CI for the $d-$vector functional $\psi(\theta) \equiv \psi_{\widetilde{\upsilon}}(\theta)$ for a fixed value $\widetilde{\upsilon}$ and define $h(O, A) \equiv s(O, a, \widetilde{\upsilon})$ and $b(a, X) \equiv E_\theta[h(O, a)|A = a, X]$. Then $\psi_{\widetilde{\upsilon}}(\theta)$ has influence function

$$IF_1(\theta) = \frac{f^*(A)}{f(A|X)}\{h(O, A) - b(A, X)\} + \int b(a, X) \, dF^*(a) - \psi(\theta).$$

Next define $p(a, X) = 1/f(a|X)$, $\psi(\theta, a) = E_{f_X}[b(a, X)]$. Then $IF_1(\theta)$ is the integral

$$IF_1(\theta) = \int dF^*(a) \, IF_1(a, \theta),$$

$$IF_1(a, \theta) = H_1(a) p(a, X) b(a, X)$$
$$+ H_2(a) b(a, X) + H_3(a) p(a, X) - \psi(\theta, a),$$
$$H_1(a) = -I(A = a), H_2(a) = 1, H_3(a) = I(A = a) h(O, a).$$

It follows that $IF_1(\theta)$ is a integral over $a \in A$ of influence functions $IF_1(a, \theta)$ for parameters $\psi(\theta, a)$ in our class with $H_4(a) = 0$. Thus we can estimate $\psi(\theta)$ by $\int dF^*(a) \widehat{\psi}(a)$, where $\widehat{\psi}(a)$ is an estimator of $\psi(\theta, a)$. If the support of $A$ is of greater cardinality than $d^*$, the model is not locally nonparametric. Different choices for $s^*(a)$ and $f^*(a)$ for which $\{\partial/\partial\upsilon^T\} E_\theta\left[s(O, A, \upsilon)\frac{f^*(A)}{f(A|X)}\right]$ is invertible may result in difference influence functions. All yield the same rate of convergence, although the constants differ. See Remark 2.5 above. Extension of our methods to continuous $A$ will be treated elsewhere.

**Example 4** (Confidence intervals for the optimal treatment strategy in a randomized clinical trial). Consider a randomized clinical trial with data $O = \{Y, Y^*, A, X\}$, $A$ is a binary treatment taking values in $\{0, 1\}$, $Y^*$ and $Y$ are univariate responses, $X$ is a vector of pretreatment covariates. In a randomized trial, the randomization probabilities $\pi_0(X) = P(A = 1|X)$ are known by design. Let $b(x) = E_\theta(Y^*|A = 1, X = x) - E_\theta(Y^*|A = 0, X = x)$ and $p(x) = E_\theta(Y|A = 1, X = x) - E_\theta(Y|A = 0, X = x)$ be the average treatment effects at level $X = x$ on $Y^*$ and $Y$. We assume $Y$ and $Y^*$ have been coded so that positive treatment effects are desirable. Let $\psi(\theta) = E[b(X) p(X)]$. Because the model is not locally nonparametric there exists more than a single first order influence function. Indeed, for any given function $c(\cdot)$,

$$IF_{1,\psi}(\theta, c) = b(X) p(X) - \psi(\theta) + [b(X)\{Y - Ap(X)\} + p(X)\{Y^* - Ab(X)\}]$$
$$\times \{A - \pi_0(X)\} \sigma_0^{-2}(X) + c(X)\{A - \pi_0(X)\},$$

with $\sigma_0^2(X) = \pi_0(X)\{1 - \pi_0(X)\}$ is an influence function in our class [provided it is square integrable] with

$$H_1 = 1 - 2A\{A - \pi_0(X)\}\sigma_0^{-2}(X),$$
$$H_2 = \{A - \pi_0(X)\}\sigma_0^{-2}(X)Y,$$
$$H_3 = \{A - \pi_0(X)\}\sigma_0^{-2}(X)Y^*,$$
$$H_4 = c(X)\{A - \pi_0(X)\}.$$

As $c(\cdot)$ is varied, one obtains all first order influence functions. We do not discuss the efficient choice of $c(\cdot)$ in this paper.

Our interest lies in the special case where $Y = Y^*$ w.p.1 (so there is but one response of interest) and thus, as in assumption $Aiiib$), $b = p$, $H_2 = H_3$ and we construct confidence interval for $\psi(\theta) = E[b^2(X)]$. In Section 5 we describe how we can use a confidence interval for $\psi(\theta) = E[b^2(X)]$ to obtain confidence intervals for the treatment effect function $b(x)$ and, most importantly, for the optimal treatment strategy $d_{opt}(x) = I[b(x) > 0]$ under which a subject with covariate value $x$ is treated if and only if the treatment effect $b(x)$ is positive ( i.e., $d_{opt}(x) = 1$).

### 3.2. Higher order influence functions for our model

#### 3.2.1. Dirac kernels, truncation bias, and a truncated parameter

In all of our examples the functions $p(\cdot)$ and $b(\cdot)$ are functions of conditional expectations given the continuous random variable $X$. It is well known that the associated point-evaluation functional $p(x)$ and $b(x)$ do not have first order influence functions. It then follows from part 5c of Theorem 2.3 and the dependence of $\mathbb{IF}_{1,\psi}(\theta) = \mathbb{V}[if_{1,\psi}(O_{i_1}; \theta)]$ on $b(\cdot)$ and $p(\cdot)$ evaluated at the point $X$ that, in none of our examples, does $\psi(\theta)$ have a second (or higher) order influence function.

As a precise understanding of the reason for the nonexistence of higher order influence functions for $\psi(\theta)$ is fundamental to our approach, we now use part 5c of Theorem 2.3 to prove that $\mathbb{IF}_{2,\psi}(\theta)$ does not exist by showing that the functional $if_{1,\psi}(o; \theta)$ does not have a first order influence function $\mathbb{V}[if_{1,if_{1,\psi}(o;\cdot)}(O; \theta)]$. Let $F_X$ and $f_X = f_X(\cdot)$ denote the marginal CDF and density of $X$. In this proof, we do not assume that $p(\cdot)$ and $b(\cdot)$ are functions of conditional expectations. Rather we only assume that our functional satisfies assumptions Ai)-Aiv)

Consider paths (parametric submodels) $\widetilde{\theta}_l(t)$ such that $\widetilde{\theta}_l(0) = \theta$ satisfying

$$p_l(t) \equiv p_l(x, t) \equiv p(x) + tc_l(x),$$
$$b_l(t) \equiv b_l(x, t) \equiv b(x) + ta_l(x),$$

where the sequences $c_l(\cdot)$ and $a_l(\cdot), l = 1, 2, \ldots$, are each dense in $L_2[F_X(x)]$. Let

$$s_l(O; \theta) = s_l(O|X; \theta) + s_l(X; \theta),$$

$s_l(O|X; \theta)$, and $s_l(X; \theta)$ denote the overall, conditional, and marginal scores

$$\partial ln f\left(O; \widetilde{\theta}_l(0)\right)/\partial t,\ \partial ln f\left(O|X; \widetilde{\theta}_l(0)\right)/\partial t,\ \partial ln f_X\left(X; \widetilde{\theta}_l(0)\right)/\partial t.$$

By linearity, $if_{1,\psi}(o; \theta)$ has an influence function only if the functionals $b(x)$ and $p(x)$ have one as well. Now by differentiating the identity

$$E_{\widetilde{\theta}_l(t)}[\{H_1 b_l(X, t) + H_3\}|X = x] = 0$$

with respect to $t$ and evaluating at $t = 0$, we have

$$-E_\theta \left[\{\{H_1 b(X) + H_3\}\} s_l(O|X) | X = x\right] = E_\theta[H_1|X = x] a_l(x).$$

However, by definition, $b(x)$ has an influence function $\mathbb{V}\left[if_{1,b(x)}(O;\theta)\right]$ at $\theta$ only if for $l = 1, 2, \ldots$, both $\partial b_l(x, t)/\partial t_{|t=0} = a_l(x)$ equals $E_\theta\left[if_{1,b(x)}(O;\theta) s_l(O;\theta)\right]$ and $E_\theta\left[if_{1,b(x)}(O;\theta)\right] = 0$. Thus if $if_{1,b(x)}(O;\theta)$ exists, it must satisfy

$$-E_\theta\left[\{H_1 b(X) + H_3\} s_l(O|X) | X = x\right]$$
$$= E_\theta[H_1|X = x] E_\theta\left[if_{1,b(x)}(O;\theta) s_l(O;\theta)\right].$$

Without loss of generality, suppose $H_1 \geq 0$ w.p.1. Now if we could find a 'kernel' $K_{f_X,\infty}(x, X)$ such that

$$r(x) = E_{f_X}\left[K_{f_X,\infty}(x, X) r(X)\right]$$

(3.7)
$$\equiv \int K_{f_X,\infty}(x, x^*) r(x^*) f_X(x^*) dx^* \text{ for all } r(\cdot) \in L_2(F_X)$$

then

$$if_{1,b(x)}(O;\theta) \equiv - \left[\begin{array}{c} \{E_\theta[H_1|X = x]\}^{-1/2} K_{f_X,\infty}(x, X) \\ \times \{E_\theta[H_1|X]\}^{-1/2} \{H_1 b(X) + H_3\} \end{array}\right]$$

would be an influence function since

$$E_\theta[H_1|X = x] E_\theta \left[\begin{array}{c} -\{E_\theta[H_1|X = x]\}^{-1/2} K_{f_X,\infty}(x, X) \times \\ \{E_\theta[H_1|X]\}^{-1/2} \{H_1 b(X) + H_3\} s_l(O;\theta) \end{array}\right]$$

$$= E[H_1|X = x]^{1/2} E_\theta \left[\begin{array}{c} -K_{f_X,\infty}(x, X) \{E_\theta[H_1|X]\}^{-1/2} \\ \times \{H_1 b(X) + H_3\} \{s_l(O|X) + s_l(X)\} \end{array}\right]$$

$$= E[H_1|X = x]^{1/2} E_{f_X} \left\{ E_\theta \left[\begin{array}{c} -K_{f_X,\infty}(x, X) \{E_\theta[H_1|X]\}^{-1/2} \times \\ \{H_1 b(X) + H_3\} s_l(O|X) | X \end{array}\right] \right\}$$

$$= -E_\theta\left[(H_1 b(X) + H_3) s_l(O|X) | X = x\right].$$

By an analogous argument

$$if_{1,p(x)}(O;\theta) = -\left[\begin{array}{c} \{E_\theta[H_1|X = x]\}^{-1/2} K_{f_X,\infty}(x, X) \\ \times \{E_\theta[H_1|X]\}^{-1/2} \{H_1 p(X) + H_2\} \end{array}\right]$$

would be an influence function.

Indeed since the sequences $\{c_l(\cdot)\}$ and $\{a_l(\cdot)\}$ are dense the existence of such a kernel is also a necessary condition for $if_{1,b(x)}(O;\theta)$ and $if_{1,p(x)}(O;\theta)$ to exist and thus for $if_{1,\psi}(o;\theta)$ to exist. A kernel satisfying Equation (3.7) is referred to as the Dirac delta function with respect to the measure $dF_X(x)$ and would clearly have to satisfy

(3.8)
$$K_{f_X,\infty}(x_{i_1}, x_{i_2}) = 0 \text{ if } x_{i_2} \neq x_{i_1}$$

were it to exist. Of course a kernel satisfying Equation (3.7) is known not to exist in $L_2[F_X] \times L_2[F_X]$. We conclude that $if_{1,\psi}(o;\theta)$ does not have an influence function and therefore $\mathbb{IF}_{2,2,\psi}(\theta)$ does not exist.

### A formal approach

To motivate how one might overcome this difficulty, we note that kernels satisfying Equation (3.7) exist as generalized functions or kernels (also known as Schwartz functions or distributions). We shall 'formally' derive higher order influence func-tions that appear to be elements of the space of generalized functions. However, we use these calculations only as motivation for statistical procedures based on ordinary kernels living in $L_2[F_X] \times L_2[F_X]$. Thus it does not matter whether these formal calculations could be made rigorous with appropriate redefinitions. Rather we can simply regard the following as results obtained by applying a "formal calculus" to part 5c of Theorem 2.3 that adds to the usual calculus additional identities licensed by Equations (3.7) and (3.8).

We will need the fact that, for any function $v(x;\theta)$, Eq. (3.8) implies that

$$v(x;\theta) K_{f_X,\infty}(x,X) = v(X;\theta) K_{f_X,\infty}(x,X).$$

We now show that

$$\mathbb{IF}_{2,2,\psi}(\theta) \equiv \mathbb{V}[IF_{2,2,\psi,i_1,i_2}(\theta)] = \Pi_{\theta,2}\left[\mathbb{V}\left[if_{1,if_{1,\psi}(O_{i_1};\cdot)}(O_{i_2};\theta)/2\right]|\mathcal{U}_1^{\perp 2,\theta}(\theta)\right]$$

would formally have U-statistic kernel

$$(3.9) \quad IF_{2,2,\psi,i_1,i_2}(\theta) = -\begin{bmatrix} \varepsilon_{b,i_1}(\theta) E_\theta[H_1|X_{i_1}]^{-\frac{1}{2}} K_{f_X,\infty}(X_{i_1},X_{i_2}) \\ E_\theta[H_1|X_{i_2}]^{-\frac{1}{2}} \varepsilon_{p,i_2}(\theta) \end{bmatrix},$$

with $\varepsilon_{b,i_1}(\theta) = \{B_{i_1} H_{1,i_1} + H_{3,i_1}\}$, $\varepsilon_{p,i_2}(\theta) = \{H_{1,i_2} P_{i_2} + H_{2,i_2}\}$.

To show Equation (3.9) note, by

$$\partial H(b,p)/\partial P = \partial\{BPH_1 + BH_2 + PH_3 + H_4\}/\partial P = BH_1 + H_3$$

and

$$\partial H(b,p)/\partial B = PH_1 + H_2,$$

we have

$$if_{1,if_{1,\psi}(O_{i_1};\cdot)}(O_{i_2};\theta) = Q_{2,b,\bar{i}_2}(\theta) + Q_{2,p,\bar{i}_2}(\theta) - IF_{1,\psi,i_2}(\theta),$$

where

$$Q_{2,p,\bar{i}_2}(\theta) \equiv \{B_{i_1} H_{1,i_1} + H_{3,i_1}\} if_{1,p(X_{i_1})}(O_{i_2};\theta)$$

$$= -\{B_{i_1} H_{1,i_1} + H_{3,i_1}\} E_\theta[H_1|X_{i_1}]^{-\frac{1}{2}}$$

$$\times K_{f_X,\infty}(X_{i_1},X_{i_2}) E_\theta[H_1|X_{i_2}]^{-\frac{1}{2}} \{P_{i_2} H_{1,i_2} + H_{2,i_2}\}$$

$$= -\varepsilon_{b,i_1}(\theta) E_\theta[H_1|X_{i_1}]^{-\frac{1}{2}} K_{f_X,\infty}(X_{i_1},X_{i_2}) E_\theta[H_1|X_{i_2}]^{-\frac{1}{2}} \varepsilon_{p,i_2}(\theta)$$

$$Q_{2,b,\bar{i}_2}(\theta) \equiv \{P_{i_1} H_{1,i_1} + H_{2,i_1}\} if_{1,b(X_{i_1})}(O_{i_2};\theta)$$

$$= -\varepsilon_{b,i_2}(\theta) E_\theta[H_1|X_{i_2}]^{-\frac{1}{2}} K_{f_X,\infty}(X_{i_2},X_{i_1}) E_\theta[H_1|X_{i_1}]^{-\frac{1}{2}} \varepsilon_{p,i_1}(\theta).$$

Thus, by part 5(c) of Theorem 2.3,

$$\mathbb{IF}_{2,2,\psi}(\theta) = \Pi_{\theta,2}\left[\frac{1}{2}\left\{\mathbb{Q}_{2,p,\bar{i}_2}(\theta) + \mathbb{Q}_{2,b,\bar{i}_2}(\theta) + \mathbb{IF}_{1,\psi,i_2}\right\}|\mathcal{U}_1^{\perp 2,\theta}(\theta)\right]$$

$$= \frac{1}{2}\left\{\mathbb{Q}_{2,p,\bar{i}_2}(\theta) + \mathbb{Q}_{2,b,\bar{i}_2}(\theta)\right\}$$

$$= \mathbb{Q}_{2,p,\bar{i}_2}(\theta) \equiv \mathbb{V}[RHS \text{ of Equation (3.9)}]$$

since $\mathbb{IF}_{1,\psi,i_2}$ is a function of only one subject's data and $Q_{2,p,\bar{i}_2}(\theta)$ and $Q_{2,b,\bar{i}_2}(\theta)$ are the same up to a permutation that exchanges $i_2$ with $i_1$.

To obtain $IF_{3,3,\psi,\bar{i}_m}(\theta)$, one must derive the influence function $if_{1,if_{2,2,\psi}(O_{i_1},O_{i_2};\cdot)}(O_{i_3};\theta)$ of $if_{2,2,\psi}(O_{i_1},O_{i_2};\theta)$. The formula for $IF_{3,3,\psi,\bar{i}_m}(\theta)$ is given in Equation (3.13). A detailed derivation is given in our technical report. Here we simply note that the only essentially new point is that we now require the influence function of $K_{f_X,\infty}(X_{i_1}, X_{i_2})$, which, as shown next, is given by

$$(3.10) \quad IF_{1,K_{f_X,\infty}}(X_{i_1},X_{i_2}) = -\left\{ \begin{array}{c} K_{f_X,\infty}(X_{i_1}, X_{i_3}) K_{f_X,\infty}(X_{i_3}, X_{i_2}) \\ -K_{f_X,\infty}(X_{i_1}, X_{i_2}) \end{array} \right\}.$$

To see that if Equation (3.7) held, Equation (3.10) would hold, note that for any path $\widetilde{\theta}(t)$ with $\widetilde{\theta}(0) = f_X(\cdot)$, $h(x) = E_{\widetilde{\theta}(t)}\left[K_{\widetilde{\theta}(t),\infty}(x, X_{i_1}) h(X_{i_1})\right]$. Differentiating with respect to $t$ and evaluating at $t = 0$, we have

$$0 = E_\theta[K_{f_X,\infty}(x, X) h(X) S(\theta)] + E_\theta\left[\left\{\frac{\partial}{\partial t}K_{\widetilde{\theta}(t),\infty}(x, X_{i_1})_{|t=0}\right\} h(X_{i_1})\right].$$

Hence it suffices to show that

$$- E_\theta[K_{f_X,\infty}(x, X) h(X) S(\theta)]$$
$$= E_\theta[\{E_\theta\{-K_{f_X,\infty}(x, X_{i_2}) K_{f_X,\infty}(X_{i_2}, X_{i_1}) S_{i_2}(\theta) | X_{i_1}\}\} h(X_{i_1})].$$

But, by Equation (3.7),

$$E_\theta[\{E_\theta\{-K_{f_X,\infty}(x, X_{i_2}) K_{f_X,\infty}(X_{i_2}, X_{i_1}) S_{i_2}(\theta) | X_{i_1}\}\} h(X_{i_1})]$$
$$= E_\theta[-K_{f_X,\infty}(x, X_{i_1}) S_{i_1}(\theta) h(X_{i_1})].$$

*Feasible estimators*

These "formal" calculations motivate a "truncated Dirac" approach to estimate $\psi(\theta)$. Let $\{z_l(\cdot)\} \equiv \{z_l(X); 1, 2, \ldots\}$ be a countable sequence of known basis functions with dense span in $L_2(F_X)$ and define $\overline{z}_k(X)^T = (z_1(X), \ldots, z_k(X))$. Define

$$K_{f_X,k}(X_{i_1}, X_{i_2}) \equiv \overline{z}_k(X_{i_1})^T \left\{E_{f_X}\left[\overline{z}_k(X)\overline{z}_k(X)^T\right]\right\}^{-1} \overline{z}_k(X_{i_2})$$

to be the projection kernel in $L_2(F_X)$ onto the subspace

$$lin\{\overline{z}_k(X)\} \equiv \{\eta^T \overline{z}_k(x); \eta \in R^k, \eta^T \overline{z}_k(x) \in L_2(F_X)\}$$

spanned by the elements of $\overline{z}_k(X)$. That is, for any $h(x)$,

$$\Pi_{f_X}[h(X) | lin\{\overline{z}_k(x)\}]$$
$$= E_{f_X}[K_{f_X,k}(x, X) h(X)]$$
$$= \overline{z}_k(x)^T \left\{E_{f_X}\left[\overline{z}_k(X)\overline{z}_k(X)^T\right]\right\}^{-1} E_{f_X}[\overline{z}_k(X) h(X)].$$

Then we can view $K_{f_X,k}(x_{i_1}, x_{i_2})$ as a truncated at $k$ approximation to $K_{f_X,\infty}(x_{i_1}, x_{i_2})$ that is in $L_2[F_X] \times L_2[F_X]$ and satisfies Equation (3.7) for all $r(x) \in lin\{\overline{z}_k(X)\}$. Then a natural idea would be to substitute

$$IF_{2,2,\psi,i_1,i_2}^{(k)}(\widehat{\theta}) \equiv \left( \begin{array}{c} -\varepsilon_{b,i_1}(\widehat{\theta}) E_{\widehat{\theta}}[H_1|X_{i_1}]^{-\frac{1}{2}} K_{\widehat{f_X},k}(X_{i_1}, X_{i_2}) \\ \times E_{\widehat{\theta}}[H_1|X_{i_2}]^{-\frac{1}{2}} \varepsilon_{p,i_2}(\widehat{\theta}) \end{array} \right)$$

with, for example,

$$\varepsilon_{b,i_1}(\widehat{\theta}) = E_{\widehat{\theta}}[H_1|X_{i_1}]^{-\frac{1}{2}}\left\{\widehat{B}_{i_1} H_{1,i_1} + H_{3,i_1}\right\}$$

for the generalized function $IF_{2,2,\psi,i_1,i_2}(\widehat{\theta})$ based on Equation (3.9) resulting in the feasible second $U$-statistic estimator

$$\widehat{\psi}_2^{(k)} = \psi\left(\widehat{\theta}\right) + \mathbb{IF}_{1,\psi}\left(\widehat{\theta}\right) + \mathbb{IF}_{2,2,\psi(\theta)}^{(k)}\left(\widehat{\theta}\right)$$

where

$$\mathbb{IF}_{2,2,\psi}^{(k)}\left(\widehat{\theta}\right) \equiv \mathbb{V}\left[IF_{2,2,\psi,i_1,i_2}^{(k)}\left(\widehat{\theta}\right)\right].$$

To avoid having to do a matrix inversion it would be convenient, when possible, to choose $\overline{z}_k(X) = \overline{\varphi}_k(X) / \left\{\widehat{f}_X(X)\right\}^{1/2}$ where $\varphi_1(X), \varphi_2(X), \ldots$ is a complete orthonormal basis with respect to Lebesgue measure in $R^d$. Then $E_{\widehat{f}_X}[\overline{z}_k(X) \times \overline{z}_k(X)^T] = I_{k \times k}$ so

$$K_{\widehat{f}_X, k}(X_{i_1}, X_{i_2}) = \overline{z}_k(X_{i_1})^T \overline{z}_k(X_{i_2}) = \frac{K_{Leb,k}(X_{i_1}, X_{i_2})}{\left\{\widehat{f}_X(X_{i_1}) \widehat{f}_X(X_{i_2})\right\}^{1/2}},$$

where

$$K_{Leb,k}(X_{i_1}, X_{i_2}) \equiv \overline{\varphi}_k(X_{i_1})^T \overline{\varphi}_k(X_{i_2}).$$

This choice corresponds to having taken

$$K_{f_X, \infty}(X_{i_1}, X_{i_2}) = K_{Leb,\infty}(X_{i_1}, X_{i_2}) / \left\{f_X(X_{i_1}) f_X(X_{i_2})\right\}^{1/2}$$

in our formal calculations where $K_{Leb,\infty}(X_{i_1}, X_{i_2})$ is the Dirac delta function with respect to Lebesgue measure. In that case with $G \equiv g(X) \equiv f_X(X) E_\theta[H_1|X]$ and $\widehat{G} \equiv \widehat{g}(X) \equiv \widehat{f}_X(X) E_{\widehat{\theta}}[H_1|X]$,

(3.11) $IF_{2,2,\psi,i_1,i_2}(\theta) = -\varepsilon_{b,i_1}(\theta) g(X_{i_1})^{-\frac{1}{2}} K_{Leb,\infty}(X_{i_1}, X_{i_2}) g(X_{i_2})^{-\frac{1}{2}} \varepsilon_{p,i_2}(\theta)$

(3.12)
$$IF_{2,2,\psi,i_1,i_2}^{(k)}\left(\widehat{\theta}\right) = -\varepsilon_{b,i_1}\left(\widehat{\theta}\right) \widehat{g}(X_{i_1})^{-\frac{1}{2}} K_{Leb,k}(X_{i_1}, X_{i_2}) \widehat{g}(X_{i_2})^{-\frac{1}{2}} \varepsilon_{p,i_2}\left(\widehat{\theta}\right)$$

Unfortunately, we show later that this choice for $\overline{z}_k$ has good statistical properties only when $f_X$ is known to lie in a Holder ball with exponent exceeding $\max\{\beta_b, \beta_p\}$. In our technical report we show one can proceed by induction to formally obtain that for $m = 3, 4, \ldots$,

$IF_{m,m,\psi,\overline{i}_m}(\theta)$
(3.13)
$$= \varepsilon_{b,i_1}(\theta) \varsigma(X_{i_1})^{-\frac{1}{2}} \left[ \begin{array}{c} \sum_{j=0}^{m-2} c(m,j) \times \\ \prod_{s=1}^{j} \frac{H_{1,i_{s+1}}}{\varsigma(X_{i_{s+1}})} K_{f_X,\infty}(X_{i_s}, X_{i_{s+1}}) \\ \times K_{f_X,\infty}(X_{i_{j+1}}, X_{i_m}) \end{array} \right] \varsigma(X_{i_m})^{-\frac{1}{2}} \varepsilon_{p,i_m}(\theta)$$

where $\varsigma(X) = E[H_1|X]$ and $c(m,j) = \binom{m-2}{j}(-1)^{(j+1)}$, which we then use to obtain $IF_{m,m,\psi,\overline{i}_m}^{(k)}\left(\widehat{\theta}\right)$ and $\widehat{\psi}_m^{(k)} = \widehat{\psi}_2^{(k)} + \sum_{j=3}^{m} \mathbb{IF}_{j,j,\psi}^{(k)}(\widehat{\theta})$.

*Statistical properties*

We shall prove below that the estimator $\widehat{\psi}_m^{(k)}$ has variance

$$\text{var}_\theta\left[\widehat{\psi}_m^{(k)}\right] \asymp \left(\frac{1}{n}\max\left[1,\left(\frac{k}{n}\right)^{m-1}\right]\right),$$

when $\{z_l(X); l = 1, 2, \ldots\}$ is a compact wavelet basis. (Robins et al. [18] proves this result for more general bases). We also prove that the bias

$$E_\theta\left[\widehat{\psi}_m^{(k)}\right] - \psi(\theta) = TB_k(\theta) + EB_m(\theta),$$

of $\widehat{\psi}_m^{(k)}$ is the sum of a truncation bias term of order

$$TB_k(\theta) = O_p\left(k^{-(\beta_b+\beta_p)/d}\right)$$

(for a basis $\{z_l(X); l = 1, 2, \ldots\}$ that provides optimal rate approximation for Hölder balls) and an estimation bias term of order

$$EB_m(\theta) = O_p\left(\{P-\widehat{P}\}\{B-\widehat{B}\}\left(\frac{G-\widehat{G}}{\widehat{G}}\right)^{m-1}\right)$$

$$= O_p\left(n^{-\frac{(m-1)\beta_g}{2\beta_g+d} - \frac{\beta_b}{2\beta_b+d} - \frac{\beta_p}{2\beta_p+d}}\right).$$

Note this estimation bias is $O_P\left(\left\|\theta-\widehat{\theta}\right\|^{m+1}\right)$. It gets its name from the fact that, unlike the truncation bias, it would be exactly zero if the initial estimator $\widehat{\theta}$ happened to equal $\theta$. Thus, the U-statistic estimator $\widehat{\psi}_m^{(k)}$ for our functional $\psi(\theta)$ (which does not admit a second order influence function) differs from the U-statistic estimators $\widehat{\psi}_m$ of Eq. (2.3) for functionals that admit second order influence functions in that, owing to truncation bias, the total bias of $\widehat{\psi}_m^{(k)}$ is not $O_p\left(\left\|\theta-\widehat{\theta}\right\|^{m+1}\right)$. The choice of $k$ determines the trade-off between the variance and truncation bias. As $k \to \infty$ with $n$ fixed, $\text{var}_\theta\left[\widehat{\psi}_m^{(k)}\right] \to \infty$ and $TB_k(\theta) \to 0$. Thus, we can heuristically view the non-existent estimator $\widehat{\psi}_m = \widehat{\psi}_m^{(k=\infty)}$ as the choice of $k$ that results in no truncation bias (and therefore a total bias of $O_p\left(\left\|\theta-\widehat{\theta}\right\|^{m+1}\right)$) at the expense of an infinite variance. Writing $k = k(n) = n^\rho$, the order of the asymptotic MSE of $\widehat{\psi}_m^{(k)}$ is minimized at the value of $\rho$ for which order of the variance equals the order of the sum of the truncation and estimation bias.

**Remark 3.6.** The models of Examples 1-4 exhibit a spectrum of different likelihood functions and therefore a spectrum of different first order and higher order scores. Nonetheless, because the first order influence functions of the functionals $\psi(\theta)$ share a common structure, we were able to use part 5c of Theorem 2.3 to formally derive $IF_{m,m,\psi,\bar{i}_m}(\theta)$ and, thus, the feasible $IF^{(k)}_{m,m,\psi,\bar{i}_m}(\widehat{\theta})$ in Examples 1-4 in a unified manner without needing to consult the full likelihood function for any of the models. See Remark (2.5) above for a closely related discussion.

*A critical non-uniqueness*

We have as yet neglected a critical non-uniqueness in our definition of $\mathbb{IF}^{(k)}_{m,m,\psi(\theta)}\left(\widehat{\theta}\right)$ and thus $\widehat{\psi}^{(k)}_m$ that poses a significant problem for our "truncated Dirac" approach. For instance, when $m = 3$, the two generalized $U$-statistic kernels $IF_{3,3,\psi,i_1,i_2,i_3}(\theta)$ of Equation (3.13) and

$$IF^*_{3,3,\psi,i_1,i_2,i_3}(\theta) \equiv \frac{\varepsilon_{b,i_1}(\theta)}{\varsigma(X_{i_1})^{\frac{1}{2}}} \left[ \begin{array}{c} \frac{H_{1,i_2}}{\varsigma(X_{i_2})} K_{f_X,\infty}(X_{i_1},X_{i_2}) \\ -E_\theta[K_{f_X,\infty}(X_{i_1},X_{i_2})|X_{i_1}] \end{array} \right]$$
$$\times K_{f_X,\infty}(X_{i_1},X_{i_3}) \frac{\varepsilon_{p,i_3}(\theta)}{\varsigma(X_{i_3})^{\frac{1}{2}}}$$

are precisely equal, by Eq. (3.8); nonetheless, upon truncation, they result in different feasible kernels;

$$IF^{(k)}_{3,3,\psi,i_1,i_2}\left(\widehat{\theta}\right)$$
$$= \frac{\widehat{\varepsilon}_{b,i_1}(\theta)}{\widehat{\varsigma}(X_{i_1})^{\frac{1}{2}}} \left[ \begin{array}{c} \frac{H_{1,i_2}}{\widehat{\varsigma}(X_{i_2})} K_{\widehat{f}_X,k}(X_{i_1},X_{i_2}) K_{\widehat{f}_X,k}(X_{i_2},X_{i_3}) \\ -K_{\widehat{f}_X,k}(X_{i_1},X_{i_3}) \end{array} \right] \times \frac{\widehat{\varepsilon}_{p,i_3}(\theta)}{\widehat{\varsigma}(X_{i_3})^{\frac{1}{2}}}$$

and

$$IF^{(k),*}_{3,3,\psi,i_1,i_2,i_3}\left(\widehat{\theta}\right) \equiv \frac{\widehat{\varepsilon}_{b,i_1}(\theta)}{\widehat{\varsigma}(X_{i_1})^{\frac{1}{2}}} \left[ \begin{array}{c} \frac{H_{1,i_2}}{\widehat{\varsigma}(X_{i_2})} K_{\widehat{f}_X,k}(X_{i_1},X_{i_2}) \\ -E_{\widehat{\theta}}\left[K_{\widehat{f}_X,k}(X_{i_1},X_{i_2})|X_{i_1}\right] \end{array} \right]$$
$$\times K_{\widehat{f}_X,k}(X_{i_1},X_{i_3}) \frac{\widehat{\varepsilon}_{p,i_3}(\theta)}{\widehat{\varsigma}(X_{i_3})^{\frac{1}{2}}}$$

with possibly different orders of bias. For simplicity, we consider the case where $H_1 = 1$ as in Examples 1b. Let $\delta B \equiv B - \widehat{B}$, $\delta P \equiv P - \widehat{P}$, $\delta f = \delta g \equiv \frac{f}{\widehat{f}} - 1$, and

$$\overline{Z}_k \equiv \overline{z}_k(X) \equiv \left\{\widehat{E}\left[\overline{\varphi}_k(X)\overline{\varphi}_k(X)^T\right]\right\}^{-1/2} \overline{\varphi}_k(X),$$

then,

$$E_\theta\left[IF^{(k),*}_{3,3,\psi,i_1,i_2,i_3}\left(\widehat{\theta}\right)\right]$$

$$= E_\theta \left[ \begin{array}{c} \delta B_{i_1} \times \\ E_{\widehat{\theta}}\left[\left(\frac{f(X_{i_2})}{\widehat{f}(X_{i_2})} - 1\right)\overline{z}_k(X_{i_2})^T\right]\overline{z}_k(X_{i_1}) \\ E_\theta\left[\delta P_{i_3}\overline{z}_k(X_{i_3})^T\right]\overline{z}_k(X_{i_1}) \end{array} \right]$$

$$= E_{\widehat{\theta}} \left[ \begin{array}{c} \left(\frac{f(X_{i_1})}{\widehat{f}(X_{i_1})} - 1 + 1\right)\delta B_{i_1} \times \\ E_{\widehat{\theta}}\left[\left(\frac{f(X_{i_2})}{\widehat{f}(X_{i_2})} - 1\right)\overline{z}_k(X_{i_2})^T\right]\overline{z}_k(X_{i_1}) \\ E_{\widehat{\theta}}\left[\left(\frac{f(X_{i_3})}{\widehat{f}(X_{i_3})} - 1 + 1\right)\delta P_{i_3}\overline{z}_k(X_{i_3})^T\right]\overline{z}_k(X_{i_1}) \end{array} \right]$$

$$= E_{\widehat{\theta}}\left[\delta B_{i_1} E_{\widehat{\theta}}\left[\delta f(X_{i_2})\overline{Z}^T_{k,i_2}\right]\overline{Z}_{k,i_1} E_{\widehat{\theta}}\left[\delta P_{i_3}\overline{Z}^T_{k,i_3}\right]\overline{Z}_{k,i_1}\right]$$
$$+ O_p\left(\{B - \widehat{B}\}\{P - \widehat{P}\}\{G - \widehat{G}\}^2\right)$$

and

$$E_\theta \left[ IF^{(k)}_{3,3,\psi,i_1,i_2,i_3}\left(\widehat{\theta}\right) \right]$$

$$= E_{\widehat{\theta}} \left[ \begin{array}{c} E_\theta \left[ \delta B_{i_1} \overline{z}_k(X_{i_1})^T \right] \overline{z}_k(X_{i_2}) \left( \frac{f(X_{i_2})}{\widehat{f}(X_{i_2})} - 1 \right) \\ \times \overline{z}_k(X_{i_2})^T E_\theta \left[ \delta P_{i_3} \overline{z}_k(X_{i_3}) \right] \end{array} \right]$$

$$= E_{\widehat{\theta}} \left[ \begin{array}{c} E_{\widehat{\theta}} \left[ (\delta f(X_{i_1}) + 1) \delta B_{i_1} \overline{Z}^T_{k,i_1} \right] \\ \times \overline{Z}_{k,i_2} \left( \frac{f(X_{i_2})}{\widehat{f}(X_{i_2})} - 1 \right) \\ \times \overline{Z}^T_{k,i_2} E_{\widehat{\theta}} \left[ (\delta f(X_{i_3}) + 1) \delta P_{i_3} \overline{Z}_{k,i_3} \right] \end{array} \right]$$

$$= E_{\widehat{\theta}} \left[ \begin{array}{c} E_{\widehat{\theta}} \left[ \delta B_{i_1} \overline{Z}^T_{k,i_1} \right] \overline{Z}_{k,i_2} \left( \frac{f(X_{i_2})}{\widehat{f}(X_{i_2})} - 1 \right) \\ \times \overline{Z}^T_{k,i_2} E_{\widehat{\theta}} \left[ \delta P_{i_3} \overline{Z}_{k,i_3} \right] \end{array} \right]$$

$$+ O_p \left( \{B - \widehat{B}\} \{P - \widehat{P}\} \{G - \widehat{G}\}^2 \right).$$

That is,

$$E_\theta \left[ IF^{(k),*}_{3,3,\psi,i_1,i_2,i_3}\left(\widehat{\theta}\right) \right] - E_\theta \left[ IF^{(k)}_{3,3,\psi,i_1,i_2,i_3}\left(\widehat{\theta}\right) \right]$$

$$= E_{\widehat{\theta}} \left[ \delta B_{i_1} E_{\widehat{\theta}} \left[ \delta f(X_{i_2}) \overline{Z}^T_{k,i_2} \right] \overline{Z}_{k,i_1} E_{\widehat{\theta}} \left[ \delta P_{i_3} \overline{Z}^T_{k,i_3} \right] \overline{Z}_{k,i_1} \right]$$

$$- E_{\widehat{\theta}} \left[ \begin{array}{c} E_{\widehat{\theta}} \left[ \delta B_{i_1} \overline{Z}^T_{k,i_1} \right] \overline{Z}_{k,i_2} \delta f(X_{i_2}) \\ \times \overline{Z}^T_{k,i_2} E_{\widehat{\theta}} \left[ \delta P_{i_3} \overline{Z}_{k,i_3} \right] \end{array} \right]$$

$$+ O_p \left( \{B - \widehat{B}\} \{P - \widehat{P}\} \{G - \widehat{G}\}^2 \right)$$

$$= E_{\widehat{\theta}} \left[ \delta B \Pi_{\widehat{\theta}} \left[ \delta P | \overline{Z}_k \right] \Pi_{\widehat{\theta}} \left[ \delta f(X) | \overline{Z}_k \right] \right]$$

$$- E_{\widehat{\theta}} \left[ \Pi_{\widehat{\theta}} \left[ \delta P | \overline{Z}_k \right] \delta f(X) \Pi_{\widehat{\theta}} \left[ \delta B | \overline{Z}_k \right] \right]$$

$$+ O_p \left( \{B - \widehat{B}\} \{P - \widehat{P}\} \{G - \widehat{G}\}^2 \right)$$

$$= E_{\widehat{\theta}} \left[ \left\{ \begin{array}{c} \Pi_{\widehat{\theta}} \left[ \delta P | \overline{Z}_k \right] \times \\ \Pi^\perp_{\widehat{\theta}} \left[ \delta B | \overline{Z}_k \right] \Pi_{\widehat{\theta}} \left[ \delta f(X) | \overline{Z}_k \right] \\ - \Pi_{\widehat{\theta}} \left[ \delta B | \overline{Z}_k \right] \Pi^\perp_{\widehat{\theta}} \left[ \delta f(X) | \overline{Z}_k \right] \end{array} \right\} \right]$$

$$+ O_p \left( \{B - \widehat{B}\} \{P - \widehat{P}\} \{G - \widehat{G}\}^2 \right)$$

where $\Pi_{\widehat{\theta}}\left[h(X)|\overline{Z}_k\right]$ and $\Pi^\perp_{\widehat{\theta}}\left[h(X)|\overline{Z}_k\right]$ respectively denote the projection under $F\left(\cdot;\widehat{\theta}\right)$ in $L_2\left(\widehat{F}\right)$ of $h(X)$ on the $k$ dimensional linear subspace $lin\left\{\overline{z}_k(X)\right\}$ spanned by the components of the vector $\overline{z}_k(X)$ and the projection on the orthocomplement of this subspace. Since the basis $\{\varphi_l(X); l = 1, 2, \ldots\}$ provides optimal rate approximation for Hölder balls, it is easy to verify that the difference is of order

$$O_p \left( \begin{array}{c} n^{-\frac{\beta_p/d}{1+2\beta_p/d} - \frac{\beta_g/d}{1+2\beta_g/d}} k^{-\beta_b/d} + n^{-\frac{\beta_p/d}{1+2\beta_p/d} - \frac{\beta_b/d}{1+2\beta_b/d}} k^{-\beta_g/d} \\ + n^{-\frac{\beta_p/d}{1+2\beta_p/d} - \frac{\beta_b/d}{1+2\beta_b/d} - \frac{2\beta_g/d}{1+2\beta_g/d}} \end{array} \right),$$

which we expect to be sharp for many bases, although not for Haar. For concreteness, we shall look at an example. Suppose $\beta_b/d = \beta_p/d = 0.3$ and $\beta_g/d = 0.1$, thus, by choosing $k = n^{\frac{5}{6}}$, $\widehat{\psi}_3^{(k)}$ converges to $\psi(\theta)$ at rate $n^{-\frac{1}{2}}$. In contrast, the order,

$$\min_k \left( n^{-\frac{\beta_p/d}{1+2\beta_p/d} - \frac{\beta_g/d}{1+2\beta_g/d}} k^{-\beta_b/d} + n^{-\frac{\beta_p/d}{1+2\beta_p/d} - \frac{\beta_b/d}{1+2\beta_b/d}} k^{-\beta_g/d} + \sqrt{\frac{1}{n} \max\left(1, \frac{k^2}{n^2}\right)} \right),$$

of the optimal root mean squares error of $\widehat{\psi}_3^{(k),*}$ that uses $IF_{3,3,\psi,\bar{i}_3}^{(k),*}\left(\widehat{\theta}\right)$ is $n^{-0.477} \gg n^{-0.5}$. Thus, for many orthonormal bases, $\widehat{\psi}_3^{(k),*}$ converges to $\psi(\theta)$ at a slower rate than $\widehat{\psi}_3^{(k)}$ which uses $IF_{3,3,\psi,\bar{i}_3}^{(k)}\left(\widehat{\theta}\right)$. Nothing in our development up to this point provides any guidance as to which of the many equivalent generalized U-statistic kernels should be selected for truncation. To provide some guidance, we introduce an alternative approach to the estimation of $\psi(\theta)$ based on truncated parameters that admit higher order influence functions. The class of estimators we derive using this alternative approach includes members algebraically identical to the estimators $\widehat{\psi}_m^{(k)}$ but does not include estimators equivalent to less efficient estimators such as $\widehat{\psi}_3^{(k),*}$.

*An approach based on truncated parameters*

We introduce a class of truncated parameters $\widetilde{\psi}_k(\theta)$ that (i) depend on the sample size through a positive integer index $k = k(n)$ (which we refer to as the truncation index and will be optimized below), (ii) have influence functions $\mathbb{IF}_{m,\widetilde{\psi}_k}(\theta)$ of all orders $m$, (iii) equal $\psi(\theta)$ on a large subset $\Theta_{sub,k}$ of $\Theta$ and (iv) the initial estimator $\widehat{\theta}$ is an element of $\Theta_{sub,k}$ so that the plug-ins $\psi\left(\widehat{\theta}\right)$ and $\widetilde{\psi}_k\left(\widehat{\theta}\right)$ are equal. To prepare we introduce a simplified notation. For functions $h(o, \cdot)$ or $r(\cdot)$ of $\theta$, we will often write $h\left(o, \widehat{\theta}\right)$ and $r\left(\widehat{\theta}\right)$ as $\widehat{h}(o)$ and $\widehat{r}$, and $E_{\widehat{\theta}}[\cdot]$ as $\widehat{E}[\cdot]$. Similarly, we often write $h(o, \theta)$ and $r(\theta)$ as $h(o)$ and $r$, and $E_\theta[\cdot]$ as $E[\cdot]$. Further we shall introduce slightly different definitions of truncation and estimation bias.

Define the estimator $\psi_{m,k}\left(\widehat{\theta}\right) \equiv \psi\left(\widehat{\theta}\right) + \mathbb{IF}_{m,\widetilde{\psi}_k}\left(\widehat{\theta}\right)$ or, equivalently, $\widehat{\psi}_{m,k} \equiv \widehat{\psi} + \widehat{\mathbb{IF}}_{m,\widetilde{\psi}_k}$. Then the conditional bias $E\left[\widehat{\psi}_{m,k}|\widehat{\theta}\right] - \psi$ of $\widehat{\psi}_{m,k}$ is $TB_k + EB_m$, where the truncation bias $TB_k = \widetilde{\psi}_k - \psi$ is zero for $\theta \in \Theta_{sub,k}$ and does not depend on $m$ and the estimation bias $EB_{m,k} = E\left[\widehat{\psi}_{m,k}|\widehat{\theta}\right] - \widetilde{\psi}_k$ is $O_P\left(\|\widehat{\theta} - \theta\|^{m+1}\right)$ by Theorem 2.2. Since, as we show later, the order of $EB_{m,k}$ does not depend on $k$, we will abbreviate $EB_{m,k}$ as $EB_m$, suppressing the dependence on $k$. Under minimal conditions, the conditional variance of $\widehat{\psi}_{m,k}$ is of the order of $\text{var}\left[\mathbb{IF}_{m,\widetilde{\psi}_k}\right]$ whenever $k \equiv k(n) \geq n$. The rate of convergence of $\widehat{\psi}_{m,k}$ to $\psi$ can depend on the choice of $\widetilde{\psi}_k$. Nevertheless, many different choices $\widetilde{\psi}_k$ result in estimators $\widehat{\psi}_{m,k}$ that achieve what we conjecture to be the optimal rate for estimators of the form $\widehat{\psi}_{m,k}$. We choose, among all such $\widetilde{\psi}_k$, the class that minimizes the computational complexity of $\widehat{\psi}_{m,k}$. Specifically for all $\widetilde{\psi}_k$ in our chosen class and all $j$, $\mathbb{IF}_{jj,\widetilde{\psi}_k}$ consists of a single term rather than a sum of many terms. We conjecture this appealing property does not hold for any $\widetilde{\psi}_k$ outside our class. We now describe this choice. The parameter $\widetilde{\psi}_k$ is defined in terms of $k(n)$−dimensional 'working' linear parametric submodels for $p(\cdot)$ and $b(\cdot)$ depending on unknown parameters $\overline{\alpha}_k$ and $\overline{\eta}_K$ through the basepoints $\widehat{p}(\cdot)$ and $\widehat{b}(\cdot)$, where $\widehat{p}(\cdot)$ and $\widehat{b}(\cdot)$ are initial estimators from the training sample. Specifically let $\dot{p}(X)$ and $\dot{b}(X)$ be arbitrary known functions chosen by the analyst satisfying Eqs. (3.14)-(3.16) below.

$$\text{(3.14)} \qquad \dot{p}(X)\dot{b}(X)E[H_1|X] \geq 0 \; w.p.1,$$

$$\text{(3.15)} \qquad \left\|\frac{\dot{p}(X)}{\dot{b}(X)}\right\|_\infty < C^*, \left\|\frac{\dot{b}(X)}{\dot{p}(X)}\right\|_\infty < C^*,$$

$$\text{(3.16)} \qquad \frac{\dot{p}(X)}{\dot{b}(X)} \text{ has at least } \lceil max\{\beta_b,\; \beta_p\} \rceil \text{ derivatives.}$$

Particular choices of $\dot{p}(X)$ and $\dot{b}(X)$ can make the form of $\mathbb{IF}_{m,\widetilde{\psi}_k}\left(\widehat{\theta}\right)$ more aesthetic. The choice has no bearing on the rate of convergence of the estimator $\widehat{\psi}_{m,k}$ to $\psi(\theta)$. Often there are fairly natural choices for $\dot{p}(\cdot)$ and $\dot{b}(\cdot)$. See Remark 3.9 below for examples. Let $\overline{\alpha}_k, \overline{\eta}_k$ be $k$−vectors of unknown parameters and consider the 'working' linear models

$$\text{(3.17)} \qquad p^*(X,\overline{\alpha}_k) \equiv \widehat{p}(X) + \dot{p}(X)\overline{\alpha}_k^T \overline{z}_k(X) \equiv \widehat{P} + \dot{P}\overline{\alpha}_k^T\overline{Z}_k,$$

$$\text{(3.18)} \qquad b^*(X,\overline{\eta}_k) = \widehat{b}(X) + \dot{b}(X)\overline{\eta}_k^T \overline{z}_k(X) = \widehat{B} + \dot{B}\overline{\eta}_k^T\overline{Z}_k.$$

We define the parameters $\widetilde{\overline{\eta}}_k(\theta)$ and $\widetilde{\overline{\alpha}}_k(\theta)$ respectively to be the solution to

$$\text{(3.19)} \qquad \begin{aligned} 0 &= E_\theta\left[\partial H\left(b^*(X,\overline{\eta}_k),p^*(X,\overline{\alpha}_k)\right)/\partial\overline{\alpha}_k\right] \\ &= E_\theta\left[\{H_1 b^*(X,\overline{\eta}_k) + H_3\}\dot{P}\overline{Z}_k\right], \end{aligned}$$

$$\text{(3.20)} \qquad \begin{aligned} 0 &= E_\theta\left[\partial H\left(b^*(X,\overline{\eta}_k),p^*(X,\overline{\alpha}_k)\right)/\partial\overline{\eta}_k\right] \\ &= E_\theta\left[\{H_1 p^*(X,\overline{\alpha}_k) + H_2\}\dot{B}\overline{Z}_k\right]. \end{aligned}$$

The solution to (3.19) and (3.20) exist in closed form as

$$\text{(3.21)} \qquad \widetilde{\overline{\eta}}_k(\theta) = -E_\theta\left[\dot{B}\dot{P}H_1\overline{Z}_k\overline{Z}_k^T\right]^{-1} E_\theta\left[\overline{Z}_k\dot{P}\left\{H_1\widehat{B} + H_3\right\}\right],$$

$$\text{(3.22)} \qquad \widetilde{\overline{\alpha}}_k(\theta) = -E_\theta\left[\dot{P}\dot{B}H_1\overline{Z}_k\overline{Z}_k^T\right]^{-1} E_\theta\left[\overline{Z}_k\dot{B}\left\{H_1\widehat{P} + H_2\right\}\right].$$

Next define $\widetilde{b}(\theta) = \widetilde{b}(\cdot,\theta) = b^*\left(\cdot,\widetilde{\overline{\eta}}_k(\theta)\right)$ and $\widetilde{p}(\theta) = \widetilde{p}(\cdot,\theta) = p^*\left(\cdot,\widetilde{\overline{\alpha}}_k(\theta)\right)$ and

$$\widetilde{\psi}_k(\theta) = E_\theta\left[H\left(\widetilde{b}(\theta),\widetilde{p}(\theta)\right)\right].$$

Note the models $p^*(\cdot,\overline{\alpha}_k)$ and $b^*(\cdot,\overline{\eta}_k)$ are used only to define the truncated parameter $\widetilde{\psi}_k(\theta)$. They are not assumed to be correctly specified. In particular, the training sample estimates $\widehat{p},\widehat{b}$ need not be based on the models $p^*(\cdot,\overline{\alpha}_k)$, $b^*(\cdot,\overline{\eta}_k)$. We now compare our truncated parameter $\widetilde{\psi}_k(\theta)$ with $\psi(\theta)$ and calculate the truncation bias. It is important to keep in mind that $b,p$ are components of the unknown $\theta$ while $\dot{p},\dot{b},\widehat{p},\widehat{b}$ are regarded as known functions.

**Theorem 3.7.** *If our model satisfies (Ai)–(Aiii) and*

$$\theta \in \Theta_{sub,k} = \{\theta; p(\cdot) = p^*(\cdot, \overline{\alpha}_k) \text{ for some } \overline{\alpha}_k \text{ or } b(\cdot) = b^*(\cdot, \overline{\eta}_k) \text{ for some } \overline{\eta}_k\} \cap \Theta$$

*then $\widetilde{\psi}_k(\theta) = \psi(\theta)$.*

*Further $TB_k(\theta) = \widetilde{\psi}_k(\theta) - \psi(\theta) = E_\theta \left[ \left\{ \widetilde{B}(\theta) - B \right\} \left\{ \widetilde{P}(\theta) - P \right\} H_1 \right].$*

*Proof.* Immediate from Theorem 3.2 and Lemma 3.3. □

We know from the above Theorem that $TB_k(\theta) = 0$ for $\theta \in \Theta_{sub,k}$. However to control the truncation bias in forming confidence intervals for $\psi(\theta)$ we will need to know how fast $\sup_{\theta \in \Theta} \{TB_k(\theta)\}$ decreases as $k$ increases. The following theorem is a key step towards determining an upper bound.

**Theorem 3.8.** *Suppose $\dot{b}(X)$ and $\dot{p}(X)$ are chosen so that $\dot{B}\dot{P}E[H_1|X] \geq 0$ w.p.1. Let*

$$Q \equiv q(X) = \left\{ \dot{B}\dot{P}E[H_1|X] \right\}^{1/2}$$

*and $\Pi\left[h(Z)|Q\overline{Z}_k\right]$ and $\Pi^\perp\left[h(X)|Q\overline{Z}_k\right]$ be, respectively, the projection in $L_2(F_X(x))$ of $h(X)$ on the $k$ dimensional linear subspace $lin\{Q\overline{Z}_k\}$ spanned by the components of the vector $Q\overline{Z}_k = q(X)\overline{z}_k(X)$ and the projection on the ortho-complement of this subspace. Then if $Ai) - Aiii)$ are satisfied,*

$$TB_k = E\left[\Pi^\perp\left[\left(\frac{P - \widehat{P}}{\dot{P}}\right)Q|Q\overline{Z}_k\right]\Pi^\perp\left[\left(\frac{B - \widehat{B}}{\dot{B}}\right)Q|Q\overline{Z}_k\right]\right].$$

**Remark 3.9.** To simplify various formulae it is often convenient and aesthetically pleasing to have $\widehat{Q} = 1$. We can choose $\dot{B}$ and $\dot{P}$ to guarantee $\widehat{Q} = 1$ w.p.1. For the functional $\psi(\theta) = E_\theta[b(X)p(X)]$ of Example 1a, $H_1 = -1$ w.p.1. Thus choosing $\dot{B}$ and $\dot{P}$ equal to 1 and $-1$, respectively, w.p.1 makes $\widehat{Q} = 1$ w.p.1. In the missing data Example 2a, the function $H_1 = -A$ so $\widehat{E}[H_1|X] = 1/\widehat{P}$ and thus the choice $\dot{B} = -1, \dot{P} = \widehat{P}$ makes $\widehat{Q} = 1$ w.p.1. Note since inference on $\psi(\theta)$ is conditional on the training sample data, we view the initial estimator $\widehat{p}(\cdot)$ of $p(\cdot)$ from the training sample as known and thus an analyst is free to choose $\dot{P}$ to be $\widehat{P}$.

**Examples continued.** In Example 1a, recall $\psi = E[BP]$. Choose $\dot{B} = -\dot{P} = 1$ w.p.1 so $\widehat{Q} = Q = 1$, and take $\widehat{B} \in lin\{\overline{Z}_k\}$. Then

$$\widetilde{B} = \widehat{B} + \Pi\left[\left(B - \widehat{B}\right)|\overline{Z}_k\right] = \Pi\left[B|\overline{Z}_k\right],$$
$$\widetilde{P} = \Pi\left[P|\overline{Z}_k\right],$$
$$TB_k = E\left\{\left[\Pi^\perp\left[B|\overline{Z}_k\right]\Pi^\perp\left[P|\overline{Z}_k\right]\right]\right\},$$
$$\widetilde{\psi}_k = \psi - TB_k = E\left\{\Pi\left[B|\overline{Z}_k\right]\Pi\left[P|\overline{Z}_k\right]\right\}.$$

Thus $\widetilde{\psi}_k$ appears to be the natural choice for a truncated parameter.

In Example 2a with $\psi = E[B], \dot{B} = -1, \dot{P} = \widehat{P} = 1/\widehat{\pi}, \widehat{Q} = 1, Q = \left\{\widehat{P}/P\right\}^{1/2} = \left\{\frac{\pi}{\widehat{\pi}}\right\}^{1/2}, \widehat{\pi} \equiv \widehat{\pi}(X), \pi \equiv \pi(X)$, we obtain

$$TB_k = E\left[\begin{array}{c}\Pi^\perp\left[\widehat{\pi}\left(\frac{1}{\pi} - \frac{1}{\widehat{\pi}}\right)\left\{\frac{\pi}{\widehat{\pi}}\right\}^{1/2} | \left\{\frac{\pi}{\widehat{\pi}}\right\}^{1/2}\overline{Z}_k\right] \\ \times \Pi^\perp\left[\left\{\frac{\pi}{\widehat{\pi}}\right\}^{1/2}\left(B - \widehat{B}\right) | \left\{\frac{\pi}{\widehat{\pi}}\right\}^{1/2}\overline{Z}_k\right]\end{array}\right].$$

Thus the truncated parameter $\widetilde{\psi}_k = \psi - TB_k$ does not seem to be a particular natural or obvious choice. The complexity of $\widetilde{\psi}_k$ is not simply due to the fact that we chose $\dot{P} = \widehat{P}$ rather than $\dot{P} = 1$ as we now demonstrate.

In Example 2a with $\dot{B} = -1, \dot{P} = 1, \widehat{Q} = \widehat{\pi}^{1/2}, Q = \pi^{1/2}$,

$$TB_k = E \begin{bmatrix} \Pi^\perp \left[ \left(\frac{1}{\pi} - \frac{1}{\widehat{\pi}}\right) \pi^{1/2} |\pi^{1/2} \overline{Z}_k \right] \times \\ \Pi^\perp \left[ \left\{\frac{\pi}{\widehat{\pi}}\right\}^{1/2} \left(B - \widehat{B}\right) |\pi^{1/2} \overline{Z}_k \right] \end{bmatrix}.$$

Nonetheless we will see that, for either choice of $\left(\dot{B}, \dot{P}\right)$, the parameter $\widetilde{\psi}_k$ will result in estimators with good properties.

**Remark 3.10.** Henceforth, given $(\beta_p, \beta_b, \beta_g)$, $\{\varphi_l(X), l = 1, 2, \ldots\}$ will always denote a complete orthonormal basis with respect to Lebesgue measure in $R^d$ or in the unit cube in $R^d$ that provides optimal rate approximation for Hölder balls $H(\beta^*, C), \beta^* \leq \lceil \max(\beta_p, \beta_b, \beta_g) \rceil$, i.e.

$$(3.23) \quad \sup_{h \in H(\beta^*, C)} \inf_{\varsigma_l} \int_{R^d} \left( h(x) - \sum_{l=1}^k \varsigma_l \varphi_l(x) \right)^2 dx = O\left(k^{-2\beta^*/d}\right).$$

The basis consisting of $d$-fold tensor products of univariate orthonormal polynomials satisfies (3.23) for all $\beta^*$. The basis consisting of $d$-fold tensor products of a univariate Daubechies compact wavelet basis with mother wavelet $\varphi_w(u)$ satisfying

$$\int_{R^1} u^m \varphi_w(u) \, du = 0, m = 0, 1, \ldots, M$$

also satisfies (3.23) for $\beta^* < M + 1$.

**Theorem 3.11.** *Suppose that* (Ai)–(Aiv) *are satisfied, that* $\dot{b}(X)$ *and* $\dot{p}(X)$ *satisfy* (3.14) – (3.16) *and that we take*

$$(3.24) \quad \overline{z}_k(X) = \overline{\varphi}_{k,\widehat{f}}(X) \left\{ \widehat{E}[H_1|X] \dot{b}(X) \dot{p}(X) \right\}^{-1/2},$$

*where*

$$\overline{\varphi}_{k,\widehat{f}}(X) \equiv \left\{ \widehat{E}\left[\overline{\varphi}_k(X) \overline{\varphi}_k(X)^T\right] \right\}^{-1/2} \overline{\varphi}_k(X).$$

*Then*

$$\sup_{\theta \in \Theta} \left\{ TB_k^2(\theta) \right\} = O_p\left(k^{-2(\beta_b + \beta_p)/d}\right)$$

**Remark 3.12.** Note if we have chosen $\dot{b}(\cdot)$ and $\dot{p}(\cdot)$ so that $\widehat{Q} = 1 \; wp1$ then $\overline{z}_k(X) = \overline{\varphi}_{k,\widehat{f}}(X)$ simplifies. The preceeding theorem does not hold with $\overline{\varphi}_k(X) / \left\{\widehat{f}(X)\right\}^{1/2}$ replacing $\overline{\varphi}_{k,\widehat{f}}(X)$ unless $f_X$ is known to lie in a Holder ball with exponent exceeding $\max\{\beta_b, \beta_p\}$.

*3.2.2. Derivation of the higher order influence functions of the truncated parameter*

We begin by proving that the first order influence functions of $\widetilde{\psi}_k$ and $\psi$ are identical except with $\widetilde{b}(\theta), \widetilde{p}(\theta), \widetilde{\psi}_k(\theta)$ replacing $b, p, \psi(\theta)$.

**Theorem 3.13.**
$$\mathbb{IF}_{1,\widetilde{\psi}_k}(\theta) = \mathbb{V}\left[IF_{1,\widetilde{\psi}_k,i_1}(\theta)\right]$$

with
$$IF_{1,\widetilde{\psi}_k}(\theta) = H\left(\widetilde{b}(\theta), \widetilde{p}(\theta)\right) - \widetilde{\psi}_k(\theta).$$

*Proof.* Since $\widetilde{\psi}_k(\theta) = E_\theta\left[H\left(\widetilde{b}(\theta), \widetilde{p}(\theta)\right)\right]$,

$$\begin{aligned}
IF_{1,\widetilde{\psi}_k}(\theta) &= H\left(\widetilde{b}(\theta), \widetilde{p}(\theta)\right) - \widetilde{\psi}_k(\theta) \\
&+ E\left[\partial H\left(b^*\left(X, \widetilde{\overline{\eta}}_k(\theta)\right), p^*\left(X, \widetilde{\overline{\alpha}}_k(\theta)\right)\right)/\partial \overline{\eta}_k^T\right] IF_{1,\widetilde{\overline{\eta}}_k(\cdot)}(\theta) \\
&+ E\left[H\left(b^*\left(X, \widetilde{\overline{\eta}}_k(\theta)\right), p^*\left(X, \widetilde{\overline{\alpha}}_k(\theta)\right)\right)/\partial \overline{\alpha}_k^T\right] IF_{1,\widetilde{\overline{\alpha}}_k(\cdot)}(\theta).
\end{aligned}$$

But, by definition of $\widetilde{\overline{\eta}}_k(\theta)$ and $\widetilde{\overline{\alpha}}_k(\theta)$, both expectations are zero. $\square$

Note that $\widetilde{\overline{\eta}}_k(\theta)$ and $\widetilde{\overline{\alpha}}_k(\theta)$ are not maximizers of the expected log-likelihood for $\overline{\alpha}_k$ and $\overline{\eta}_k$. This choice was deliberate. Had we defined $\widetilde{\overline{\eta}}_k(\theta)$ and $\widetilde{\overline{\alpha}}_k(\theta)$ as the maximizers of the expected log-likelihood, then $\mathbb{IF}_{1,\widetilde{\psi}_k}(\theta)$ would have had additional terms since the expectations in the preceding proof would not be zero. The existence of these extra terms would translate to many extra terms in $\mathbb{IF}_{m,\widetilde{\psi}_k}(\theta)$ for large $m$ leading to computational difficulties. Similarly had we chosen models $p^*(X, \overline{\alpha}_k) \equiv \Phi\left(\widehat{P} + \dot{P}\overline{\alpha}_k^T \overline{Z}_k\right)$ and $b^*(X, \overline{\eta}_k) = \Phi\left(\widehat{B} + \dot{B}\overline{\eta}_k^T \overline{Z}_k\right)$ with $\Phi(\cdot)$ a non-linear inverse-link function, $\mathbb{IF}_{m,\widetilde{\psi}_k}(\theta)$ would also have had many extra terms without an improvement in the rate of convergence.

The following is proved in the Appendix.

**Theorem 3.14.** $\mathbb{IF}_{m,\widetilde{\psi}_k} = \mathbb{IF}_{1,\widetilde{\psi}_k} + \sum_{j=2}^m \mathbb{IF}_{jj,\widetilde{\psi}_k}$ where $\mathbb{IF}_{jj,\widetilde{\psi}_k} = \mathbb{V}\left[IF_{jj,\widetilde{\psi}_k,\overline{i}_j}\right]$ is a $j$th order degenerate $U-$statistic given by

$$IF_{22,\widetilde{\psi}_k,\overline{i}_2} = -\left\{\begin{array}{c}\left[\left(H_1\widetilde{P} + H_2\right)\dot{B}\overline{Z}_k^T\right]_{i_1}\left\{E\left[\dot{P}\dot{B}H_1\overline{Z}_k\overline{Z}_k^T\right]\right\}^{-1} \\ \times \left[\overline{Z}_k\left(H_1\widetilde{B} + H_3\right)\dot{P}\right]_{i_2}\end{array}\right\},$$

$$\begin{aligned}
IF_{jj,\widetilde{\psi}_k,\overline{i}_j} &= (-1)^{j-1}\left[\left(H_1\widetilde{P} + H_2\right)\dot{B}\overline{Z}_k^T\right]_{i_1} \\
&\times \left[\begin{array}{c}\prod_{s=3}^j\left\{E\left[\dot{P}\dot{B}H_1\overline{Z}_k\overline{Z}_k^T\right]\right\}^{-1} \\ \left\{\left(\dot{P}\dot{B}H_1\overline{Z}_k\overline{Z}_k^T\right)_{i_s} - E\left[\dot{P}\dot{B}H_1\overline{Z}_k\overline{Z}_k^T\right]\right\}\end{array}\right] \\
&\times \left\{E\left[\dot{P}\dot{B}H_1\overline{Z}_k\overline{Z}_k^T\right]\right\}^{-1}\left[\overline{Z}_k\left(H_1\widetilde{B} + H_3\right)\dot{P}\right]_{i_2}.
\end{aligned}$$

*3.2.3. The Estimator $\widehat{\psi}_{m,k} \equiv \widehat{\psi} + \widehat{\mathbb{IF}}_{m,\widetilde{\psi}_k}$ and its Estimation Bias*

We can now calculate the estimator $\widehat{\psi}_{m,k} \equiv \widehat{\psi} + \widehat{\mathbb{IF}}_{m,\widetilde{\psi}_k}$ by substitution of $\widehat{\theta}$ for $\theta$ in $\widehat{\mathbb{IF}}_{m,\widetilde{\psi}_k} \equiv \mathbb{IF}_{m,\widetilde{\psi}_k}(\theta)$ to obtain the following.

**Theorem 3.15.** *Suppose* (3.24) *holds and define* $\widehat{\varsigma}(X) = \widehat{E}[H_1|X]$. *Then* $\widehat{\psi}_{m,k} = \widehat{\psi} + \widehat{\mathbb{IF}}_{1,\widetilde{\psi}_k} + \sum_{j=2}^{m} \widehat{\mathbb{IF}}_{jj,\widetilde{\psi}_k}$ *where*

$$\widehat{\psi} + \widehat{\mathbb{IF}}_{1,\widetilde{\psi}_k} = \widehat{B}\widehat{P}H_1 + \widehat{B}H_2 + \widehat{P}H_3 + H_4,$$

$$\widehat{IF}_{22,\widetilde{\psi}_k,\overline{i}_2} = -\left[\left(H_1\widehat{P} + H_2\right)\dot{B}\overline{Z}_k^T\right]_{i_1}\left[\overline{Z}_k\left(H_1\widehat{B} + H_3\right)\dot{P}\right]_{i_2}$$

$$= -\left\{ \begin{array}{l} \left[\left(H_1\widehat{P} + H_2\right)\left(\frac{\dot{B}}{\dot{P}}\right)^{1/2}\frac{\overline{\varphi}_{k,\widehat{f}}^T(X)}{\widehat{\varsigma}(X)^{1/2}}\right]_{i_1} \\ \times \left[\left(H_1\widehat{B} + H_3\right)\left\{\frac{\dot{B}}{\dot{P}}\right\}^{-1/2}\frac{\overline{\varphi}_{k,\widehat{f}}(X)}{\widehat{\varsigma}(X)^{1/2}}\right]_{i_2} \end{array} \right\},$$

$$\widehat{IF}_{jj,\widetilde{\psi}_k,\overline{i}_j} = (-1)^{j-1}\left\{\left[\left(H_1\widehat{P} + H_2\right)\dot{B}\overline{Z}_k^T\right]_{i_1}\left[\prod_{s=3}^{j}\left\{\begin{array}{c}\left(\dot{P}\dot{B}H_1\overline{Z}_k\overline{Z}_k^T\right)_{i_s} \\ -I_{k\times k}\end{array}\right\}\right] \right. \\ \left. \times \left[\overline{Z}_k\left(H_1\widehat{B} + H_3\right)\dot{P}\right]_{i_2}\right\}$$

$$= (-1)^{j-1}\left\{\begin{array}{l}\left[\left(H_1\widehat{P} + H_2\right)\left(\frac{\dot{B}}{\dot{P}}\right)^{1/2}\frac{\overline{\varphi}_{k,\widehat{f}}^T(X)}{\widehat{\varsigma}(X)^{1/2}}\right]_{i_1} \\ \times \left[\prod_{s=3}^{j}\left\{\begin{array}{c}\left(H_1\frac{\overline{\varphi}_{k,\widehat{f}}(X)\overline{\varphi}_{k,\widehat{f}}^T(X)}{\widehat{\varsigma}(X)}\right)_{i_s} \\ -I_{k\times k}\end{array}\right\}\right] \\ \times \left[\left(H_1\widehat{B} + H_3\right)\left\{\frac{\dot{B}}{\dot{P}}\right\}^{-1/2}\frac{\overline{\varphi}_{k,\widehat{f}}(X)}{\widehat{\varsigma}(X)^{1/2}}\right]_{i_2}\end{array}\right\}.$$

*Proof.* By Lemma 3.3 $E_{\widehat{\theta}}\left[\left\{H_1\widehat{B} + H_3\right\}\dot{P}\overline{Z}_k\right] = E_{\widehat{\theta}}\left[\left\{H_1\widehat{P} + H_2\right\}\dot{B}\overline{Z}_k\right] = 0$. Thus by Eqs. (3.21) and (3.22) $\widetilde{\overline{\eta}}_k\left(\widehat{\theta}\right) = \widetilde{\overline{\alpha}}_k\left(\widehat{\theta}\right) = 0$ so $\widetilde{B}\left(\widehat{\theta}\right) = \widehat{B}$ and $\widetilde{P}\left(\widehat{\theta}\right) = \widehat{P}$. Further, by Eq. (3.24), $\widehat{E}\left[\dot{P}\dot{B}H_1\overline{Z}_k\overline{Z}_k^T\right] = \widehat{E}\left[\dot{P}\dot{B}\widehat{E}\left[H_1|X\right]\overline{Z}_k\overline{Z}_k^T\right] = \widehat{E}\left[\widehat{Q}^2\overline{Z}_k\times\overline{Z}_k^T\right] = I_{k\times k}$. □

**Remark 3.16.** The reader can easily check that when $\dot{B} = \dot{P} = 1$ and $H_1 \geq 0$ wp 1, $\widehat{IF}_{22,\widetilde{\psi}_k,\overline{i}_2}$ is precisely the same as $IF_{2,2,\psi,i_1,i_2}^{(k)}\left(\widehat{\theta}\right)$ of equation (3.12) in Section 3.2.1. By expanding the product $\prod_{s=3}^{m}\left\{\left(H_1\frac{\overline{\varphi}_{k,\widehat{f}}\overline{\varphi}_{k,\widehat{f}}^T}{\widehat{\varsigma}(Z)}\right)_{i_s} - I_{k\times k}\right\}$, the equality of $\widehat{IF}_{mm,\widetilde{\psi}_k,\overline{i}_m}$ and $IF_{m,m,\psi,\overline{i}_m}^{(k)}\left(\widehat{\theta}\right)$ can be demonstrated for all $m$. Note that $\widehat{IF}_{jj,\widetilde{\psi}_k,\overline{i}_j}$ depends on $\dot{P}$ and $\dot{B}$ only through their ratio.

**Example 1a** (Continued). $\psi = E[BP], \dot{B} = -\dot{P} = 1$ w.p.1, $\widehat{Q} = 1, H_1 = -1, \overline{Z}_k = \overline{\varphi}_{k,\widehat{f}}(X)$ so $\widehat{E}\left[\overline{Z}_k\overline{Z}_{ki_s}^T\right] = I_{k\times k}$. Then

$$\dot{P}\left\{H_1\widehat{B} + H_3\right\} = -\left(Y - \widehat{B}\right), \dot{B}\left\{H_1\widehat{P} + H_2\right\} = A - \widehat{P}$$

and thus

$$\widehat{IF}_{22,\widetilde{\psi}_k,\bar{i}_2} = \left[\left(A - \widehat{P}\right)\overline{Z}_k^T\right]_{i_1} \left[\overline{Z}_k\left(Y - \widehat{B}\right)\right]_{i_2},$$

$$\widehat{IF}_{jj,\widetilde{\psi}_k,\bar{i}_j} = (-1)^j \left[\left(A - \widehat{P}\right)\overline{Z}_k^T\right]_{i_1} \left[\prod_{s=3}^{j}\left\{\overline{Z}_{ki_s}\overline{Z}_{ki_s}^T - I_{k\times k}\right\}\right] \left[\overline{Z}_k\left(Y - \widehat{B}\right)\right]_{i_2}.$$

**Example 2a** (Continued). $H_1 = -A$, $\dot{B} = -1$, $\dot{P} = \widehat{P} = 1/\widehat{\pi}$, $\widehat{Q} = 1$, $\psi = E[B]$, $Q = \left\{\widehat{P}/P\right\}^{1/2} = \left\{\frac{\pi}{\widehat{\pi}}\right\}^{1/2}$ and $\overline{Z}_k = \overline{\varphi}_{k,\widehat{f}}(X)$, so $\widehat{E}\left[\overline{Z}_k\overline{Z}_{ki_s}^T\right] = I_{k\times k}$. Then $\dot{P}\left\{H_1\widehat{B} + H_3\right\} = \frac{A}{\widehat{\pi}}\left(Y - \widehat{B}\right)$, $\dot{B}\left\{H_1\widehat{P} + H_2\right\} = \left(\frac{A}{\widehat{\pi}} - 1\right)$, so

$$\widehat{IF}_{22,\widetilde{\psi}_k,\bar{i}_2} = -\left[\left(\frac{A}{\widehat{\pi}} - 1\right)\overline{Z}_k^T\right]_{i_1} \left[\overline{Z}_k\frac{A}{\widehat{\pi}}\left(Y - \widehat{B}\right)\right]_{i_2},$$

$$\widehat{IF}_{jj,\widetilde{\psi}_k,\bar{i}_j} = (-1)^{j-1}\left[\left(\frac{A}{\widehat{\pi}} - 1\right)\overline{Z}_k^T\right]_{i_1}$$
$$\times \left[\prod_{s=3}^{j}\left\{\frac{A}{\widehat{\pi}}\overline{Z}_k\overline{Z}_k^T - I_{k\times k}\right\}_{i_s}\right] \left[\overline{Z}_k\frac{A}{\widehat{\pi}}\left(Y - \widehat{B}\right)\right]_{i_2}.$$

Consider Example 2a with $\dot{B} = -1, \dot{P} = 1, \widehat{Q} = \widehat{\pi}^{1/2}$, and $\overline{Z}_k = \overline{\varphi}_{k,\widehat{f}}(X)$, $\widehat{E}\left[\widehat{Q}^2\overline{Z}_k\overline{Z}_{ki_s}^T\right] = I_{k\times k}$, $\dot{P}\left\{H_1\widehat{B} + H_3\right\} = \left[A\left(Y - \widehat{B}\right)\right]$, $\dot{B}\left\{H_1\widehat{P} + H_2\right\} = \left(\frac{A}{\widehat{\pi}} - 1\right)$, so

$$\widehat{IF}_{22,\widetilde{\psi}_k,\bar{i}_2} = -\left[\left(\frac{A}{\widehat{\pi}} - 1\right)\overline{Z}_k^T\right]_{i_1} \left[\overline{Z}_kA\left(Y - \widehat{B}\right)\right]_{i_2},$$

$$\widehat{IF}_{jj,\widetilde{\psi}_k,\bar{i}_j} = (-1)^{j-1}\left[\left(\frac{A}{\widehat{\pi}} - 1\right)\overline{Z}_k^T\right]_{i_1}$$
$$\times \left[\prod_{s=3}^{j}\left\{A\overline{Z}_k\overline{Z}_k^T - I_{k\times k}\right\}_{i_s}\right] \left[\overline{Z}_kA\left(Y - \widehat{B}\right)\right]_{i_2}.$$

Our next theorem, proved in the Appendix of our technical report, derives the estimation bias $EB_m = E\left[\widehat{\psi}_{m,k}\right] - \widetilde{\psi}_k$.

**Theorem 3.17.** *Suppose* (3.14)–(3.16) *and* (Ai)–(Aiv) *hold then*

$$(3.25) \quad EB_m = (-1)^{m-1}\left\{\begin{array}{l} E\left[Q^2\left(\frac{B-\widehat{B}}{\dot{B}}\right)\overline{Z}_k^T\right]\left\{E\left[Q^2\overline{Z}_k\overline{Z}_k^T\right] - I_{k\times k}\right\}^{m-1} \\ \times \left\{E\left[Q^2\overline{Z}_k\overline{Z}_k^T\right]\right\}^{-1} E\left[\overline{Z}_kQ^2\left(\frac{P-\widehat{P}}{\dot{P}}\right)\right] \end{array}\right\}$$

$$|EB_m|$$

$$(3.26) \quad \leq \left\{\begin{array}{l} \left\|\left\{\frac{\dot{B}}{\dot{P}}G\right\}^{1/2}\right\|_\infty \left\|\left\{\frac{\dot{P}}{\dot{B}}G\right\}^{1/2}\right\|_\infty \|\delta g\|_\infty^{m-1}(1 + o_p(1)) \times \\ \left\{\int (p(X) - \widehat{p}(X))^2 dX\right\}^{1/2} \left\{\int \left(b(X) - \widehat{b}(X)\right)^2 dX\right\}^{1/2} \end{array}\right\}$$

$$(3.27) \quad = O_P\left(\left(\frac{\log n}{n}\right)^{\frac{(m-1)\beta_g}{d+2\beta_g}} n^{-\left(\frac{\beta_b}{d+2\beta_b} + \frac{\beta_p}{d+2\beta_p}\right)}\right)$$

for $m \geq 1$, where $\delta g = \frac{g(X) - \widehat{g}(X)}{\widehat{g}(X)}$.

**Remark 3.18.** At the cost of a longer proof we could have used Hölder's inequality repeatedly to control $\delta g$ in the $L_p$ norm $||\delta g||_{m+1}$ with $p = m+1$ to show that $|EB_m| = O_P\left(||\delta g||_{m+1}^{m-1} \left\|b(\cdot) - \widehat{b}(\cdot)\right\|_{m+1} ||p(\cdot) - \widehat{p}(\cdot)||_{m+1}\right)$. Thus, $|EB_m|$ is $O_P\left(\left\|\theta - \widehat{\theta}\right\|^{m+1}\right)$, consistent with the form of the bias given in our fundamental Theorem 2.2.

**Remark 3.19** (An alternate derivation of $\widehat{\psi}_{m,k}$). The above derivation of $\widehat{\psi}_{m,k}$ required that one have facility in calculating higher order influence functions $\mathbb{IF}_{m,\widetilde{\psi}_k}$, as done in the proof of Theorem 3.14 in the Appendix. However, there exists an alternate derivation of $\widehat{\psi}_{m,k}$ that does not require one learn how to calculate influence functions. Specifically, we know from Theorems 2.2 and 2.3 that in a (locally) nonparametric model $\widehat{\mathbb{IF}}_{jj,\widetilde{\psi}_k}, j \geq 2$ is the unique $j^{th}$ order U-statistic that is degenerate under $\widehat{\theta}$ and satisfies

$$(3.28) \qquad EB_{j-1} + E\left[\widehat{\mathbb{IF}}_{jj,\widetilde{\psi}_k}|\widehat{\theta}\right] \equiv EB_j = O_p\left(\left\|\widehat{\theta} - \theta\right\|^{j+1}\right)$$

with $EB_1 = E\left[\widehat{\psi}_1|\widehat{\theta}\right] - \widetilde{\psi}_k$. In fact, we first derived $\widehat{\psi}_{m,k}$ by beginning with $\widehat{\psi}_1 = \widehat{\psi} + \widehat{\mathbb{IF}}_{1,\widetilde{\psi}_k}$, calculating $EB_1 = E\left[\widehat{\psi}_1|\widehat{\theta}\right] - \widetilde{\psi}_k$, and then, recursively for $j = 2, \ldots$ finding $\widehat{\mathbb{IF}}_{jj,\widetilde{\psi}_k}$ satisfying the above equation. In fact if one did not even know how to derive $\mathbb{IF}_{1,\widetilde{\psi}_k}$, one could begin the recursion by obtaining $\widehat{\mathbb{IF}}_{1,\widetilde{\psi}_k}$ as the unique first order U-statistic with mean zero under $\widehat{\theta}$ satisfying $\widehat{\psi} - \widetilde{\psi}_k + E\left[\widehat{\mathbb{IF}}_{1,\widetilde{\psi}_k}|\widehat{\theta}\right] = O_p\left(\left\|\widehat{\theta} - \theta\right\|^2\right)$.

### 3.2.4. The variance of $\widehat{\psi}_{m,k} \equiv \widehat{\psi} + \widehat{\mathbb{IF}}_{m,\widetilde{\psi}_k}$ using compact wavelets

In this section, we derive the order of the variance of $\widehat{\psi}_{m,k}$ when the orthonormal system $\{\varphi_j(X)\}$ used to construct our U-statistics are a compact wavelet basis. First consider the case where $X$ is univariate; without loss of generality, assume that $X \sim Uniform[0,1]$. Because we are primarily interested in convergence rates, the fact that $X$ may not follow the uniform distribution will not affect the rate results given below, but can influence the size of the constants. We use $\phi_j(X)$ in place of $\varphi_j(X)$ to indicate univariate basis functions.

Let $k^*$, $k$ be integer powers of two with $k > k^*$. Denote by $\overline{\phi}(X) \equiv \overline{\phi}_1^k(X)$ the $k-$dimensional basis vector whose first $k^*$ components $\overline{\phi}_1^{k^*}(X)$ are the $k^*-$vector of level $\log_2 k^*$ scaled and translated versions of a compactly supported 'father' wavelet (Mallat [10]) and whose last $k - k^*$ components $\overline{\phi}_{k^*+1}^{k}(X)$ are the associated compact mother wavelets between levels $\log_2 k^*$ and $\log_2 k$. In particular, one may use periodic wavelets, folded wavelets or Daubechies' boundary wavelets with enough vanishing moments to obtain the optimal approximation rate of $O\left(k^{-2\beta/d}\right)$ for $\beta = \max(\beta_g, \beta_p, \beta_b)$. The multiresolution analysis (MRA) property of wavelets allows us to decompose the vector space spanned by the $\log_2(k)$-level father wavelets $\mathcal{V}_{\log_2(k)}$ into the direct sum of the subspace spanned by $\log_2(k^*)$-level father wavelets $\mathcal{V}_{\log_2(k^*)} = \left\{ a^T \overline{\phi}_1^{k^*}(X) : a \in R^{k^*} \right\}$ and the span of mother wavelets for each level between $\log_2(k^*)$ and $\log_2(k) - 1$ which we respectively write as

$$\mathcal{W}_{\log_2(k^*)} = \left\{ a^T \overline{\phi}_{k^*+1}^{2k^*}(X) : a \in R^{k^*} \right\},$$

$$\mathcal{W}_{\log_2(k_0)+1} = \left\{ a^T \overline{\phi}_{2k^*+1}^{4k^*}(X) : a \in R^{2k^*} \right\},$$

$$\vdots$$

$$\mathcal{W}_{\log_2(k)-1} = \left\{ a^T \overline{\phi}_{\frac{k}{2}+1}^{k}(X) : a \in R^{\frac{k}{2}} \right\}.$$

Then for any integer $s$ with $\log_2(k^*) + 1 \leq s$, we have

$$\mathcal{V}_s = \mathcal{V}_{\log_2(k^*)} \oplus \left( \bigoplus_{v=\log_2(k^*)}^{s-1} \mathcal{W}_v \right).$$

As $s \to \infty$, the resulting basis system is dense in $L_2(X)$ (Mallat [10]). Since, in fact, $X$ is $d-$dimensional we require a generalization that allows for multivariate tensor wavelet basis functions. In fact, suppose $X^T = (X^1, \ldots, X^d)$ is now multivariate, and we again assume $X \sim Uniform$ on $[0,1]^d$. Given $d$ univariate vector spaces

$$\mathcal{V}_{1,\log_2(k)}, \mathcal{V}_{2,\log_2(k)}, \ldots, \mathcal{V}_{d,\log_2(k)}$$

respectively spanned by vectors $\overline{\phi}_1^k(X^1), \overline{\phi}_1^k(X^2), \ldots, \overline{\phi}_1^k(X^d)$, so that for $1 \leq r \leq d$,

$$\mathcal{V}_{r,\log_2(k^*)} \subset \mathcal{V}_{r,\log_2(k^*)+1} \subset \ldots \subset \mathcal{V}_{r,\log_2(k)-1} \subset \mathcal{V}_{r,\log_2(k)}$$

and

$$\mathcal{V}_{r,\log_2(k)} = \mathcal{V}_{r,\log_2(k^*)} \oplus \left( \bigoplus_{v=\log_2(k^*)}^{\log_2(k)-1} \mathcal{W}_{r,v} \right).$$

One may define $d$ dimensional tensor vector spaces

$$\mathcal{Y}_{d,\log_2(k^*)}, \mathcal{Y}_{d,\log_2(k^*)+1}, \ldots, \mathcal{Y}_{d,\log_2(k)}$$

such that

$$\mathcal{Y}_{d,\log_2(k^*)} \subset \mathcal{Y}_{d,\log_2(k^*)+1} \subset \ldots \subset \mathcal{Y}_{d,\log_2(k)},$$

where for $s \geq 0$,

$$\mathcal{Y}_{d,\log_2(k_0)+s} = \bigotimes_{1 \leq r \leq d} \mathcal{V}_{r,\log_2(k_0)+s}.$$

As $s \to \infty$, the resulting tensor basis system is dense in $L_2(X)$ (Mallat [10]).

Next, suppose that we have a set of multivariate basis functions

$$\left\{ \overline{\varphi}_1^{k_j}(X), j = 0, 1, \ldots, 2m \right\}$$

such that for each $k_j$, $\overline{\varphi}_1^{k_j}(X)$ spans $\bigotimes_{1 \leq r \leq d} \mathcal{V}_{r,\log_2(k_{j,r})}$, where $\prod_{r=1}^{d} k_{j,r} = k_j$. Define $||\cdot||_2$ as the $L_2$ norm with respect to the Lebesgue measure. The following theorem is key to our derivation of the order of the variance of $\widehat{\psi}_{m,k}$.

**Theorem 3.20.** For $m \geq 0$,

$$\left\| \overline{\varphi}_1^{k_1}(X_{i_1})^T \prod_{j=1}^{m} \left\{ \overline{\varphi}_1^{k_j}(X_{i_{j+1}}) \overline{\varphi}_1^{k_{j+1}}(X_{i_{j+1}})^T \right\} \overline{\varphi}_1^{k_{m+1}}(X_{i_{m+2}}) \right\|_2^2$$

$$= E \left( \prod_{j=1}^{m+1} K_{(1,k_j)}(X_{i_j}, X_{i_{j+1}}) \right)^2$$

$$\asymp \prod_{j=1}^{m+1} k_j.$$

The following theorem is an immediate consequence of Theorem 3.20 obtained by taking $k_j = k^* = k$ (which implies we use the father wavelets at level $\log_2(k)$ but no mother wavelets.)

**Theorem 3.21.** For all $\theta \in \Theta$,

$$var_\theta \left[ \mathbb{IF}_{1,\widetilde{\psi}_k}(\theta) \right] \asymp \frac{1}{n},$$

$$var_\theta \left[ \mathbb{IF}_{jj,\widetilde{\psi}_k}(\theta) \right] \asymp \left( \frac{1}{n} \max \left\{ 1, \left( \frac{k}{n} \right)^{j-1} \right\} \right),$$

$$var_\theta \left[ \mathbb{IF}_{m,\widetilde{\psi}_k}(\theta) \right] \approx var_{\widehat{\theta}} \left[ \widehat{\mathbb{IF}}_{m,\widetilde{\psi}_k} | \widehat{\theta} \right] \asymp \frac{1}{n} \max \left\{ 1, \left( \frac{k}{n} \right)^{m-1} \right\}.$$

We now use Theorem 3.21 to derive the order of the conditional variance of $\widehat{\psi}_{m,k}$ given $\widehat{\theta}$.

**Theorem 3.22.** If $\sup_{o \in \mathcal{O}} \left| f(o; \widehat{\theta}) - f(o; \theta) \right| \to 0$ as $\|\widehat{\theta} - \theta\| \to 0$, then for a fixed $m$,

$$var_\theta \left[ \widehat{\psi}_{m,\widetilde{\psi}_k} | \widehat{\theta} \right] = var_\theta \left[ \widehat{\mathbb{IF}}_{m,\widetilde{\psi}_k} | \widehat{\theta} \right]$$

$$= var_{\widehat{\theta}} \left[ \widehat{\mathbb{IF}}_{m,\widetilde{\psi}_k} | \widehat{\theta} \right] (1 + o_p(1))$$

$$\asymp \left( \frac{1}{n} \max \left\{ 1, \left( \frac{k}{n} \right)^{m-1} \right\} \right).$$

The proof is in our technical report.

For a given $m$, the estimator $\widehat{\psi}_{m,k_{opt}(m)}$ that minimizes the maximum asymptotic MSE over the model $\mathcal{M}(\Theta)$ defined by (Ai)–(Aiv) among the candidates $\widehat{\psi}_{m,k}$ uses the value $k_{opt}(m) \equiv k_{opt}(m,n)$ of $k$ that equates the order $\frac{1}{n} \max \left\{ 1, \left( \frac{k}{n} \right)^{m-1} \right\}$ of $var \left[ \widehat{\psi}_{m,\widetilde{\psi}_k} | \widehat{\theta} \right]$ to the order

$$\max \left[ \{TB_k\}^2, \{EB_m(\theta)\}^2 \right] =$$

$$\max \left[ \begin{array}{c} \left( \frac{\log n}{n} \right)^{\frac{2(m-1)\beta_g}{d+2\beta_g}} n^{-\left( \frac{2\beta_b}{d+2\beta_b} + \frac{2\beta_p}{d+2\beta_p} \right)}, \\ k^{-\frac{2(\beta_b+\beta_p)}{d}} \end{array} \right]$$

of the maximal squared bias. The estimator $\widehat{\psi}_{m_{opt},k_{opt}} \equiv \widehat{\psi}_{m_{opt},k_{opt}(m_{opt})}$ that minimizes the maximum asymptotic MSE over the model $\mathcal{M}(\Theta)$ among all candidates $\widehat{\psi}_{m,k}$ is the estimator $\widehat{\psi}_{m,k_{opt}(m,n)}$ which minimizes $\frac{1}{n}\max\left(1, \left(\frac{k_{opt}(m,n)}{n}\right)^{m-1}\right)$.

*3.2.5. Distribution theory and confidence interval construction*

We derive a consistent estimator of the variance and give the asymptotic distribution of $\widehat{\psi}_{m,k}$ for any model and functional satisfying (Ai)–(Aiv). Let $z_\alpha$ be the upper $\alpha$-quantile of a standard normal distribution, i.e. a $N(0,1)$.

**Theorem 3.23.** (a) *Let* $\widehat{\mathbb{W}}^2_{1,\widetilde{\psi}_k} = n^{-1}\mathbb{V}\left[\left\{\widehat{IF}_{1,\widetilde{\psi}_k,i_1}\right\}^2\right]$,

$$\widehat{\mathbb{W}}^2_{jj,\widetilde{\psi}_k} = \binom{n}{j}^{-1}\mathbb{V}\left[\left(\widehat{IF}^{(s)}_{j,j,\widetilde{\psi}_k(\cdot)}\right)^2\right],$$

*for* $j \geq 2$, *and*

$$\widehat{\mathbb{W}}^2_{m,\widetilde{\psi}_k} = \widehat{\mathbb{W}}^2_{1,\widetilde{\psi}_k} + \sum_{j=2}^{m}\widehat{\mathbb{W}}^2_{jj,\widetilde{\psi}_k},$$

*where* $\widehat{IF}^{(s)}_{j,j,\widetilde{\psi}_k(\cdot)}$ *is the symmetric kernel of* $\widehat{\mathbb{IF}}_{jj,\widetilde{\psi}_k(\cdot)}$. *We have,*

$$\widehat{E}\left[\widehat{\mathbb{W}}^2_{1,\widetilde{\psi}_k}\right] = \widehat{var}\left[\widehat{\mathbb{IF}}_{1,\widetilde{\psi}_k}|\widehat{\theta}\right],$$
$$\widehat{E}\left[\widehat{\mathbb{W}}^2_{jj,\widetilde{\psi}_k}\right] = \widehat{var}\left[\widehat{\mathbb{IF}}_{jj,\widetilde{\psi}_k}|\widehat{\theta}\right],$$
$$\widehat{E}\left[\widehat{\mathbb{W}}^2_{m,\widetilde{\psi}_k}\right] = \widehat{var}\left[\widehat{\mathbb{IF}}_{m,\widetilde{\psi}_k}|\widehat{\theta}\right],$$

*where* $\widehat{var}[\cdot] = var_{\widehat{\theta}}[\cdot]$.
(b) *Conditional on the training sample,*

$$\left\{\frac{1}{n}\max\left\{1,\left(\frac{k_{opt}(m,n)}{n}\right)^{m-1}\right\}\right\}^{-1/2}\left\{\widehat{\psi}_{m,k_{opt}(m,n)} - E\left[\widehat{\psi}_{m,k_{opt}(m,n)}|\widehat{\theta}\right]\right\}$$

*converges uniformly for* $\theta \in \Theta$ *to a normal distribution with finite variance as* $n \to \infty$. *The asymptotic variance is uniformly consistently estimated by*

$$\left\{\frac{1}{n}\max\left\{1,\left(\frac{k}{n}\right)^{m-1}\right\}\right\}^{-1}\widehat{\mathbb{W}}^2_{m,\widetilde{\psi}_{k_{opt}(m,n)}}.$$

*Thus*

$$\left\{\widehat{\psi}_{m,k_{opt}(m,n)} - E\left[\widehat{\psi}_{m,k_{opt}(m,n)}|\widehat{\theta}\right]\right\}/\widehat{\mathbb{W}}_{m,\widetilde{\psi}_{k_{opt}(m,n)}}$$

*is converging in distribution to a standard normal distribution.*
(c) *Define the interval* $C_{m,k} = \widehat{\psi}_{m,k} \pm z_\alpha \widehat{\mathbb{W}}_{m,\widetilde{\psi}_k}$. *Suppose* $k_{opt}(m,n) = n^{\rho_{opt}(m,n)}$. *Then for* $k^* = n^{\rho^*}$, $\rho^* > \rho_{opt}(m,n)$,

$$\sup_{\theta \in \Theta}\left[\frac{E_\theta\left[\widehat{\psi}_{2,k^*}|\widehat{\theta}\right]}{\sqrt{var_\theta\left[\widehat{\psi}_{2,k^*}|\widehat{\theta}\right]}}\right] = o_p(1)$$

and $\left\{\widehat{\psi}_{m,k^*} - \psi(\theta)\right\}/\widehat{\mathbb{W}}_{m,\widetilde{\psi}_{k^*}}$ converges uniformly in $\theta \in \Theta$ to a $N(0,1)$. Moreover, $C_{m,k^*}$ is a conservative uniform asymptotic $(1-\alpha)$ confidence interval for $\psi(\theta)$.

(d) Suppose we could derive a constant $C_{bias}$ and a constant $N^*$ such that

$$\sup_\theta \left| E_\theta \left[ \left\{ \widehat{\psi}_{m,k_{opt}(m,n)} - \psi(\theta) \right\} \right] \right|$$
$$= \sup_\theta \left| \left\{ TB_{k_{opt}(m,n)}(\theta) + EB_m(\theta) \right\} \right|$$
$$\leq C_{bias} \left\{ \frac{1}{n} \max \left\{ 1, \left( \frac{n^{\rho_{opt}(m,n)}}{n} \right)^{m-1} \right\} \right\}^{1/2}$$

for $n > N^*$. Then

$$BC_{m,k_{opt}(m,n)}$$
$$= \widehat{\psi}_{m,k_{opt}(m,n)} \pm \left\{ z_\alpha \widehat{\mathbb{W}}_{m,\widetilde{\psi}_{k^*}} + C_{bias} \left\{ \frac{1}{n} \max \left\{ 1, \left( \frac{n^{\rho_{opt}(m,n)}}{n} \right)^{m-1} \right\} \right\}^{1/2} \right\}$$

is a conservative uniform asymptotic $(1-\alpha)$ confidence interval for $\psi(\theta)$.

Part (a) of the theorem is an easy calculation. The asymptotic normality of $\widehat{\psi}_{m,k_{opt}(m,n)}$ is based on new results on the asymptotic distribution of higher order $U$-statistics with kernels depending on $n$ to be published elsewhere (Robins et al. [18]).

Part (c) of the theorem implies we obtain a conservative uniform asymptotic $(1-\alpha)$ confidence interval for any value of $\rho^*$ exceeding $\rho_{opt(m,n)}$. However, for the actual fixed sample size of our study, say $n = 5000$, there is no guarantee the interval of part (c) based on given difference $\rho^* - \rho_{opt(m,n)}$, say 0.3, will provide conservative finite sample coverage.

Because of this difficulty, a better approach, described in part (d), would be to determine a constant $C_{bias}$ that can be used to bound the maximal bias under the model at a sample sizes exceeding $N^*$, with $N^*$ no greater than the actual fixed sample size $n$ of the study. Then the interval $BC_{m,k_{opt}(m,n)}$ will be a honest conservative finite sample $1-\alpha$ confidence interval, provided that $\widehat{\psi}_{m,k_{opt}(m,n)}$ has nearly converged to its normal limit at sample size $n$. Unfortunately, as yet we do not know how to determine the constants $C_{bias}$ and $N^*$ of part (d) as a function of our model and of our initial estimator $\widehat{\theta}$. This is an important open problem.

### 3.2.6. Models of increasing dimension and multi-robustness

**A model of increasing dimension.** The previous results can also be used for the analysis of models whose dimension increases with sample size. In fact, consider the $\mathcal{M}(\Theta_{n^\eta})$, $\eta$ known, that differs from model $\mathcal{M}(\Theta)$ in that, rather than assuming $b(x)$ and $p(x)$ live in particular Hölder balls, we instead assume the working models of Eqs. (3.17) and (3.18) are precisely true for $k = n^\eta$, so $\psi(\theta) \equiv \widetilde{\psi}_{n^\eta}(\theta)$ and the dimensions of $b(x)$ and $p(x)$ increase as $n^\eta$. Valid point and interval estimation for $\widetilde{\psi}_{n^\eta}(\theta)$ can still be based on the estimators $\widehat{\psi}_{m,k}$ except now (i) there is truncation bias only when $k < n^\eta$, (ii) the variance remains of the order of $\frac{1}{n} \max\left(1, \left(\frac{k}{n}\right)^{m-1}\right)$, and (iii) the estimation and truncation bias (when it exists) orders will be determined by any additional complexity reducing restrictions placed on the fraction of non-zero components or on the rate of decay of the components of the vectors $\widetilde{\overline{\eta}}_{n^\eta}(\theta)$ and $\widetilde{\overline{\alpha}}_{n^\eta}(\theta)$, and, for estimation bias, by $\beta_g$ as well. As a consequence, $m_{opt}$ and $k_{opt}$ under model $\mathcal{M}(\Theta_{n^\eta})$ will differ from their values under model $\mathcal{M}(\Theta)$. Note we need not take $k = n^\eta$ as we did in the heuristic discussion following Remark 2.8. Indeed $\widehat{\psi}_{m_{best}}$ in that discussion corresponds to the estimator in the class $\widehat{\psi}_{m,k=n^\eta}$ with the fastest rate of convergence. In general, $\widehat{\psi}_{m_{best}}$ will have convergence rate slower than $\widehat{\psi}_{m_{opt},k_{opt}}$. Furthermore, the discussion in Section 4.1.1 implies that, when $n^\eta \gg n$ and the minimax rate for estimation of $\psi(\theta)$ is slower than $n^{-1/2}$, even $\widehat{\psi}_{m_{opt},k_{opt}}$ will typically fail to converge at the minimax rate when complexity reducing restrictions have been imposed on $\widetilde{\overline{\eta}}_{n^\eta}(\theta)$ and $\widetilde{\overline{\alpha}}_{n^\eta}(\theta)$.

**Multi-robustness and a practical data analysis strategy.** Conditional on $\widehat{\theta}$, for $m \geq 2$, $EB_m$ is zero and thus estimator $\widehat{\psi}_{m,k}$ is unbiased for $\widetilde{\psi}_k$ if $\widehat{p}(\cdot) = p(\cdot), \widehat{b}(\cdot) = b(\cdot)$, or $\widehat{g}(\cdot) = g(\cdot)$. We refer to $\widehat{\psi}_{m,k}$ as triply-robust for $\widetilde{\psi}_k$, generalizing Robins and Rotnitzky [17] and van der Laan and Robins [24] who referred to $\widehat{\psi}_1$ as doubly-robust because of its being unbiased for $\widetilde{\psi}_k$ if either $\widehat{p}(\cdot) = p(\cdot)$ or $\widehat{b}(\cdot) = b(\cdot)$. In fact, for $m \geq 3$, we can construct a modified estimator $\widehat{\psi}_{m,k}^{\text{mod}}$ that is $m+1$-fold robust as follows. Let $\widehat{g}_s(\cdot), s = 3, \ldots, m$, denote $m-2$ additional initial estimators of $g(\cdot)$ that differ from one another and from $\widehat{g}(\cdot)$. Define $\widehat{\psi}_{m,k}^{\text{mod}} = \widehat{\psi} + \widehat{\mathbb{IF}}_{1,\widetilde{\psi}_k} + \widehat{\mathbb{IF}}_{22,\widetilde{\psi}_k,\widetilde{i}_j} + \sum_{j=3}^{m} \widehat{\mathbb{IF}}_{jj,\widetilde{\psi}_k}^{\text{mod}}$, where

$$\widehat{IF}_{jj,\widetilde{\psi}_k,\widetilde{i}_j}^{\text{mod}} = (-1)^{j-1} \left[\left(H_1 \widehat{P} + H_2\right) \dot{B} \overline{Z}_k^T\right]_{i_1} \left\{\left(\dot{P} \dot{B} H_1 \overline{Z}_k \overline{Z}_k^T\right)_{i_2} - I_{k \times k}\right\}$$

$$\times \begin{bmatrix} \prod_{s=3}^{j-1} \left\{\widehat{E}_s\left[\dot{P}\dot{B}H_1\overline{Z}_k\overline{Z}_k^T\right]\right\}^{-1} \\ \left\{\left(\dot{P}\dot{B}H_1\overline{Z}_k\overline{Z}_k^T\right)_{i_s} - \widehat{E}_s\left[\dot{P}\dot{B}H_1\overline{Z}_k\overline{Z}_k^T\right]\right\} \end{bmatrix}$$

$$\times \left\{\widehat{E}_j\left[\dot{P}\dot{B}H_1\overline{Z}_k\overline{Z}_k^T\right]\right\}^{-1} \times \left[\overline{Z}_k\left(H_1\widehat{B} + H_3\right)\dot{P}\right]_{i_j}$$

with $\widehat{E}_s$ defined like $\widehat{E}$, except with $\widehat{g}_s(\cdot)$ replacing $\widehat{g}(\cdot)$. In the Appendix of our technical report, we prove that $EB_m^{\text{mod}} = E\left[\widehat{\psi}_{m,k}^{\text{mod}}\right] - \widetilde{\psi}_k$ is

$$(-1)^{m-1} \left\{ \begin{array}{c} E\left[\dot{B}\dot{P}H_1\left(\frac{P-\widehat{P}}{\dot{P}}\right)\overline{Z}_k^T\right]\left\{E\left[\dot{B}\dot{P}H_1\overline{Z}_k\overline{Z}_k^T\right] - I_{k\times k}\right\} \\ \times \prod_{s=3}^{m} \left\{\widehat{E}_s\left[\dot{B}\dot{P}H_1\overline{Z}_k\overline{Z}_k^T\right]\right\}^{-1} \\ \times \left\{E\left[\dot{B}\dot{P}H_1\overline{Z}_k\overline{Z}_k^T\right] - \widehat{E}_s\left[\dot{B}\dot{P}H_1\overline{Z}_k\overline{Z}_k^T\right]\right\} \\ \times \left\{E\left[\dot{B}\dot{P}H_1\overline{Z}_k\overline{Z}_k^T\right]\right\}^{-1} E\left[\dot{B}\dot{P}H_1\left(\frac{B-\widehat{B}}{\dot{B}}\right)\right] \end{array} \right\}$$

which is zero if $\widehat{p}(\cdot) = p(\cdot), \widehat{b}(\cdot) = b(\cdot), \widehat{g}(\cdot) = g(\cdot)$, or if any of the $m - 2$ $\widehat{g}_s(\cdot)$ equals $g(\cdot)$. (We note that if $\widehat{p}(\cdot) = p(\cdot)$ or $\widehat{b}(\cdot) = b(\cdot)$, $\psi = \widetilde{\psi}_k$ and thus $\widehat{\psi}_{m,k}^{\text{mod}}$ and $\widehat{\psi}_{m,k}$ are unbiased for $\psi$.)

In settings where the dimension $d$ of $X$ is so large (e.g., $30 - 100$) that the above asymptotic results fail as a guide to the finite sample performance of our procedures at the moderate sample sizes, say $n = 500$–$5000$, commonly found in practice, one might consider, as a practical data analysis strategy, using the $m + 1$-fold robust estimator $\widehat{\psi}_{m,k}^{\mathrm{mod}}$ with $\widehat{p}(\cdot), \widehat{b}(\cdot), \widehat{g}(\cdot)$, and the $\widehat{g}_s(\cdot)$ selected by cross-validation as in van der Laan and Dudoit [23]. Specifically, the training sample is split into two random subsamples – a candidate estimator subsample of size $n_c$ and a validation subsample of size $n_v$, where both $n_c/n$ and $n_v/n$ are bounded away from 0 as $n \to \infty$. A large number (e.g., $n^3$) candidate parametric models of various dimensions and functional forms for $p$, $b$, and $g$ are fit to the candidate estimator subsample and the validation sample is used to find the candidate estimators $\widehat{p}(\cdot)$ and $\widehat{b}(\cdot)$ for $p$ and $b$ and the $m - 1$ candidate estimators $\widehat{g}(\cdot)$ and $\widehat{g}_s(\cdot), s = 3, \ldots, m$, for $g$ with the smallest estimated risks (with respect to an appropriate risk function such as squared error or Kullback-Leibler).

In the setting of very high dimensional $X$, current practice is to use a doubly robust estimator, say $\widehat{\psi}_1$ with $\widehat{p}$ and $\widehat{b}$ selected by cross validation (van der Laan and Dudoit [23]). An $m + 1$-robust estimator $\widehat{\psi}_{m,k}^{\mathrm{mod}}$ with $k \gg n$ and $m$ rather large may be preferable to a doubly robust estimator for two reasons. First, if one uses an $m + 1$-robust estimator of $\psi$ rather than a DR estimator, it may be more likely that the estimator will have very small bias, as it is more likely that at least one of $m + 1$, rather than one of two, models is very nearly correct. Second, because $k \gg n$, nominal $1 - \alpha$ Wald confidence intervals centered at $\widehat{\psi}_{m,k}^{\mathrm{mod}}$ will be wider than the interval of length proportional to $n^{-1/2}$ centered at $\widehat{\psi}_1$. A wide interval is a more appropriate measure of the actual uncertainty about $\psi$. However, it is also the case that the bias of $\widehat{\psi}_{m,k}^{\mathrm{mod}}$ can exceed that of $\widehat{\psi}_1$ when all of the $m + 1$ models selected by cross-validation are far from correct, owing to the product structure of the estimation bias. The product of $m + 1$ errors, each greater than 1, will exceed the product of just 2 of the errors. We therefore plan to compare through simulations the finite sample performances of $\widehat{\psi}_{m,k}^{\mathrm{mod}}$ and $\widehat{\psi}_1$ in the setting of very high dimensional $X$.

## 4. Rates of convergence and minimaxity

We consider a generic version in which we only assume a model and functional satisfying $Ai) - Aiv)$. To examine efficiency issues, we first consider the estimator $\widehat{\psi}_1$ based on the first order influence function and sample splitting. Without loss of generality we assume $\beta_p \geq \beta_b$. (Otherwise simply interchange $\beta_p$ and $\beta_b$ in what follows.) It will be useful to consider the alternative parametrization

$$\beta = \frac{\beta_p + \beta_b}{2},$$

$$\Delta = \left(\frac{\beta_p}{\beta_b} - 1\right) \geq 0.$$

The (conditional) variance of $\widehat{\psi}_1$ is of the order of $1/n$ and the (conditional) bias of $\widehat{\psi}_1$ in estimating $\psi$ is $O_p\left(n^{-\left(\frac{\beta_b}{d+2\beta_b} + \frac{\beta_p}{d+2\beta_p}\right)}\right)$. If $\Delta = 0$ and thus $\beta_p = \beta_b$, the bias of $\widehat{\psi}_1$ is $n^{-\frac{2\beta}{d+2\beta}}$ and $\widehat{\psi}_1$ is not $n^{1/2}-$ consistent for $\psi$ when $\beta/d < 1/2$. At the other extreme, as $\Delta \to \infty$, i.e. $\beta_b \to 0$, the bias of $\widehat{\psi}_1$ is $n^{-\frac{2\beta}{d+4\beta}}$ which fails to be $n^{1/2}-$consistent for any finite $\beta$.

**Minimaxity with $g$ known.**    To further examine efficiency issues, it is instructive to first consider the estimation of $\psi$ with $g(\cdot)$ known. If $g(\cdot)$ were known, we could set $\widehat{g}(X) = g(X)$ when calculating $\widehat{\psi}_{m,k}$. Then $EB_2 = 0$ and $\widehat{\psi}_{2,k}$ would therefore be an unbiased estimator of $\widetilde{\psi}_k$. Letting a superscript $g$ denote the model with $g$ known, it is easy to see that $\widehat{\psi}_{m^g_{opt}, k^g_{opt}}(m^g_{opt})$ would be $\widehat{\psi}_{2,k^g_{opt}(2)}$ where $k^g_{opt}(2)$ satisfies $\max(1/n, k/n^2) \asymp \operatorname{var}\left(\widehat{\psi}_{2,k}\right) = TB_k^2 = k^{-2(\beta_b + \beta_p)/d} = k^{-4\beta/d}$. Solving this, we find that when $\beta/d$ is greater than or equal to $1/4$, we can take $k^g_{opt} = n^{\frac{1}{4\beta/d}} \leq n$ and $\left|\widehat{\psi}_{2,k^g_{opt}(2)} - \psi\right| = O_p\left(n^{-\frac{1}{2}}\right)$ regardless of $\Delta$, which is, of course, the minimax rate.

In contrast if $\beta/d < 1/4$, $k^g_{opt}(2) = n^{\frac{2}{1+4\beta/d}}$ and $\left(\widehat{\psi}_{2,k^g_{opt}(2)} - \psi\right) = n^{-\frac{4\beta}{4\beta+d}}$. In an unpublished paper, we have proved that this is the minimax rate when $g(\cdot)$ is known.

This raises the question of whether the lower bounds of rate $n^{-\frac{1}{2}}$ for $\beta/d \geq 1/4$ and/or rate $n^{-\frac{4\beta}{4\beta+d}}$ for $\beta/d < 1/4$ are still achievable when $g$ is unknown, without restrictions on the smoothness of $g$.

Before addressing this question, we take the opportunity to compare the relative efficiencies of competing rate-optimal unbiased estimators in the case of $g$ known. This discussion will provide further insight into the results given in Remark 2.6 for models which are not locally nonparametric.

**Relative efficiency of various unbiased estimators with $g$ known.**    For simplicity, we restrict the following discussion to the truncated version of the parameter $\psi = E\left[\{b(X)\}^2\right]$, with $b(X) = E[Y|X]$, $g(\cdot)$ known, and $Y$ Bernoulli. For this choice of $\psi$, $g(\cdot)$ is the marginal density of $X$. In this subsection, we assume $\widehat{g}(X)$ is chosen equal to the known $g(X)$ so $E\left[\overline{Z}_k \overline{Z}_k^T\right] = I_{k \times k}$. Also we choose $\dot{B} = -\dot{P} = 1$ and take $\widehat{B} = \widehat{b}(X) \in lin\{\overline{Z}_k\}$, so $\widetilde{B} = \Pi\left[B|\overline{Z}_k\right] = E\left[B\overline{Z}_k^T\right]\overline{Z}_k$ and $\widetilde{\psi}_k = E\left[\left\{\Pi\left[B|\overline{Z}_k\right]\right\}^2\right]$ do not depend on $\widehat{B}$. Further we only concern ourselves with efficiency relative to the $n$ observations in the estimation sample. We thus ignore any efficiency loss from using $N - n$ observations to construct $\widehat{b}$.

Let $\Theta^g = \{b : x \mapsto b(x) \in [0,1]\} \subset \Theta$ denote the subset of $\Theta$ corresponding to the known $g$, which consists of all functions from the unit cube in $R^d$ to the unit interval. The model $\mathcal{M}(\Theta^g)$ is not locally nonparametric. For example, the first order tangent space $\Gamma_1(\theta)$ does not include first order scores for $g$. Its second order tangent space $\Gamma_2(b)$ does not contain second order scores for $g$ or mixed scores for $g$ and $b$. Rather, $\Gamma_2(b)$ is the closed linear span of the first and second order scores for $b$. Thus

$$\Gamma_2(b) = \{\mathbb{S}(a,c) ; \operatorname{var}_b[\mathbb{S}(a,c)]\} < \infty; a \in \mathcal{A}, c \in \mathcal{C} \}$$

where

$$S_{ij}(a,c) = \{(Y-B)a(X)\}_i + \left[(Y-B)_i c(X_i, X_j)(Y-B)_j\right],$$

and $\mathcal{A}$ and $\mathcal{C}$ are the set of one and two dimensional functions of $x$. Since, for

$\widehat{b} \in \lin \{\overline{z}_k(x)\}$,

$$\widehat{\psi}_{2,k}(\widehat{b}) \equiv \widehat{\psi}_{2,k}$$
$$= \mathbb{V}\left\{ \begin{array}{c} \left[\widehat{B}^2 + 2\widehat{B}\left(Y - \widehat{B}\right)\right]_i \\ + \left[\left(Y - \widehat{B}\right)\overline{Z}_k^T\right]_i \left[\overline{Z}_k\left(Y - \widehat{B}\right)\right]_j \end{array} \right\}$$

is unbiased for $\widetilde{\psi}_k(b) = E\left[\{\Pi [B|\overline{Z}_k]\}^2\right]$ in model $\mathcal{M}(\Theta^g)$, we know, by Remark 2.6, that $\mathbb{IF}^{eff}_{2,\widetilde{\psi}_k}(b)$ for $\widetilde{\psi}_k(b)$ is the projection $\Pi_b\left[\widehat{\psi}_{2,k} - \widetilde{\psi}_k(b) | \Gamma_2(\theta)\right]$ of the second order influence function $\widehat{\psi}_{2,k} - \widetilde{\psi}_k(b)$ onto $\Gamma_2(b)$. Now if $\widehat{\psi}_{2,k} - \widetilde{\psi}_k(b)$ was an element of $\Gamma_2(b)$, $\widehat{\psi}_{2,k} - \widetilde{\psi}_k(b)$ would equal $\mathbb{IF}^{eff}_{2,\widetilde{\psi}_k}(b)$ and thus be second order 'unbiased locally efficient', at $b \in \Theta^g$, as defined earlier in Remark 2.6. However we show below that, when $\widehat{b}(X) = c$ for some $c$ $w.p.1$ does not hold, $\widehat{\psi}_{2,k} - \widetilde{\psi}_k(b)$ is not an element of $\Gamma_2(b)$ for any $b$. Rather, a straightforward calculation gives

$$\mathbb{IF}^{eff}_{2,\widetilde{\psi}_k}(b) = \mathbb{V}\left\{ \begin{array}{c} \left[2E\left[B\overline{Z}_k^T\right]\overline{Z}_k(Y - B)\right]_i \\ + \left[(Y - B)\overline{Z}_k^T\right]_i \left[\overline{Z}_k(Y - B)\right]_j \end{array} \right\}.$$

Now one can check that $\widetilde{\psi}_k(\widehat{b}) + \mathbb{IF}^{eff}_{2,\widetilde{\psi}_k}(\widehat{b})$ is a function of $\widehat{b}$, so by Theorem 2.7 of Remark 2.6, we conclude no unbiased globally efficient estimator exists. However, we prove below that $\widetilde{\psi}_k(\widehat{b}) + \mathbb{IF}^{eff}_{2,\widetilde{\psi}_k}(\widehat{b})$ and $\widehat{\psi}_{2,k}$ have identical means. It follows that $\widetilde{\psi}_k(\widehat{b}) + \mathbb{IF}^{eff}_{2,\widetilde{\psi}_k}(\widehat{b})$ is an unbiased estimator of $\widetilde{\psi}_k(b) = E\left[(\Pi[B|\overline{Z}_k])^2\right]$ for any $\widehat{b} \in lin\{\overline{z}_k(x)\}$. Thus, for a given choice of $\widehat{b} \in lin\{\overline{z}_k(x)\}$, $\widetilde{\psi}_k(\widehat{b}) + \mathbb{IF}^{eff}_{2,\widetilde{\psi}_k}(\widehat{b})$ is second order unbiased locally efficient at $b = \widehat{b}$. However, one can show using a proof analogous to that in Theorem 3.22 that for $k \ll n^2$

$$\var_b\left[\widetilde{\psi}_k(\widehat{b}) + \mathbb{IF}^{eff}_{2,\widetilde{\psi}_k}(\widehat{b})\right] / \var_b\left[\mathbb{IF}^{eff}_{2,\widetilde{\psi}_k}(b)\right]$$
$$= 1 + o_P\left(\left\|\widehat{b} - b\right\|_\infty\right).$$

Henceforth we assume that $b$ lies in a Hölder ball $H(\beta_b, C_b)$. That is, we consider the submodel $b \in \Theta^g \cap H(\beta_b, C_b)$ and assume $\widehat{b}(x) \in lin\{\overline{z}_k(x)\}$ converges to $b$ in sup norm at the optimal rate of $\left(\frac{n}{\log n}\right)^{-\beta_b/(2\beta_b+d)}$ uniformly over $\Theta^g \cap H(\beta_b, C_b)$. The submodel and the model $\Theta^g$ have identical tangent spaces. For all $b \in \Theta^g \cap H(\beta_b, C_b)$, $\left(max\left(n^{-1}, k/n^2\right)\right)^{-1/2}\left\{\widetilde{\psi}_k(\widehat{b}) + \mathbb{IF}^{eff}_{2,\widetilde{\psi}_k}(\widehat{b}) - \widetilde{\psi}_k(b)\right\}$ has an asymptotic distribution with mean zero and variance equal to

$$\lim_{n\to\infty}\left(max\left(n^{-1}, k/n^2\right)\right)^{-1}\var_b\left[\mathbb{IF}^{eff}_{2,\widetilde{\psi}_k}(b)\right],$$

for all $b \in \Theta^g \cap H(\beta_b, C_b)$. In a slight abuse of language, we shall refer to $\var_b\left[\mathbb{IF}^{eff}_{2,\widetilde{\psi}_k}(b)\right]$ as the asymptotic variance of $\left\{\widetilde{\psi}_k(\widehat{b}) + \mathbb{IF}^{eff}_{2,\widetilde{\psi}_k}(\widehat{b}) - \widetilde{\psi}_k(b)\right\}$. Thus, as with standard first order theory, even when no unbiased estimator has finite sample variance that attains the Bhattacharyya bound for all $b \in \Theta^g \cap H(\beta_b, C_b)$, there can exist an unbiased estimator sequence whose asymptotic variance does attain the bound globally.

We next compare the means and variances of $\widetilde{\psi}_k\left(\widehat{b}\right) + \mathbb{IF}^{eff}_{2,\widetilde{\psi}_k}\left(\widehat{b}\right)$ and $\widehat{\psi}_{2,k}$. Now the two estimators are algebraically related by

$$\widehat{\psi}_{2,k} = \left\{\widetilde{\psi}_k\left(\widehat{b}\right) + \mathbb{IF}^{eff}_{2,\widetilde{\psi}_k}\left(\widehat{b}\right)\right\} + \left\{\mathbb{V}\left[\widehat{B}^2\right] - E\left[\widehat{B}^2\right]\right\}.$$

Since $\mathbb{V}\left[\widehat{B}^2\right] - E\left[\widehat{B}^2\right]$ is unbiased for zero, we conclude that $\widehat{\psi}_{2,k}$ and $\widetilde{\psi}_k\left(\widehat{b}\right) + \mathbb{IF}^{eff}_{2,\widetilde{\psi}_k}\left(\widehat{b}\right)$ have the same mean but $\text{var}_b\left[\widehat{\psi}_{2,k}\right]/\text{var}_b\left[\mathbb{IF}^{eff}_{2,\widetilde{\psi}_k}(b)\right] > 1$ except when $\widehat{b}(X) = b(X) = c$ w.p.1 for some $c$. Thus, since $\widetilde{\psi}_k\left(\widehat{b}\right) + \mathbb{IF}^{eff}_{2,\widetilde{\psi}_k}\left(\widehat{b}\right)$ has asymptotic variance $\text{var}_b\left[\mathbb{IF}^{eff}_{2,\widetilde{\psi}_k}(b)\right]$ and, except when $\widehat{b}(X) = c + o_p(1)$, $\text{var}\left(\mathbb{V}\left[\widehat{B}^2\right]\right) \asymp n^{-1}$, we conclude the asymptotic variance of $\widehat{\psi}_{2,k}$ attains the bound $\text{var}_b\left[\mathbb{IF}^{eff}_{2,\widetilde{\psi}_k}(b)\right]$ when $k \gg n$, but exceeds the bound when $k \leq n$, except when $\widehat{b}(X) = c + o_p(1)$.

Finally, for completeness, Robins and van der Vaart [19] considered an alternative particularly simple rate-optimal unbiased estimator of $\widetilde{\psi}_k(b) = E\left[\{\Pi\left[B|\overline{Z}_k\right]\}^2\right]$ given by $\widehat{\psi}_{RV} = \mathbb{V}\left\{\left[Y\overline{Z}_k^T\right]_i\left[\overline{Z}_kY\right]_j\right\}$. The Hoeffding decomposition of $\widehat{\psi}_{RV} - \widetilde{\psi}_k(b)$ is

$$\mathbb{V}\left[E\left[B\overline{Z}_k^T\right]\overline{Z}_kY - \widetilde{\psi}_k(b)\right] + \mathbb{V}\left\{\left[Y\overline{Z}_k^T - E\left[B\overline{Z}_k^T\right]\right]_i\left[\overline{Z}_kY - E\left[B\overline{Z}_k\right]\right]_j\right\}$$
$$= \mathbb{IF}^{eff}_{2,\widetilde{\psi}_k}(b) + Q + T$$

with

$$Q = \mathbb{V}\left[\{\Pi\left[B|\overline{Z}_k\right]B - \psi\}\right]$$
$$T = \mathbb{V}\left\{\begin{array}{c} 2\left(B_i\overline{Z}_{k,i}^T\overline{Z}_{k,j} - \Pi\left[B|\overline{Z}_k\right]_j\right)(Y-B)_j \\ +B_i\overline{Z}_{k,i}^T\overline{Z}_{k,j}B_j - \Pi\left[B|\overline{Z}_k\right]_i B_i - \Pi\left[B|\overline{Z}_k\right]_j B_j + \psi \end{array}\right\}.$$

Since, except when $B = c$ w.p.1, $\text{var}_b(Q) \asymp n^{-1}$ and $\text{var}_b(T) \asymp k/n^2$, we conclude that the asymptotic variance of $\widehat{\psi}_{RV}$ exceeds the bound $\text{var}_b\left[\mathbb{IF}^{eff}_{2,\widetilde{\psi}_k}(b)\right]$ regardless of whether $k \gg n$ does or does not hold except when $b(X) = c$ w.p.1.

**Minimaxity with unknown $g$ and $\beta/d \geq 1/4$.** We now show that the bound $n^{-\frac{1}{2}}$ for $\beta/d \geq 1/4$ is achievable for each $\beta_g > 0$. Consider the estimator $\widehat{\psi}_{m,k}$ with $n^{\frac{2}{1+4\beta/d}} \leq k \leq n$ and

$$m \geq 1 + \left\{\frac{1}{2} - \frac{\beta_b}{d+2\beta_b} - \frac{\beta_p}{d+2\beta_p}\right\}\frac{2\beta_g + d}{\beta_g}$$

so that $EB_m = O_p\left(n^{-\left(\frac{(m-1)\beta_g}{2\beta_g+d} + \frac{\beta_b}{d+2\beta_b} + \frac{\beta_p}{d+2\beta_p}\right)}\right)$ is $O_p\left(n^{-1/2}\right)$. Then $\text{var}\left(\widehat{\psi}_{m,k}\right) \asymp 1/n$, $TB_k^2 = O_p(1/n)$ and $EB_m^2 = O_p(1/n)$ so $\widehat{\psi}_{m,k}$ will be $n^{\frac{1}{2}}$−consistent for $\psi$. If $\Delta = 0$ and $\beta < 1/2$, the above expression implies that $m \geq \frac{d-2\beta}{2(2\beta+d)}/\frac{\beta_g}{(2\beta_g+d)} + 1$ for $n^{\frac{1}{2}}$-consistency. Similarly, if $\Delta \to \infty$, i.e. $\beta_b \to 0$, it is necessary that $m \geq \frac{d}{2(4\beta + d)} / \frac{\beta_g}{(2\beta_g + d)} + 1$ for $n^{\frac{1}{2}}$-consistency. These results imply that estimators $\widehat{\psi}_{m,k}$ in our class can always achieve $n^{\frac{1}{2}}$-consistency whenever $\beta_g > 0$, but for fixed $\beta < d/2$, the order $m$ of the required $U$-statistic increases without bound as the smoothness $\beta_g$ of $g$ approaches zero.

**Efficiency.** We now show that when $\beta/d$ is strictly greater than $1/4$, we can construct an unconditional asymptotically linear estimator based on all $N$ subjects with influence function $N^{-1} \sum_{i=1}^{N} IF_{1,\psi,i}(\theta)$ by having the number of the $N$ subjects allotted to the validation sample and analysis sample be $N^{1-\epsilon}$ and $n = n(\epsilon) = N - N^{1-\epsilon}$, respectively, for $1 > \epsilon > 0$. It then follows from van der Vaart [26] that the estimator is regular and semiparametric efficient. Specifically, suppose $\beta/d = 1/4 + \delta$, $\delta > 0$. Consider the estimator $\widehat{\psi}_{m^*,k}$ with $m^* > 1 + \left\{\frac{1}{2(1-\epsilon)} - \frac{\beta_b}{d+2\beta_b} - \frac{\beta_p}{d+2\beta_p}\right\} \frac{2\beta_g + d}{\beta_g}$ so that $EB_{m^*} = O_p\left(N^{-(1-\epsilon)\left(\frac{(m^*-1)\beta_g}{2\beta_g+d} + \frac{\beta_b}{d+2\beta_b} + \frac{\beta_p}{d+2\beta_p}\right)}\right)$ is $o_p\left(N^{-1/2}\right)$ and $k = n(\epsilon)^{\frac{1}{1+2\delta}} < n(\epsilon)$ so that $TB_k^2 = o_p(1/N)$ and $\text{var}\left[\widehat{IF}_{jj,\widetilde{\psi}_k}\right] = o_p(1/N)$ for $j \geq 2$. Then, by our previous results,

$$\widehat{\psi}_{m^*,k} - \psi(\theta) = n(\epsilon)^{-1} \sum_{i=1}^{n(\epsilon)} IF_{1,\psi,i}(\theta) + o_p\left(N^{-1/2}\right).$$

It remains to show that

$$N^{-1} \sum_{i=1}^{N} IF_{1,\psi,i}(\theta) - n(\epsilon)^{-1} \sum_{i=1}^{n(\epsilon)} IF_{1,\psi,i}(\theta) = o_p\left(N^{-1/2}\right).$$

But the LHS is

$$n(\epsilon)^{-1} \sum_{i=1}^{n(\epsilon)} IF_{1,\psi,i}(\theta) \left\{\frac{n(\epsilon)}{N} - 1\right\} + N^{-1} \sum_{i=n(\epsilon)+1}^{N} IF_{1,\psi,i}(\theta)$$
$$= O_p\left(n(\epsilon)^{-1/2} N^{-\epsilon}\right) + O_p\left(N^{(1-\epsilon)/2} N^{-1}\right)$$
$$= O_p\left(N^{-1/2} N^{-\epsilon}\right) + O_p\left(N^{-1/2} N^{-\epsilon/2}\right)$$
$$= o_p\left(N^{-1/2}\right).$$

**Adaptivity when $\beta/d > 1/4$.** We next prove that if we let $n \equiv n(\epsilon) = N - N^{1-\epsilon}$, $m \equiv m(N) = o(N)$ with $\ln(N) = O(m(N))$ and $k = n(\epsilon)/\ln(n)$, $\widehat{\psi}_{m,k}$ will be semiparametric efficient for each $\beta > 1/4$, provided $\|\widehat{g}(\cdot) - g(\cdot)\|_\infty = o_p\left(m\left(N^{(1-\epsilon)}\right)^{-2}\right)$. Clearly, the truncation bias is $o\left(N^{-1/2}\right)$. The estimation bias $EB_{m(N)}$ is $O_p\left(m\left(N^{(1-\epsilon)}\right)^{-2[m(N)-1]} N^{-(1-\epsilon)\left\{\frac{\beta_b}{d+2\beta_b} + \frac{\beta_p}{d+2\beta_p}\right\}}\right)$. Thus $EB_{m(N)} = o_p\left(N^{-1/2}\right)$ if $m\left(N^{(1-\epsilon)}\right)^{-2[m(N)-1]} = o\left(N^{-\frac{1}{2}+(1-\epsilon)\left\{\frac{\beta_b}{d+2\beta_b} + \frac{\beta_p}{d+2\beta_p}\right\}}\right)$. So we require $2[m(N)-1]\ln\left\{m\left(N^{(1-\epsilon)}\right)\right\}/\left[\frac{1}{2} - (1-\epsilon)\left\{\frac{\beta_b}{d+2\beta_b} + \frac{\beta_p}{d+2\beta_p}\right\}\right]\ln(N) \to \infty$, which is satisfied if $\ln(N) = O(m(N))$. In the appendix of our technical report we prove that $\mathrm{var}_\theta \left[ \widehat{\psi}_{m,k} \right] = \mathrm{var}_{\widehat{\theta}} \left[ \widehat{\psi}_{m,k} \right] \{1 + o_p(1)\}$ provided $||\widehat{g}(\cdot) - g(\cdot)||_\infty = o_p\left( m \left( N^{(1-\epsilon)} \right)^{-2} \right)$. Now

$$\mathrm{var}_{\widehat{\theta}} \left[ \widehat{\psi}_{m,k} \right] = n^{-1} \mathrm{var}_{\widehat{\theta}} \left\{ IF_{1,\psi,i}\left(\widehat{\theta}\right) \right\} \left[ O\left( \sum_{l=0}^{m(N)} \{\ln n\}^{-l} \right) \right].$$

But

$$\sum_{l=0}^{m(N)} \{\ln n\}^{-l} = O\left( \frac{1 - (\ln n)^{-[m(N)+1]}}{1 - \{\ln n\}^{-1}} \right) = O\left( 1 + \{\ln n\}^{-1} \right),$$

so $\mathrm{var}_\theta \left[ \widehat{\psi}_{m,k} \right]$ is $n^{-1} \mathrm{var}_{\widehat{\theta}} \left\{ IF_{1,\psi,i}\left(\widehat{\theta}\right) \right\} \{1 + o_p(1)\} = n^{-1} \mathrm{var}\{\mathbb{IF}_{1,\psi}\} \{1 + o_p(1)\}$. The proof of efficiency now proceeds as above.

**Alternative estimators when $\beta/d > 1/4$.** When $\beta/d > 1/4$, there actually exist, at least for certain functionals in our class, $n^{\frac{1}{2}}$-consistent estimators of $\psi$ that are much simpler than our very high order $U$-statistic estimators. For example consider the expected conditional covariance $\psi = E[\mathrm{cov}\{A, Y|X\}]$ of Example 1b with $d = 1$.

**Example 1b (Continued).** Number the study subjects $i = 0, \ldots, N-1$ ordered by their realized values $X_i$, where we have not split the sample. Following Wang et al. [27], consider the difference-based estimator

$$\widehat{\psi}_d = N^{-1} \sum_{i=0}^{N/2-1} \{Y_{2i} A_{2i} + Y_{2i+1} A_{2i+1} - Y_{2i+1} A_{2i} - Y_{2i} A_{2i+1}\},$$

which has conditional mean given $\{X_1, \ldots, X_N\}$ of

$$N^{-1} \sum_{i=0}^{N/2-1} \mathrm{cov}\{A, Y|X_{2i}\} + \mathrm{cov}\{A, Y|X_{2i+1}\}$$

$$+ N^{-1} \sum_{i=0}^{N/2-1} \left( \{b(X_{i+1}) - b(X_i)\} \{p(X_{i+1}) - p(X_i)\} \right).$$

Hence

$$E\left[ \widehat{\psi}_d - \psi \right]$$

$$= N^{-1} E\left[ \sum_{i=0}^{N/2-1} \{b(X_{i+1}) - b(X_i)\} \{p(X_{i+1}) - p(X_i)\} \right]$$

$$= O_p\left( N^{-1} \sum_{i=0}^{\frac{N}{2}-1} E\{X_{i+1} - X_i\}^{2\beta} \right) = O\left( N^{-2\beta} \right)$$

by the theory of spacings (Pyke [13]). But $O\left(N^{-2\beta}\right)$ is $o_p\left(N^{-1/2}\right)$ when $\beta > 1/4$. The variance of $\widehat{\psi}_d$ is $O\left(N^{-1}\right)$ so $\widehat{\psi}_d$ is $N^{1/2}$−consistent. However, $\frac{\mathrm{var}_\theta\left(\widehat{\psi}_d\right)}{\mathrm{var}_\theta(\mathbb{IF}_{1,\psi}(\theta))} \neq 1 + o_p(1)$ so $\widehat{\psi}_d$ is not (semiparametric) efficient. As discussed by Arellano [1], by using a $m^{\mathrm{th}}$ order rather than a second order difference operator and letting $m \to \infty$ at an appropriate rate as $N \to \infty$, the $m^{\mathrm{th}}$ order estimator $\widehat{\psi}_d$ can be made efficient.

**Minimaxity with unknown $g$ and $\beta/d < 1/4$.** Consider next whether the lower bound of $n^{-\frac{4\beta}{4\beta+d}}$ for $\beta/d < 1/4$ is achievable when $g$ is unknown but $\beta_g > 0$. We will show the next section that the bound $n^{-\frac{4\beta}{4\beta+d}}$ is achievable provided

$$(4.1) \qquad \frac{2\beta_g/d}{2\beta_g/d + 1} > \frac{4\beta/d^{\frac{1-4\beta/d}{1+4\beta/d}}(\Delta+1)}{(\Delta+2)},$$

i.e., $\beta_g > \frac{2\beta(\Delta+1)(1-4\beta/d)}{(\Delta+2)(1+4\beta/d)-4(\beta/d)(1-4\beta/d)(\Delta+1)}$. To attain the bound $n^{-\frac{4\beta}{4\beta+d}}$ whenever Equation (4.1) holds, we introduce new more efficient estimators, owing to the fact that an estimator $\widehat{\psi}_{m,k}$ in our class can attain the bound $n^{-\frac{4\beta}{4\beta+d}}$ only in the special case where the second order estimation bias $EB_2 = O_p\left(n^{-\left(\frac{\beta_g}{2\beta_g+d} + \frac{\beta_b}{d+2\beta_b} + \frac{\beta_p}{d+2\beta_p}\right)}\right)$ is less than $n^{-\frac{4\beta}{4\beta+d}}$.

For a fixed $\beta = (\beta_p + \beta_b)/2$, the right hand side of Equation (4.1) is minimized over $\Delta \geq 0$ at $\Delta = 0$. At $\Delta = 0$, Equation (4.1) reduces to

$$(4.2) \qquad \frac{\beta_g/d}{2\beta_g/d + 1} > \frac{1 - 4\beta/d}{1 + 4\beta/d}\beta/d \Rightarrow$$

$$(4.3) \qquad \beta_g > \frac{\beta(1-4\beta/d)}{1 + 2\beta/d + 8(\beta/d)^2}.$$

The right hand side of Equation (4.1) increases with $\Delta$ with asymptote equal to twice the RHS of Equation (4.2) as $\Delta \to \infty$. Hence, in order to attain the optimal rate $n^{-\frac{4\beta}{4\beta+d}}$ when $\beta_p = 2\beta$ and $\beta_b = 0$, the quantity $\frac{\beta_g}{2\beta_g+d}$ must be twice as large as when $\beta_p = \beta_b = \beta$.

In the next section, we construct an estimator with a convergence rate of $\log(n) n^{-\frac{4\beta}{4\beta+d}}$ at the cut-point $\frac{\beta_g}{1+2\beta_g} = \frac{(1-4\beta)\beta}{1+4\beta}$. In this paper we do not consider the construction of estimators that are rate optimal below this cutpoint.

However, for the special case $\Delta = 0$, in an unpublished paper Li et. al. [9] have constructed estimators which converge at a rate given in Eq. (1.3), whenever inequality (4.1) fails to hold. We conjecture that this rate is minimax, possibly only up to log factors, when inequality (4.1) fails to hold and $\Delta = 0$. At the cut-point $\frac{\beta_g}{1+2\beta_g} = \frac{(1-4\beta)\beta}{1+4\beta}$, we obtain $m^* = 0$ and thus Equation (1.3) becomes $\log(n) n^{-\frac{4\beta}{4\beta+d}}$, in agreement with the rate of the estimator of Section 4.1.2 below. In the extreme case in which $\beta_g \to 0$ with $\beta$ remaining fixed, $\log(n) n^{-\frac{1}{2} + \frac{\beta_g/d}{1+2\beta_g/d} \frac{(m^*+1)^2}{2\beta/d}} \to \log(n) n^{-\frac{1}{2} + \frac{\beta_g}{1+2\beta_g} \frac{1}{\beta}\beta(1-4\beta/d)\frac{1+2\beta_g}{2\beta_g}} = \log(n) n^{-2\beta/d}$, which agrees (up to a log factor) with the rate of $n^{-2\beta/d}$ given by the simple estimator of Wang et al. [27] analyzed above under "Example 1b continued." Based on the arguments given in the Appendix of our technical report, we conjecture that when $\beta < d/4$ and $r(n)||\widehat{g} - g||_\infty = o_p(1)$ for some $r(n) \to \infty$ as $n \to \infty$, $\left|\widehat{\psi}_{m,k} - \psi\right| = O_p\left(n^{-2\beta/d}(\log n)^t\right)$ for some natural number $t$ and any $m = m(n) \geq \frac{4(\beta/d)^2 \log(n)}{(1+2\beta/d)\log(r(n))}$ and $k = k(n) = n^{\frac{m(n)-4\beta/d}{m(n)-1}}$.

**Remark 4.1.** Suppose $b(\cdot) = E[Y|X = \cdot]$ is known to be contained in a Hölder ball $H(\beta, C)$. We provide a heuristic argument as to why the minimax rate for the linear functional $b(x)$ does not depend on a priori knowledge of the smoothness of $f_X(x)$ but the minimax rate for the functional $\psi = E[b^2(X)]$ may when $\beta/d < 1/4$. First consider the case where $f_X(x)$ is known. Let $\{\phi_l(x); l = 1, \ldots,\}$ be a complete linear independent basis for $L_2(F_X)$. Define

$$\eta_k^T = E_{F_X}\left[b(X)\overline{\phi}_k(X)^T\right] E_{F_X}\left[\overline{\phi}_k(X)\overline{\phi}_k(X)^T\right]^{-1}.$$

Then $\widetilde{b}_k(x) = \sum_{l=1}^k \Pi\left[b(x)|\overline{\phi}_k(x)\right] = \eta_k^T \overline{\phi}_k(x)$ is the projection in $L_2(F_X)$ of $b(x)$ onto linear span of the first $k$ basis functions. With $f_X(x)$ known and $\widetilde{\eta}_k^T \equiv P_n\left[Y\overline{\phi}_k(X)^T\right] E_{F_X}\left[\overline{\phi}_k(X)\overline{\phi}_k(X)^T\right]^{-1}$, unbiased estimators $\sum_{l=1}^k \widetilde{\eta}_k^T \overline{\phi}_k(x)$ and $\widetilde{\eta}_{1,k}^T E\left[\overline{\phi}_k(X)\overline{\phi}_k(X)^T\right]\widetilde{\eta}_{2,k}$ of the the truncated functionals $\widetilde{b}_k(x) = \eta_k^T \overline{\phi}_k(x)$ and $\widetilde{\psi}_k = E\left[\widetilde{b}_k^2(X)\right] = E\left[\eta_k^T \overline{\phi}_k(X)\overline{\phi}_k(X)^T \eta_k\right]$ are, respectively, rate minimax for $b(x)$ and $\psi$, when $k \equiv k_{opt}$ is chosen to equate the order of the respective variances $k/n$ and $max\left(\frac{1}{n}, k/n^2\right)$ with the order of the respective squared truncation biases $|b(x) - b_k(x)|^2 = k^{-2\beta/d}$ and $\left(\psi - \widetilde{\psi}_k\right)^2 = k^{-4\beta/d}$. For $b(x)$, $k_{opt} = n^{1/(1+2\beta/d)} \ll n$ and the rate is $n^{-\beta/(d+2\beta)}$. For $\psi$, $k_{opt} = n^{2/(1+4\beta/d)} \gg n$ and the rate is $n^{-4(\beta/d)/(1+4\beta/d)}$ when $\beta < 1/4$. The minimax rate for $b(x)$ with $f_X(x)$ unknown and without smoothness assumptions imposed is the same as with $f_X(x)$ known, since, subject to some regularity conditions, for $k_{opt} < n$,

$$\left|P_n\left[Y\overline{\phi}_{k_{opt}}(X)^T\right] P_n\left[\overline{\phi}_{k_{opt}}(X)\overline{\phi}_{k_{opt}}(X)^T\right]^{-1}\overline{\phi}_{k_{opt}}(x) - \eta_{k_{opt}}^T \overline{\phi}_{k_{opt}}(x)\right|$$

remains of order $O_p(k_{opt}/n)$. In contrast, with $k > n$, $P_n\left[\overline{\phi}_k(X)\overline{\phi}_k(X)^T\right]$ is not invertible. It is for this reason that the minimax rate for $\psi$ with $f_X(x)$ unknown is slower than $n^{-4(\beta/d)/(1+4\beta/d)}$ unless the model places sufficient restrictions on the complexity of $f_X(x)$.

**Improved rates of convergence with $X$ random in a semiparametric model.** We now, as promised in the Introduction, construct an estimator of $\sigma^2$ under the homoscedastic model $E[Y|X] = b(X)$, $var[Y|X] = \sigma^2$ with $X$ random with unknown density that, whenever $\beta < \min\{1, d/4\}$ and, regardless of the smoothness of $f_X(x)$, converges at the rate $n^{-\frac{4\beta/d}{4\beta/d+1}}$, which is faster than equal-spaced non-random minimax rate of $n^{-2\beta/d}$. Specifically we divide the support of $X$, i.e. the unit cube in $R^d$, into $k = k(n) = n^\gamma, \gamma > 1$ identical subcubes with edge length $k^{-1/d}$. We continue to assume the unknown density $f_X(x)$ is absolutely continuous with respect to Lebesgue measure and both it and its inverse are bounded in sup-norm. Then it is a standard probability calculation that the number of subcubes containing at least two observations is $O_p(n^2/k)$. We estimate $\sigma^2$ in each such subcube by $(Y_i - Y_j)^2/2$, where, for any subcube with 3 or more observations, $i$ and $j$ are chosen randomly, without replacement. Our final estimator $\widehat{\sigma}^2$ of $\sigma^2$ is the average of our subcube-specific estimates $(Y_i - Y_j)^2/2$ over the $O_p(n^2/k)$ subcubes with at least two observations. The rate of convergence of the estimator is minimized at $n^{-\frac{4\beta/d}{4\beta/d+1}}$ by taking $k = n^{\frac{2}{1+4\beta/d}}$, as we now show.

We note that

$$E\left[(Y_i - Y_j)^2/2|X_i, X_j\right] = \sigma^2 + \{b(X_i) - b(X_j)\}^2/2,$$

$|b(X_i) - b(X_j)| = O\|X_i - X_j\|^\beta$ by $\beta < 1$, and $\|X_i - X_j\| = d^{1/2} O\left(k^{-1/d}\right)$ when $X_i$ and $X_j$ are in the same subcube. It follows that the estimator has variance $O_p(k/n^2)$ and bias of $O(k^{-2\beta/d})$. To minimize the convergence rate we equate the orders of the variance and the squared bias by solving $k/n^2 = k^{-4\beta/d}$ which gives $k = n^{\frac{2}{1+4\beta/d}}$. Our random design estimator has better bias control and hence converges faster than the optimal equal-spaced fixed $X$ estimator, because the random design estimator exploits the $O_p\left(n^2/n^{\frac{2}{1+4\beta/d}}\right)$ random fluctuations for which $X's$ corresponding to two different observations are a distance of $O\left(\left\{n^{\frac{2}{1+4\beta/d}}\right\}^{-1/d}\right)$ apart. Our estimator will not converge at rate $n^{-\frac{4\beta/d}{4\beta/d+1}}$ to $E\left[\text{var}(Y|X)\right]$ in our nonparametric model, because it then no longer suffices to average estimates of var $(Y|X)$ only over subcubes containing 2 or more observations. Indeed, when var $[Y|X]$ depends on $X$, the estimator $\widehat{\sigma}^2 = \widehat{\sigma}_n^2$ satisfies $\left\{\widehat{\sigma}_n^2 - E\left[f_{k(n)}(X)\text{var}\{Y|X\}\right]/E[f_{k(n)}(X)]\right\} = O_p\left(n^{-\frac{4\beta/d}{4\beta/d+1}}\right)$, $k(n) = n^{\frac{2}{1+4\beta/d}}$, $f_{k(n)}(X)$ is $1/k(n)$ times the integral of $f_X(x)$ with respect to Lesbegue measure over the subcube containing $X$.

**Remark 4.2.** Consider again Example 1c with $\tau(\theta)$ being the variance weighted average treatment effect. We impose no smoothness assumptions on $f_X(x)$. The argument in the previous two paragraphs implies that if $\beta = (\beta_p + \beta_b)/2 < d/4$ and $\max(\beta_p, \beta_b) < 1$, we can construct an estimator $\widehat{\tau}$ of $\tau(\theta)$ that converges at rate $n^{-\frac{4\beta/d}{4\beta/d+1}}$ when the semiparametric model (3.5) holds, which is faster than our conjectured minimax rate of $n^{-2\beta/d}$ for $\tau(\theta)$ when (3.5) is not assumed. Specifically, again create $k = n^{\frac{2}{1+4\beta/d}}$ subcubes. Let $\widehat{\tau}$ make the sum $\widehat{\psi}(\tau)$ over subcubes containing at least 2 observations of $\{Y_i^*(\tau) - Y_j^*(\tau)\}\{A_i - A_j\}/2$ equal to 0 (treating subcubes with greater than 2 observations as above), where $Y^*(\tau) = Y^* - \tau A$. When (3.5) holds, cov $\{Y^*(\tau(\theta)), A|X\} = 0$. Thus, an argument analogous to that above implies that $\widehat{\psi}(\tau(\theta))$ converges to cov $\{Y^*(\tau(\theta)), A|X\} = 0$ at rate $n^{-\frac{4\beta/d}{4\beta/d+1}}$. That $\widehat{\tau}$ converges to $\tau(\theta)$ at rate $n^{-\frac{4\beta/d}{4\beta/d+1}}$ is then proved by a Taylor expansion of $0 = \widehat{\psi}(\widehat{\tau})$ around $\tau(\theta)$.

### *4.1. More efficient estimators*

#### *4.1.1. Case 1: The estimation bias of the third order estimator is less than the optimal rate*

In a (locally) nonparametric model $\mathcal{M}(\Theta)$, the estimator $\widehat{\psi}_{m,k} = \widehat{\psi} + \widehat{\mathbb{IF}}_{m,\widetilde{\psi}_k}$ is essentially the unique $m-th$ order U-statistic estimator of the truncated parameter $\widetilde{\psi}_k$ for which the leading term in the bias is $O\left(\left\|\widehat{\theta} - \theta\right\|^{m+1}\right)$. However, when the minimax rate of convergence for $\psi$ is less than $n^{-1/2}$, other $m^{th}$ order U-statistics estimators will often converge to $\widetilde{\psi}_k$ (and thus $\psi$) at a faster rate uniformly over the model than does any estimator $\widehat{\psi}_{m,k}$ (constructed from an estimated higher order influence function $\widehat{\mathbb{IF}}_{m,\widetilde{\psi}_k}$ for $\widetilde{\psi}_k$) by tolerating bias at orders less than $m+1$ in exchange for a savings in variance.

**Remark 4.3.** A heuristic understanding as to why this is so can be gained from the following considerations. The theory of higher order influence functions as developed in Theorems 2.2 and 2.3 is a theory of score functions (derivatives). Thus it can directly incorporate the restriction that a function, say $b(x)$, has an expansion $b(x) = \sum_{l=1}^{\infty} \eta_l z_l(x)$ for which $\eta_l = 0$ for $l > k$, as the restriction is equivalent to various scores being equal to zero. However the theory cannot directly incorporate restrictions such as $\sum_{l=k}^{\infty} \eta_l^2 = k^{-2\beta_b}$ or $\eta_l \propto l^{-(\beta_b+\frac{1}{2})}$ that do not imply any restrictions on score functions. Thus to find an optimal estimator, one must perform additional "side calculations" to quantify the estimation and truncation bias of various candidate estimators under these restrictions. As the assumption that $b(x)$ lies in a Hölder ball can be expressed in terms of such restrictions, this remark is relevant to a search for an optimal rate estimator.

We now construct more efficient estimators. We first consider the case where $\beta_b, \beta_b$, and $\beta_g$ are such that the estimation bias $O\left(n^{-\left(\frac{\beta_g}{2\beta_g+d}+\frac{\beta_b}{d+2\beta_b}+\frac{\beta_p}{d+2\beta_p}\right)}\right)$ of the second order estimator is greater than $O\left(n^{-\frac{4\beta}{4\beta+d}}\right)$ but the estimation bias $O\left(n^{-\left(\frac{2\beta_g}{2\beta_g+d}+\frac{\beta_b}{d+2\beta_b}+\frac{\beta_p}{d+2\beta_p}\right)}\right)$ of the third order estimator is less than $O\left(n^{-\frac{4\beta}{4\beta+d}}\right)$. That is,

$$(4.4) \qquad n^{-\left(\frac{2\beta_g}{2\beta_g+d}+\frac{\beta_b}{d+2\beta_b}+\frac{\beta_p}{d+2\beta_p}\right)} < n^{-\frac{4\beta}{4\beta+d}} < n^{-\left(\frac{\beta_g}{2\beta_g+d}+\frac{\beta_b}{d+2\beta_b}+\frac{\beta_p}{d+2\beta_p}\right)}.$$

Then the most efficient estimator $\widehat{\psi}_{m,k}$ in our class has rate of convergence slower than $n^{-\frac{4\beta}{4\beta+d}}$ because $\widehat{\psi}_{2,k_{opt}(2)}$ converges at rate $n^{-\left(\frac{\beta_g}{2\beta_g+d}+\frac{\beta_b}{d+2\beta_b}+\frac{\beta_p}{d+2\beta_p}\right)}$ determined by the second order estimation bias and, for $m > 3$, $\widehat{\psi}_{m,k_{opt}(m)}$ converges at a rate no faster than $n^{-\frac{6\beta}{(d+2\beta)}} = n^{-4\frac{\beta}{d}3/\left((3-1)+4\frac{\beta}{d}\right)} = \min_{\{m;m>3\}} n^{-4\frac{\beta}{d}m/\left((m-1)+4\frac{\beta}{d}\right)}$. We obtained $n^{-4\frac{\beta}{d}m/\left((m-1)+4\frac{\beta}{d}\right)}$ as $\left(k^{-4\beta/d}\right)^{1/2}$, where $k$ solves the equation $k^m/n^{m+1} = k^{-4\beta/d}$ that equates the variance $k^m/n^{m+1}$ of $\mathbb{IF}_m$ to the squared truncation bias $k^{-4\beta/d}$.

First, for the remainder of the paper, we redefine $\bar{Z}_k \equiv \bar{z}_k(X) \equiv \bar{\varphi}_{k,\hat{f}}(X)\{\hat{E}[H_1|X]\dot{b}(X)\dot{p}(x)\}^{-1/2}$ by redefining $\bar{\varphi}_{n^2,\hat{f}}(X)$ to have orthonormal components under $L_2(\hat{F}_X)$ such that the linear spans of $\{\varphi_{t,\hat{f}}(X),\ldots,\varphi_{n^2,\hat{f}}(X)\}$ and $\{\varphi_t(X),\ldots,\varphi_{n^2}(X)\}$ agree for $t < n^2$. This can be accomplished by Gram-Schmidt orthogonalization of $\bar{\varphi}_{n^2}(X)$ in $L_2(\hat{F}_X)$ beginning with $\varphi_{n^2}(X)$ and working backwards.

To describe our more efficient estimator, define for nonnegative integers $k(0), k(1), k^*(0), k^*(1)$ with $k(0) < k(1)$ and $k^*(0) < k^*(1)$ the $U$-statistic

$$\widehat{\mathbb{U}}_3 \begin{pmatrix} k(1), k^*(1) \\ k(0), k^*(0) \end{pmatrix} = \mathbb{V}\left(\widehat{U}_3 \begin{pmatrix} k(1), k^*(1) \\ k(0), k^*(0) \end{pmatrix}\right)$$

with

$$\widehat{\mathbb{U}}_3 \begin{pmatrix} k(1),k^*(1) \\ k(0),k^*(0) \end{pmatrix}$$
$$= \widehat{\epsilon}_{i_1} \overline{Z}_{k(0),i_1}^{k(1),T} \left(\left[\dot{P}\dot{B}H_1 \overline{Z}_{k(0)}^{k(1)} \overline{Z}_{k^*(0)}^{k^*(1),T}\right]_{i_2} - I_{\{k(1)-k(0)\}\times\{k^*(1)-k^*(0)\}}\right) \overline{Z}_{k^*(0),i_3}^{k^*(1)} \widehat{\Delta}_{i_3}$$
$$= \sum_{s_1=k(0)+1}^{k(1)} \sum_{s_2=k^*(0)+1}^{k^*(1)} \left\{ \begin{array}{c} \widehat{\epsilon}_{i_1} z_{s_1}(X_{i_1}) \times \\ \left\{\left[\dot{B}\dot{P}H_1\right]_{i_2} z_{s_1}(X_{i_2}) z_{s_2}(X_{i_2}) - I\left[s_1 = s_2\right]\right\} \\ \times z_{s_2}(X_{i_3}) \widehat{\Delta}_{i_3} \end{array} \right\},$$

where $\overline{Z}_{k(0)}^{k(1)} = (Z_{k(0)+1},\ldots,Z_{k(1)})^T$, $\widehat{\epsilon} = \left(H_1\widehat{P}+H_2\right)\dot{B}$, $\widehat{\Delta} = \left(H_1\widehat{B}+H_3\right)\dot{P}$, $I_{r\times v} = (I_{ij})_{r\times v}$ with $I_{ij} = I(i=j)$.

As an example, $\widehat{\mathbb{IF}}_{33,\widetilde{\psi}_k} = \widehat{\mathbb{U}}_3 \binom{k,\,k}{0,\,0}$. We can identify $\binom{k(1),\,k^*(1)}{k(0),\,k^*(0)}$ with the rectangle in $R^2$ defined by $\{(r_1, r_2)\,;\, k(0)+1 \leq r_1 \leq k(1), k^*(0)+1 \leq r_1 \leq k^*(1)\}$ with $(k(0)+1, k^*(0)+1)$ and $(k(1)+1, k^*(1)+1)$, respectively, the vertices closest and furthest from the origin. Thus $\widehat{\mathbb{IF}}_{33,\widetilde{\psi}_k} = \widehat{\mathbb{U}}_3 \binom{k,\,k}{0,\,0}$ is identified with the rectangle $\binom{k,\,k}{0,\,0}$. Indeed we can write

$$\widehat{\mathbb{U}}_3 \binom{k(1),\,k^*(1)}{k(0),\,k^*(0)}$$

$$= \sum_{(s_1,s_2) \in \binom{k(1),\,k^*(1)}{k(0),\,k^*(0)}} \left\{ \begin{array}{c} \widehat{\epsilon}_{i_1} z_{s_1}(X_{i_1}) \times \\ \left\{ \left[\dot{B}\dot{P}H_1\right]_{i_2} z_{s_1}(X_{i_2}) z_{s_2}(X_{i_2}) - I[s_1=s_2] \right\} \\ \times z_{s_2}(X_{i_3}) \widehat{\Delta}_{i_3} \end{array} \right\}$$

where, here and below, $s_1$ and $s_2$ are restricted to be integers, so $(s_1, s_2) \in \binom{k(1),\,k^*(1)}{k(0),\,k^*(0)}$ are the lattice points of the rectangle.

We next study the variance of $\widehat{\mathbb{U}}_3 \binom{k(1),\,k^*(1)}{k(0),\,k^*(0)}$. It follows from Theorem 3.20 above that the number of lattice points in $\binom{k(1),\,k^*(1)}{k(0),\,k^*(0)}$ is proportional to the variance of $\widehat{\mathbb{U}}_3 \binom{k(1),\,k^*(1)}{k(0),\,k^*(0)}$ so if $k(0) \ll k(1)$ and $k^*(0) \ll k^*(1)$ then $\mathrm{var}\left[\widehat{\mathbb{U}}_3 \binom{k(1),\,k^*(1)}{k(0),\,k^*(0)}\right]$ and $\mathrm{var}\left[\widehat{\mathbb{U}}_3 \binom{k(1),\,k^*(1)}{0,\,0}\right]$ are both of order $k(1)k^*(1)/n^3$. Hence the order of the variance of $\widehat{\mathbb{U}}_3 \binom{k(1),\,k^*(1)}{k(0),\,k^*(0)}$ is determined by the vertex of the rectangle $\binom{k(1),\,k^*(1)}{k(0),\,k^*(0)}$ furthest from the origin.

In contrast by a theorem in the Appendix of our technical report, the mean $E\left[\widehat{\mathbb{U}}_3 \binom{k(1),\,k^*(1)}{k(0),\,k^*(0)}\right]$ is

$$\widehat{E}\left(\widehat{\Pi}\left[\delta b | \overline{Z}_{k(0)}^{k(1)}\right] \delta g \widehat{Q}^2 \widehat{\Pi}\left[\delta p | \overline{Z}_{k^*(0)}^{k^*(1)}\right]\right)(1+o_p(1))$$

with $\delta b = \dot{P}\widehat{E}(H_1|X)\left(\widehat{B} - B\right)$, $\delta p = \dot{B}\widehat{E}(H_1|X)\left(\widehat{P} - P\right)$, $\delta g = \frac{g(X) - \widehat{g}(X)}{g(X)}$ and $\widehat{Q}^2 = \dot{B}\dot{P}\widehat{E}(H_1|X)$. It follows that if $k(0) \ll k(1)$ and $k^*(0) \ll k^*(1)$ then $E\left[\widehat{\mathbb{U}}_3 \binom{k(1),\,k^*(1)}{k(0),\,k^*(0)}\right]$ and $E\left[\widehat{\mathbb{U}}_3 \binom{\infty,\,\infty}{k(0),\,k^*(0)}\right]$ are both of order

$$O_p\left[k(0)^{-\beta_b} k^*(0)^{-\beta_p} (n/\log n)^{\frac{-\beta_g}{2\beta_g + 1}}\right].$$

This "bias" is a "product mixture" of truncation bias through the term $k(0)^{-\beta_b} k^*(0)^{-\beta_p}$ and estimaton bias through the term $(n/\log n)^{\frac{-\beta_g}{2\beta_g + 1}}$. To see this for $E\left[\widehat{\mathbb{U}}_3 \binom{k(1),\,k^*(1)}{k(0),\,k^*(0)}\right]$, we 'sup out' $\left|\delta g \widehat{Q}^2\right|$ from

$$\widehat{E}\left(\left|\widehat{\Pi}\left[\delta b | \overline{Z}_{k(0)}^{k(1)}\right] \delta g \widehat{Q}^2 \widehat{\Pi}\left[\delta p | \overline{Z}_{k^*(0)}^{k^*(1)}\right]\right|\right)$$

which is

$$O_p\left[(n/\log n)^{\frac{-\beta_g}{2\beta_g+1}}\right] \widehat{E}\left(\left|\widehat{\Pi}\left[\delta b | \overline{Z}_{k(0)}^{k(1)}\right] \widehat{\Pi}\left[\delta p | \overline{Z}_{k^*(0)}^{k^*(1)}\right]\right|\right).$$

We then apply Cauchy Schwartz to $\widehat{E}\left(\left|\widehat{\Pi}\left[\delta b | \overline{Z}_{k(0)}^{k(1)}\right] \widehat{\Pi}\left[\delta p | \overline{Z}_{k^*(0)}^{k^*(1)}\right]\right|\right)$, noting that $\widehat{E}\left(\left\{\widehat{\Pi}\left[\delta b | \overline{Z}_{k(0)}^{k(1)}\right]\right\}^2\right)^{1/2} = O\left(k(0)^{-\beta_b}\right)$. Again a more careful argument using

Hölder's inequality would show the log factor is unnecessary. Hence the order of the mean of $\widehat{\mathbb{U}}_3 \begin{pmatrix} k(1), k^*(1) \\ k(0), k^*(0) \end{pmatrix}$ is determined by the vertex of the rectangle $\begin{pmatrix} k(1), k^*(1) \\ k(0), k^*(0) \end{pmatrix}$ closest to the origin.

**Motivation.** With this background we are ready to motivate our new estimator. Recall from Section 3.2.5, that with $g$ known, the choice $k_{opt}^g(2) = n^{\frac{2}{1+4\beta/d}}$ gives $\left(\widehat{\psi}_{2,k_{opt}^g(2)} - \psi\right) = O_p\left(n^{-\frac{4\beta}{4\beta+d}}\right)$ because the truncation bias $\left|\widetilde{\psi}_{k_{opt}^g(2)} - \psi\right|$ and variance are of order $n^{-\frac{4\beta}{4\beta+d}}$ and the estimation bias is zero. Any choice of $k$ larger than $k_{opt}^g(2)$ will result in a slower rate of convergence.

However, when $g$ is unknown and thus estimated, $\widehat{\psi}_{2,k_{opt}^g(2)} - \psi$ does not attain the optimal rate of convergence because the estimation bias $n^{-\left(\frac{\beta_g}{2\beta_g+d}+\frac{\beta_b}{d+2\beta_b}+\frac{\beta_p}{d+2\beta_p}\right)}$ exceeds $n^{-\frac{4\beta}{4\beta+d}}$. The estimator $\widehat{\psi}_{3,k_{opt}^g(2)} = \widehat{\psi}_{2,k_{opt}^g(2)} + \widehat{\mathbb{U}}_3 \begin{pmatrix} k_{opt}^g(2), k_{opt}^g(2) \\ 0, 0 \end{pmatrix}$ also fails to attain the rate $n^{-\frac{4\beta}{4\beta+d}}$ because it has variance of the order of

$$\frac{k_{opt}^g(2)}{n}\frac{k_{opt}^g(2)}{n^2} = O\left(\frac{n^{\frac{2}{1+4\beta/d}}}{n}n^{-\frac{8\beta}{4\beta+d}}\right),$$

which exceeds $O\left(n^{-\frac{8\beta}{4\beta+d}}\right)$. On the other hand, $\widehat{\psi}_{3,k_{opt}^g(2)}$ has bias of $O_p\left(n^{-\frac{4\beta}{4\beta+d}}\right)$ because the truncation bias is $O_p\left(n^{-\frac{4\beta}{4\beta+d}}\right)$ and the estimation bias

$$O_p\left(n^{-\left(\frac{2\beta_g}{2\beta_g+d}+\frac{\beta_b}{d+2\beta_b}+\frac{\beta_p}{d+2\beta_p}\right)}\right)$$

is also $O_p\left(n^{-\frac{4\beta}{4\beta+d}}\right)$ under our assumption (4.4). Our strategy will be to try to replace the term $\widehat{\mathbb{U}}_3 \begin{pmatrix} k_{opt}^g(2), k_{opt}^g(2) \\ 0, 0 \end{pmatrix}$ in the estimator $\widehat{\psi}_{3,k_{opt}^g(2)} = \widehat{\psi}_{2,k_{opt}^g(2)} + \widehat{\mathbb{U}}_3 \begin{pmatrix} k_{opt}^g(2), k_{opt}^g(2) \\ 0, 0 \end{pmatrix}$ by

$$\widehat{U}_3(\Omega) = \sum_{(s_1,s_2)\in\Omega} \left\{ \begin{array}{c} \widehat{\epsilon}_{i_1} z_{s_1}(X_{i_1}) z_{s_2}(X_{i_3}) \widehat{\Delta}_{i_3} \times \\ \left\{ \left[\dot{B}\dot{P}H_1\right]_{i_2} z_{s_1}(X_{i_2}) z_{s_2}(X_{i_2}) - I[s_1 = s_2] \right\} \end{array} \right\}$$

where $\Omega$ is a subset of the rectangle $\begin{pmatrix} k_{opt}^g(2), k_{opt}^g(2) \\ 0, 0 \end{pmatrix}$ such that $\text{var}\left(\widehat{U}_3(\Omega)\right) \asymp n^{-\frac{8\beta}{4\beta+d}}$ but the additional bias

$$E\left[\widehat{\mathbb{U}}_3 \begin{pmatrix} k_{opt}^g(2), k_{opt}^g(2) \\ 0, 0 \end{pmatrix} - \widehat{U}_3(\Omega)\right]$$

$$= E\left[\widehat{\mathbb{U}}_3\left(\begin{pmatrix} k_{opt}^g(2), k_{opt}^g(2) \\ 0, 0 \end{pmatrix}\backslash\Omega\right)\right]$$

$$\equiv E\left[\sum_{(s_1,s_2)\in\begin{pmatrix} k_{opt}^g(2), k_{opt}^g(2) \\ 0, 0 \end{pmatrix}\backslash\Omega} \left\{\widehat{\epsilon}_{i_1} z_{s_1}(X_{i_1}) \left\{\begin{array}{c} \left[\dot{B}\dot{P}H_1\right]_{i_2} z_{s_1}(X_{i_2}) \\ \times z_{s_2}(X_{i_2}) - I[s_1 = s_2] \\ \times z_{s_2}(X_{i_3}) \widehat{\Delta}_{i_3} \end{array}\right\}\right\}\right]$$

is $O_p\left(n^{-\frac{4\beta}{4\beta+d}}\right)$. This approach will succeed if we can chose $\Omega$ and thus

$$\begin{pmatrix} k_{opt}^g(2\ ), k_{opt}^g(2\ ) \\ 0, \qquad\qquad 0 \end{pmatrix} \setminus \Omega$$

to be sums of rectangles (whose number does not increase with $n$) such that (i) each rectangle in $\begin{pmatrix} k_{opt}^g(2\ ), k_{opt}^g(2\ ) \\ 0, \qquad\qquad 0 \end{pmatrix} \setminus \Omega$ has its closest vertex to the origin, say $(k(0), k^*(0))$, satisfying $O_p\left[k(0)^{-\beta_b} k^*(0)^{-\beta_p} n^{\frac{-\beta_g}{2\beta_g+1}}\right] \leq n^{-\frac{4\beta}{4\beta+d}}$ and (ii) simultaneously each rectangle in $\Omega$ has its furthest vertex from the origin, say $(k(1), k^*(1))$, satisfying $O\left(k(1) k^*(1)/n^3\right) = O\left(n^{-\frac{8\beta}{4\beta+d}}\right)$.

We index the vertices of our set of rectangles as follows. Consider a natural number $J$ and a set of non-negative integers $\mathcal{K}_{J,tot} = \{k_{-2}, k_{-1}, k_0, \ldots, k_{2J}, k_{2J+1}, k_{2J+2}\}$ satisfying $0 = k_{-2} < k_0 < k_2 < \cdots < k_{2J-2} < k_{2J} < k_{2J+2} = k_{2J+1} < k_{2J-1} < \cdots < k_1 < k_{-1}$.

Note the elements with even subscripts increase from $0$ to $2J+2$ while elements with odd subscripts decrease from $-1$ to $2J-1$. Further the smallest element with odd subscript equals the largest element with even subscript. We will use two such sets $\mathcal{K}_{b,J,tot}$ and $\mathcal{K}_{p,J,tot}$ with corresponding elements $k_{bl}$ and $k_{pl}$ with $k_{b,-1} = k_{p,-1}$.

Set for $s \in \{-1, 0, \ldots, J\}$

(4.5) $$k_{b,2s+1} = n^{\frac{3d+4\beta}{(d+4\beta)}}/k_{p,2s+2},$$

(4.6) $$k_{p,2s+1} = n^{\frac{3d+4\beta}{(d+4\beta)}}/k_{b,2s+2}, \text{ so}$$

$$\frac{k_{p,2s+1}k_{b,2s+2}}{n^3} = \frac{k_{b,2s+1}k_{p,2s+2}}{n^3} = n^{-\frac{8\beta}{4\beta+d}}.$$

We leave $J$, $\mathcal{K}_{p,J} = \{k_{p,2s}, s = 0, \ldots, J+1\}$, and $\mathcal{K}_{b,J} = \{k_{b,2s}, s = 0, \ldots, J+1\}$ unspecified for now but derive optimal values below.

Let $\Omega = \Omega(\mathcal{K}_{pJ}, \mathcal{K}_{bJ})$ be the union of rectangles

$$\Omega(\mathcal{K}_{pJ}, \mathcal{K}_{bJ}) = \left\{\bigcup_{s=0}^{J} \begin{pmatrix} k_{p,2s-1}, k_{b,2s} \\ k_{p,2s-2}, k_{b,2s-2} \end{pmatrix} \cup \begin{pmatrix} k_{p,2s}, k_{b,2s-1} \\ k_{p,2s-2}\ \ k_{b,2s} \end{pmatrix}\right\} \cup \begin{pmatrix} k_{p,2J+1}k_{b,2J+1} \\ k_{p,2J}\ \ k_{b,2J} \end{pmatrix}.$$

The points $(k_{p,2s+1}, k_{b,2s+2}), (k_{p,2s+2}, k_{b,2s+1})$ for $s = -1, 0, \ldots, J+1$ lie on a hyperbola $Hy$ in $R^2$ defined by $Hy = \left\{(r_1, r_2); r_1 r_2 = n^{\frac{3d+4\beta}{(d+4\beta)}}\right\}$ shown in Figure 1 for $J = 2$. The set $\Omega(\mathcal{K}_{pJ}, \mathcal{K}_{bJ}) \subset \begin{pmatrix} k_{opt}^g(2\ ), k_{opt}^g(2\ ) \\ 0, \qquad\qquad 0 \end{pmatrix}$ lies below $Hy$.

Define

$$\widehat{\widetilde{\psi}}_{3,(\mathcal{K}_{pJ},\mathcal{K}_{bJ})} = \widehat{\psi}_{2,k_{-1}} + \widehat{\mathbb{U}}_3\left(\Omega(\mathcal{K}_{pJ},\mathcal{K}_{bJ})\right).$$

We then have

**Theorem 4.4.** (i) *The estimator* $\widehat{\widetilde{\psi}}_{3,(\mathcal{K}_{pJ},\mathcal{K}_{bJ})}$ *has variance of the order of*

$$\frac{k_{-1}}{n^2} + (2J+1) n^{-\frac{8\beta}{4\beta+d}}$$

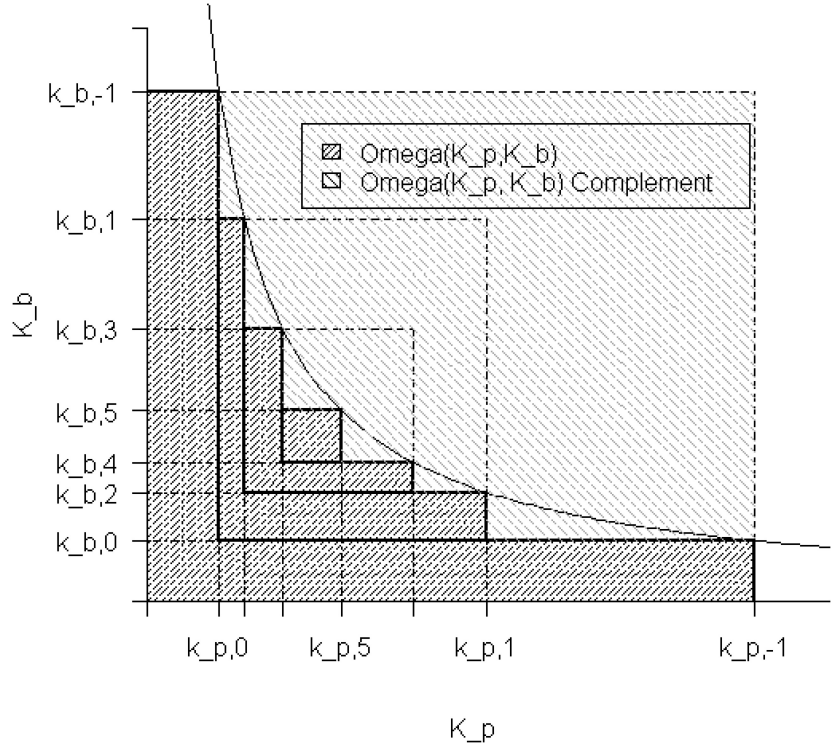

Fig 1. *Hyperbola Hy and Associated Rectangles.*

and bias $E\left(\widehat{\overline{\psi}}_{3,(\mathcal{K}_{pJ},\mathcal{K}_{bJ})}\right) - \psi$ of order

$$O_p\left\{n^{-\frac{\beta_g}{2\beta_g+d}}\left(\sum_{s=0}^{J}\left(k_{p,2s+1}^{-\beta_p/d}k_{b,2s}^{-\beta_b/d} + k_{b,2s+1}^{-\beta_b/d}k_{p,2s}^{-\beta_p/d}\right)\right)\right\}$$
$$+O_p\left(n^{-\left(\frac{2\beta_g}{2\beta_g+d}+\frac{\beta_b}{d+2\beta_b}+\frac{\beta_p}{d+2\beta_p}\right)}\right) + O_p\left(k_{-1}^{-(\beta_p+\beta_b)/d}\right).$$

*Proof.* Each of the $2J+1$ rectangles whose union is $\Omega\left(\mathcal{K}_{pJ},\mathcal{K}_{bJ}\right)$ has $(k_{p,2s+1}, k_{b,2s+2})$ or $(k_{p,2s+2}, k_{b,2s+1})$ for some $s \in \{-1,0,\ldots,J\}$ as the vertex furthest from the origin and thus contributes $\frac{k_{p,2s+1}k_{b,2s+2}}{n^3} = n^{-\frac{8\beta}{4\beta+d}}$ to the variance of $\widehat{\overline{\psi}}_{3,(\mathcal{K}_{pJ},\mathcal{K}_{bJ})}$. The variance of $\widehat{\psi}_{2,k_{-1}} \asymp \frac{k_{-1}}{n^2}$. Now

$$E\left(\widehat{\overline{\psi}}_{3,(\mathcal{K}_{pJ},\mathcal{K}_{bJ},)}\right) - \psi$$
$$= \left\{E\left(\widehat{\psi}_{3,k_{-1}}\right) - \psi\right\} + \left\{E\left[\widehat{\mathbb{U}}_3\left\{\Omega\left((\mathcal{K}_{pJ},\mathcal{K}_{bJ})\right)\right\}\right] - E\left[\widehat{\mathbb{U}}_3\left\{\begin{pmatrix}k_{-1},k_{-1}\\0, &0\end{pmatrix}\right\}\right]\right\}$$
$$= O_p\left(k_{-1}^{-(\beta_p+\beta_b)/d}\right) + O_p\left(n^{-\left(\frac{2\beta_g}{2\beta_g+d}+\frac{\beta_b}{d+2\beta_b}+\frac{\beta_p}{d+2\beta_p}\right)}\right)$$
$$+ E\left[\widehat{\mathbb{U}}_3\left\{\begin{pmatrix}k_{-1},k_{-1}\\0, &0\end{pmatrix}\setminus\Omega\left((\mathcal{K}_{pJ},\mathcal{K}_{bJ})\right)\right\}\right].$$

As is evident from Figure 1, $\Omega^c(\mathcal{K}_{pJ},\mathcal{K}_{bJ}) \equiv \begin{pmatrix}k_{-1},k_{-1}\\0, &0\end{pmatrix}\setminus\Omega\left((\mathcal{K}_{pJ},\mathcal{K}_{bJ})\right)$ is the union of rectangles $\cup_{s=0}^{J} \left\{ \begin{pmatrix} k_{p,2s-1}, & k_{b,2s-1} \\ k_{p,2s}, & k_{b,2s+1} \end{pmatrix} \cup \begin{pmatrix} k_{p,2s-1}, & k_{b,2s+1} \\ k_{p,2s+1}, & k_{b,2s} \end{pmatrix} \right\}$ which have

$$\{(k_{p,2s}, k_{b,2s+1}), (k_{p,2s+1}, k_{b,2s}); s \in \{-1, 0, \ldots, J\}\}$$

as the set of vertices closest to the origin, leading to the expression for the bias given in the theorem. □

**Theorem 4.5.** *Given $(\beta_b, \beta_p, \beta_g)$ with $\beta_p \geq \beta_b$ so $\Delta \geq 0$, Equation (4.1) holds if and only if there exists $J, \mathcal{K}_{pJ}, \mathcal{K}_{bJ}$ such that $\widehat{\overline{\psi}}_{3,(\mathcal{K}_{pJ}, \mathcal{K}_{bJ},)} - \psi = O_p\left(n^{-\frac{4\beta}{4\beta+d}}\right)$.*

*If Equation (4.1) holds, $E\left[\widehat{\mathbb{U}}_3 \left\{ \begin{pmatrix} k_{-1}, k_{-1} \\ 0, 0 \end{pmatrix} \setminus \Omega\left((\mathcal{K}_{pJ}, \mathcal{K}_{bJ})\right) \right\} \right] = O_p\left(n^{-\frac{4\beta}{4\beta+d}}\right)$ and thus $\widehat{\overline{\psi}}_{3,(\mathcal{K}_{pJ}, \mathcal{K}_{bJ})} - \psi = O_p\left(n^{-\frac{4\beta}{4\beta+d}}\right)$, when we choose $J$ to be the smallest integer such that*

$$(1+\Delta)(J+1) + c^*(\beta_g, \beta, \Delta) \sum_{l=1}^{J+1} (1+\Delta)^{l-1} > \frac{3+4\beta/d}{2(1+4\beta/d)} \text{ with}$$

$$c^*(\beta_g, \beta, \Delta) = \left(\frac{2\beta_g/d}{2\beta_g/d+1}\right) \frac{(\Delta+2)}{4\beta/d} - \frac{2(\Delta+2)}{4\beta/d+1} + \frac{3+4\beta/d}{(1+4\beta/d)},$$

$k_{b,0} = k_{p,0} = n$, $k_{b,2s} = k_{p,2s} = n^{(1+\Delta)s} n^{q \sum_{l=1}^{s}(1+\Delta)^{l-1}}$ *for $s = 1, \ldots, J+1$, with* $q = \left\{ \frac{3+4\beta/d}{2(1+4\beta/d)} - (1+\Delta)(J+1) \right\} / \sum_{l=1}^{J+1} (1+\Delta)^{l-1}$.

*Note $J$ does not depend on the sample size $n$.*

*Proof.* From Theorem 4.4, for the variance of $\widehat{\overline{\psi}}_{3,(\mathcal{K}_{pJ}, \mathcal{K}_{bJ})}$ to be $O_p\left(n^{-\frac{8\beta}{4\beta+d}}\right)$, $J$ cannot increase with $n$. Further for the second order truncation bias $O_p\left(k_{-1}^{-(\beta_p+\beta_b)/d}\right)$ and the square root of the variance $\frac{k_{-1}}{n^2}$ of $\widehat{\psi}_{2,k_{-1}}$ both to be $O_p\left(n^{-\frac{4\beta}{4\beta+d}}\right)$, we must have $k_{-1} = k_{opt}^g(2) = n^{\frac{2}{1+4\beta/d}}$. It then follows from Equations (4.5) and (4.6) that $k_{p,0} = k_{b,0} = n$.

In order for $E\left[\widehat{\mathbb{U}}_3 \left\{ \begin{pmatrix} k_{-1}, k_{-1} \\ 0, 0 \end{pmatrix} \setminus \Omega\left((\mathcal{K}_{pJ}, \mathcal{K}_{bJ})\right) \right\} \right] = O_p\left(n^{-\frac{4\beta}{4\beta+d}}\right)$, we require for $s = 0, \ldots, J$,

(4.7) $\qquad n^{-\frac{2\beta_g}{2\beta_g+d}} \left\{ k_{b,2s}^{-2\beta_b/d} k_{p,2s+1}^{-2\beta_p/d} \right\} \leq n^{-\frac{8\beta/d}{4\beta/d+1}},$

(4.8) $\qquad n^{-\frac{2\beta_g}{2\beta_g+d}} \left\{ k_{p,2s}^{-2\beta_p/d} k_{b,2s+1}^{-2\beta_b/d} \right\} \leq n^{-\frac{8\beta/d}{4\beta/d+1}}.$

Substituting for $k_{b,2s+1}$ in Equation (4.8) using Equation (4.5) and recalling that

$\beta_p \geq \beta_b$ so $\Delta \geq 0$, we obtain

$$
(4.9) \quad n^{-\frac{2\beta_g}{2\beta_g+d}} k_{p,2s}^{-2\beta_p/d} \left\{ \frac{n^{\frac{3d+4\beta}{(d+4\beta)}}}{k_{p,2s+2}} \right\}^{-2\beta_b/d} \leq n^{-\frac{8\beta}{4\beta+d}}
$$

$$
\Leftrightarrow k_{p,2s+2}^{2\beta_b/d} \leq n^{\frac{2\beta_g/d}{2\beta_g/d+1}} n^{-\frac{8\beta/d}{4\beta/d+1}} k_{p,2s}^{2\beta_p/d} \left( n^{\frac{3+4\beta/d}{(1+4\beta/d)}} \right)^{2\beta_b/d}
$$

$$
\Leftrightarrow k_{p,2s+2} \leq n^{\left(\frac{2\beta_g/d}{2\beta_g/d+1}\right)\frac{1}{2\beta_b/d}} n^{-\frac{8\beta/d}{4\beta/d+1}\frac{1}{2\beta_b/d}} k_{p,2s}^{\frac{\beta_p}{\beta_b}} \left( n^{\frac{3+4\beta/d}{(1+4\beta/d)}} \right)
$$

$$
\Leftrightarrow 1 \leq \frac{k_{p,2s+2}}{k_{p,2s}} \leq n^{\left(\frac{2\beta_g/d}{2\beta_g/d+1}\right)\frac{1}{2\beta_b/d}} n^{-\frac{8\beta/d}{4\beta/d+1}\frac{1}{2\beta_b/d}} k_{p,2s}^{\Delta} \left( n^{\frac{3+4\beta/d}{(1+4\beta/d)}} \right)
$$

$$
(4.10) \quad \Leftrightarrow 1 \leq \frac{k_{p,2s+2}}{k_{p,2s}} \leq n^{c^*(\beta_g,\beta,\Delta)} k_{p,2s}^{\Delta}
$$

$$
\Leftrightarrow 1 \leq n^{c^*(\beta_g,\beta,\Delta)} n^{\Delta}
$$

$$
\Leftrightarrow 0 \leq c^*(\beta_g, \beta, \Delta) + \Delta
$$

since $n = k_0 \leq k_{p,2s} \leq k_{p,2s+2}$.

Solving the last expression for $\frac{2\beta_g/d}{2\beta_g/d+1}$, we obtain

$$
(4.11) \quad \frac{2\beta_g/d}{2\beta_g/d+1} \geq \frac{\frac{1-4\beta/d}{1+4\beta/d} + \Delta\left\{\frac{2}{4\beta/d+1} - 1\right\}}{\frac{(\Delta+2)}{4\beta/d}} = \left\{\frac{4\beta/d}{(\Delta+2)}\right\}(\Delta+1)\frac{1-4\beta/d}{1+4\beta/d},
$$

which is Equation (4.1), except with a nonstrict inequality. We have just deduced that the constraint (4.11) was due to restriction (4.8). We have not yet considered whether the restriction (4.7) implies additional constraints. We now show that it does not. Specifically if we set $k_{p,2l} = k_{b,2l}$ for all $l \in \{1, 2, \ldots, J+1\}$, then equation (4.7) is true whenever Equation (4.8) holds because of our assumption that $\Delta \geq 0$. Thus we can set $\mathcal{K}_{pJ} = \mathcal{K}_{bJ}$.

Thus we have shown that if $\widehat{\overline{\psi}}_{3,(\mathcal{K}_{pJ},\mathcal{K}_{bJ},\,)} - \psi = O_p\left(n^{-\frac{4\beta}{4\beta+d}}\right)$, then $k_{-1} = n^{\frac{2}{1+4\beta/d}}$, (4.11) holds, and $J$ must not increase with $n$.

We next show that when the inequality is strict in (4.11) and Equation (4.4) holds, we can find $\mathcal{K}_J = \mathcal{K}_{pJ} = \mathcal{K}_{bJ}$ for which $\widehat{\overline{\psi}}_{3,\mathcal{K}_J} - \psi = O_p\left(n^{-\frac{4\beta}{4\beta+d}}\right)$. We then complete the proof of the theorem by showing that when (4.11) holds with an equality, there is no choice of $\mathcal{K}_J$ for which $\widehat{\overline{\psi}}_{3,\mathcal{K}_J}$ converges at a rate better than $O_p\left((\log n) n^{-\frac{4\beta}{4\beta+d}}\right)$.

Suppose the inequality is strict in (4.11). Since $k_0 = n$. Equation (4.10) applied recursively suggests we define $k_{2s} = n^{(1+\Delta)^s} n^{c^*(\beta_g,\beta,\Delta)\sum_{l=1}^{s}(1+\Delta)^{l-1}}$ for $s = 1, \ldots, J+1$ and take $k_{2s+1} = \frac{n^{\frac{3d+4\beta}{(d+4\beta)}}}{k_{2s+2}}$. However, this will not generally give $k_{2J+1} = k_{2J+2} = n^{\left\{\frac{3d+4\beta}{(d+4\beta)}\right\}\frac{1}{2}}$ as required when $\mathcal{K}_{pJ} = \mathcal{K}_{bJ}$. Instead we use the modified algorithm given in the statement of the theorem which insures that $k_{2J+1} = k_{2J+2} = n^{\frac{3+4\beta/d}{2(1+4\beta/d)}}$, as required. Since $J$ is not a function of $n$, in order to show $\widehat{\overline{\psi}}_{3,\mathcal{K}_J}$ converges at rate $n^{-\frac{4\beta}{4\beta+d}}$, we only need to check the bias.

Now $\frac{k_{2s+2}}{k_{2s}} = n^{(1+\Delta)} n^{q(1+\Delta)^{s-1}} = k_0^{(1+\Delta)} n^{q(1+\Delta)^{s-1}} \leq k_0^{(1+\Delta)} n^{c^*(\beta_g,\beta,\Delta)(1+\Delta)^{s-1}}$ since $q \leq c^*(\beta_g, \beta, \Delta)$ so the bias of $\widehat{\overline{\psi}}_{3,\mathcal{K}_J}$ is $O_P\left(n^{-\frac{4\beta}{4\beta+d}}\right)$, as required.

Suppose now the equality holds in Equation (4.11) so $c^*(\beta_g, \beta, \Delta) + \Delta = 0$ and continue to assume Eq. (4.4) holds. We now construct an estimator $\widehat{\overline{\psi}}_{3,\mathcal{K}_J}$ that converges at rate $O_P\left(n^{-\frac{4\beta}{4\beta+d}}\ln(n)\right)$ and show that no estimator in our class $\widehat{\overline{\psi}}_{3,\mathcal{K}_J}$ converges at a faster rate. We conjecture this rate is minimax when the equality in Equation (4.11) holds. Again $k_{2s+1} = \frac{n^{\frac{3d+4\beta}{(d+4\beta)}}}{k_{2s+2}}$ and by the previous arguments, $k_0 = n, k_{-1} = n^{\frac{2}{(1+4\beta/d)}}, k_{2J+1} = k_{2J+2} = \left\{n^{\frac{3d+4\beta}{(d+4\beta)}}\right\}^{1/2}$. We can suppose that $k_{2s} = n\{v(n)\}^s$. It remains to determine $v(n)$ and $J = J(n)$. We know $J(n)$ must satisfy

$$k_{2J(n)+2} = \left\{n^{\frac{3d+4\beta}{(d+4\beta)}}\right\}^{1/2} = n\{v(n)\}^{J(n)+1} \text{ so}$$
$$v(n) = n^{\left(\frac{3d+4\beta}{2(d+4\beta)}-1\right)\frac{1}{J(n)+1}}.$$

The variance of $\widehat{\overline{\psi}}_{3,\mathcal{K}_J}$ is of order $n^{-\frac{8\beta}{4\beta+d}}J(n)$. Thus the order of the bias will still equal that of the variance provided we multiply the RHS of Eq. (4.9) by $J(n)$. Then Equation (4.10) becomes $1 \leq \frac{k_{p,2s+2}}{k_{p,2s}} \leq n^{c^*(\beta_g,\beta,\Delta)}k_{p,2s}^\Delta J(n)^{\frac{1}{2\beta/d}}$. Since, $\frac{k_{p,2s+2}}{k_{p,2s}} = v(n)$ and $n = k_0 \leq k_{p,2s}$, we substitute $n^\Delta = k_0^\Delta$ for $k_{p,2s}^\Delta$ in the modified Equation (4.10) which gives $v(n) = J(n)^{\frac{1}{2\beta/d}}$. Hence $n^{\left(\frac{3d+4\beta}{2(d+4\beta)}-1\right)\frac{1}{J(n)+1}} = J(n)^{\frac{1}{2\beta/d}}$ which implies that.

(4.12) $$\frac{\ln(n)}{J(n)} = O(\ln[J(n)]).$$

To minimize the variance, we want the slowest growing function of $n$ that satisfies Equation (4.12), which is $J(n) = \ln(n)$, as claimed. $\square$

### 4.1.2. *Case 2: The estimation bias of the third order estimator exceeds the optimal rate*

In this section we no longer assume that the estimation bias $n^{-\left(\frac{2\beta_g}{2\beta_g+d}+\frac{\beta_b}{d+2\beta_b}+\frac{\beta_p}{d+2\beta_p}\right)}$ of a third order estimator is less than $n^{-\frac{4\beta}{4\beta+d}}$. Then even when Equation (4.11) holds with a strict inequality, $\widehat{\overline{\psi}}_{3,\mathcal{K}_J}$ does not achieve a $n^{-\frac{4\beta}{4\beta+d}}$ rate of convergence because the fourth order bias $n^{-\left(\frac{2\beta_g}{2\beta_g+d}+\frac{\beta_b}{d+2\beta_b}+\frac{\beta_p}{d+2\beta_p}\right)}$ exceeds $n^{-\frac{4\beta}{4\beta+d}}$. However, we will now construct an estimator $\widehat{\psi}_{\mathcal{K}_J}^{eff} \equiv \widehat{\psi}_{\mathcal{K}_J}^{eff}(\beta_g, \beta_b, \beta_p)$ that under our assumptions (Ai)–(Aiv) does converge at rate $n^{-\frac{4\beta}{4\beta+d}}$ whenever $(\beta_g, \beta_b, \beta_p)$ given in assumption (Aiv) satisfy Equation (4.11) with a strict inequality. Because the estimator is very complicated, we have chosen to only define the estimator and give its properties in the text. The motivating ideas for and the formal proofs of these properties are provided in the appendix of our technical report.

To define the estimator, we need some additional notation. Define

$$\widehat{\mathbb{U}}_m\left((l)_{k(l,0)}^{k(l,1)}, 1 \leq l \leq m-1\right)$$
$$= \mathbb{V}_m\left(\widehat{\epsilon}_{i_1}\overline{Z}_{k(1,0),i_1}^{k(1,1)T}\prod_{u=2}^{m-1}\left(\dot{B}\dot{P}H_1\overline{Z}_{k(u-1,0)}^{k(u-1,1)}\overline{Z}_{k(u,0)}^{k(u,1)T} - I_{k_{u-1}\times k_u}\right)\overline{Z}_{k(m-1,0)}^{k(m-1,1)}\widehat{\Delta}_{i_m}\right)$$

where, $k_u = k(u, 1) - k(u, 0)$, $I_{k_{u-1} \times k_u} = (I_{ij})_{k_{u-1} \times k_u}$ with $I_{ij} = I(i = j)$.

Then define $\widehat{\mathbb{U}}_m \begin{pmatrix} k(1) \\ k(0) \end{pmatrix}$ as $\widehat{\mathbb{U}}_m \left( (l)_{k(0)}^{k(1)}, 1 \leq l \leq m-1 \right)$. $\widehat{\mathbb{U}}_m^{(u)} \begin{pmatrix} k^*(1) \, k(1) \\ k^*(0), k(0) \end{pmatrix}$ is defined as $\widehat{\mathbb{U}}_m \left( (l)_{k(l,0)}^{k(l,1)}, 1 \leq l \leq m-1 \right)$ with $k(l, 1) = k(1), k(l, 0) = k(0)$ for $l \neq u$, and $k(u, 1) = k^*(1), k(u, 0) = k^*(0)$. Next $\widehat{\mathbb{U}}_m^{(u,u+1)} \begin{pmatrix} k^*(1) \, k^{**}(1) \, k(1) \\ k^*(0), k^{**}(0), k(0) \end{pmatrix}$ is defined as $\widehat{\mathbb{U}}_m \left( (l)_{k(l,0)}^{k(l,1)}, 1 \leq l \leq m-1 \right)$ with $k(l, 1) = k(1), k(l, 0) = k(0)$ for $l \neq u$ and $l \neq u+1$, $k(u, 1) = k^*(1), k(u, 0) = k^*(0), k(u+1, 1) = k^{**}(1), k(u+1, 0) = k^{**}(0)$. We will use this notation for $m = 3$, even though $\widehat{\mathbb{U}}_3^{(1,2)} \begin{pmatrix} k^*(1) \, k^{**}(1) \, k(1) \\ k^*(0), k^{**}(0), k(0) \end{pmatrix}$ does not depend on $k(0), k(1)$ and is equal to $\widehat{\mathbb{U}}_3 \begin{pmatrix} k^*(1) \, k^{**}(1) \\ k^*(0), k^{**}(0) \end{pmatrix}$ of the previous subsection.

Finally, define

$$\mathbb{H}_v^* = \widehat{\mathbb{U}}_v \begin{pmatrix} k_0 \\ 0 \end{pmatrix} + \sum_{u=1}^{v-1} \widehat{\mathbb{U}}_v^{(u)} \begin{pmatrix} k_{-1} \, k_0 \\ k_0 \, , 0 \end{pmatrix},$$

$$\mathbb{G}(s, v) = \sum_{u=1}^{v-2} \left\{ \widehat{\mathbb{U}}_v^{(u,u+1)} \begin{pmatrix} k_{2s-1} \, k_{2s} \, k_0 \\ k_{2s-2}, k_{2s-2}, 0 \end{pmatrix} + \widehat{\mathbb{U}}_v^{(u,u+1)} \begin{pmatrix} k_{2s} \, k_{2s-1} \, k_0 \\ k_{2s-2}, k_{2s}, 0 \end{pmatrix} \right\},$$

$$\mathbb{Q}_v = \sum_{u=1}^{v-2} \widehat{\mathbb{U}}_v^{(u,u+1)} \begin{pmatrix} k_{2J+1} \, k_{2J+1} \, k_0 \\ k_{2J}, k_{2J}, 0 \end{pmatrix}.$$

**Theorem 4.6.** *Given $(\beta_g, \beta_b, \beta_p)$ satisfying Equation (4.11) with a strict inequality, define*

$$(4.13) \quad m(\beta_g, \beta_b, \beta_p) = int \left\{ \left( \frac{4\beta}{d+4\beta} - \frac{\beta_b}{d+2\beta_b} - \frac{\beta_p}{d+2\beta_p} \right) \left( 2 + \frac{d}{\beta_g} \right) + 1 \right\} + 1$$

*to be the smallest integer such that $\left( \frac{\log n}{n} \right)^{\frac{(m-1)\beta_g}{d+2\beta_g}} n^{-\frac{\beta_b}{d+2\beta_b} - \frac{\beta_p}{d+2\beta_p}} < n^{-\frac{4\beta}{d+4\beta}}$, where $\beta = \frac{\beta_b + \beta_p}{2}$. Let $\mathcal{K}_J$, $J$, $\widehat{\overline{\psi}}_{3,\mathcal{K}_J}$ be as in Theorem 4.5 and define*

$$\widehat{\psi}_{\mathcal{K}_J}^{eff}(\beta_g, \beta_b, \beta_p)$$

$$= \widehat{\overline{\psi}}_{3,\mathcal{K}_J} + \sum_{v=4}^{m(\beta_g,\beta_b,\beta_p)} (-1)^{v-1} \mathbb{H}_v^* + \sum_{s=1}^{J} \sum_{v=4}^{m(\beta_g,\beta_b,\beta_p)} (-1)^{v-1} \mathbb{G}(s, v)$$

$$+ \sum_{v=4}^{m(\beta_g,\beta_b,\beta_p)} (-1)^{v-1} \mathbb{Q}_v$$

$$= \mathbb{V}_{n,1} \left( H_1 \widehat{B} \widehat{P} + H_2 \widehat{B} + H_3 \widehat{P} + H_4 \right) - \mathbb{H}_2^*$$

$$+ \sum_{v=3}^{m(\beta_g,\beta_b,\beta_p)} (-1)^{v-1} \mathbb{H}_v^* + \sum_{s=1}^{J} \sum_{v=3}^{m(\beta_g,\beta_b,\beta_p)} (-1)^{v-1} \mathbb{G}(s, v) + \sum_{v=3}^{m(\beta_g,\beta_b,\beta_p)} (-1)^{v-1} \mathbb{Q}_v$$

Then

$$E\left(\widehat{\psi}^{eff}_{\mathcal{K}_J}(\beta_g, \beta_b, \beta_p)\right) - \psi(\theta)$$

$$= O_p\left(\max\left[\begin{array}{c} k_{-1}^{2\beta/d}, \left(\frac{\log n}{n}\right)^{-\frac{\beta_g}{d+2\beta_g}} k_{2s}^{-\beta_b/d} k_{2s+1}^{-\beta_p/d}, \left(\frac{\log n}{n}\right)^{-\frac{\beta_g}{d+2\beta_g}} k_{2s+1}^{-\beta_b/d} k_{2s}^{-\beta_p/d}, \\ \left(\frac{\log n}{n}\right)^{-\frac{2\beta_g}{d+2\beta_g}} k_0^{-2\beta/d}, \left(\frac{\log n}{n}\right)^{-\frac{(m-1)\beta_g}{d+2\beta_g}} n^{-\frac{\beta_b}{d+2\beta_b} - \frac{\beta_p}{d+2\beta_p}}, \\ \left(\frac{\log n}{n}\right)^{-\frac{2\beta_g}{d+2\beta_g}} \max_{1\leq s\leq J}\left(k_{2s}^{-\beta_b/d} k_0^{-\beta_p/d}, k_0^{-\beta_b/d} k_{2s}^{-\beta_p/d}\right) \end{array}\right]\right)$$

$$= O_p\left(\max\left[\begin{array}{c} k_{-1}^{2\beta/d}, \left(\frac{\log n}{n}\right)^{-\frac{\beta_g}{d+2\beta_g}} k_{2s}^{-\beta_b/d} k_{2s+1}^{-\beta_p/d}, \left(\frac{\log n}{n}\right)^{-\frac{\beta_g}{d+2\beta_g}} k_{2s+1}^{-\beta_b/d} k_{2s}^{-\beta_p/d}, \\ \left(\frac{\log n}{n}\right)^{-\frac{2\beta_g}{d+2\beta_g}} k_0^{-2\beta/d}, \left(\frac{\log n}{n}\right)^{-\frac{(m-1)\beta_g}{d+2\beta_g}} n^{-\frac{\beta_b}{d+2\beta_b} - \frac{\beta_p}{d+2\beta_p}} \end{array}\right]\right)$$

$$= O_p\left(n^{-\frac{4\beta}{d+4\beta}}\right)$$

and

$$\operatorname{var}\left(\widehat{\psi}^{eff}_{\mathcal{K}_J}(\beta_g, \beta_b, \beta_p)\right)$$
$$\asymp \frac{k_{-1}}{n^2} + \sum_{s=0}^{J} \frac{k_{2s}k_{2s-1}}{n^3} + \frac{k_{2J+1}^2}{n^3} \asymp n^{-\frac{8\beta}{d+4\beta}}.$$

**Inference.** Elsewhere, we prove that $\widehat{\psi}^{eff}_{\mathcal{K}_J}(\beta_g, \beta_b, \beta_p)$ is asymptotically normal. Here, to avoid the problem of unknown 'constants' for confidence interval construction that we discussed in Section 3.2.5, we will construct nearly optimal rather than optimal confidence intervals. We suppose that Equation (4.11) holds with strict equality for the $(\beta_g, \beta_b, \beta_p)$ associated with the parameter space $\Theta$. Then there exists $\epsilon > 0$ such that for all $0 < \sigma < \epsilon$, $(\beta_g, \beta_b - \sigma, \beta_p - \sigma)$ satisfies Equation (4.11) with strict equality,

$$\sup_{\theta\in\Theta}\left[\frac{E_\theta\left[\widehat{\psi}^{eff}_{\mathcal{K}_J}(\beta_g, \beta_b - \sigma, \beta_p - \sigma)|\widehat{\theta}\right]}{\operatorname{var}_\theta\left[\widehat{\psi}^{eff}_{\mathcal{K}_J}(\beta_g, \beta_b - \sigma, \beta_p - \sigma)|\widehat{\theta}\right]}\right] = o_p(1)$$

and

$$\sup_{\theta\in\Theta}\left\{\operatorname{var}_\theta\left[\widehat{\psi}^{eff}_{\mathcal{K}_J}(\beta_g, \beta_b - \sigma, \beta_p - \sigma)|\widehat{\theta}\right]\right\} \asymp n^{-\frac{8(\beta-\sigma)}{d+4(\beta-\sigma)}}.$$

Let $\widehat{\mathbb{W}}\left[\widehat{\psi}^{eff}_{\mathcal{K}_J}(\beta_g, \beta_b, \beta_p)\right]$ be a uniformly consistent estimator of (the properly standardized) $\operatorname{var}_\theta\left[\widehat{\psi}^{eff}_{\mathcal{K}_J}(\beta_g, \beta_b, \beta_p)|\widehat{\theta}\right]$ constructed in the same manner as in Section 3.2.5. Then, for all $\sigma < \epsilon$,

$$\left\{\widehat{\psi}^{eff}_{\mathcal{K}_J}(\beta_g, \beta_b - \sigma, \beta_p - \sigma) - \psi(\theta)\right\}\left(\widehat{\mathbb{W}}\left[\widehat{\psi}^{eff}_{\mathcal{K}_J}(\beta_g, \beta_b - \sigma, \beta_p - \sigma)\right]\right)^{-1}$$

converges uniformly in $\theta \in \Theta$ to a $N(0,1)$. Moreover,

$$\widehat{\psi}^{eff}_{\mathcal{K}_J}(\beta_g, \beta_b - \sigma, \beta_p - \sigma) \pm z_\alpha \widehat{\mathbb{W}}\left[\widehat{\psi}^{eff}_{\mathcal{K}_J}(\beta_g, \beta_b - \sigma, \beta_p - \sigma)\right]$$

is a conservative uniform asymptotic $(1-\alpha)$ confidence interval for $\psi(\theta)$ with diameter of the order of $n^{-\frac{4(\beta-\sigma)}{d+4(\beta-\sigma)}}$.

**Remark 4.7.** If Equation (4.11) holds with an equality and $\mathcal{K}_J$, $J$, $\widehat{\overline{\psi}}_{3,\mathcal{K}_J}$ are as in the final paragraph of the preceding subsection then the proof of Theorem 4.6 in the appendix of our technical report implies $\widehat{\psi}_{\mathcal{K}_J}^{eff}(\beta_g, \beta_b, \beta_p) - \psi(\theta) = O_p\left((\log n) n^{-\frac{4\beta}{d+4\beta}}\right)$

## 5. Adaptive confidence intervals for regression and treatment effect functions with unknown marginal of $X$

In this section we describe how to construct adaptive confidence intervals (i) for a regression function $b(X) = E[Y|X]$ when the marginal of $X$ is unknown and (ii) for the treatment effect function and optimal treatment regime in a randomized clinical trial.

### 5.1. Regression functions

**Example 1a** (Continued). Consider the case $b = p$, $O = (Y, X)$ with $b(X) = E(Y|X)$. As usual, we assume for all $\theta \in \Theta$, $b(\cdot)$ and the density $g(\cdot)$ of $X$ are contained in known Hölder balls $H(\beta_b, C_b)$ and $H(\beta_g, C_g)$. Redefine $\psi(\theta) \equiv E_\theta\left[\left(b(X) - \widehat{b}(X)\right)^2\right]$ where $\widehat{b}(\cdot)$ is an adaptive estimate of $b(\cdot)$ from the training sample and expectations and probabilities remain conditional on the training sample. Adaptivity of $\widehat{b}(\cdot)$ implies that if $b(\cdot) \in \theta$ is also contained in a smaller Hölder ball $H(\beta^*, C)$, $\beta^* > \beta_b, C < C_b$, then $\widehat{b}(\cdot)$ will converge to $b(\cdot)$ under $F(\cdot, \theta)$ at rate $O_p\left(n^{-\frac{\beta^*/d}{1+2\beta^*/d}}\right)$. Robins and van der Vaart [19] showed that, when the marginal density $g(x)$ of $X$ is known, the key to constructing optimal (rate) adaptive confidence balls for $b(X)$ was to find a rate optimal estimator of $E_\theta\left[\left(b(X) - \widehat{b}(X)\right)^2\right]$. We shall show that their approach fails when the marginal of $X$ is unknown, but that a modification described below succeeds. Specifically, if $b(\cdot) \in \theta$ lies in a smaller Hölder ball $H(\beta^*, C)$, $\beta^* > \beta_b, C < C_b$, our modification results in honest asymptotic confidence balls under $F(\cdot, \theta)$, $\theta \in \Theta$, whose diameter is (essentially) of the same order $O_p\left(\max\left\{n^{-\frac{\beta^*/d}{1+2\beta^*/d}}, n^{-\frac{2\beta_b}{d+4\beta_b}}\right\}\right)$ as the diameter of Robins and van der Vaart's optimal adaptive region or ball, provided either (i) $\beta_b/d > 1/4$ and $\beta_g/d > 0$ or (ii) $\beta_b/d \leq 1/4$ and Equation (4.1) holds with $\beta = \beta_b$. This order is the maximum of the minimax rate $n^{-\frac{\beta^*/d}{1+2\beta^*/d}}$ of convergence of $\widehat{b}(X)$ to $b(X)$ were $b(X)$ known to lie in $H(\beta^*, C)$ and the square root of the minimax rate of convergence of an estimator of $E_\theta\left[\left(b(X) - \widehat{b}(X)\right)^2\right]$ in the larger model $\mathcal{M}(\Theta)$ with $b(\cdot)$ and $g(\cdot)$ only known to lie in $H(\beta_b, C_b)$ and $H(\beta_g, C_g)$.

The case where $\beta_b/d \leq 1/4$ and Equation (4.1) does not hold will be considered elsewhere.

Now, since $E_\theta\left[\widehat{b}(X) b(X)\right] = E_\theta\left[\widehat{b}(X) Y\right]$,

$$\begin{aligned}\psi(\theta) &\equiv E_\theta\left[\left(b(X) - \widehat{b}(X)\right)^2\right] \\ &= E_\theta\left[\{b(X)\}^2\right] - 2E_\theta\left[\widehat{b}(X) b(X)\right] + E_\theta\left[\{\widehat{b}(X)\}^2\right]\end{aligned}$$

has first order influence function $\mathbb{IF}_{1,\psi}(\theta) = \mathbb{V}[H(b,b) - \psi(\theta)]$ where

$$H(b,b) = b^2(X) + 2b(X)[Y - b(X)] - 2\widehat{b}(X)Y + \widehat{b}^2(X),$$

so $H_1 = -1, H_2 = H_3 = Y, H_4 = -2\widehat{b}(X)Y + \widehat{b}^2(X)$. Thus $H(b,b)$ for $E_\theta\left[b(X)^2\right]$ differs from $H(b,b)$ for $E_\theta\left[\left(b(X) - \widehat{b}(X)\right)^2\right]$ only in $H_4$. Since the truncation bias $\widetilde{\psi}_k(\theta) - \psi(\theta)$, higher order influence functions of $\widetilde{\psi}_k(\theta)$ and estimation bias do not depend on $H_4$, it follows that $TB_k(\theta), \mathbb{IF}_{jj,\widetilde{\psi}_k}(\theta), \widehat{\mathbb{W}}^2_{jj,\widetilde{\psi}_k}$, and $EB_m(\theta)$ are identical for $\psi(\theta) \equiv E_\theta\left[\left(b(X) - \widehat{b}(X)\right)^2\right]$ and $\psi(\theta) \equiv E_\theta\left[b(X)^2\right]$. In contrast, $IF_{1,\psi}\left(\widehat{\theta}\right)$ is identically zero for $\psi(\theta) \equiv E_\theta\left[\left(b(X) - \widehat{b}(X)\right)^2\right]$ but not for $\psi(\theta) \equiv E_\theta\left[b(X)^2\right]$. Thus, by Theorem 3.21, for $\psi(\theta) \equiv E_\theta\left[\left(b(X) - \widehat{b}(X)\right)^2\right]$, $\mathrm{var}_\theta\left[\widehat{\psi}_{m,\widetilde{\psi}_k}\right] \asymp \frac{1}{n}\left(\frac{k}{n}\right)^{m-1}$ if $k > n$ and $m > 1$, and $\mathrm{var}_\theta\left[\widehat{\psi}_{m,\widetilde{\psi}_k}\right] = 0$ if $k \leq n$ and $m = 1$. In the case when $k \leq n$ and $m > 1$, by the Hoeffding decomposition,

$$\mathrm{var}_\theta\left[\widehat{\psi}_{m,\widetilde{\psi}_k}\right] = \mathrm{var}_\theta\left(\sum_{s=1}^{m} \mathbb{D}_s^{\left(\widehat{\psi}_{m,\widetilde{\psi}_k}\right)}(\theta)\right)$$

where $\mathbb{D}_s^{\left(\widehat{\psi}_{m,\widetilde{\psi}_k}\right)}$ is a $s$th order degenerate U-statistic. Further by Theorem 3.21, we have

$$\mathrm{var}_\theta\left[\widehat{\psi}_{m,\widetilde{\psi}_k}\right] \asymp \max\left(\mathrm{var}_\theta\left(\mathbb{D}_1^{\left(\widehat{\psi}_{m,\widetilde{\psi}_k}\right)}\right), \mathrm{var}_\theta\left(\mathbb{D}_2^{\left(\widehat{\psi}_{m,\widetilde{\psi}_k}\right)}\right)\right)$$

as $\mathrm{var}_\theta\left(\mathbb{D}_s^{\left(\widehat{\psi}_{m,\widetilde{\psi}_k}\right)}\right) \asymp \frac{1}{n}\left(\frac{k}{n}\right)^{s-1} = o\left(\frac{k}{n^2}\right)$ for any $s > 2$. Moreover,

$$\mathrm{var}_\theta\left(\mathbb{D}_1^{\left(\widehat{\psi}_{m,\widetilde{\psi}_k}\right)}\right) \asymp \frac{\left\|b(X) - \widehat{b}(X)\right\|_2^2}{n}$$

since the kernel of $\mathbb{D}_1^{\left(\widehat{\psi}_{m,\widetilde{\psi}_k}\right)}$ is of order $O_p\left(\left\|b(X) - \widehat{b}(X)\right\|_2\right)$. In summary

$$\mathrm{var}_\theta\left[\widehat{\psi}_{m,\widetilde{\psi}_k}\right] \asymp \max\left(\frac{\left\|b(X) - \widehat{b}(X)\right\|_2^2}{n}, \frac{k}{n^2}\right)$$

$$= \max\left(n^{-\frac{2\beta_b/d}{1+2\beta_b/d}-1}, \frac{k}{n^2}\right)$$

if $k \leq n$ and $m > 1$. (In contrast, for $\psi(\theta) \equiv E_\theta\left[b(X)^2\right]$, $\mathrm{var}_\theta\left[\widehat{\psi}_{m,\widetilde{\psi}_k}\right] \asymp \frac{1}{n}$ if $k \leq n$). Thus, if $\beta_b/d > 1/4$, (i) $\widehat{\psi}_{m_{opt},k_{opt}(m_{opt})}$ has $k_{opt}(m_{opt})$ of $O\left(n^{\frac{2}{1+4\beta_b/d}}\right)$,

where $n^{\frac{2}{1+4\beta/d}} < n$ comes from equating the order $k^{-4\beta_b/d}$ of $TB_k^2(\theta)$ to the order $k/n^2 = n^{-\frac{8\beta_b/d}{1+4\beta_b/d}} \ll n^{-1}$ of the variance and (ii) $m_{opt}$ is the smallest integer $m$ such that the order $n^{-\left(\frac{(m-1)\beta_g}{2\beta_g+d} + \frac{2\beta_b}{d+2\beta_b}\right)}$ of

$$EB_m = O_p\left(n^{-\left(\frac{(m-1)\beta_g}{2\beta_g+d} + \frac{2\beta_b}{d+2\beta_b}\right)}\right)$$

is less than the order $n^{-\frac{4\beta_b/d}{1+4\beta_b/d}}$ of the standard error. It follows that, for $\beta_b/d > 1/4$, in contrast to $\psi(\theta) \equiv E_\theta\left[b(X)^2\right]$, we can estimate $\psi(\theta) \equiv E_\theta\left[\left(b(X) - \widehat{b}(X)\right)^2\right]$ at (the minimax) rate $n^{-\frac{4\beta_b/d}{1+4\beta_b/d}}$ which is faster (i.e., less ) than the usual parametric rate of $n^{-1/2}$.

When $\beta_b/d < 1/4$, the minimax rates for $\psi(\theta) \equiv E_\theta\left[\left(b(X) - \widehat{b}(X)\right)^2\right]$ and $\psi(\theta) \equiv E_\theta\left[b(X)^2\right]$ are identical and, when Eq. (4.1) holds, it follows from Theorem 4.6 that $\widehat{\psi}_{KJ}^{eff}(\beta_g, \beta_b, \beta_b)$ achieves the minimax rate of $n^{-\frac{4\beta_b/d}{1+4\beta_b/d}} \geq n^{-1/2}$.

Henceforth assume either $(i)$ $\beta_b/d > 1/4$ or $(ii)$ $\beta_b/d < 1/4$ and Equation (4.1) holds. Pick an $\epsilon$ so that Equation (4.1) holds for $(\beta_g, \beta_b - \epsilon, \beta_p - \epsilon)$. Let $0 < \sigma < \epsilon$ and define $\widehat{\psi}^* \equiv \widehat{\psi}(\sigma) = \widehat{\psi}_{m_{opt},\{k_{opt}(m_{opt})\}^{1+\sigma}}$ and

$$\widehat{\mathbb{W}}^* \equiv \widehat{\mathbb{W}}^*(\sigma) = \widehat{\mathbb{W}}_{m_{opt},\widetilde{\psi}_{\{k_{opt}(m_{opt})\}^{1+\sigma}}}$$

if $\beta_b/d > 1/4$ and $\widehat{\psi}^* = \widehat{\psi}_{KJ}^{eff}(\beta_g, \beta_b - \sigma, \beta_p - \sigma)$ and
$\widehat{\mathbb{W}}^* = \widehat{\mathbb{W}}\left[\widehat{\psi}_{KJ}^{eff}(\beta_g, \beta_b - \sigma, \beta_p - \sigma)\right]$ if $\beta_b/d < 1/4$. Note $\widehat{\mathbb{W}}^*$ is $O_p\left(n^{-\frac{4(\beta_b-\sigma)}{d+4(\beta_b-\sigma)}}\right)$ uniformly over $\Theta$, where $\Theta$ is the parameter space with smoothness parameters $(\beta_g, \beta_b)$. Then, by Equation (4.1) and results in Section 4.1.2, as $n \to \infty$,

$$\inf_{\theta \in \Theta} \Pr_\theta\left[\left\{\widehat{\psi}^* - \psi(\theta)\right\} \geq -z_\alpha \widehat{\mathbb{W}}^*\right] \geq 1 - \alpha.$$

Thus, if $\psi(\theta)$ were a function of $\theta$ only through $b(\cdot)$ so $\psi(\theta) = \psi(b)$, the set

(5.1) $$\left\{b^*(\cdot); \psi(\theta) \leq \widehat{\psi}^* + z_\alpha \widehat{\mathbb{W}}^*\right\}$$

would be an uniform asymptotic $(1 - \alpha)$ confidence region for $b(\cdot)$. However, for $\psi(\theta) = E_\theta\left[\left(b(X) - \widehat{b}(X)\right)^2\right]$, this approach fails because $\psi(\theta)$ also depends on $\theta$ through the unknown density $g(x)$ of $X$. This approach succeeded in Robins and van der Vaart [19] because $g(x)$ was assumed known.

We consider two solutions. The first gives (near) optimal adaptive honest intervals. The second would give honest, but non-optimal, intervals. The first solution is to replace $\psi(\theta)$ with its empirical mean $\psi_{emp}(b) \equiv \mathbb{V}\left[\left\{b(X) - \widehat{b}(X)\right\}^2\right]$ in Equation (5.1),

$$\psi_{emp}(b) - \psi(\theta) = O_p\left(\left[\left\{b(X) - \widehat{b}(X)\right\}^2\right] n^{-1/2}\right) = O_p\left(n^{-\left(\frac{2\beta_b}{d+2\beta_b} + \frac{1}{2}\right)}\right)$$

uniformly in $\theta \in \Theta$. It is straightforward to check that for all $\beta_b > 0$, $n^{-\left(\frac{2\beta_b}{d+2\beta_b}+\frac{1}{2}\right)} \ll n^{-\frac{4\beta_b/d}{1+4\beta_b/d}}$. Thus, for $\sigma < \epsilon$, $\left\{\widehat{\psi}^* - \psi_{emp}(b)\right\} / \left\{\widehat{\psi}^* - \psi(\theta)\right\} = 1 + o_p(1)$ uniformly over $\theta \in \Theta$, so $\inf_{\theta \in \Theta} \text{pr}_\theta \left[\left\{\widehat{\psi}^* - \psi_{emp}(b)\right\} \geq -z_\alpha \widehat{\mathbb{W}}^*\right] \geq 1 - \alpha$ and

$$\text{(5.2)} \qquad \left\{b^*(\cdot); \mathbb{V}\left[\left\{b^*(X) - \widehat{b}(X)\right\}^2\right] \leq \widehat{\psi}^* + z_\alpha \widehat{\mathbb{W}}^*\right\}$$

is a uniform asymptotic $(1 - \alpha)$ confidence region for $b(\cdot)$. Moreover, if $b(\cdot) \in \theta$ lies in a smaller Hölder ball $H(\beta^*, C)$, $\beta^* > \beta_b, C < C_b$, then, under $F(\cdot, \theta)$, the diameter

$$\left\{\widehat{\psi}^* + z_\alpha \widehat{\mathbb{W}}^*\right\}^{1/2} = \left\{\psi(\theta) + O_p\left(n^{-\frac{4(\beta_b-\sigma)}{d+4(\beta_b-\sigma)}}\right)\right\}^{1/2}$$

$$= O_p\left(\max\left\{n^{-\frac{2\beta^*/d}{1+2\beta^*/d}}, n^{-\frac{4(\beta_b-\sigma)}{d+4(\beta_b-\sigma)}}\right\}\right)^{1/2}$$

$$= O_p\left(\max\left\{n^{-\frac{\beta^*/d}{1+2\beta^*/d}}, n^{-\frac{2(\beta_b-\sigma)}{d+4(\beta_b-\sigma)}}\right\}\right)$$

since $\psi(\theta) = O_p\left(n^{-\frac{2\beta^*/d}{1+2\beta^*/d}}\right)$ and $\widehat{\psi}^* - \psi(\theta)$ and $\widehat{\mathbb{W}}^*$ are $O_p\left(n^{-\frac{4(\beta_b-\sigma)}{d+4(\beta_b-\sigma)}}\right)$.

The second, non-optimal, solution would be to replace the functional $\psi(\theta) \equiv E_\theta\left[\left(b(X) - \widehat{b}(X)\right)^2\right]$ with $\psi(b) = \int \left\{b(x) - \widehat{b}(x)\right\}^2 dx$. The functional $\psi(b)$ is the first functional we have considered that is not in our doubly robust class of functionals. Arguing as above, if we can construct an asymptotically normal higher order $U$-statistic estimator $\widehat{\psi}^*$ that converges to $\psi(b)$ at rate $n^{-\omega}$ on $\mathcal{M}(\Theta)$ and a consistent estimator $\widehat{\mathbb{W}}^*$ of its standard error, then

$$\left\{b^*(\cdot); \int \left\{b(x) - \widehat{b}(x)\right\}^2 dx \leq \widehat{\psi}^* + z_\alpha \widehat{\mathbb{W}}^*\right\}$$

would be an honest adaptive confidence interval of diameter $O_p\left(\max\left\{n^{-\frac{\beta^*/d}{1+2\beta^*/d}}, n^{-\omega/2}\right\}\right)$. We conjecture, based on arguments given elsewhere, that the minimax rate for estimation of $\psi(b) = \int \left\{b(x) - \widehat{b}(x)\right\}^2 dx$ exceeds $O_p\left[n^{-\frac{4\beta_b}{d+4\beta_b}}\right]$ whenever $\frac{\beta_g/d}{2\beta_g/d+1} < \frac{\beta/d}{(1+4\beta/d)(1+2\beta/d)}$. Since $\frac{\beta/d}{(1+4\beta/d)(1+2\beta/d)} > \frac{1-4\beta/d}{1+4\beta/d}\beta/d$ for all $\beta > 0$, it follows that, when the marginal of $X$ is unknown and $\frac{\beta/d}{(1+4\beta/d)(1+2\beta/d)} > \frac{\beta_g/d}{2\beta_g/d+1} > \frac{1-4\beta/d}{1+4\beta/d}\beta/d$, intervals based on $\mathbb{V}\left[\left\{b^*(X) - \widehat{b}(X)\right\}^2\right]$ will, but intervals based on $\int \left\{b(x) - \widehat{b}(x)\right\}^2 dx$ will not, have diameter of the same order as the optimal interval with the marginal of $X$ known.

### 5.2. Treatment effect functions in a randomized trial

**Example 4** (Continued). Consider the case $b = p$, $Y = Y^*$ w.p.1 so we have data $O = \{Y, A, X\}$, where $A$ is a binary treatment, $Y$ is the response, and $X$ is a vector of prerandomization covariates. The randomization probabilities $\pi_0(X) = P(A = 1|X)$ are known by design and $b(x) = E_\theta(Y|A = 1, X = x) - E_\theta(Y|A = 0, X = x)$ is the average treatment effects function. For $\theta \in \Theta$, $b(\cdot)$ and the density $g(\cdot)$ of $X$ are contained in known Hölder balls $H(\beta_b, C_b)$ and $H(\beta_g, C_g)$. Suppose we have an adaptive estimator $\widehat{b}(\cdot)$ of $b(\cdot)$ based on the training sample constructed as described below. Now, since $E_\theta\left[\widehat{b}(X)b(X)\right] = E_\theta\left[\widehat{b}(X)Y|A=1\right] - E_\theta\left[\widehat{b}(X)Y|A=0\right]$ has influence function $\frac{A}{\pi_0(X)}Y\widehat{b}(X) - \frac{1-A}{1-\pi_0(X)}Y\widehat{b}(X) - E_\theta\left[\widehat{b}(X)b(X)\right] = (A - \pi_0(X))\sigma_0^{-2}(X)Y\widehat{b}(X) - E_\theta[\widehat{b}(X) \times b(X)]$, where $\sigma_0^2(X) = \pi_0(X)\{1-\pi_0(X)\}$, $\psi(\theta) \equiv E_\theta\left[\left(b(X) - \widehat{b}(X)\right)^2\right]$ has first order influence functions, indexed by arbitrary functions $c(x)$, $\mathbb{IF}_{1,\psi}(\theta, c) \equiv \mathbb{IF}_{1,\psi}(\theta) = \mathbb{V}[H(b,b) - \psi(\theta)]$ with

$$H_1 = 1 - 2A\{A - \pi_0(X)\}\sigma_0^{-2}(X),$$
$$H_2 = H_3 = \{A - \pi_0(X)\}\sigma_0^{-2}(X)Y,$$
$$H_4 = \{A - \pi_0(X)\}c(X) - 2(A - \pi_0(X))\sigma_0^{-2}(X)Y\widehat{b}(X) + \widehat{b}^2(X).$$

Thus $H(b,b)$ for $E_\theta\left[\left(b(X) - \widehat{b}(X)\right)^2\right]$ differs from $H(b,b)$ for $\psi(\theta) \equiv E_\theta\left[b(X)^2\right]$ only in $H_4$. It follows that all the properties of the confidence ball 5.2 for $b(\cdot) = E_\theta(Y|X = \cdot)$ in the setting of the last subsection remain true for $b(\cdot) = E_\theta(Y|A = 1, X = \cdot) - E_\theta(Y|A = 0, X = \cdot)$ in the setting of this subsection.

Now define $d_{b^*}(x) = I[b^*(x) > 0]$. Then it then follows that an honest $1 - \alpha$ uniform asymptotic confidence set for the optimal treatment regime $d_{opt}(\cdot) = I[b(\cdot) > 0]$ is given by $\left\{d_{b^*}(\cdot); \mathbb{V}\left[\left\{b^*(X) - \widehat{b}(X)\right\}^2\right] \leq \widehat{\psi}^* + z_\alpha \widehat{\mathbb{W}}^*\right\}$.

**Adaptive estimator of the treatment effect function.** One among many approaches to constructing a rate-adaptive estimator of $b(\cdot)$ is as follows. Split the training sample into two random subsamples - a candidate estimator subsample of size $n_c$ and a validation subsample of size $n_v$, where both $n_c/n$ and $n_v/n$ are bounded away from 0 as $n \to \infty$. Noting that

$$0 = \mathbb{E}_\theta\left[\{Y - Ab(X)\}q(X)\{A - \pi_0(X)\}\right]$$

for all $q(\cdot)$, we construct candidate estimators of $b(\cdot)$ as follows. For $s = 1, 2, \ldots, n-1$, let $\widehat{\overline{\varkappa}}_s$ be the solution, if any, to the $s$ equations

$$0 = \mathbb{P}_c\left[\{Y - A\overline{\varkappa}_s^T\overline{\varphi}_s(X)\}\overline{\varphi}_s(X)\{A - \pi_0(X)\}\right],$$

where $\varphi_1(X), \varphi_2(X), \ldots$ is a complete basis with respect to Lebesgue measure in $R^d$ that provides optimal rate approximation for Hölder balls and $\mathbb{P}_c$ is the empirical measure for the candidate estimator subsample. Our candidates for $b(X)$ are the $\widehat{b}^{(s)}(X) = \overline{\varphi}_s(X)^T\widehat{\overline{\varkappa}}_s$. Robins [16] proved that $b(\cdot)$ is the unique function $b^*(\cdot)$ minimizing $\text{Risk}(b^*) \equiv E_\theta\left[\sigma_0^{-2}(X)\{Y - [A - \pi_0(X)]b^*(X)\}^2\right]$. In fact, the candidate $\widehat{b}^{(s)}(X)$ in our set for which $\text{Risk}\left(\widehat{b}^{(s)}\right)$ is smallest is also the candidate that minimizes $E\left[\left(b(X) - \widehat{b}^{(s)}(X)\right)^2\right]$ since $\text{Risk}\left(\widehat{b}^{(s)}\right) - \text{Risk}(b) =$

$E\left[\left(b\left(X\right)-\widehat{b}^{(s)}\left(X\right)\right)^{2}\right]$. Specifically,

$$E\left[\begin{array}{c}\sigma_{0}^{-2}\left(X\right)\left\{Y-\left[A-\pi_{0}\left(X\right)\right]\widehat{b}^{(s)}\left(X\right)\right\}^{2}\\ -\sigma_{0}^{-2}\left(X\right)\left\{Y-\left[A-\pi_{0}\left(X\right)\right]b\left(X\right)\right\}^{2}\end{array}\right]$$

$$=E\left[\begin{array}{c}\sigma_{0}^{-2}\left(X\right)\left(A-\pi_{0}\left(X\right)\right)\left(b\left(X\right)-\widehat{b}^{(s)}\left(X\right)\right)\times\\ \left(2\left(Ab\left(X\right)-E\left(Y|A=0,X\right)\right)-\left(A-\pi_{0}\left(X\right)\right)\left(b\left(X\right)+\widehat{b}^{(s)}\left(X\right)\right)\right)\end{array}\right]$$

$$=E\left(\sigma_{0}^{-2}\left(X\right)\left(A-\pi_{0}\left(X\right)\right)A\left(b\left(X\right)-\widehat{b}^{(s)}\left(X\right)\right)^{2}\right)$$

$$=E\left[\left(b\left(X\right)-\widehat{b}^{(s)}\left(X\right)\right)^{2}\right].$$

We use these results to select among our candidates by cross-validation. Let $\widehat{b}\left(\cdot\right)$ be the $\widehat{b}^{(s)}\left(\cdot\right)$ minimizing $\mathbb{P}_{v}\left[\sigma_{0}^{-2}\left(X\right)\left\{Y-\left[A-\pi_{0}\left(X\right)\right]\widehat{b}^{(s)}\left(X\right)\right\}^{2}\right]$ over $s=1,2,\ldots,n-1$, where $\mathbb{P}_{v}$ is the validation subsample empirical measure. If $b\left(\cdot\right)$ were known to lie in a Hölder ball $H\left(\beta,C\right)$, it is easy to check that the candidate $\widehat{b}^{(s)}\left(\cdot\right)$ with $s=\lfloor n^{\frac{1}{2\beta+1}}\rfloor$ obtains the optimal rate of $n^{\frac{-\beta}{2\beta+1}}$ for estimating $E\left[\left(b\left(X\right)-\widehat{b}^{(s)}\left(X\right)\right)^{2}\right]$. Since the number of candidates at sample size $n$ is less than $n$, it then follows at once from van der Laan and Dudoit's [23] results on model selection by cross validation that $\widehat{b}\left(\cdot\right)$ is adaptive over Hölder balls.

## 6. Testing, confidence sets, and implicitly defined functionals

In Example 1c of Section 3.1, we considered the following problem. We were given a functional $\psi\left(\tau,\theta\right)$ indexed by a real number $\tau$ and the parameter $\theta\in\Theta$. The implicitly defined-functional $\tau\left(\theta\right)$ was the assumed unique solution to $0=\psi\left(\tau,\theta\right)$. We noted that a $(1-\alpha)$ confidence set for $\tau\left(\theta\right)$ is the set of $\tau$ such that a $(1-\alpha)$ CI interval for $\psi\left(\tau,\theta\right)$ contains 0. In the following subsection we derive the width of the confidence set for $\tau\left(\theta\right)$. We then generalize the problem in the second subsection by introducing the notions of the testing tangent space, a testing influence function, and the higher order efficient testing score. In the final subsection, we show how the two earlier subsections are related.

### 6.1. Confidence intervals for implicitly defined functionals

To derive the order of the length of the confidence interval for the parameter $\tau\left(\theta\right)$ in Example 1c, we can use the next theorem as follows. Assume Equation (4.1) holds and $\beta\leq1/4$. Then we can take the estimator $\widetilde{\psi}\left(\tau\right)$ and rate $n^{-\gamma}$ in the theorem to be the estimator $\widehat{\psi}_{\mathcal{K},J}^{eff}$ and rate $n^{-\frac{4\beta}{4\beta+1}+\sigma}$ for a very small positive $\sigma$ and conclude that the length of the confidence interval for $\tau\left(\theta\right)$ in Example 1c to be $O_{p}\left(n^{-\frac{4\beta}{4\beta+1}+\sigma}\right)$.

**Theorem 6.1.** *Suppose for an estimator $\widehat{\psi}\left(\tau\right)$ and functional $\psi\left(\tau,\theta\right)$, there is a scale estimator $\widehat{\mathbb{W}}\left(\tau\right)$ such that $n^{\gamma}\widehat{\mathbb{W}}\left(\tau\right)\rightarrow w\left(\tau,\theta\right)$ in $\theta-$probability, $w\left(\tau,\theta\right)>$* $c^* > 0$ and $\left(\widehat{\psi}(\tau) - \psi(\tau,\theta)\right)/\widehat{\mathbb{W}}(\tau)$ converges in law to $N(0,1)$ uniformly for $\theta \in \Theta$, $\tau \in \{\tau(\theta); \theta \in \Theta\}$. Then, (i) with $z_\alpha$ the $\alpha$-quantile and $\Phi(\cdot)$ the CDF of a $N(0,1)$, the confidence set $\mathcal{C}_n = \left\{\tau; -z_{1-\alpha/2} < \frac{\widehat{\psi}(\tau)}{\widehat{\mathbb{W}}(\tau)} < z_{1-\alpha/2}\right\}$ is a uniform asymptotic $1-\alpha$ confidence set for the (assumed) unique solution $\tau(\theta)$ to $\psi(\tau,\theta) = 0$; (ii) the probability under $\theta$ that a sequence $\tau = \tau_n$ satisfying $\psi(\tau_n,\theta) = a_n n^{-\rho}, a_n \to a \neq 0$ is contained in $\mathcal{C}_n$ converges to 1 when $\rho > \gamma$, is $o(1)$ when $\rho < \gamma$, and converges to $\Phi\left(z_{1-\alpha/2} - \frac{a}{w(\tau(\theta),\theta)}\right) - \Phi\left(-z_{1-\alpha/2} - \frac{a}{w(\tau(\theta),\theta)}\right)$ when $\rho = \gamma$. (iii) If $\psi(\tau,\theta)$ is uniformly twice continuously differentiable in $\tau$ and $0 < \sigma < |\psi_\tau(\tau(\theta),\theta)| < c$ and $|\psi_{\tau^2}(\tau(\theta),\theta)| < c$ for constants $(\sigma, c)$, then $(ii)$ holds for a sequence $\tau = \tau_n$ satisfying $\tau_n - \tau(\theta) = \{\psi_\tau(\tau(\theta),\theta)\}^{-1} a_n n^{-\rho}, a_n \to a \neq 0, \rho > 0$.

*Proof.* (i) That $\mathcal{C}_n$ is a uniform asymptotic $1-\alpha$ confidence set is immediate. (ii) Now

$$Pr_\theta\left\{z_{1-\alpha/2} > \frac{\widehat{\psi}(\tau_n)}{\widehat{\mathbb{W}}(\tau_n)} > -z_{1-\alpha/2}\right\}$$

$$= Pr_\theta\left\{z_{1-\alpha/2} - \frac{\psi(\tau_n,\theta)}{\widehat{\mathbb{W}}(\tau_n)} > \frac{\widehat{\psi}(\tau_n) - \psi(\tau_n,\theta)}{\widehat{\mathbb{W}}(\tau_n)} > -z_{1-\alpha/2} - \frac{\psi(\tau_n,\theta)}{\widehat{\mathbb{W}}(\tau_n)}\right\}$$

$$\xrightarrow[n\to\infty]{} \Phi\left(z_{1-\alpha/2} - \lim_{n\to\infty}\frac{n^\gamma \psi(\tau_n,\theta)}{n^\gamma \widehat{\mathbb{W}}(\tau_n)}\right) - \Phi\left(-z_{1-\alpha/2} - \lim_{n\to\infty}\frac{n^\gamma \psi(\tau_n,\theta)}{n^\gamma \widehat{\mathbb{W}}(\tau_n)}\right)$$

$$= \Phi\left(z_{1-\alpha/2} - \frac{a \lim_{n\to\infty} n^{\gamma-\rho}}{w(\tau(\theta),\theta)}\right) - \Phi\left(-z_{1-\alpha/2} - \frac{a \lim_{n\to\infty} n^{\gamma-\rho}}{w(\tau(\theta),\theta)}\right).$$

(iii) Since $\psi(\tau_n,\theta) = \psi_\tau(\tau(\theta),\theta)(\tau_n - \tau(\theta)) + \frac{1}{2}\psi_{\tau^2}(\tau^*(\theta),\theta)(\tau_n - \tau(\theta))^2$ for some $\tau^*(\theta)$ between $\tau(\theta)$ and $\tau$, we have that $\psi(\tau_n,\theta) = a_n n^{-\rho} + o_p(a_n n^{-\rho}) = a_n(1+o_p(1))n^{-\rho}$ satisfies the assumption in (ii). □

**Remark 6.2.** Under some further regularity conditions, the solution $\widetilde{\tau}$ to $0 = \widetilde{\psi}(\tau)$ is asymptotically normal with mean $\tau(\theta)$ and variance $\psi_\tau^{-2}(\tau,\theta)\left[\{w(\tau(\theta),\theta)\}^2\right]$ uniformly over $\theta \in \Theta$, $\tau \in \{\tau(\theta); \theta \in \Theta\}$.

### 6.2. Testing influence functions and a higher order efficient score

In the following, we repeatedly use definitions from Section 2, which might usefully be reviewed at this point.

**Definition 6.3.** $m^{\text{th}}$ order testing nuisance tangent space, testing tangent space, testing influence functions, efficient score, efficient information, and efficient testing variance: Given a model $\mathcal{M}(\Theta)$ with parameter space $\Theta$ and a functional $\tau(\theta)$, define $\mathcal{M}(\Theta(\tau^\dagger))$ to be the submodel with parameter space $\Theta(\tau^\dagger) \equiv \Theta \cap \{\theta; \tau(\theta) = \tau^\dagger\}$). Thus $\mathcal{M}(\Theta(\tau^\dagger))$ is the submodel with $\tau(\theta)$ equal to $\tau^\dagger$. Define, for $\theta \in \Theta(\tau^\dagger)$, the $m^{\text{th}}$ order (i) testing nuisance tangent space $\Gamma_m^{nuis,test}(\theta,\tau^\dagger)$ to be the $m^{th}$ order tangent space for the submodel $\mathcal{M}(\Theta(\tau^\dagger))$, (ii) testing tangent space $\Gamma_m^{test}(\theta,\tau^\dagger)$ to be the closed linear span of $\mathbb{IF}_{1,\tau(\cdot)}(\theta) \cup \Gamma_m^{nuis,test}(\theta,\tau^\dagger)$, (iiia) set $\Gamma_m^{nuis,test,\perp}(\theta,\tau^\dagger) \equiv \left\{\mathbb{IF}_{m,\tau(\cdot)}^{test}\right\}$ of testing influence functions to be the orthocomplement of $\Gamma_m^{nuis,test}(\theta,\tau^\dagger)$ in $\mathcal{U}_m(\theta)$, (iiib) set $\Gamma_m^{std,nuis,test,\perp}(\theta,\tau^\dagger) \equiv \left\{\mathbb{IF}_{m,\tau(\cdot)}^{std,test}\right\}$

of standardized testing influence functions to be

$$\left\{ \mathbb{IF}_{m,\tau(\cdot)}^{std,test} \in \Gamma_{m}^{nuis,test,\perp}\left(\theta,\tau^{\dagger}\right); \; E_{\theta}\left[\mathbb{IF}_{m,\tau(\cdot)}^{std,test}\mathbb{IF}_{1,\tau(\cdot)}^{eff}\left(\theta\right)\right] = \text{var}_{\theta}\left[\mathbb{IF}_{1,\tau(\cdot)}^{eff}\left(\theta\right)\right]\right\},$$

(iv) efficient testing score $\mathbb{ES}_{m}^{test}\left(\theta\right) \equiv \mathbb{ES}_{m,\tau(\cdot)}^{test}\left(\theta\right) \in \Gamma_{m}^{test}\left(\theta,\tau^{\dagger}\right)$ to be

$$\mathbb{ES}_{m,\tau(\cdot)}^{test}\left(\theta\right) = \mathbb{ES}_{1}^{test} - \Pi_{\theta}\left[\mathbb{ES}_{1}^{test}|\Gamma_{m}^{nuis,test}\left(\theta,\tau^{\dagger}\right)\right]$$
$$\equiv \Pi_{\theta}\left[\mathbb{ES}_{1}^{test}\left(\theta\right)|\Gamma_{m}^{nuis,test,\perp}\left(\theta,\tau^{\dagger}\right)\right]$$

where $\mathbb{ES}_{1}^{test}\left(\theta\right) \equiv \mathbb{ES}_{1,\tau(\cdot)}^{test}\left(\theta\right) \equiv \text{var}_{\theta}\left\{\mathbb{IF}_{1,\tau(\cdot)}^{eff}\left(\theta\right)\right\}^{-1}\mathbb{IF}_{1,\tau(\cdot)}^{eff}\left(\theta\right)$, (v) efficient testing information to be $\text{var}_{\theta}\left\{\mathbb{ES}_{m}^{test}\left(\theta\right)\right\}$, and (vi) the efficient testing variance to be $\left[\text{var}_{\theta}\left\{\mathbb{ES}_{m}^{test}\left(\theta\right)\right\}\right]^{-1}$.

Further define, for $\theta \in \Theta$, the $m^{\text{th}}$ order (i) estimation nuisance tangent space $\Gamma_{m}^{nuis}\left(\theta\right)$ to be $\Gamma_{m}^{nuis}\left(\theta\right) \equiv \left\{\mathbb{A}_{m} \in \Gamma_{m}\left(\theta\right); E\left[\mathbb{A}_{m}\mathbb{IF}_{m,\tau(\cdot)}^{eff}\left(\theta\right)\right] = 0\right\}$, and (ii) efficient estimation variance to be $\text{var}_{\theta}\left[\mathbb{IF}_{m,\tau(\cdot)}^{eff}\left(\theta\right)\right]$.

**Remark 6.4.** For $m = 1$, the testing and estimation nuisance tangent spaces $\Gamma_{m}^{nuis,test}\left(\theta,\tau^{\dagger}\right)$ and $\Gamma_{m}^{nuis}\left(\theta\right)$ are identical. However for $m > 1$, $\Gamma_{m}^{nuis,test}\left(\theta,\tau^{\dagger}\right)$ is generally a strict subset of $\Gamma_{m}^{nuis}\left(\theta\right)$. For example, if the model can be parametrized as $\theta = \left(\tau,\rho\right)$ and $\Theta$ is the product of the parameter spaces for $\tau$ and $\rho$, the $\Gamma_{m}^{nuis,test}\left(\theta,\tau^{\dagger}\right)$ is the space of $m^{\text{th}}$ order scores for $\rho$; however, $\Gamma_{m}^{nuis}\left(\theta\right)$ also includes the mixed scores that have $s$ derivatives in the direction $\tau$ and $m - s \geq 1$ derivatives in $\rho$ directions. It is this strict inclusion that gives rise to higher order phenomena that do not occur in the first order theory.

**Theorem 6.5.** *Suppose $\mathbb{ES}_{m}^{test}\left(\theta\right)$ exists in $\mathcal{U}_{m}\left(\theta\right)$. Then for $\theta \in \Theta\left(\tau^{\dagger}\right)$,*

*(i) the set of estimation nuisance scores $\Gamma_{m}^{nuis}\left(\theta\right)$ includes the set of testing nuisance scores $\Gamma_{m}^{nuis,test}\left(\theta,\tau^{\dagger}\right)$ with equality of the sets when $m = 1$,*

*(ii) $\mathbb{IF}_{m,\tau(\cdot)}^{test}\left(\theta\right), \theta \in \Theta\left(\tau^{\dagger}\right)$ is standardized if and only if $E\left[\mathbb{IF}_{m,\tau(\cdot)}^{test}\left(\theta\right) \times \mathbb{ES}_{m}^{test}\left(\theta\right)\right] = 1$ if and only if $E\left[\mathbb{IF}_{m,\tau(\cdot)}^{test}\left(\theta\right)\mathbb{ES}_{1}^{test}\left(\theta\right)\right] = 1$,*

*(iii) $\left\{\mathbb{IF}_{m,\tau(\cdot)}^{std,test}\right\} = \left\{E_{\theta}\left[\mathbb{IF}_{m,\tau(\cdot)}^{test}\mathbb{ES}_{1}^{test}\left(\theta\right)\right]^{-1}\mathbb{IF}_{m,\tau(\cdot)}^{test}; \mathbb{IF}_{m,\tau(\cdot)}^{test} \in \left\{\mathbb{IF}_{m,\tau(\cdot)}^{test}\right\}\right\}$,*

*(iv) the set $\left\{\mathbb{IF}_{m,\tau(\cdot)}\left(\theta\right)\right\}$ of all $m^{\text{th}}$ order estimation influence functions is contained in $\left\{\mathbb{IF}_{m,\tau(\cdot)}^{std,test}\right\}$ with equality of the sets when $m = 1$,*

*(v) $\Pi_{\theta}\left[\mathbb{IF}_{m,\tau(\cdot)}^{std,test}\left(\theta\right)|\Gamma_{m}^{test}\left(\theta,\tau^{\dagger}\right)\right] = \left\{\text{var}\left[\mathbb{ES}_{m}^{test}\left(\theta\right)\right]\right\}^{-1}\mathbb{ES}_{m}^{test}\left(\theta\right)$,*

*(vi) $\left\{\text{var}_{\theta}\left[\mathbb{ES}_{m}^{test}\left(\theta\right)\right]\right\}^{-1}\mathbb{ES}_{m}^{test}\left(\theta\right) \in \left\{\mathbb{IF}_{m,\tau(\cdot)}^{std,test}\right\}$ and has the minimum variance $\left\{\text{var}_{\theta}\left[\mathbb{ES}_{m}^{test}\left(\theta\right)\right]\right\}^{-1}$ among members of $\left\{\mathbb{IF}_{m,\tau(\cdot)}^{std,test}\right\}$. In particular $\left\{\text{var}_{\theta}\left[\mathbb{ES}_{m}^{test}\left(\theta\right)\right]\right\}^{-1} \leq \text{var}_{\theta}\left[\mathbb{IF}_{m,\tau(\cdot)}^{eff}\left(\theta\right)\right]$ with equality when $m = 1$,*

*(vii) Given $\mathbb{IF}_{m,\tau(\cdot)}^{test}\left(\cdot\right) \in \left\{\mathbb{IF}_{m,\tau(\cdot)}^{test}\left(\cdot\right)\right\}$, any smooth submodel $\widetilde{\theta}\left(\zeta\right)$ with range containing $\theta$ and contained in $\Theta\left(\tau^{\dagger}\right)$, and $s \leq m$, we have*

$$\partial^{s} E_{\theta}\left[\mathbb{IF}_{m,\tau(\cdot)}^{test}\left(\widetilde{\theta}\left(\zeta\right)\right)\right]/\partial\zeta_{l_{1}}\ldots\partial\zeta_{l_{s}}\Big|_{\zeta=\widetilde{\theta}^{-1}\{\theta\}} = 0.$$

Thus, if $E_\theta \left[ \mathbb{IF}^{test}_{m,\tau(\cdot)}(\theta^*) \right]$ is Fréchet differentiable with respect to $\theta^*$ to order $m+1$ for a norm $||\cdot||$, $E_\theta \left[ \mathbb{IF}^{test}_{m,\tau(\cdot)}(\theta + \delta\theta) \right] = O\left( ||\delta\theta||^{m+1} \right)$ for $\theta$ and $\theta + \delta\theta$ in an open neighborhood contained in $\Theta(\tau^\dagger)$, since the Taylor expansion of $E_\theta \left[ \mathbb{IF}^{test}_{m,\tau(\cdot)}(\theta^*) \right]$ around $\theta$ through order $m$ is identically zero.

The proof of the Theorem will use the following two lemmas:

**Lemma 6.6.** *For any* $\mathbb{IF}^{test}_{m,\tau(\cdot)}(\theta), \theta \in \Theta(\tau^\dagger)$

$$E_\theta \left[ \mathbb{IF}^{test}_{m,\tau(\cdot)}(\theta) \, \mathbb{ES}^{test}_1(\theta) \right] = E_\theta \left[ \mathbb{IF}^{test}_{m,\tau(\cdot)}(\theta) \, \mathbb{ES}^{test}_m(\theta) \right].$$

*Proof.*

$$E_\theta \left[ \mathbb{IF}^{test}_{m,\tau(\cdot)} \mathbb{ES}^{test}_m(\theta) \right]$$
$$= E_\theta \left[ \mathbb{IF}^{test}_{m,\tau(\cdot)} \Pi_\theta \left[ \mathbb{ES}^{test}_1(\theta) \, | \, \Gamma^{nuis,test,\perp}_m(\theta,\tau^\dagger) \right] \right] = E_\theta \left[ \mathbb{IF}^{test}_{m,\tau(\cdot)} \mathbb{ES}^{test}_1(\theta) \right],$$

where the last equality holds by $\mathbb{IF}^{test}_{m,\tau(\cdot)} \in \Gamma^{nuis,test,\perp}_m(\theta,\tau^\dagger)$. □

**Lemma 6.7.** *For any* $\mathbb{IF}^{test}_{m,\tau(\cdot)}(\theta), \theta \in \Theta(\tau^\dagger)$,

$$\Pi_\theta \left[ \mathbb{IF}^{test}_{m,\tau(\cdot)}(\theta) \, | \, \Gamma^{test}_m(\theta,\tau^\dagger) \right]$$
$$= E \left[ \mathbb{IF}^{test}_{m,\tau(\cdot)}(\theta) \, \mathbb{ES}^{test}_m(\theta) \right] \left\{ var \left[ \mathbb{ES}^{test}_m(\theta) \right] \right\}^{-1} \mathbb{ES}^{test}_m(\theta)$$
$$= E \left[ \mathbb{IF}^{test}_{m,\tau(\cdot)}(\theta) \, \mathbb{ES}^{test}_1(\theta) \right] \left\{ var \left[ \mathbb{ES}^{test}_m(\theta) \right] \right\}^{-1} \mathbb{ES}^{test}_m(\theta).$$

*Proof.* $\Gamma^{test}_m(\theta,\tau^\dagger) = \left\{ c\mathbb{ES}^{test}_m(\theta) \,;\, c \in R^1 \right\} \oplus \Gamma^{nuis,test}_m(\theta,\tau^\dagger)$. Thus, by $\mathbb{IF}^{test}_{m,\tau(\cdot)}(\theta) \in \Gamma^{nuis,test,\perp}_m(\theta,\tau^\dagger)$,

$$\Pi_\theta \left[ \mathbb{IF}^{test}_{m,\tau(\cdot)}(\theta) \, | \, \Gamma^{test}_m(\theta,\tau^\dagger) \right] = \Pi_\theta \left[ \mathbb{IF}^{test}_{m,\tau(\cdot)}(\theta) \, | \, \left\{ c\mathbb{ES}^{test}_m(\theta);\, c \in R^1 \right\} \right]$$
$$= E \left[ \mathbb{IF}^{test}_{m,\tau(\cdot)}(\theta) \, \mathbb{ES}^{test}_m(\theta) \right] \left\{ var \left[ \mathbb{ES}^{test}_m(\theta) \right] \right\}^{-1} \mathbb{ES}^{test}_m(\theta).$$

Now apply Lemma 6.6. □

*Proof of Theorem 6.5.* (i) is immediate from the definitions. (ii) and (iiii) follow from

$$E \left[ \mathbb{IF}^{test}_{m,\tau(\cdot)}(\theta) \, \mathbb{ES}^{test}_m(\theta) \right] = 1 \Leftrightarrow E \left[ \mathbb{IF}^{test}_{m,\tau(\cdot)}(\theta) \, \mathbb{ES}^{test}_1(\theta) \right] = 1$$
$$\Leftrightarrow E_\theta \left[ \mathbb{IF}^{test}_{m,\tau(\cdot)} \mathbb{IF}^{eff}_{1,\tau(\cdot)}(\theta) \right] = var_\theta \left[ \mathbb{IF}^{eff}_{1,\tau(\cdot)}(\theta) \right],$$

where we have used Lemma 6.6. For (iv), note $\left\{ \mathbb{IF}_{m,\tau(\cdot)}(\theta) \right\} \subset \left\{ \mathbb{IF}^{test}_{m,\tau(\cdot)} \right\}$ follows from the fact that every smooth submodel through $\theta$ in model $\mathcal{M}(\Theta(\tau^\dagger))$ is a smooth submodel through $\theta$ in model $\mathcal{M}(\Theta)$. Thus it remains to prove that $\mathbb{IF}_{m,\tau(\cdot)}(\theta)$ is standardized. But, by Part 4 of Theorem 2.3,

$$E_\theta \left[ \mathbb{IF}_{m,\tau(\cdot)}(\theta) \, \mathbb{IF}^{eff}_{1,\tau(\cdot)}(\theta) \right] = var_\theta \left[ \mathbb{IF}^{eff}_{1,\tau(\cdot)}(\theta) \right].$$

(v) follows at once from Lemma 6.6 and Part (ii). For (vi), note that

$$\left\{\text{var}_\theta\left[\mathbb{ES}_m^{test}(\theta)\right]\right\}^{-1}\mathbb{ES}_m^{test}(\theta) \in \left\{\mathbb{IF}_{m,\tau(\cdot)}^{std,test}\right\}$$

by definition. Thus

$$\text{var}_\theta\left\{E_\theta\left[\mathbb{IF}_{m,\tau(\cdot)}^{test}\mathbb{ES}_m^{test}(\theta)\right]^{-1}\mathbb{IF}_{m,\tau(\cdot)}^{test}\right\} \geq \left\{\text{var}_\theta\left[\mathbb{ES}_m^{test}(\theta)\right]\right\}^{-1}$$

follows from (v). The result then follows from part (iii). Part (vii) is proved analogously to Theorem 2.2 except now all scores lie in $\Gamma_m^{nuis}(\theta)$ by range $\widetilde{\theta}(\zeta)$ in $\Theta(\tau^\dagger)$. □

In the case of (locally) nonparametric models, we can explicitly characterize $\Gamma_m^{test,\perp}(\theta,\tau^\dagger)$. Let $\left\{\mathbb{U}_{j,j}^{test,\perp}(\theta,\tau^\dagger)\right\}$ be the set of all

$$\mathbb{U}_{j,j}^{test,\perp}(\theta,\tau^\dagger) = \mathbb{V}\left[U_{j,j}^{test,\perp}(\theta,\tau^\dagger)\right]$$

with the $U_{j,j}^{test,\perp}(\theta,\tau^\dagger) = \sum_{l=1}^\infty c_l IF_{1,\tau(\cdot),i_1}^{eff}(\theta)\prod_{s=2}^j h_{l,s}(O_{i_s};\theta) \in \mathcal{U}_j(\theta)$, indexed by constants $c_l \in R^1$, and functions $h_{l,s}(O_{i_s};\theta)$ satisfying $E_\theta[h_{l,s}(O_{i_s};\theta)] = 0$. We remark that the subset of $\mathcal{U}_j(\theta)$ comprised of all $j$th order degenerate $U$-statistics can be written $\left\{\mathbb{V}\left[\sum_{l=1}^\infty \prod_{s=1}^j h_{l,s}(O_{i_s};\theta)\right]\right\}$. Thus $\left\{\mathbb{U}_{j,j}^{test,\perp}(\theta,\tau^\dagger)\right\}$ simply restricts one of the functions $h_{l,s}(O;\theta)$ to be $c_l IF_{1,\tau(\cdot)}^{eff}$.

**Theorem 6.8.** *If the model $\mathcal{M}(\Theta)$ is (locally) nonparametric, then $\Gamma_m^{test,\perp}(\theta,\tau^\dagger) = \left\{\sum_{j=2}^m \mathbb{U}_{j,j}^{test,\perp}(\theta,\tau^\dagger); \mathbb{U}_{j,j}^{test,\perp}(\theta,\tau^\dagger) \in \left\{\mathbb{U}_{j,j}^{test,\perp}(\theta,\tau^\dagger)\right\}\right\}$.*

*Proof.* Since the model is locally nonparametric $\Gamma_m^{test}(\theta,\tau^\dagger)$ includes the set of all mean zero first order $U$-statistics $\mathcal{U}_1(\theta)$ and thus any element of $\Gamma_m^{test,\perp}(\theta,\tau^\dagger)$ must be a sum of degenerate $U$-statistics of orders 2 through $m$. We continue by induction. First we prove the theorem for $m = 2$. Now, $\Gamma_2^{test}(\theta,\tau^\dagger) = \mathcal{U}_1(\theta) + \mathcal{U}_{2,2}^{nuis,test}(\theta)$ where $\mathcal{U}_{2,2}^{nuis,test}$ is the closed linear span of the second order degenerate part $\sum_{s\neq j} S_{l_1,j}S_{l_2,s}$ of second order scores $\widetilde{\mathbb{S}}_{2,\overline{l}_2} = \sum_j S_{l_1 l_2,j} + \sum_{s\neq j} S_{l_1,j}S_{l_2,s}$ in model $\mathcal{M}(\Theta(\tau^\dagger))$, where $\sum_{s\neq j} S_{l_1,j}S_{l_2,s}$ is a sum of products $S_{l_1,j}S_{l_2,s}$ of first order scores in model $\mathcal{M}(\Theta(\tau^\dagger))$ for two different subjects. By model $\mathcal{M}(\Theta)$ being (locally) nonparametric, the set of first order scores in model $\mathcal{M}(\Theta(\tau^\dagger))$ is precisely the set of random variables $\Gamma_1^{nuis,test}(\theta,\tau^\dagger)$ orthogonal to $IF_{1,\tau(\cdot)}^{eff}(\theta)$. But the set of degenerate $U$-statistics of order 2 orthogonal to the product of two scores in $\Gamma_1^{nuis,test}(\theta,\tau^\dagger)$ is clearly $\left\{\mathbb{U}_{2,2}^{test,\perp}(\theta,\tau^\dagger)\right\}$. Suppose now the theorem is true for $m$, $m \geq 2$, we show it is true for $m + 1$. By $\mathcal{M}(\Theta)$ (locally) nonparametric and the induction assumption, $\Gamma_{m+1}^{test}(\theta,\tau^\dagger) = \Gamma_m^{test}(\theta,\tau^\dagger) + \mathcal{U}_{m+1,m+1}^{test}(\theta)$ where $\mathcal{U}_{m+1,m+1}^{nuis,test}(\theta)$ is the closed linear span of the sum of products of first order scores in model $\mathcal{M}(\Theta(\tau^\dagger))$ for $m+1$ different subjects. But $\left\{\mathbb{U}_{m+1,m+1}^{test,\perp}(\theta,\tau^\dagger)\right\}$ is the set of set of degenerate $U$-statistics of order $m+1$ orthogonal to $\mathcal{U}_{m+1,m+1}^{nuis,test}(\theta)$. □

### 6.3. Implicitly defined functionals and testing influence functions

In the following theorem we show that estimation influence functions $\mathbb{IF}_{m,\psi(\tau,\cdot)}(\theta)$ for the parameter $\psi(\tau,\cdot)$ evaluated at the solution $\tau(\theta)$ to $0 = \psi(\tau,\theta)$ is contained in the set $\left\{\mathbb{IF}^{test}_{m,\tau(\cdot)}(\theta)\right\}$ of testing influence functions for $\tau(\theta)$. We also derive the estimation influence functions $\mathbb{IF}_{m,\tau(\cdot)}(\theta) = \sum_{s=1}^{m}\mathbb{IF}_{s,s,\tau(\cdot)}(\theta)$ for $\tau(\theta)$ in terms of the estimation influence functions $\mathbb{IF}_{m,\psi(\tau,\cdot)}(\theta)$ for $\psi(\tau,\cdot)$ and their derivatives with respect to $\tau$.

**Theorem 6.9.** *Let $\tau(\theta)$ be the assumed unique functional defined by $0 = \psi(\tau(\theta), \theta)$, $\theta \in \Theta$. Then, for $\theta \in \Theta(\tau^\dagger)$, whenever $\mathbb{IF}_{m,\psi(\tau^\dagger,\cdot)}(\theta)$ and $\mathbb{IF}_{m,\tau(\cdot)}(\theta)$ exist,*

(i) $\mathbb{IF}_{m,\psi(\tau^\dagger,\cdot)}(\theta) \in \left\{\mathbb{IF}^{test}_{m,\tau(\cdot)}(\theta)\right\}$,

(ii) $\mathbb{IF}_{1,\tau(\cdot)}(\theta) = -\psi_\tau^{-1}\mathbb{IF}_{1,\psi(\tau^\dagger,\cdot)}(\theta) \in \left\{\mathbb{IF}^{std,test}_{1,\tau(\cdot)}(\theta)\right\}$ *where $\psi_\tau \equiv \partial\psi(\tau,\theta)/\partial\tau_{|\tau=\tau^\dagger}$,*

(iii) $\mathbb{IF}_{m,m,\tau(\cdot)}(\theta) = -\psi_\tau^{-1}\left\{\mathbb{IF}_{m,m,\psi(\tau^\dagger,\cdot)}(\theta) + \mathbb{Q}_{m,m}(\theta)\right\}$, *where $\mathbb{Q}_{m,m}(\theta) \equiv \mathbb{Q}_{m,m,\tau(\cdot)}(\theta) = \mathbb{V}\{Q_{m,m}(\theta)\} \in \left\{\mathbb{U}^{test,\perp}_{m,m}(\theta,\tau^\dagger)\right\}$. For $m = 2$,*

$$
\begin{aligned}
Q_{2,2}(\theta) &= \frac{1}{2}\psi_{\backslash\tau^2}IF_{1,\tau(\cdot),i_1}(\theta)\ IF_{1,\tau(\cdot),i_2}(\theta) \\
&\quad + \frac{1}{2}\left[\begin{array}{l}\left(\begin{array}{l}\frac{\partial IF_{1,\psi(\tau^\dagger,\cdot),i_1}(\theta)}{\partial\tau} \\ -E_\theta\left[\frac{\partial IF_{1,\psi(\tau^\dagger,\cdot),i_1}(\theta)}{\partial\tau}\right]\end{array}\right)IF_{1,\tau(\cdot),i_2}(\theta) \\ + \left(\begin{array}{l}\frac{\partial IF_{1,\psi(\tau^\dagger,\cdot),i_2}(\theta)}{\partial\tau} \\ -E_\theta\left[\frac{\partial IF_{1,\psi(\tau^\dagger,\cdot),i_2}(\theta)}{\partial\tau}\right]\end{array}\right)IF_{1,\tau,i_1}(\theta)\end{array}\right]
\end{aligned}
$$
(6.1)

*where $\frac{\partial IF_{1,\psi(\tau^\dagger,\cdot),i_1}(\theta)}{\partial\tau} = \partial IF_{1,\psi(\tau,\cdot),i_1}(\theta)/\partial\tau_{|\tau=\tau^\dagger}$. $Q_{m,m}(\theta)$ is given in the appendix of our technical report as well as the general formula.*

*Proof.* (i) For $r \leq m$, consider any suitably smooth $r$ dimensional parametric submodel $\widetilde{\theta}(\zeta)$ with range containing $\theta$ and contained in $\Theta(\tau^\dagger)$. Let $\widetilde{\mathbb{S}}_{s\backslash\bar{l}_s}(\theta)$ be any associated $s^{\text{th}}$-order score $s \leq m$. By definition of $\tau(\theta)$, $\psi(\tau(\theta(\zeta)),\theta(\zeta)) = 0$. Hence, $0 = \partial^s\psi(\tau(\theta(\zeta)),\theta(\zeta))/\partial\zeta_{l_1}\cdots\partial\zeta_{l_s}|_{\zeta=\widetilde{\theta}^{-1}(\theta)}$. Now we expand the RHS using the chain rule and note that the only non-zero term is the term $\psi_{\backslash\bar{l}_s}(\tau^\dagger,\theta)$ in which all $s$−derivatives are taken with respect to the second $\theta(\zeta)$ in $\psi(\tau(\theta(\zeta)),\theta(\zeta))$; all other terms include derivatives of $\tau(\theta(\zeta))$, which are zero by range $\widetilde{\theta}(\zeta) \subset \Theta(\tau^\dagger)$. Further $\psi_{\backslash\bar{l}_s}(\tau^\dagger,\theta) = E_\theta\left[\mathbb{IF}_{m,\psi(\tau^\dagger,\cdot)}(\theta)\widetilde{\mathbb{S}}_{s\backslash\bar{l}_s}(\theta)\right]$ by the definition of the estimation influence function $\mathbb{IF}_{m,\psi(\tau^\dagger,\cdot)}(\theta)$. We conclude that $\mathbb{IF}_{m,\psi(\tau^\dagger,\cdot)}(\theta)$ is in $\Gamma^{nuis,test}_m(\theta,\tau^\dagger)^\perp$. (ii) $\mathbb{IF}_{1,\tau(\cdot)} = -\psi_\tau^{-1}\mathbb{IF}_{1,\psi(\tau^\dagger,\cdot)}$ is straightforward. That $\mathbb{IF}_{1,\tau(\cdot)}$ is contained in $\left\{\mathbb{IF}^{std,test}_{1,\tau(\cdot)}\right\}$ follows by Part (iv) of Theorem 6.5. (iii) See Appendix of our technical report for proof. □

### 6.4. "Inefficiency" of the efficient score

We now provide an example to show that, contrary to what one might expect based on Part (vi) of Theorem 6.5, inference concerning $\tau(\theta)$ may be more efficient when based on an 'inefficient' member of the set $\left\{\mathbb{IF}_{m,\tau(\cdot)}^{test}(\theta)\right\}$ such as $\mathbb{IF}_{m,\psi(\tau^\dagger,\cdot)}(\theta)$ than when based on the efficient score $\mathbb{ES}_{m,\tau(\cdot)}^{test}(\theta)$. Without loss of generality, it is sufficient to consider the case $m = 2$. In the following example it is $\widetilde{\tau}_k(\theta)$ and $\widetilde{\psi}_k(\tau^\dagger,\theta)$ that play the role of $\tau(\theta)$ and $\psi(\tau^\dagger,\theta)$ in the preceding theorem, because $\widetilde{\tau}_k(\theta)$ and $\widetilde{\psi}_k(\tau^\dagger,\theta)$ have, but $\tau(\theta)$ and $\psi(\tau^\dagger,\theta)$ do not have, higher order estimation and testing influence functions.

**Example 1c** (Continued). In this example, with $Y^*(\tau) \equiv Y^* - \tau A$, $A$ and $Y^*$ binary,
$$\psi(\tau,\theta) = E_\theta\left[\{Y^*(\tau) - E_\theta(Y^*(\tau)|X)\}\{A - E_\theta(A|X)\}\right]$$
and $\tau(\theta)$ satisfies $\psi(\tau(\theta),\theta) = 0$. Let $\widetilde{\tau}_k(\theta)$ satisfy $\widetilde{\psi}_k(\widetilde{\tau}_k(\theta),\theta) = 0$ where $\widetilde{\psi}_k(\tau,\theta) = E_\theta[Y^*(\tau)A] - E_\theta\left\{\left[\Pi_\theta\left[B(\tau)|\overline{Z}_k\right]\Pi_\theta\left[P|\overline{Z}_k\right]\right]\right\}$ is defined in Section 3.1 with $\tau$ a real-valued index and $B(\tau) = b(X,\tau) = E_\theta(Y^*(\tau)|X)$. Note $\widetilde{\psi}_{k,\tau}(\tau,\theta) \equiv \partial\widetilde{\psi}_k(\tau,\theta)/\partial\tau = -\left\{E_\theta[A^2] - E_\theta\left[\left\{\Pi_\theta\left[P|\overline{Z}_k\right]\right\}^2\right]\right\}$, $\psi_\tau(\tau,\theta) = -E_\theta[\text{var}_\theta(A|X)]$, $\widetilde{\psi}_{k,\tau^2}(\tau,\theta) = \psi_{\tau^2}(\tau,\theta) = 0$. Below we freely use results of Theorems 3.11, 3.14, and 3.17. We suppose that $0 < \sigma < \text{var}_\theta(A|X)$ and $E_\theta[A^2] < c$ for some $(\sigma,c)$, $\beta = \frac{\beta_p + \beta_b}{2} < 1/4$. Choose $k = k_{opt}(2)n^{2\sigma} = n^{\frac{2}{1+4\beta}+2\sigma}$, $\sigma > 0$ so the truncation bias of $\widehat{\psi}_{2,k}(\tau) \equiv \psi_{2,k}\left(\tau,\widehat{\theta}\right)$ is $O_p\left(n^{-\frac{4\beta}{4\beta+d}}\right)$ and $n^{-\frac{4\beta}{4\beta+d}} \ll \text{var}_\theta\left[\widehat{\psi}_{2,k}(\tau)\right] \asymp k/n^2 = n^{-2\left(\frac{4\beta}{4\beta+d}+\sigma\right)}$. We assume the given $(\beta_g,\beta_b,\beta_p)$ are such that the order $O_p\left[n^{-\left(\frac{\beta_g}{2\beta_g+d}+\frac{\beta_b}{d+2\beta_b}+\frac{\beta_p}{d+2\beta_p}\right)}\right]$ of the estimation bias of $\widehat{\psi}_{2,k}(\tau)$ is $O_p\left(n^{-\frac{4\beta}{4\beta+d}}\right)$. Then $\left|\widehat{\psi}_{2,k}(\tau) - \widetilde{\psi}_k(\tau,\theta)\right|$ and $\left|\widehat{\psi}_{2,k}(\tau) - \psi(\tau,\theta)\right|$ are $O_p\left(n^{-\frac{4\beta}{4\beta+d}+\sigma}\right)$ which just exceeds the minimax rate $O_p\left(n^{-\frac{4\beta}{4\beta+d}}\right)$ for $\sigma$ very small.

Our goal is to compare the coverage and length of confidence intervals for $\widetilde{\tau}_k(\theta)$ and $\tau(\theta)$ based on

$$C_{1-\alpha,\widetilde{\psi}_k(\tau)} \equiv \left\{\tau; -z_{1-\alpha/2} < \frac{\psi_{2,k}\left(\tau,\widehat{\theta}\right)}{\mathbb{W}_{2,\widetilde{\psi}_k(\tau)}\left(\widehat{\theta}\right)} < z_{1-\alpha/2}\right\},$$

$$C_{1-\alpha,2,\widetilde{\tau}_k} \equiv \left\{\tau; -z_{1-\alpha/2} < \frac{\tau_{2,k}\left(\widehat{\theta}\right) - \tau}{\mathbb{W}_{2,\widetilde{\tau}_k}\left(\widehat{\theta}\right)} < z_{1-\alpha/2}\right\},$$

$$C_{1-\alpha,2,ES} \equiv \left\{\tau; -z_{1-\alpha/2} < \frac{\mathbb{ES}_{2,\widetilde{\tau}_k}^{test}\left(\widehat{\theta}(\tau)\right)}{\mathbb{W}_{2,\widetilde{\tau}_k}^{ES}\left(\widehat{\theta}(\tau)\right)} < z_{1-\alpha/2}\right\},$$

where $\mathbb{W}_{2,\widetilde{\psi}_k(\tau)}\left(\widehat{\theta}\right), \mathbb{W}_{2,\widetilde{\tau}_k}\left(\widehat{\theta}\right), \mathbb{W}_{2,\widetilde{\tau}_k}^{ES}\left(\widehat{\theta}(\tau)\right)$ are appropriate variance estimators, $\widehat{\theta}$ is our usual split sample initial estimator, and $\widehat{\theta}(\tau^\dagger)$ is an initial split sample estimator depending on $\tau^\dagger$ that satisfies $\widetilde{\psi}_k\left(\tau,\widehat{\theta}(\tau^\dagger)\right) = 0$ if $\tau = \tau^\dagger$, i.e., $\tau\left[\widehat{\theta}(\tau^\dagger)\right] = \tau^\dagger$. We assume that if $\tau(\theta) = \tau^\dagger$ then the convergence rate under $\theta$ of our estimator of $b(X,\tau^*)$ for any $\tau^*$ remains $n^{-\frac{\beta_b}{d+2\beta_b}}$.

We shall see that the interval based on $C_{1-\alpha,\widetilde{\psi}_k(\tau)}$ outperforms the other two interval estimators. The next theorem gives explicit formulae for $\psi_{2,k}\left(\tau,\widehat{\theta}\right)$,

$\mathbb{ES}_{2,\widetilde{\tau}_k}^{test}\left(\widehat{\theta}(\tau)\right)$, and $\tau_{2,k}\left(\widehat{\theta}\right)$. Using these formulae we calculate the biases and variances necessary to compare the coverage of the three intervals. Before proceeding, note the assumption $0 < \sigma < E_\theta\left[\text{var}_\theta(A|X)\right]$, $E_\theta\left[A^2\right] < c$ implies

$$\left|\widetilde{\tau}_k(\theta) - \tau(\theta)\right| / \left|\widetilde{\psi}_k(\tau,\theta) - \psi(\tau,\theta)\right|$$

is uniformly bounded away from zero and infinity. It then follows from earlier results on $\widehat{\psi}_{2,k}(\tau)$, the assumption $0 < \sigma < E_\theta\left[\text{var}_\theta(A|X)\right], E_\theta\left[A^2\right] < c$, and Theorem 6.9 that $C_{1-\alpha,\widetilde{\psi}_k(\tau)}$ is a uniform asymptotic $1-\alpha$ confidence interval for both $\tau(\theta)$ and $\widetilde{\tau}_k(\theta)$ of length $O_p\left(n^{-\frac{4\beta}{4\beta+d}+\sigma}\right)$.

Our comparison requires each of our three candidate procedures to be on the same scale. Therefore we used standardized versions of the relevant statistics.

**Theorem 6.10.** *Suppose the assumptions described in the preceding example hold. Then*

(i)
$$\psi_{2,k}\left(\tau,\widehat{\theta}\right) = \widetilde{\psi}_{k,\tau}\left(\tau,\widehat{\theta}\right) + \mathbb{IF}_{2,\widetilde{\psi}_k(\tau,\cdot)}\left(\widehat{\theta}\right)$$
$$= \widetilde{\psi}_{k,\tau}\left(\tau,\widehat{\theta}\right) + \mathbb{V}\left[\left(Y^*(\tau) - \widehat{b}(X,\tau)\right)\{A - \widehat{p}(X)\}\right]$$
$$+ \mathbb{V}\left[\left\{\left[Y^*(\tau) - \widehat{b}(X,\tau)\right]\overline{Z}_k^T\right\}_{i_1}\left\{\overline{Z}_k\left[A - \widehat{p}(X)\right]\right\}_{i_2}\right],$$

where $\widehat{b}(X,\tau) = \widehat{B}(\tau) = E_{\widehat{\theta}}(Y^*(\tau)|X), \widehat{p}(X) = \widehat{P} = E_{\widehat{\theta}}(A|X);$

(ii) Let $\widehat{\epsilon}$ denote $Y - \widehat{b}(X)$, and $\widehat{\Delta}$ denote $A - \widehat{p}(X)$. Thus,

$\mathbb{ES}_{2,\widetilde{\tau}_k}^{test}\left(\widehat{\theta}(\tau^\dagger)\right)$

$= v\left(\widehat{\theta}(\tau^\dagger)\right)\left\{var_{\widehat{\theta}(\tau^\dagger)}\left[\mathbb{IF}_{2,\widetilde{\psi}_k(\tau^\dagger,\cdot)}\left(\widehat{\theta}(\tau^\dagger)\right)\right] - var\left[\mathbb{U}_{2,2,\widetilde{\tau}_k(\cdot)}^{*,test,\perp}\left(\widehat{\theta}(\tau^\dagger),\tau^\dagger\right)\right]\right\}^{-1}$

$\times\left\{\mathbb{IF}_{2,\widetilde{\psi}_k(\tau^\dagger,\cdot)}\left(\widehat{\theta}(\tau^\dagger)\right) - \mathbb{U}_{2,2,\widetilde{\tau}_k(\cdot)}^{*,test,\perp}\left(\widehat{\theta}(\tau^\dagger),\tau^\dagger\right)\right\},$

where

$$\mathbb{U}_{2,2,\widetilde{\tau}_k(\cdot),ij}^{*,test,\perp}\left(\widehat{\theta}(\tau^\dagger),\tau^\dagger\right)$$
$$= \left(E_{\widehat{\theta}}\left[\widehat{\epsilon}_i^2\widehat{\Delta}_i^2\right]\right)^{-1} \times \widehat{\epsilon}_i\widehat{\Delta}_i$$
$$\left\{\begin{array}{c} -\left\{\begin{array}{c}\left(E_{\widehat{\theta}}\left[\widehat{\epsilon}_i^2\widehat{\Delta}_i^2\right]\right)^{-1}E_{\widehat{\theta}}\left[\widehat{\epsilon}\widehat{\Delta}^2\overline{Z}_k^T\right] \\ \times E_{\widehat{\theta}}\left[\widehat{\epsilon}^2\widehat{\Delta}\overline{Z}_k^T\right]\widehat{\epsilon}_j\widehat{\Delta}_j \end{array}\right\} \\ +E_{\widehat{\theta}}\left[\widehat{\epsilon}_i^2\widehat{\Delta}_i\overline{Z}_{k,i}^T\right]\overline{Z}_{k,j}\widehat{\Delta}_j \\ +E_{\widehat{\theta}}\left[\widehat{\epsilon}_i\widehat{\Delta}_i^2\overline{Z}_{k,i}^T\right]\overline{Z}_{k,j}\widehat{\epsilon}_j \end{array}\right\}$$

and
$$v(\theta) = E_\theta\left[var_\theta(A|X)\right].$$

Also,

$$var_{\widehat{\theta}(\tau^\dagger)}\left\{\mathbb{ES}_{2,\widetilde{\tau}_k(\cdot)}^{test}\left(\widehat{\theta}(\tau^\dagger)\right)\right\}^{-1}\mathbb{ES}_{2,\widetilde{\tau}_k(\cdot)}^{test}\left(\widehat{\theta}(\tau^\dagger)\right)$$

(6.2) $= v\left(\widehat{\theta}(\tau^\dagger)\right)^{-1}\left\{\mathbb{IF}_{2,\widetilde{\psi}_k(\tau^\dagger,\cdot)}\left(\widehat{\theta}(\tau^\dagger)\right) - \mathbb{U}_{2,2,\widetilde{\tau}_k(\cdot)}^{*,test,\perp}\left(\widehat{\theta}(\tau^\dagger),\tau^\dagger\right)\right\}.$

(iii)
$$\tau_{2,k}\left(\widehat{\theta}\right) \equiv \widetilde{\tau}_k\left(\widehat{\theta}\right) + \mathbb{IF}_{1,\widetilde{\tau}_k(\cdot)}\left(\widehat{\theta}\right) + \mathbb{IF}_{2,2,\widetilde{\tau}_k(\cdot)}\left(\widehat{\theta}\right),$$

where

$$\mathbb{IF}_{1,\widetilde{\tau}_k(\cdot)}\left(\widehat{\theta}\right) = \mathbb{IF}_{1,\widetilde{\tau}_k(\cdot)}\left(\widehat{\theta}\right) = \mathbb{V}\left\{v\left(\widehat{\theta}\right)^{-1}\left[\left\{Y - \widehat{b}(X)\right\}\left\{A - \widehat{p}(X)\right\}\right]\right\}$$

with $Y = Y^*\left(\tau\left(\widehat{\theta}\right)\right), \widehat{b}(X) = \widehat{b}\left(X, \tau\left(\widehat{\theta}\right)\right),$
$$\mathbb{IF}_{2,2,\widetilde{\tau}_k(\cdot)}\left(\widehat{\theta}\right) = v\left(\widehat{\theta}\right)^{-1}\left[\mathbb{IF}_{2,2,\widetilde{\psi}_k\left(\tau\left(\widehat{\theta}\right),\cdot\right)}\left(\widehat{\theta}\right) + \mathbb{Q}_{2,2,\widetilde{\tau}_k(\cdot)}\left(\widehat{\theta}\right)\right] \text{ where}$$

$$\mathbb{Q}_{2,2,\widetilde{\tau}_k(\cdot),\overline{i}_2}\left(\widehat{\theta}\right)$$
$$= -\frac{1}{2}v\left(\widehat{\theta}\right)^{-1}\left[\begin{array}{l}\left[\{A-\widehat{p}(X)\}_{i_1}^2 - v\left(\widehat{\theta}\right)\right]\left[\{Y - \widehat{b}(X)\}\{A - \widehat{p}(X)\}\right]_{i_2} + \\ \left[\{A-\widehat{p}(X)\}_{i_2}^2 - v\left(\widehat{\theta}\right)\right]\left[\{Y - \widehat{b}(X)\}\{A - \widehat{p}(X)\}\right]_{i_1}\end{array}\right].$$

*Proof.* The proof of (i) was given earlier. The proofs of (ii) and (iii) are in the Appendix of our technical report. □

**Theorem 6.11.** *Suppose $\widetilde{\tau}_k(\theta) = \tau^\dagger$ and the assumptions of the preceding theorem hold. Then*

(i) $var_\theta\left[\mathbb{U}_{2,2,\widetilde{\tau}_k(\cdot)}^{*,test,\perp}\left(\widehat{\theta}\left(\tau^\dagger\right), \tau^\dagger\right)\right] = o\left(\frac{1}{n}\right),$

$$var_\theta\left[v\left(\widehat{\theta}\right)^{-1}\psi_{2,k}\left(\tau,\widehat{\theta}\right)\right]$$
$$\times\left[var_\theta\left\{var_{\widehat{\theta}(\tau^\dagger)}\left\{\mathbb{ES}_{2,\widetilde{\tau}_k(\cdot)}^{test}\left(\widehat{\theta}\left(\tau^\dagger\right)\right)\right\}^{-1}\mathbb{ES}_{2,\widetilde{\tau}_k(\cdot)}^{test}\left(\widehat{\theta}\left(\tau^\dagger\right)\right)\right\}\right]^{-1}$$
$$= 1 + o_p(1)$$

(ii)
$$var_\theta\left[\mathbb{Q}_{2,2,\widetilde{\tau}_k(\cdot)}\left(\widehat{\theta}\right)\right] = o\left(\frac{1}{n}\right),$$
$$var_\theta\left[v\left(\widehat{\theta}\right)^{-1}\psi_{2,k}\left(\tau,\widehat{\theta}\right)\right] / var_\theta\left\{\tau_{2,k}\left(\widehat{\theta}\right) - \tau^\dagger\right\} = 1 + o_p(1)$$

(iii)
$$v\left(\widehat{\theta}\right)^{-1}E_\theta\left[\psi_{2,k}\left(\tau^\dagger,\widehat{\theta}\right)\right] = O_p\left\{\left(P - \widehat{P}\right)\left(B\left(\tau^\dagger\right) - \widehat{B}\left(\tau^\dagger\right)\right)\left(\frac{g(X)}{\widehat{g}(X)} - 1\right)\right\}$$
$$= O_p\left(n^{-\left(\frac{\beta_g}{2\beta_g+d} + \frac{\beta_b}{d+2\beta_b} + \frac{\beta_p}{d+2\beta_p}\right)}\right)$$

(iv)

$$E_\theta \left[ var_{\widehat{\theta}(\tau^\dagger)} \left\{ \mathbb{ES}^{test}_{2,\widetilde{\tau}_k(\cdot)} \left( \widehat{\theta}(\tau^\dagger) \right) \right\}^{-1} \mathbb{ES}^{test}_{2,\widetilde{\tau}_k} \left( \widehat{\theta}(\tau^\dagger) \right) \right]$$

$$= O_p \left\{ \left( P - \widehat{P} \right) \left( B(\tau^\dagger) - \widehat{B}(\tau^\dagger) \right) \right.$$
$$\left. \times \left[ \left( \frac{g(X)}{\widehat{g}(X)} - 1 \right) + \left( P - \widehat{P} \right) + \left( B(\tau^\dagger) - \widehat{B}(\tau^\dagger) \right) \right] \right\}$$

$$= O_p \left[ \max \left\{ n^{-\left( \frac{\beta_g}{2\beta_g+d} + \frac{\beta_b}{d+2\beta_b} + \frac{\beta_p}{d+2\beta_p} \right)}, n^{-\left( \frac{\beta_b}{d+2\beta_b} + \frac{2\beta_p}{d+2\beta_p} \right)}, n^{-\left( \frac{2\beta_b}{d+2\beta_b} + \frac{\beta_p}{d+2\beta_p} \right)} \right\} \right]$$

(v)

$$E_\theta \left[ \tau_{2,k}\left(\widehat{\theta}\right) - \tau^\dagger \right]$$

$$= O_p \left\{ \begin{array}{c} \left(P - \widehat{P}\right)\left(\frac{g(X)}{\widehat{g}(X)} - 1\right)\left(B - \widehat{B}\right) + \\ \left(P - \widehat{P}\right)^2 \left(P - \widehat{P}\right)\left(B - \widehat{B}\right) \\ + \left(\frac{g(X)}{\widehat{g}(X)} - 1\right)^2 \left[\left(P - \widehat{P}\right) + \left(\frac{g(X)}{\widehat{g}(X)} - 1\right) + \left(B - \widehat{B}\right)\right] \end{array} \right\}$$

$$= O_p \left\{ \max \left\{ \begin{array}{c} n^{-\left(\frac{\beta_g}{2\beta_g+d} + \frac{\beta_p}{d+2\beta_p} + \frac{\beta_b}{d+2\beta_b}\right)}, \\ n^{-\left(\frac{2\beta_p}{d+2\beta_p}\right)} n^{-\frac{\beta_p}{d+2\beta_p}} n^{-\frac{\beta_b}{d+2\beta_b}}, \\ n^{-\frac{2\beta_g}{2\beta_g+d}} \left\{ n^{-\frac{\beta_b}{d+2\beta_b}} + n^{-\frac{\beta_p}{d+2\beta_p}} + + n^{-\frac{\beta_g}{2\beta_g+d}} \right\} \end{array} \right\} \right\}.$$

*Proof.* The proof of part (iii) was given earlier. The remaining parts are proved in the Appendix of our technical report. □

We conclude from this theorem that the savings in variance that comes with using $\mathbb{ES}^{test}_{2,\widetilde{\tau}_k(\cdot)}\left(\widehat{\theta}(\tau^\dagger)\right)$ rather than $\psi_{2,k}\left(\tau,\widehat{\theta}\right)$ is asymptotically negligible even in regard to constants. Similarly, we conclude that the difference in variance that comes with using $\psi_{2,k}\left(\tau,\widehat{\theta}\right)$ rather than $\mathbb{IF}_{2,2,\widetilde{\tau}_k(\cdot)}\left(\widehat{\theta}\right)$ is asymptotically negligible, again even in regard to constants. Further, because $var_\theta \left[ \mathbb{U}^{*,test,\perp}_{2,2,\widetilde{\tau}_k(\cdot)}\left(\widehat{\theta}(\tau^\dagger), \tau^\dagger\right) \right]$ and $var_\theta \left[ \mathbb{Q}_{2,2,\widetilde{\tau}_k(\cdot)}\left(\widehat{\theta}\right) \right]$ are of the order of $o\left(\frac{1}{n}\right)$ as their first order degenerate kernels are both of order $o_p(1)$, and $n^{\frac{4\beta}{4\beta+d}-\sigma}\left\{\psi_{2,k}\left(\tau,\widehat{\theta}\right) - E_\theta\left[\psi_{2,k}\left(\tau,\widehat{\theta}\right)\right]\right\}$ is asymptotically normal, we conclude that

$$n^{\frac{4\beta}{4\beta+d}-\sigma}\left\{\tau_{2,k}\left(\widehat{\theta}\right) - E_\theta\left[\tau_{2,k}\left(\widehat{\theta}\right)\right]\right\},$$

$$n^{\frac{4\beta}{4\beta+d}-\sigma}\left\{\mathbb{ES}^{test}_{2,\widetilde{\tau}_k(\cdot)}\left(\widehat{\theta}(\tau^\dagger)\right)\right\}^{-1}\left[\mathbb{ES}^{test}_{2,\widetilde{\tau}_k(\cdot)}\left(\widehat{\theta}(\tau^\dagger)\right) - E_\theta\left[\mathbb{ES}^{test}_{2,\widetilde{\tau}_k(\cdot)}\left(\widehat{\theta}(\tau^\dagger)\right)\right]\right]$$

and $n^{\frac{4\beta}{4\beta+d}-\sigma}v\left(\widehat{\theta}\right)^{-1}\left\{\psi_{2,k}\left(\tau,\widehat{\theta}\right) - E_\theta\left[\psi_{2,k}\left(\tau,\widehat{\theta}\right)\right]\right\}$ are all asymptotically normal with the same asymptotic variance.

It then follows that a necessary condition for the intervals based on $\psi_{2,k}\left(\tau^\dagger,\widehat{\theta}\right)$, $\mathbb{ES}^{test}_{2,\widetilde{\tau}_k(\cdot)}\left(\widehat{\theta}(\tau)\right)$, and $\tau_{2,k}\left(\widehat{\theta}\right) - \tau$ to cover $\widetilde{\tau}_k(\theta) = \tau^\dagger$ at the nominal $1-\alpha$ level as

$n \to \infty$ is that

$$v\left(\widehat{\theta}\right)^{-1} E_\theta \left[\psi_{2,k}\left(\tau^\dagger, \widehat{\theta}\right)\right],$$

$$\operatorname{var}_{\widehat{\theta}(\tau^\dagger)} \left\{\mathbb{ES}^{test}_{2,\widetilde{\tau}_k(\cdot)}\left(\widehat{\theta}\left(\tau^\dagger\right)\right)\right\}^{-1} E_\theta \left[\mathbb{ES}^{test}_{2,\widetilde{\tau}_k}\left(\widehat{\theta}\left(\tau^\dagger\right)\right)\right]$$

and $E_\theta \left[\tau_{2,k}\left(\widehat{\theta}\right) - \tau^\dagger\right]$ are $O_p\left(n^{-\frac{4\beta}{4\beta+d}+\sigma}\right)$.

Now we know under the assumptions of Theorem 6.11 that this necessary condition holds for $v\left(\widehat{\theta}\right)^{-1} E_\theta \left[\psi_{2,k}\left(\tau^\dagger, \widehat{\theta}\right)\right]$ since $v\left(\widehat{\theta}\right)$ is bounded away from zero and one and, by assumption, $n^{-\left(\frac{\beta_g}{2\beta_g+d}+\frac{\beta_b}{d+2\beta_b}+\frac{\beta_p}{d+2\beta_p}\right)} = O_p\left(n^{-\frac{4\beta}{4\beta+d}}\right)$. However, this necessary condition need not hold for either

$$\operatorname{var}_{\widehat{\theta}(\tau^\dagger)} \left\{\mathbb{ES}^{test}_{2,\widetilde{\tau}_k(\cdot)}\left(\widehat{\theta}\left(\tau^\dagger\right)\right)\right\}^{-1} E_\theta \left[\mathbb{ES}^{test}_{2,\widetilde{\tau}_k}\left(\widehat{\theta}\left(\tau^\dagger\right)\right)\right],$$

or $E_\theta \left[\tau_{2,k}\left(\widehat{\theta}\right) - \tau^\dagger\right]$. For example, consider the following specification consistent with our assumptions: $\beta_p/d = 0$, $\beta_b/d = \beta_g/d = 1/4$. Then $\beta/d = 1/8$, so $n^{-\left(\frac{\beta_g}{2\beta_g+d}+\frac{\beta_b}{d+2\beta_b}+\frac{\beta_p}{d+2\beta_p}\right)} = n^{-\frac{4\beta}{4\beta+d}} = n^{-1/3}$. However, $E_\theta \left[\tau_{2,k}\left(\widehat{\theta}\right) - \tau^\dagger\right]$ converges to zero at rate $n^{-\frac{\beta_b}{d+2\beta_b}} = n^{-\frac{1}{6}}$. Next

$$\operatorname{var}_{\widehat{\theta}(\tau^\dagger)} \left\{\mathbb{ES}^{test}_{2,\widetilde{\tau}_k(\cdot)}\left(\widehat{\theta}\left(\tau^\dagger\right)\right)\right\}^{-1} E_\theta \left[\mathbb{ES}^{test}_{2,\widetilde{\tau}_k}\left(\widehat{\theta}\left(\tau^\dagger\right)\right)\right]$$
$$= O_p\left(n^{-\left(\frac{\beta_b}{d+2\beta_b}+\frac{2\beta_p}{d+2\beta_p}\right)}\right) = n^{-1/6} \gg O_p\left(n^{-\frac{4\beta}{4\beta+d}+\sigma}\right) = n^{-1/3+\sigma},$$

for small $\sigma$. We conclude that the intervals based on $\mathbb{ES}^{test}_{2,\widetilde{\tau}_k(\cdot)}\left(\widehat{\theta}(\tau)\right)$ and $\tau_{2,k}\left(\widehat{\theta}\right) - \tau$ fail to cover $\widetilde{\tau}_k(\theta) = \tau^\dagger$ at the nominal $1-\alpha$ level uniformly over $\Theta$ as $n \to \infty$. We reach the identical conclusion with regard to the parameter $\tau(\theta)$ because under our assumptions $|\tau(\theta) - \widetilde{\tau}_k(\theta)| = O_p\left(n^{-\frac{4\beta}{4\beta+d}+\sigma}\right)$.

Furthermore, by the argument used in the proof of Theorem 6.9, it is easy to see that the length of each interval is $O_p\left(k/n^2\right) = O_p\left(n^{-\frac{4\beta}{4\beta+d}+\sigma}\right)$. It follows that if we try to improve the coverage of the intervals based on $\mathbb{ES}^{test}_{2,\widetilde{\tau}_k(\cdot)}\left(\widehat{\theta}(\tau)\right)$ and $\tau_{2,k}\left(\widehat{\theta}\right) - \tau$ by further increasing $k$, the length of the intervals will increase beyond $O_p\left(n^{-\frac{4\beta}{4\beta+d}+\sigma}\right)$. We conclude that the interval based on $\psi_{2,k}\left(\tau, \widehat{\theta}\right)$ is strictly preferred to the other two intervals when $\beta_p/d = 0$, $\beta_b/d = \beta_g/d = 1/4$ and is never worse in terms of shrinkage rate and coverage than the other two intervals whatever be $\beta_p$, $\beta_b$, and $\beta_g$. We reach the identical conclusion with regard to the coverage of the parameter $\tau(\theta)$ because, under our assumptions including our choice of $k$, $|\tau(\theta) - \widetilde{\tau}_k(\theta)| = O_p\left(n^{-\frac{4\beta}{4\beta+d}}\right)$ and $n^{-\frac{4\beta}{4\beta+d}} \ll n^{-\frac{4\beta}{4\beta+d}+\sigma}$, the order of the interval lengths.

These results translate directly into analogous results concerning the associated estimators. Under our assumptions the estimator solving $\psi_{2,k}\left(\tau, \widehat{\theta}\right) = 0$ converges to both $\tau(\theta)$ and $\widetilde{\tau}_k(\theta)$ at rate $O_p\left(n^{-\frac{4\beta}{4\beta+d}+\sigma}\right)$. In contrast the rate of convergence of $\tau_{2,k}\left(\widehat{\theta}\right)$ and the estimator solving $\mathbb{ES}^{test}_{2,\widetilde{\tau}_k(\cdot)}\left(\widehat{\theta}(\tau)\right) = 0$ converge to $\tau(\theta)$ and $\widetilde{\tau}_k(\theta)$ at the rates given in (iv) and (v) of Theorem 6.11.

What is the intuition behind the above findings? First note that, as promised by Theorem 2.2 and part (vii) of the theorem in the last subsection, the bias away from zero of $\mathrm{var}_{\widehat{\theta}(\tau^\dagger)}\left\{\mathbb{ES}^{test}_{2,\widetilde{\tau}_k(\cdot)}\left(\widehat{\theta}(\tau^\dagger)\right)\right\}^{-1} E_\theta\left[\mathbb{ES}^{test}_{2,\widetilde{\tau}_k}\left(\widehat{\theta}(\tau^\dagger)\right)\right]$, $E_\theta\left[\tau_{2,k}\left(\widehat{\theta}\right) - \widetilde{\tau}_k(\theta)\right]$, and $v\left(\widehat{\theta}\right)^{-1} E_\theta\left[\psi_{2,k}\left(\tau^\dagger,\widehat{\theta}\right)\right]$ are all $O_p\left(\left\|\widehat{\theta}-\theta\right\|^3\right)$. However the nature and convergence rate of the $O_p\left(\left\|\widehat{\theta}-\theta\right\|^3\right)$ term can vary markedly between estimators, attaining a minimum for $E_\theta\left[\psi_{2,k}\left(\tau^\dagger,\widehat{\theta}\right)\right]$. Now it is not surprising that, for the same order of variance, the order of $E_\theta\left[\tau_{2,k}\left(\widehat{\theta}\right) - \widetilde{\tau}_k(\theta)\right]$ often exceeds that of $E_\theta\left[\psi_{2,k}\left(\tau^\dagger,\widehat{\theta}\right)\right]$. Confidence intervals for $\widetilde{\tau}_k(\theta)$ based on $\tau_{2,k}\left(\widehat{\theta}\right)$ are centered at (i.e are symmetric around) $\tau_{2,k}\left(\widehat{\theta}\right)$, which is a quite stringent constraint on the form of the interval. In that sense, intervals based on $\tau_{2,k}\left(\widehat{\theta}\right)$ are a higher order generalization of the first order asymptotic Wald intervals for $\widetilde{\tau}_k(\theta)$. It is well known that when $\widetilde{\tau}_k(\theta)$ is an implicit parameter that sets a functional such as $\widetilde{\psi}_k(\tau,\theta)$ to zero, first-order Wald confidence intervals are often outperformed in finite samples by confidence sets obtained by inverting a 'score-like' test based on first order 'estimating functions' for the functional that depend on the parameter $\widetilde{\tau}_k$ and, frequently, on estimated nuisance parameters as well, although this fact is not reflected in the first order asymptotics. Our example is higher order version of this phenomenon, where the benefit of the interval $C_{1-\alpha,\widetilde{\psi}_k(\tau)}$ obtained by inverting tests based on the estimating function $\psi_{2,k}\left(\tau,\widehat{\theta}\right)$ for the functional $\widetilde{\psi}_k(\tau,\theta)$ is clearly and quantitatively revealed by the asymptotics. Note that, like first order Wald intervals, the interval based on $\tau_{2,k}\left(\widehat{\theta}\right)$ will differ from the interval for $\widetilde{\tau}_k(\theta)$ based on applying an inverse nonlinear monotone transform $h^{-1}(\cdot)$ to the end points of a Wald interval for the transformed parameter $h\{\widetilde{\tau}_k(\theta)\}$ that is centered on $h(\tau)_{2,k}\left(\widehat{\theta}\right) \equiv h\left(\widetilde{\tau}_k\left(\widehat{\theta}\right)\right) + \mathbb{IF}_{2,h(\widetilde{\tau}_k(\cdot))}\left(\widehat{\theta}\right)$. In contrast, like first order score-based intervals, the intervals based on $\psi_{2,k}\left(\tau,\widehat{\theta}\right)$ and $\mathbb{ES}^{test}_{2,\widetilde{\tau}_k(\cdot)}\left(\widehat{\theta}(\tau^\dagger)\right)$ are invariant to monotone transformations of the parameter $\widetilde{\tau}_k(\theta)$.

More interesting and perhaps more surprising is that, for the same order of variance, the order of $E_\theta\left[\mathrm{var}_{\widehat{\theta}(\tau^\dagger)}\left\{\mathbb{ES}^{test}_{2,\widetilde{\tau}_k(\cdot)}\left(\widehat{\theta}(\tau^\dagger)\right)\right\}^{-1} \mathbb{ES}^{test}_{2,\widetilde{\tau}_k}\left(\widehat{\theta}(\tau^\dagger)\right)\right]$ exceeds that of $E_\theta\left[\psi_{2,k}\left(\tau^\dagger,\widehat{\theta}\right)\right]$. The surprise derives from a failure to recognize that Theorem 6.5 is simply too general to help select among competing procedures. For example, this theorem implies that under law $\widehat{\theta}(\tau^\dagger)$, (a) the variance

$$\mathrm{var}_{\widehat{\theta}(\tau^\dagger)}\left\{\mathbb{ES}^{test}_{2,\widetilde{\tau}_k(\cdot)}\left(\widehat{\theta}(\tau^\dagger)\right)\right\}^{-1}$$

of

$$\left[\mathrm{var}_{\widehat{\theta}(\tau^\dagger)}\left\{\mathbb{ES}^{test}_{2,\widetilde{\tau}_k(\cdot)}\left(\widehat{\theta}(\tau^\dagger)\right)\right\}^{-1} E_\theta\left[\mathbb{ES}^{test}_{2,\widetilde{\tau}_k}\left(\widehat{\theta}(\tau^\dagger)\right)\right]\right]$$

is less (and generally strictly less ) than the variance of

$$v\left(\widehat{\theta}\left(\tau^{\dagger}\right)\right)^{-1} E_{\theta}\left[\psi_{2,k}\left(\tau^{\dagger},\widehat{\theta}\left(\tau^{\dagger}\right)\right)\right],$$

while (b) both have bias of $O_p\left(\left\|\widehat{\theta}\left(\tau^{\dagger}\right)-\theta\right\|^3\right)$. At first blush, this might suggest that the estimator solving $\mathbb{ES}_{2,\widetilde{\tau}_k}^{test}\left(\widehat{\theta}\left(\tau\right)\right)=0$ would likely have the same bias but smaller variance than the estimator solving $\psi_{2,k}\left(\tau,\widehat{\theta}\right)=0$. But we have seen that just the opposite is true. The reason is that the difference between the variances in (a) is negligible in the sense that their ratio is $1+o_p\left(1\right)$, while the $O_p\left(\left\|\widehat{\theta}\left(\tau^{\dagger}\right)-\theta\right\|^3\right)$ biases are often of quite different orders with that of $v\left(\widehat{\theta}\left(\tau^{\dagger}\right)\right)^{-1} E_{\theta}\left[\psi_{2,k}\left(\tau^{\dagger},\widehat{\theta}\left(\tau^{\dagger}\right)\right)\right]$ always a minimum. Furthermore, the theory of higher order estimation and testing influence functions, as a theory of score functions, is, in itself, insufficient to order these biases. Rather side calculations were required. See Remark 4.3 above for further discussion.

More generally, whenever the functional $\psi(\tau,\theta)$ is in our doubly robust class, Equation (4.1) holds so $\widehat{\psi}_{\mathcal{K}_J}^{eff}$ is rate minimax (or near minimax if $\sigma$ is chosen positive), and the suppositions of Theorem 6.1 hold for $\widetilde{\psi}(\tau) = \widehat{\psi}_{\mathcal{K}_J}^{eff}(\tau)$, Theorem 6.1 then implies the width of the interval estimator for $\tau(\theta)$ based on $\widehat{\psi}_{\mathcal{K}_J}^{eff}(\tau)$ converges to zero at the convergence rate of $\widehat{\psi}_{\mathcal{K}_J}^{eff}(\tau)$ to $\psi(\tau,\theta)$.

**Appendix**

In the following, we assume all parametric submodels are sufficiently smooth and regular that expectation and differentiation operators commute as needed. We also define $\mathbb{IF}_{1,1}$ to be $\mathbb{IF}_1$.

*Proof of Theorem 2.2.* Define the bias function $B_m\left[\theta^{\dagger},\theta\right]$ of $\mathbb{IF}_m(\theta)$ to be $E_{\theta^{\dagger}}\left[\mathbb{IF}_m(\theta)\right]$. Define

$$B_{m,l_1^*\ldots l_j^* l_{j+1},\ldots l_s}[\theta,\theta] =$$
$$\partial^s B_m\left[\widetilde{\theta}(\varsigma^*),\widetilde{\theta}(\varsigma)\right]/\partial\varsigma_{l_1^*}^*\ldots\partial\varsigma_{l_j^*}^*\partial\varsigma_{l_{j+1}}\ldots\partial\varsigma_{l_s}\big|_{\varsigma^*=\widetilde{\theta}^{-1}\{\theta\},\varsigma=\widetilde{\theta}^{-1}\{\theta\}},$$

where we reserve $*$ for differentiation with respect to the first argument of $B_m[\cdot,\cdot]$. Thus for $s \leq m$,
$$\psi_{\backslash l_1\ldots l_s}(\theta) = B_{m,l_1^*\ldots l_s^*}[\theta,\theta].$$

To prove the theorem we will first need to show that:

(A.1) $\qquad B_{m,l_1^*\ldots l_j^* l_{j+1},\ldots l_s}[\theta,\theta] = 0$ for $m \geq s > j > 0$

To this end note that for $j < m$,

$$\psi_{\backslash l_1\ldots l_{j+1}}(\theta) = \partial\psi_{\backslash l_1\ldots l_j}(\theta)/\partial\varsigma_{l_{j+1}} = \partial B_{m,l_1^*\ldots l_j^*}[\theta,\theta]/\partial\varsigma_{l_{j+1}}$$
$$= B_{m,l_1^*\ldots l_j^* l_{j+1}^*}[\theta,\theta] + B_{m,l_1^*\ldots l_j^* l_{j+1}}[\theta,\theta]$$
$$= \psi_{\backslash l_1\ldots l_{j+1}}(\theta) + B_{m,l_1^*\ldots l_j^* l_{j+1}}[\theta,\theta],$$

where the second equality is by the definition of $\mathbb{IF}_m(\theta)$, the third is by the chain rule, and the fourth is again by the definition of $\mathbb{IF}_m(\theta)$. Hence $B_{m,l_1^*...l_j^*l_{j+1}}[\theta,\theta] = 0$. Hence for $j \leq m-2$,

$$\begin{aligned}
0 &= \partial B_{m,l_1^*...l_j^*l_{j+1}}[\theta,\theta]/\partial \varsigma_{l_{j+2}...} \\
&= B_{m,l_1^*...l_j^*l_{j+2}^*l_{j+1}}[\theta,\theta] + B_{m,l_1^*...l_j^*l_{j+1}l_{j+2}}[\theta,\theta] \\
&= 0 + B_{m,l_1^*...l_j^*l_{j+1}l_{j+2}}[\theta,\theta],
\end{aligned}$$

where the last equality holds because we just proved $B_{m,l_1^*...l_j^*l_{j+1}}[\theta,\theta] = 0$ for arbitrary indices. Iterating this argument proves (A.1). We complete the proof by induction on $s$ for some $s < m$. Given a $s = 1$ dimensional regular parametric submodel $\widetilde{\theta}(\varsigma)$, $E_{\theta(\varsigma)}[\mathbb{IF}_m(\theta(\varsigma))] = 0$ by assumption. Hence, by regularity of the model, $0 = B_{m,l_1^*.}[\theta,\theta] + B_{m,l_1.}[\theta,\theta]$. Therefore $B_{m,l_1.}[\theta,\theta] = -\psi_{\setminus l_1}(\theta)$. Now suppose the theorem is true for $s$. Then

$$\begin{aligned}
-\psi_{\setminus l_1...l_{s+1}}(\theta) &= -\partial \psi_{\setminus l_1...l_s}(\theta)/\partial \varsigma_{l_{s+1}} \\
&= \partial B_{m,l_1...l_s}[\theta,\theta]/\partial \varsigma_{l_{s+1}} \\
&= B_{m,l_{s+1}^* l_1...l_s}[\theta,\theta] + B_{m,l_1...l_{s+1}}[\theta,\theta] \\
&= 0 + B_{m,l_1...l_{s+1}}[\theta,\theta],
\end{aligned}$$

where the second equality is by the induction assumption, the third by the chain rule, and the last by Equation (A.1). $\square$

*Proof of Theorem 2.3.* (1) Consider two influence functions $\mathbb{IF}_m^{(1)}(\theta)$ and $\mathbb{IF}_m^{(2)}(\theta)$ for $\psi(\theta)$. Then $E_\theta\left[\left\{\mathbb{IF}_m^{(1)}(\theta) - \mathbb{IF}_m^{(2)}(\theta)\right\}\widetilde{\mathbb{S}}_{s,\bar{l}_s}(\theta)\right] = \psi_{\setminus \bar{l}_s}(\theta) - \psi_{\setminus \bar{l}_s}(\theta) = 0$ for any score $\widetilde{\mathbb{S}}_{s,\bar{l}_s}(\theta), s \leq m$ and hence for any linear combination of scores. But, by definition, linear combinations of scores are dense in $\Gamma_m(\theta)$. Thus $\mathbb{IF}_m^{(1)}(\theta)$ and $\mathbb{IF}_m^{(2)}(\theta)$ have the same projection on $\Gamma_m(\theta)$. (2-3): Essentially immediate from the definitions. (4): For $t \leq s$,

$$\psi_{\setminus \bar{l}_t}(\theta) = E_\theta\left[\mathbb{IF}_m(\theta)\widetilde{\mathbb{S}}_{t,\bar{l}_t}(\theta)\right] = E_\theta\left[\Pi_{m,\theta}\left[\mathbb{IF}_m(\theta)|\mathcal{U}_t(\theta)\right]\widetilde{\mathbb{S}}_{t,\bar{l}_t}(\theta)\right]$$

for any $\widetilde{\mathbb{S}}_{t,\bar{l}_t}(\theta)$. (5.a): follows from (1). (5.b): follows from (4). Degeneracy of $\mathbb{IF}_{mm}(\theta)$ follows at once from the fact that $\mathbb{IF}_{mm}(\theta) \in \mathcal{U}_{m-1}(\theta)^\perp$ in $\mathcal{U}_m(\theta)$. Proof of part (5.c) requires the following. $\square$

**Lemma A.1.** *Suppose, for $m \geq 1$, $\mathbb{IF}_{m,m}(\theta)$ and $if_{1,if_{m,m}^{sym}(O_{i_1},...,O_{i_m};\cdot)}(O_{i_{m+1}};\theta)$ exist w.p.1 for a kernel $IF_{m,m}(\theta)$. Let $f\left(O;\widetilde{\theta}(\zeta)\right)$, $\zeta^T = (\zeta_1,...,\zeta_s)$, denote an arbitrary smooth s-dimensional parametric submodel. Let $l_t \in \{1,2,...,s\}$, and $S_{l_t}(O)$ be the score for $\zeta_{l_t}$ evaluated at $\theta$. Then,*

(i) $-if_{1,if_{m,m}^{sym}(O_{i_1},...,O_{i_m};\cdot)}(O_{i_{m+1}};\theta)s_{l_t}(O_{i_{m+1}})$, $-if_{m,m,\setminus l_t}(O_{i_1},...,O_{i_m};\theta)$, *and $if_{m,m}(O_{i_1},...,O_{i_m};\theta)s_{l_t}(O_{i_m})$ each have the same mean given $O_{i_1},...,O_{i_{m-1}}$,*

(ii) $E\left[if_{m,m,\setminus l_t}(O_{i_1},...,O_{i_m};\theta)|O_{i_1},...,O_{i_{m-2}}\right] = 0$,

(iii) $E_\theta\left[if_{1,if_{m,m}^{sym}(O_{i_1},...,O_{i_m};\cdot)}(O_{i_{m+1}};\theta)|O_{i_1},...,O_{i_{m-2}},O_{i_{m+1}}\right] = 0$, *so*

$$\Pi\left[\mathbb{V}\left[if_{1,if_{m,m}^{sym}(O_{i_1},...,O_{i_m};\cdot)}(O_{i_{m+1}};\theta)\right]|\mathcal{U}_m(\theta)\right]$$
$$= \Pi\left[\mathbb{V}\left[if_{1,if_{m,m}^{sym}(O_{i_1},...,O_{i_m};\cdot)}(O_{i_{m+1}};\theta)\right]|\mathcal{U}_m(\theta) \cap \mathcal{U}_{m-2}^\perp(\theta)\right],$$

(iv) $\mathbb{IF}_{m,m,\setminus l_t}(\theta)$ *satisfies* $\Pi_\theta\left[\mathbb{IF}_{m,m,\setminus l_t}(\theta)|\mathcal{U}_{m-2}(\theta)\right]=0$ *and*

$$\Pi_\theta\left[\mathbb{IF}_{m,m,\setminus l_t}(\theta)|\mathcal{U}_{m-1}(\theta)\right]=-\mathbb{V}\left[mE_\theta\left[IF^{sym}_{m,m,\setminus l_t,\bar{i}_m}(\theta)|O_{i_1},\ldots,O_{i_{m-1}}\right]\right].$$

*Proof.* (i) By $IF_{m,m}(\theta)$ degenerate,

$$E_\theta\left[IF_{m,m,\setminus l_t,\bar{i}_m}(\theta)|O_{i_1},\ldots,O_{i_{m-1}}\right]$$
$$=-E_\theta\left[IF_{m,m,\bar{i}_m}(\theta)s_{l_t}(O_{i_m})|O_{i_1},\ldots,O_{i_{m-1}}\right].$$

Further, by definition,

$$E_\theta\left[if_{1,if^{sym}_{m,m}(O_{i_1},\ldots,O_{i_m};\cdot)}(O_{i_{m+1}};\theta)s_{l_t}(O_{i_{m+1}})|O_{i_1},\ldots,O_{i_m}\right]$$
$$=E_\theta\left[IF_{m,m,\setminus l_t,\bar{i}_m}(\theta)|O_{i_1},\ldots,O_{i_m}\right].$$

(ii) By $IF_{m,m}(\theta)$ degenerate $0=E_\theta\left[IF_{m,m,\bar{i}_m}(\theta)s_{l_t}(O_{i_m})|O_{i_1},\ldots,O_{i_{m-2}}\right]$ w.p.1 and so (ii) follows from (i).

(iii) (i) and (ii) imply

$$0=E_\theta\left[if_{1,if^{sym}_{m,m}(O_{i_1},\ldots,O_{i_m};\cdot)}(O_{i_{m+1}};\theta)s_{l_t}(O_{i_{m+1}})|O_{i_1},\ldots,O_{i_{m-2}}\right]$$
$$=E_\theta\left[E_\theta\left[if_{1,if^{sym}_{m,m}(O_{i_1},\ldots,O_{i_m};\cdot)}(O_{i_{m+1}};\theta)|O_{i_{m+1}},O_{i_1},\ldots,O_{i_{m-2}}\right]\right.$$
$$\left.\times s_{l_t}(O_{i_{m+1}})|O_{i_1},\ldots,O_{i_{m-2}}\right].$$

But, by $s_{l_t}(O_{i_{m+1}})$ an arbitrary mean zero function,

$$E_\theta\left\{if_{1,if^{sym}_{m,m}(O_{i_1},\ldots,O_{i_m};\cdot)}(O_{i_{m+1}};\theta)|O_{i_{m+1}},O_{i_1},\ldots,O_{i_{m-2}}\right\}$$
$$=E_\theta\left\{if_{1,if^{sym}_{m,m}(O_{i_1},\ldots,O_{i_m};\cdot)}(O_{i_{m+1}};\theta)|O_{i_1},\ldots,O_{i_{m-2}}\right\}=0.$$

(iv) By definition,

$$\Pi_\theta\left[\mathbb{IF}_{m,m,\setminus l_t}(\theta)|\mathcal{U}_{m-1}(\theta)\right]=\mathbb{V}\left[\{I-d_{m,\theta}\}\left\{IF_{m,m,\setminus l_t,\bar{i}_m}(\theta)\right\}\right].$$

The result follows by Equation (2.1) and part (ii). □

*Proof of Theorem 5(c)(ii).* Consider a $m$-dimensional parametric submodel

$$f\left(O;\widetilde{\theta}(\zeta)\right)=f(O;\theta)\left\{1+\sum_{l=1}^m \zeta_j a_j(O)\right\},\quad \zeta^T=(\zeta_1,\ldots,\zeta_m),$$

with $E_\theta[a_l(O)]=0$. Since this model is linear in the $\zeta_j$, $f_{\setminus l_1\ldots l_m}(O_j;\theta)=0$ for $m>1$. Hence $\widetilde{\mathbb{S}}_{m,\bar{l}_m}(\theta)$ is degenerate of order $m$, i.e., $\widetilde{\mathbb{S}}_{m,\bar{l}_m}(\theta)\in\mathcal{U}^\perp_{m-1}(\theta)$. Since $\mathbb{IF}_{m-1}(\theta)$ exists, on setting $l_s=s$ for $s=1,\ldots,m$,

$$\partial^{m-1}\psi\left(\widetilde{\theta}(\zeta)\right)/\prod_{j=1}^{m-1}\partial\zeta_{j|\zeta=0}\equiv\psi_{\setminus\bar{l}_{m-1}}(\theta)=E_\theta\left[\mathbb{IF}_{m-1}(\theta)\widetilde{\mathbb{S}}_{m-1,\bar{l}_{m-1}}(\theta)\right].$$

Differentiating the last display with respect to $\zeta_m$ and evaluating at $\zeta = 0$, we obtain

$$\psi_{\backslash \bar{l}_m}(\theta) = E_\theta \left[ \mathbb{IF}_{m-1}(\theta) \widetilde{\mathbb{S}}_{m,\bar{l}_m}(\theta) \right] + E_\theta \left[ \mathbb{IF}_{m-1,\backslash l_m}(\theta) \widetilde{\mathbb{S}}_{m-1,\bar{l}_{m-1}}(\theta) \right]$$
$$= E_\theta \left[ \mathbb{IF}_{m-1,\backslash l_m}(\theta) \widetilde{\mathbb{S}}_{m-1,\bar{l}_{m-1}}(\theta) \right].$$

Now

$$E_\theta \left[ \mathbb{IF}_{m-1,\backslash l_m}(\theta) \widetilde{\mathbb{S}}_{m-1,\bar{l}_{m-1}}(\theta) \right]$$
$$= E_\theta \left[ \mathbb{IF}_{m-2,\backslash l_m}(\theta) \widetilde{\mathbb{S}}_{m-1,\bar{l}_{m-1}}(\theta) \right] + E_\theta \left[ \mathbb{IF}_{m-1,m-1,\backslash l_m}(\theta) \widetilde{\mathbb{S}}_{m-1,\bar{l}_{m-1}}(\theta) \right].$$

Setting $s_{l_r}(O_{i_r},\theta) = a_r(O_{i_r})$, $\widetilde{\mathbb{S}}_{m-1,\bar{l}_{m-1}}(\theta) = \sum_{i_1 \neq \cdots \neq i_{m-1}} \prod_{r=1}^{m-1} a_r(O_{i_r};\theta)$ is degenerate of order $m-1$ so

$$E_\theta \left[ \mathbb{IF}_{m-1,m-1,\backslash l_m}(\theta) \widetilde{\mathbb{S}}_{m-1,\bar{l}_{m-1}}(\theta) \right]$$
$$= (m-1)! E_\theta \left( \left[ i f_{m-1,m-1,\backslash l_m}^{sym}(O_{i_1},\ldots,O_{i_{m-1}};\theta) \right] \prod_{r=1}^{m-1} a_r(O_{i_r};\theta) \right),$$

and $E_\theta \left[ \mathbb{IF}_{m-2,\backslash l_m}(\theta) \widetilde{\mathbb{S}}_{m-1,\bar{l}_{m-1}}(\theta) \right] = 0$. Hence

$$\psi_{\backslash \bar{l}_m}(\theta) = (m-1)! E_\theta \left( i f_{m-1,m-1,\backslash l_m}^{sym}(O_{i_1},\ldots,O_{i_{m-1}};\theta) \prod_{r=1}^{m-1} a_r(O_{i_r};\theta) \right).$$

Now, by the assumed existence of $\mathbb{IF}_m(\theta)$, we also have $\psi_{\backslash \bar{l}_m}(\theta) = E_\theta [\mathbb{IF}_m(\theta) \times \widetilde{\mathbb{S}}_m(\theta)] = m! E_\theta \left( i f_{m,m}^{sym}(O_{i_1},\ldots,O_{i_m};\theta) \prod_{r=1}^{m} a_r(O_{i_r};\theta) \right)$. It follows that, for any choice of $m-1$ mean zero functions $a_r(O)$ under $\theta$,

$$0 = E_\theta \left( \left\{ \begin{array}{c} i f_{m-1,m-1,\backslash l_m}^{sym}(O_{i_1},\ldots,O_{i_{m-1}};\theta) \\ -m E_\theta \left[ i f_{m,m}^{sym}(O_{i_1},\ldots,O_{i_m};\theta) a_m(O_{i_m};\theta) | O_{i_1},\ldots,O_{i_{m-1}} \right] \end{array} \right\} \right.$$
$$\left. \times \prod_{r=1}^{m-1} a_r(O_{i_r};\theta) \right)$$
$$= E_\theta \left( r(O_{i_1},\ldots,O_{i_{m-1}};\theta) \prod_{r=1}^{m-1} a_r(O_{i_r};\theta) \right),$$

where

$$r(O_{i_1},\ldots,O_{i_{m-1}};\theta)$$
$$\equiv d_{m-1,\theta} \left[ i f_{m-1,m-1,\backslash l_m}^{sym}(O_{i_1},\ldots,O_{i_{m-1}};\theta) \right]$$
$$- m E_\theta \left[ i f_{m,m}^{sym}(O_{i_1},\ldots,O_{i_m};\theta) a_m(O_{i_m};\theta) | O_{i_1},\ldots,O_{i_{m-1}} \right].$$

The last equality follows from

$$i f_{m-1,m-1,\backslash l_m}^{sym}(O_{i_1},\ldots,O_{i_{m-1}};\theta) - d_{m-1,\theta} \left[ i f_{m-1,m-1,\backslash l_m}^{sym}(O_{i_1},\ldots,O_{i_{m-1}};\theta) \right]$$

orthogonal to $\prod_{r=1}^{m-1} a_r(O_{i_r}; \theta)$. We conclude $r(O_{i_1}, \ldots, O_{i_{m-1}}; \theta) = 0$ with probability 1 because $r(O_{i_1}, \ldots, O_{i_{m-1}}; \theta)$ is a degenerate U-statistic kernel of order $m-1$ and all degenerate U-statistics of order $m-1$ have kernels that are the (possibly infinite) sum of products of $m-1$ mean zero functions. It follows that, on a set $\mathcal{O}_{m-1}$ which has probability 1 under $F^{(m-1)}(\cdot; \theta)$,

$$if^{sym}_{m-1,m-1,\backslash l_m}(o_{i_1}, \ldots, o_{i_{m-1}}; \theta)$$
$$= E_\theta \left[ \{ m \times if^{sym}_{m,m}(o_{i_1}, \ldots, o_{i_{m-1}}, O_{i_m}; \theta) a_m(O_{i_m}; \theta) \} \right]$$
$$+ \{I - d_{m-1,\theta}\} \left[ if^{sym}_{m-1,m-1,\backslash l_m}(o_{i_1}, \ldots, o_{i_{m-1}}; \theta) \right]$$
$$= E_\theta \left[ \left\{ \begin{array}{c} m \times if^{sym}_{m,m}(o_{i_1}, \ldots, o_{i_{m-1}}, O; \theta) \\ -\sum_{j=1}^{m-1} if^{sym}_{m-1,m-1}(o_{i_1}, \ldots, o_{i_{j-1}}, O, o_{i_{j+1}}, \ldots, o_{i_{m-1}}; \theta) \end{array} \right\} a_m(O; \theta) \right]$$

since, by parts (i) and (ii) of the Lemma A.1 and Equation (2.1),

$$\{I - d_{m-1,\theta}\} \left[ if^{sym}_{m-1,m-1,\backslash l_m}(o_{i_1}, \ldots, o_{i_{m-1}}; \theta) \right]$$
$$= -E_\theta \left[ \sum_{j=1}^{m-1} if^{sym}_{m-1,m-1}(o_{i_1}, \ldots, o_{i_{j-1}}, O, o_{i_{j+1}}, \ldots, o_{i_{m-1}}; \theta) a_m(O; \theta) \right].$$

Here $I$ is the identity operator. Now since the model $f\left(O; \widetilde{\theta}(\zeta)\right) = f(O; \theta)\{1 + \zeta_m a_m(O)\}$ with $\zeta_s = 0$ for $s < m$ has score $a_m(O)$ and such scores are dense in the subspace of $L_2(F(\cdot; \theta))$ with mean zero, it follows that $if^{sym}_{m-1,m-1}(o_{i_1}, \ldots, o_{i_{m-1}}; \theta)$ has influence function

$$m \times if^{sym}_{m,m}(o_{i_1}, \ldots, o_{i_{m-1}}, O; \theta)$$
$$- \sum_{j=1}^{m-1} if^{sym}_{m-1,m-1}(o_{i_1}, \ldots, o_{i_{j-1}}, O, o_{i_{j+1}}, \ldots, o_{i_{m-1}}; \theta)$$

on the set $\mathcal{O}_{m-1}$. Thus

$$m \times if^{sym}_{m,m}(o_{i_1}, \ldots, o_{i_{m-1}}, O_{i_m}; \theta) = d_{m,\theta} \left[ if_{1, if^{sym}_{m-1,m-1}(o_{i_1}, \ldots, o_{i_{m-1}};\cdot)}(O_{i_m}; \theta) \right].$$

$\square$

**Corollary A.2.** *For $m \geq 2$,*

(A.2) $\qquad \Pi_\theta \left[ \mathbb{IF}_{m-1,m-1,\backslash l_t}(\theta) | \mathcal{U}^\perp_{m-2}(\theta) \right] = -\Pi_\theta \left[ \mathbb{IF}_{m,m,\backslash l_t}(\theta) | \mathcal{U}_{m-1}(\theta) \right]$

(A.3) $\qquad \mathbb{IF}_{m,\backslash l_t}(\theta) = \Pi_\theta \left[ \mathbb{IF}_{m,m,\backslash l_t}(\theta) | \mathcal{U}^\perp_{m-1}(\theta) \right]$

$$E_\theta \left[ \mathbb{IF}_{m,\backslash l_{m+1}}(\theta) \widetilde{\mathbb{S}}_{m,\bar{l}_m}(\theta) \right] =$$

(A.4) $\qquad m! E_\theta \left( if^{sym}_{m,m,\backslash l_{m+1}}(O_{i_1}, \ldots, O_{i_{l_m}}; \theta) \prod_{r=1}^{m} S_{l_r}(O_{i_r}; \theta) \right).$

*Proof of Equation (A.2).* By Lemma A.1 and Theorem 5c(ii),

$$\Pi_\theta \left[ \mathbb{IF}_{m,m,\backslash l_t} (\theta) | \mathcal{U}_{m-1} (\theta) \right]$$
$$= \mathbb{V} \left[ m E_\theta \left( IF^{sym}_{m,m,\bar{i}_m} (\theta) s_{l_t} (O_{i_m}) | O_{i_1}, \ldots, O_{i_{m-1}} \right) \right]$$
$$= \mathbb{V} \left[ m E_\theta \left( m^{-1} d_{m,\theta} \left\{ if_{1, if^{sym}_{m-1,m-1}(O_{i_1}, \ldots, O_{i_{m-1}};\cdot)} (O_{i_m}; \theta) \right\} s_{l_t} (O_{i_m}) | \right. \right.$$
$$\left. \left. \times O_{i_1}, \ldots, O_{i_{m-1}} \right) \right].$$

Now, by part (iii) of Lemma A.1 and Equation (2.1), the RHS is

$$\mathbb{V} \left[ E_\theta \left( if_{1, if^{sym}_{m-1,m-1}(O_{i_1}, \ldots, O_{i_{m-1}};\cdot)} (O_{i_m}; \theta) s_{l_t} (O_{i_m}) | O_{i_1}, \ldots, O_{i_{m-1}} \right) \right]$$
$$- \mathbb{V} \left\{ E \left[ E \left[ (m-1) E \left[ if_{1, if^{sym}_{m-1,m-1}(O_{i_1}, \ldots, O_{i_{m-1}};\cdot)} (O_{i_m}; \theta) | O_{i_m}, O_{i_1}, \ldots, O_{i_{m-2}} \right] \right] \right. \right.$$
$$\left. \left. \times s_{l_t} (O_{i_m}) | O_{i_1}, \ldots, O_{i_{m-1}} \right] \right\}$$
$$= \mathbb{V} \left[ IF^{sym}_{m-1,m-1,\backslash l_t} (\theta) \right] - \mathbb{V} \left\{ (m-1) E_\theta \left[ IF^{sym}_{m-1,m-1,\backslash l_t} (\theta) | O_{i_1}, \ldots, O_{i_{m-2}} \right] \right\}$$

On the other hand, by part (iv) of the Lemma A.1,

$$\Pi_\theta \left[ \mathbb{IF}_{m-1,m-1,\backslash l_t} (\theta) | \mathcal{U}^\perp_{m-2} (\theta) \right]$$
$$= \mathbb{V} \left[ IF_{m-1,m-1,\backslash l_t} (\theta) \right]$$
$$- \mathbb{V} \left[ (m-1) E_\theta \left[ IF^{sym}_{m-1,m-1,\backslash l_t, \bar{i}_{m-1}} (\theta) | O_{i_1}, \ldots, O_{i_{m-2}} \right] \right]. \quad \square$$

*Proof of Equation (A.3).* Write

$$\mathbb{IF}_{m,\backslash l_t} (\theta) = \Pi_\theta \left[ \mathbb{IF}_{m,m,\backslash l_t} (\theta) | \mathcal{U}^\perp_{m-1} (\theta) \right] + \left\{ \Pi \left[ \mathbb{IF}_{2,2,\backslash l_t} (\theta) | \mathcal{U}_1 (\theta) \right] + \mathbb{IF}_{1,\backslash l_t} (\theta) \right\}$$
$$+ \sum_{j=2}^{m-1} \left\{ \Pi \left[ \mathbb{IF}_{j+1,j+1,\backslash l_t} (\theta) | \mathcal{U}_j (\theta) \right] + \Pi \left[ \mathbb{IF}_{jj,\backslash l_t} (\theta) | \mathcal{U}^\perp_{j-1} (\theta) \right] \right\}$$

The RHS is $\Pi_\theta \left[ \mathbb{IF}_{m,m,\backslash l_t} (\theta) | \mathcal{U}^\perp_{m-1} (\theta) \right]$ by Equation (A.2). $\quad \square$

*Proof of Equation (A.4).*

$$E_\theta \left[ \mathbb{IF}_{m,\backslash l_{m+1}} (\theta) \widetilde{\mathbb{S}}_{m,\bar{l}_m} (\theta) \right] = E_\theta \left[ \Pi \left[ \mathbb{IF}_{m,m,\backslash l_{m+1}} (\theta) | \mathcal{U}^\perp_{m-1} (\theta) \right] \widetilde{\mathbb{S}}_{m,\bar{l}_m} (\theta) \right]$$

by Equation (A.3). But the RHS of this equation is the RHS of Equation (A.4). $\quad \square$

*Proof of Theorem 5c(i).* By assumption

$$\psi_{\backslash \bar{l}_{m-1}} (\theta) = E_\theta \left( \mathbb{IF}_{m-1} (\theta) \widetilde{\mathbb{S}}_{m-1,\bar{l}_{m-1}} (\theta) \right).$$

Hence

$$\psi_{\backslash \bar{l}_m} (\theta) = E_\theta \left( \mathbb{IF}_{m-1} (\theta) \widetilde{\mathbb{S}}_{m,\bar{l}_m} (\theta) \right) + E_\theta \left[ \mathbb{IF}_{m-1,\backslash l_m} (\theta) \widetilde{\mathbb{S}}_{m-1,\bar{l}_{m-1}} (\theta) \right].$$

By Equation (A.4), and the assumption $if^{sym}_{m-1,m-1} (O_{i_1}, \ldots, O_{i_m}; \theta)$ has an influence function, we obtain

$$E_\theta \left[ \mathbb{IF}_{m-1,\backslash l_m} (\theta) \widetilde{\mathbb{S}}_{m-1,\bar{l}_{m-1}} (\theta) \right]$$
$$= (m-1)! E_\theta \left( if_{1, if^{sym}_{m-1,m-1}(O_{i_1}, \ldots, O_{i_{m-1}};\cdot)} (O_{i_m}; \theta) S_{l_m} (O_{i_m}; \theta) \prod_{r=1}^{m-1} S_{l_r} (O_{i_r}; \theta) \right).$$

We conclude that $\mathbb{IF}_{m,m}$ exists and equals

$$\mathbb{V}\left[m^{-1}d_{m,\theta}\left\{if_{1,if^{sym}_{m-1,m-1}(O_{i_1},\ldots,O_{i_{m-1}};\cdot)}(O_{i_m};\theta)\right\}\right]. \qquad \square$$

*Proof of Theorem 3.14.* By Theorem 3.13,

$$\begin{aligned}IF_{1,\widetilde{\psi}_k,i_1}(\theta) &= if_{1,\widetilde{\psi}_k}(O_{i_1};\theta) \\ &= H_{i_1}\left(\widetilde{b}(\theta),\widetilde{p}(\theta)\right) - \widetilde{\psi}_k(\theta) \\ &= h\left(O_{i_1},\widetilde{b}(X_{i_1};\theta),\widetilde{p}(X_{i_1};\theta)\right) - \widetilde{\psi}_k(\theta),\end{aligned}$$

and by part 5.c of Theorem 2.3,

$$\mathbb{V}\left[IF_{22,\widetilde{\psi}_k,\overline{i}_2}\right] = \frac{1}{2}\left\{\Pi_\theta\left[\mathbb{V}\left[IF_{1,if_{1,\widetilde{\psi}_k}(O_{i_1},\cdot),i_2}(\theta)\right]|\mathcal{U}_1^{\perp\theta,2}(\theta)\right]\right\}.$$

Now

$$\begin{aligned}IF_{1,if_{1,\widetilde{\psi}_k}(O_{i_1},\cdot),i_2}(\theta) &= h_{\widetilde{b}}\left(O_{i_1},\widetilde{b}(X_{i_1};\theta),\widetilde{p}(X_{i_1};\theta)\right)IF_{1,\widetilde{b}(X_{i_1};\cdot),i_2}(\theta) \\ &\quad + h_{\widetilde{p}}\left(O_{i_1},\widetilde{b}(X_{i_1};\theta),\widetilde{p}(X_{i_1};\theta)\right)IF_{1,\widetilde{p}(X_{i_1};\cdot),i_2}(\theta),\end{aligned}$$

where

$$\begin{aligned}h_{\widetilde{b}}\left(O_{i_1},\widetilde{b}(X_{i_1};\theta),\widetilde{p}(X_{i_1};\theta)\right) &= H_{1,i_1}\widetilde{p}(X_{i_1};\theta) + H_{2,i_1} \\ h_{\widetilde{p}}\left(O_{i_1},\widetilde{b}(X_{i_1};\theta),\widetilde{p}(X_{i_1};\theta)\right) &= H_{1,i_1}\widetilde{b}(X_{i_1};\theta) + H_{3,i_1}.\end{aligned}$$

$$\begin{aligned}IF_{1,\widetilde{b}(X_{i_1};\cdot),i_2}(\theta) &= IF_{1,b^*(X_{i_1},\widetilde{\eta}_k(\cdot)),i_2}(\theta) \\ &= \dot{B}_{i_1}\overline{Z}^T_{ki_1}IF_{1,\widetilde{\eta}_k(\cdot),i_2}(\theta) \\ &= -\dot{B}_{i_1}\overline{Z}^T_{ki_1}\left\{E_\theta\left[\dot{P}\dot{B}H_1\overline{Z}_k\overline{Z}^T_k\right]\right\}^{-1}\left[\left\{H_1\widetilde{b}(X;\theta) + H_3\right\}\dot{P}\overline{Z}_k\right]_{i_2},\end{aligned}$$

and

$$IF_{1,\widetilde{p}(X_{i_1},\cdot),i_2}(\theta) = -\dot{P}_{i_1}\overline{Z}^T_{ki_1}\left\{E_\theta\left[\dot{P}\dot{B}H_1\overline{Z}_k\overline{Z}^T_k\right]\right\}^{-1}\left[\left\{H_1\widetilde{p}(X;\theta) + H_2\right\}\dot{B}\overline{Z}_k\right]_{i_2}$$

$$\begin{aligned}&IF_{1,if_{1,\widetilde{\psi}_k}(O_{i_1};\cdot),i_2}(\theta) \\ &= -\left\{H_1\widetilde{p}(X;\theta) + H_2\right\}_{i_1}\dot{B}_{i_1}\overline{Z}^T_{ki_1}\left\{E_\theta\left[\dot{P}\dot{B}H_1\overline{Z}_k\overline{Z}^T_k\right]\right\}^{-1} \\ &\quad\times\left[\left\{H_1\widetilde{b}(X;\theta) + H_3\right\}\dot{P}\overline{Z}_k\right]_{i_2} \\ &\quad -\left\{H_1\widetilde{b}(X;\theta) + H_3\right\}_{i_1}\dot{P}_{i_1}\overline{Z}^T_{ki_1}\left\{E_\theta\left[\dot{P}\dot{B}H_1\overline{Z}_k\overline{Z}^T_k\right]\right\}^{-1} \\ &\quad\times\left[\left\{H_1\widetilde{p}(X;\theta) + H_2\right\}\dot{B}\overline{Z}_k\right]_{i_2},\end{aligned}$$

and further

$$\Pi_\theta\left[\mathbb{V}\left[IF_{1,if_{1,\widetilde{\psi}_k}(O_{i_1};\cdot),i_2}(\theta)\right]|\mathcal{U}_1(\theta)\right] = 0,$$

since

$$E_\theta \left[ \{H_1 \widetilde{p}(X;\theta) + H_2\} \dot{B} \overline{Z}_k \right] = E_\theta \left[ \{H_1 \widetilde{b}(X;\theta) + H_3\} \dot{P} \overline{Z}_k \right] = 0$$

and thus $IF_{1,if_{1,\widetilde{\psi}_k}(O_{i_1};\cdot),i_2}(\theta)$ is degenerate. Because $IF_{1,if_{1,\widetilde{\psi}_k}(O_{i_1};\cdot),i_2}(\theta)$ has two terms, it appears that $IF_{22,\widetilde{\psi}_k,\bar{i}_2}$ will consist of two terms. However by the symmetry upon interchange of $i_2$ and $i_1$, and the permutation invariance of the operator $\mathbb{V}$

$$\mathbb{V}\left[IF_{1,if_{1,\widetilde{\psi}_k}(O_{i_1};\theta),i_2}(\theta)\right]$$
$$= \mathbb{V}\left[ \begin{array}{c} -2\{H_1\widetilde{p}(X,\theta) + H_2\}_{i_1} \dot{B}_{i_1} \overline{Z}^T_{ki_1} \left\{E_\theta\left[\dot{P}\dot{B}H_1\overline{Z}_k\overline{Z}^T_k\right]\right\}^{-1} \\ \times \left[\overline{Z}_k\{H_1\widetilde{b}(X,\theta)+H_3\}\dot{P}\right]_{i_2} \end{array} \right].$$

Thus we can take

$$IF_{22,\widetilde{\psi}_k,\bar{i}_2} = -\{H_1\widetilde{p}(X,\theta)+H_2\}_{i_1}\dot{B}_{i_1}\overline{Z}^T_{ki_1}\left\{E_\theta\left[\dot{P}\dot{B}H_1\overline{Z}_k\overline{Z}^T_k\right]\right\}^{-1}$$
$$\times \left[\overline{Z}_k\{H_1\widetilde{b}(X,\theta)+H_3\}\dot{P}\right]_{i_2}$$

as was to be proved. We now complete the proof of the Theorem by induction. We assume it is true for $IF_{mm,\widetilde{\psi}_k,\bar{i}_m}$ and prove it is true for $IF_{(m+1)(m+1),\widetilde{\psi}_k,\bar{i}_{m+1}}$. Now

$$\mathbb{V}\left[IF_{(m+1),(m+1),\widetilde{\psi}_k,\bar{i}_{m+1}}(\theta)\right] =$$
$$\frac{1}{m}\mathbb{V}\left[\Pi_\theta\left[IF_{1,if_{m,m,\widetilde{\psi}_k}(O_{\bar{i}_m},\cdot),i_{m+1}}(\theta)\,|\,\mathcal{U}^{\perp\theta,m+1}_m(\theta)\right]\right].$$

Now by the induction hypothesis,

$$if_{m,m,\widetilde{\psi}_k}(O_{\bar{i}_m},\theta)$$
$$= (-1)^{m-1}\left[\left(H_1\widetilde{P}(\theta)+H_2\right)\dot{B}\overline{Z}^T_k\right]_{i_1}$$
$$\times \left[\prod_{s=3}^m \left\{E_\theta\left[\dot{P}\dot{B}H_1\overline{Z}_k\overline{Z}^T_k\right]\right\}^{-1}\left\{\begin{array}{c}\left(\dot{P}\dot{B}H_1\overline{Z}_k\overline{Z}^T_k\right)_{i_s}\\ -E_\theta\left[\dot{P}\dot{B}H_1\overline{Z}_k\overline{Z}^T_k\right]\end{array}\right\}\right]$$
$$\times \left\{E_\theta\left[\dot{P}\dot{B}H_1\overline{Z}_k\overline{Z}^T_k\right]\right\}^{-1}\left[\overline{Z}_k\left(H_1\widetilde{B}(\theta)+H_3\right)\dot{P}\right]_{i_2}.$$

The derivatives with respect to the $\theta's$ in $\widetilde{P}(\theta), \widetilde{B}(\theta)$ and in the $m-1$ terms $\left\{E_\theta\left[\dot{P}\dot{B}H_1\overline{Z}_k\overline{Z}^T_k\right]\right\}^{-1}$ will each contribute a term to $\mathbb{V}\left[IF_{(m+1)(m+1),\widetilde{\psi}_k,\bar{i}_{m+1}}(\theta)\right]$. However differentiating with respect to the $\theta$ in the $m-2$ terms $E_\theta\left[\dot{P}\dot{B}H_1\overline{Z}_k\overline{Z}^T_k\right]$ will not contribute to $\mathbb{V}\left[IF_{(m+1),(m+1),\widetilde{\psi}_k,\bar{i}_{m+1}}(\theta)\right]$ as the contribution from each of these $m-2$ terms to $IF_{1,if_{m,m,\widetilde{\psi}_k}(O_{\bar{i}_m},\theta),i_{m+1}}(\theta)$ is only a function of $m$ units' data and is thus an element of $\mathcal{U}_m(\theta)$ which is orthogonal to the space $\mathcal{U}^{\perp\theta,m+1}_m(\theta)$

that is projected on. Now

$$IF_{1,\{E_\theta[\dot{P}\dot{B}H_1\overline{Z}_k\overline{Z}_k^T]\}^{-1},i_{m+1}}(\theta)$$
$$= -\left\{E_\theta\left[\dot{P}\dot{B}H_1\overline{Z}_k\overline{Z}_k^T\right]\right\}^{-1}\left\{\begin{array}{c}\left(\dot{P}\dot{B}H_1\overline{Z}_k\overline{Z}_k^T\right)_{i_{m+1}}\\ -E_\theta\left[\dot{P}\dot{B}H_1\overline{Z}_k\overline{Z}_k^T\right]\end{array}\right\}\left\{E_\theta\left[\dot{P}\dot{B}H_1\overline{Z}_k\overline{Z}_k^T\right]\right\}^{-1}.$$

so upon permuting the unit indices, the contribution of each of these $m-1$ terms to $IF_{1,if_{m,m,\widetilde{\psi}_k}(O_{\bar{i}_m},\theta),i_{m+1}}(\theta)$ is

(A.5)
$$-(-1)^{m-1}\left[\left(H_1\widetilde{P}(\theta)+H_2\right)\dot{B}\overline{Z}_k^T\right]_{i_1}$$
$$\times\left[\prod_{s=3}^{m+1}\left\{E_\theta\left[\dot{P}\dot{B}H_1\overline{Z}_k\overline{Z}_k^T\right]\right\}^{-1}\left\{\begin{array}{c}\left(\dot{P}\dot{B}H_1\overline{Z}_k\overline{Z}_k^T\right)_{i_s}\\ -E_\theta\left[\dot{P}\dot{B}H_1\overline{Z}_k\overline{Z}_k^T\right]\end{array}\right\}\right]$$
$$\times\left\{E_\theta\left[\dot{P}\dot{B}H_1\overline{Z}_k\overline{Z}_k^T\right]\right\}^{-1}\left[\overline{Z}_k\left(H_1\widetilde{B}(\theta)+H_3\right)\dot{P}\right]_{i_2},$$

which is already degenerate ( i.e., orthogonal to $\mathcal{U}_m(\theta)$). Differentiating with respect to the $\theta's$ of $\widetilde{P}(\theta), \widetilde{B}(\theta)$ in $IF_{1,if_{m,m,\widetilde{\psi}_k}(O_{\bar{i}_m},\theta),i_{m+1}}(\theta)$ we obtain

$$=(-1)^{m-1}IF_{1,\widetilde{b}(X_{i_1},\cdot),i_{m+1}}(\theta)\left[H_1\dot{B}\overline{Z}_k^T\right]_{i_1}$$
$$\left[\prod_{s=3}^{m}\left\{E_\theta\left[\dot{P}\dot{B}H_1\overline{Z}_k\overline{Z}_k^T\right]\right\}^{-1}\left\{\left(\dot{P}\dot{B}H_1\overline{Z}_k\overline{Z}_k^T\right)_{i_s}-E_\theta\left[\dot{P}\dot{B}H_1\overline{Z}_k\overline{Z}_k^T\right]\right\}\right]\times$$
$$\left\{E_\theta\left[\dot{P}\dot{B}H_1\overline{Z}_k\overline{Z}_k^T\right]\right\}^{-1}\left[\overline{Z}_k\left(H_1\widetilde{B}(\theta)+H_3\right)\dot{P}\right]_{i_2}$$
$$+(-1)^{m-1}\left[\left(H_1\widetilde{P}(\theta)+H_2\right)\dot{B}\overline{Z}_k^T\right]_{i_1}\times$$
$$\left[\prod_{s=3}^{m}\left\{E_\theta\left[\dot{P}\dot{B}H_1\overline{Z}_k\overline{Z}_k^T\right]\right\}^{-1}\left\{\left(\dot{P}\dot{B}H_1\overline{Z}_k\overline{Z}_k^T\right)_{i_s}-E_\theta\left[\dot{P}\dot{B}H_1\overline{Z}_k\overline{Z}_k^T\right]\right\}\right]\times$$
$$\left[\overline{Z}_kH_1\dot{P}\right]IF_{1,\widetilde{p}(X_{i_2},\cdot),i_{m+1}}(\theta)$$

Substituting in the above expressions for $IF_{1,\widetilde{b}(X_{i_1},\cdot),i_{m+1}}(\theta)$ and $IF_{1,\widetilde{p}(X_{i_2},\cdot),i_{m+1}}(\theta)$, then projecting on $\mathcal{U}_m^{\perp\theta,m+1}(\theta)$, and again permuting unit indices, we obtain two identical terms both equal to Equation (A.5). Thus we obtain $m+1$ identical terms in all. Upon dividing by $m+1$, we conclude that $\mathbb{V}\left[IF_{(m+1),(m+1),\widetilde{\psi}_k,\bar{i}_{m+1}}(\theta)\right]$ equals $\mathbb{V}$ operating on (A.5), proving the theorem. □

**Acknowledgments.** A proper acknowledgment to David Freedman on the occasion of this Festschrift is to say with Kipling: you are a better man than I. Who else has David's passion to get things exactly right, whether in an elementary textbook or an Annals article, no matter what the cost in time or effort? Who else in our community has stood up and taken as strong a public stand: statistical models and regression are not a magic cure for inadequate data, plagued by uncontrolled confounding and measurement error? How many quantitative social scientists can now write their papers without imagining the withering attack of a David Freedman critique should they take unwarranted shortcuts? How many other established statisticians would invite an unpublished, unsung, nascent want-to-be statistician into their home, treat them as an equal, and give them the gift of 23 years of loyal, if sometimes testy, friendship?

*Higher order influence functions* 421[20] Small, D. and McLeish, C. (1994). *Hilbert Space Methods in Probability and Statistical Inference.* Wiley, New York. MR1269321
[21] Tchetgen, E., Li, L., van der Vaart, A. W. and Robins, J. M. (2006). Robust inference with higher order inference functions: Part I. In *2005 JSM Proceedings* 2644–2651. American Statistical Association, Alexandria.
[22] Tchetgen, E., Li, L., van der Vaart, A. W. and Robins, J. M. (2007). Higher Order U-statistics estimators for longitudinal missing data and causal inference models. Working paper.
[23] van der Laan, M. and Dudoit, S. (2005). Asymptotics of cross-validated risk estimation in estimator selection and performance assessment. *Stat. Methodol.* **2** 131–154. MR2161394
[24] van der Laan, M. and Robins, J. M. (2003). *Unified Methods for Censored Longitudinal Data and Causality.* Springer, New York. MR1958123
[25] van der Vaart, A. W. (1991). On differentiable functionals. *Ann. Statist.* **19** 178–204. MR1091845
[26] van der Vaart, A. W. (1998). *Asymptotic Statistics.* Cambridge Series in Statistical and Probabilistic Mathematics, Cambridge.
[27] Wang, L., Brown, L. D., Cai, T. and Levine, M. (2006). Effect of mean on variance function estimation in nonparametric regression. Technical report.